\documentclass[reqno,11pt]{amsart}

\usepackage[top=2.0cm,bottom=2.8cm,left=2.8cm,right=2cm]{geometry}
\usepackage{amsthm,amsmath,amssymb,dsfont}
\usepackage{mathrsfs,amsfonts,functan,extarrows,mathtools}
\usepackage[colorlinks,urlcolor=black]{hyperref} 
\usepackage{marginnote}
\usepackage{xcolor}
\usepackage{stmaryrd}
\usepackage{esint}
\usepackage{graphicx}
\usepackage{bbm}
\usepackage{enumerate}
\usepackage{tikz}
\usetikzlibrary{positioning}

\numberwithin{equation}{section}
\newtheorem{theorem}{Theorem}[section]
\newtheorem{proposition}[theorem]{Proposition}
\newtheorem{definition}[theorem]{Definition}

\newtheorem{remark}[theorem]{Remark}
\newtheorem{lemma}[theorem]{Lemma}
\newtheorem{reduct}[theorem]{Reduction}
\newtheorem{assumption}[theorem]{Assumption}
\allowdisplaybreaks

\def\p{\partial}
\def\d{\mathrm{d}} 
\def\eps{\varepsilon} 

\newcounter{wronumber}

\newcommand{\newterm}[1]{{\textbf{#1}}}

\newcommand{\pb}{{\boldsymbol{p}}}

\newcommand{\supp}{\mathrm{supp}\,}

\newcommand{\Dirac}{\boldsymbol{\Delta}}
\newcommand{\dirac}{\boldsymbol{\delta}}

\newcommand{\Cb}{\mathbb C}

\newcommand{\Eb}{\mathbb E}

\newcommand{\Lb}{\mathbb L}
\newcommand{\Mb}{\mathbb M}
\newcommand{\Nb}{\mathbb N}

\newcommand{\Pb}{\mathbb P}
\newcommand{\Qb}{\mathbb Q}
\newcommand{\Rb}{\mathbb R}
\newcommand{\Sb}{\mathbb S}

\newcommand{\Vb}{\mathbb V}

\newcommand{\Xb}{\mathbb X}
\newcommand{\Yb}{\mathbb Y}
\newcommand{\Zb}{\mathbb Z}

\newcommand{\oneb}{\mathbbm{1}}

\newcommand{\Ac}{\mathcal A}
\newcommand{\Bc}{\mathcal B}
\newcommand{\Cc}{\mathcal C}
\newcommand{\Dc}{\mathcal D}
\newcommand{\Ec}{\mathcal E}
\newcommand{\Fc}{\mathcal F}
\newcommand{\Gc}{\mathcal G}
\newcommand{\Hc}{\mathcal H}
\newcommand{\Ic}{\mathcal I}
\newcommand{\Jc}{\mathcal J}

\newcommand{\Lc}{\mathcal L}
\renewcommand{\Mc}{\mathcal M}
\newcommand{\Nc}{\mathcal N}
\newcommand{\Oc}{\mathcal O}
\newcommand{\Pc}{\mathcal P}
\newcommand{\Qc}{\mathcal Q}

\newcommand{\Sc}{\mathcal S}
\newcommand{\Tc}{\mathcal T}

\newcommand{\Xc}{\mathcal X}
\newcommand{\Yc}{\mathcal Y}
\newcommand{\Zc}{\mathcal Z}
\newcommand{\Df}{\mathfrak D}

\newcommand{\Sf}{\mathfrak S}
\newcommand{\Bf}{\mathfrak B}

\newcommand{\Rf}{\mathfrak R}
\newcommand{\af}{\mathfrak a}
\newcommand{\mf}{\mathfrak m}
\newcommand{\pf}{\mathfrak p}

\newcommand{\rf}{\mathfrak r}
\newcommand{\of}{\mathfrak o}
\newcommand{\Af}{\mathfrak A}
\newcommand{\lf}{\mathfrak l}
\newcommand{\tf}{\mathfrak t}
\newcommand{\qf}{\mathfrak q}

\newcommand{\Cf}{\mathfrak C} 
\newcommand{\Nf}{\mathfrak N}
\newcommand{\Lf}{\mathfrak L}
\newcommand{\Pf}{\mathsf{P}}

\newcommand{\Xf}{\mathfrak X} 

\newcommand{\Mf}{\mathfrak M} 
\newcommand{\Of}{\mathfrak{O}}

\newcommand{\nf}{\mathfrak n}

\newcommand{\vz}{\boldsymbol{z}}

\newcommand{\vLambda}{\boldsymbol{\Lambda}}

\newcommand{\vt}{\boldsymbol{t}}

\newcommand{\vf}{\mathsf{v}}

\begin{document}
\title {Long time derivation of the Boltzmann equation from dynamics of short-range potentials}

\author{Yu Deng}
\address[Yu Deng]{{\newline Department of Mathematics, University of Chicago, Chicago, IL, USA}}
\email{yudeng@uchicago.edu}

\author[Ning Jiang]{Ning Jiang}
\address[Ning Jiang]{\newline School of Mathematics and Statistics, Wuhan University, Wuhan, 430072, P. R. China}
\email{njiang@whu.edu.cn}

\author[Yulong Wu]{Yulong Wu}
\address[Yulong Wu]{\newline Department of Mathematics, University of Michigan, Ann Arbor, MI, USA}
\email{yulongw@umich.edu}

\begin{abstract}
We derive the Boltzmann equation from the Newtonian dynamics of particles interacting through a short-range repulsive potential for arbitrarily long times, in the Boltzmann--Grad limit. This extends the recent hard-sphere result of Deng, Hani and Ma \cite{DHM24} to short-range potentials. Compared with hard spheres, two new difficulties arise: (i) in the absence of a Burago--Ferleger--Kononenko-type \cite{BFK98} bound, a bounded number of particles may undergo an unbounded number of collisions within a short time; and (ii) since collisions last a positive time, three or more particles may interact simultaneously. To overcome (i), we treat such clusters of recollisions as single objects, called boxes, keeping track only of their entry and exit structure, so that internal collisions do not enter the combinatorial counting. The multi-particle collisions in (ii), together with the boxes, lead to more complicated molecules, which requires adapting the counting estimates and the cutting algorithm.

\end{abstract}
\maketitle
\section{Introduction} 
 Hilbert's sixth problem, proposed by Hilbert at the International Congress of Mathematics in 1900, aimed to rigorously derive fluid equations, such as Euler and Navier-Stokes equations, from the atomistic level. Boltzmann's kinetic theory, which describes the gas at the mesoscopic level, is an intermediate step of this process. The problem thus splits into two limits: (i) the kinetic limit, in which one passes from the Newtonian dynamics of $N$ particles to the Boltzmann equation for the one-particle density, as $N\to\infty$ and the particle size $\varepsilon\to0$ under the Boltzmann--Grad scaling $N\varepsilon^{d-1}\sim1$; and (ii) the hydrodynamic limit, in which one passes from the Boltzmann equation to the equations of fluid mechanics as the Knudsen number tends to zero. 
 
The first rigorous result on the kinetic limit was established by Lanford \cite{Lan75} in 1975 for the hard-sphere case, where he derived the Boltzmann equation from hard sphere dynamics over a sufficiently short time. The starting point was the series expansion (BBGKY hierarchy) describing the time evolution of the statistical states of a hard-sphere system. Later, King \cite{King75} extended Lanford's result to the short-range case, again for short times. Gallagher, Saint-Raymond and Texier \cite{GST14} revisited and completed these arguments for both hard spheres and short-range potentials, under suitable assumptions on the potential, and Pulvirenti, Saffirio and Simonella \cite{PSS14} treated a more general class of potentials and obtained explicit convergence rates. Beyond short times, the available results concern perturbative regimes: Illner and Pulvirenti \cite{IP86,IP89} proved global validity for a rare cloud of gas expanding in vacuum, and a series of works by Bodineau, Gallagher, Saint-Raymond and Simonella studied fluctuations around global equilibrium for long times \cite{BGS16,BGS17,BGS18,BGSS20,BGSS20_2,BGSS22,BGSS22_2,BGSS23}. We also mention OpenAI's work \cite{OpenAI26}, see Remark \ref{rem.oai}.

Recently, Deng, Hani, and Ma \cite{DHM24} and \cite{DHM25} completed the Hilbert's sixth problem for the hard sphere case. They rigorously derived the Boltzmann equation from hard sphere dynamics over long time, which were the first results on arbitrarily long time derivation. Their approach is based on a cumulant expansion that keeps track of the whole collision history. The physical trajectories are reduced to diagrams called molecules, and a series of cutting algorithms decompose each molecule into elementary molecules that can be estimated individually. The idea of using such a cumulant expansion for long-time derivations was first introduced in Deng and Hani \cite{DH23_2} in the context of wave turbulence.

In this paper, we aim to rigorously derive the Boltzmann equation from short-range potentials under the Boltzmann-Grad scaling for a rarefied gas over arbitrarily long times, building on the framework of \cite{DHM24,DHM25}. We use the method in \cite{DHM24} to obtain a similar result, under suitable assumptions on the potential function. One of the novelties of this paper is to deal with the pathological collisions which contain only a small number of particles but with lots of collisions in short time. This scenario does not appear in hard sphere case by the main theorem in \cite{BFK98}. The theorem establishes that, for hard sphere model, the number of collisions has an upper bound depending only on the number of particles. To overcome such difficulty, we collect such pathological collision as a whole, and view such interaction as a single interaction, see Section \ref{sec.truncated_dynamic}. The key idea is that: if we can gain some smallness from part of the trajectory, there is no need to further expand the whole collision history. Besides that, we also deal with collisions involving more than two particles, this leads to a more complicated particle trajectory. The corresponding molecule and cutting algorithm should be modified accordingly. For more details on the novelties, see Section \ref{sec.idea_of_proof}.
\subsection{Setup and the main result}
\subsubsection{The  dynamics of short-range potentials}
Consider a system of $N$ identical  particles interacting via a short-range pair potential $\Phi$ with interaction range $\varepsilon$.
\begin{definition}[The short-range Dynamics]\label{def.short_range} We define the short-range dynamics for the $N$-particle short-range system in dimension $d\geq 2$. This will be referred to as the \textbf{Original} or \textbf{O-dynamics} when necessary, to distinguish it from the various modified dynamics in Definitions \ref{def.modified_dynamics}, \ref{def.T_dynamics} below.
\begin{enumerate}
\item\label{it.O1}\emph{State vectors $z_j$ and $\vz_N$}. We define $z_j=(x_j,v_j)$ and $\vz_N=(z_j)_{1\leq j\leq N}$, where $x_j$ and $v_j$ are the center of mass and velocity of particle $j$. These $z_j$ and $\vz_N$ are called the \textbf{state vector} of the $j$-th particle and the collection of all $N$ particles, respectively.
\item\label{it.O2}\emph{The domain $\Dc_N$}. We define the domain $\Dc_N$ as 
  \begin{equation}
    \Dc_N : = \{\vz_N = (z_1, \cdots , z_N) \in \Rb^{2dN}: x_i \neq x_j \quad (\forall i \neq j)\}.
  \end{equation}
\item\label{it.O3} \emph{{The dynamics of short-range potential $\vz_N(t)$}}. Given initial configuration $\vz_N^0=(z_1^0,\cdots,z_N^0)\in\Dc_N$, we define the \textbf{dynamics of short-range potential} $\vz_N(t)=(z_j(t))_{1\leq j\leq N}$ as the time evolution of the following system:
\begin{enumerate}
\item\label{it.O31} We have $\vz_N(0)=\vz_N^0$.
\item\label{it.O32} The Hamiltonian system of motion are
\begin{equation}\label{eq.shortrange}
\frac{\d x_i}{\d t} = v_i, \quad \frac{\d v_i}{\d t} = - \frac{1}{\eps} \sum_{j \neq i} \nabla \Phi \bigg(\frac{x_i - x_j}{\eps}\bigg), 
\end{equation}
where $\Phi$ is a repulsive function compactly supported in the unit ball of $\Rb^d$ and will be specified in subsection \ref{sec.main}.
\item\label{it.energy_Hamiltonian} The energy of the Hamiltonian system is 
\begin{equation}\label{eq.energy_Hamiltonian}
    E_N := \frac{1}{2} \sum^N_{i=1} |v_i(t)|^2 + \sum_{1\leq i< j\leq N} \Phi(\frac{x_i(t) - x_j(t)}{\eps}).
\end{equation}
\end{enumerate}
\item \label{it.O4}\emph{The flow map $\Hc_N(t)$ and flow operator $\Sc_N(t)$}. We define the \textbf{flow map} $\Hc_N(t)$ by
\begin{equation}
    \label{eq.defH^O}\Hc_N(t):\Dc_N \to \Dc_N, \quad\Hc_N(t)(\vz_N^0)=\vz_N(t)
\end{equation}
 where $\vz_N(t)$ is defined by (\ref{it.O3}). Define also the \textbf{flow operator} $\Sc_N(t)$, for functions $f=f(\vz_N)$, by
\begin{equation}\label{eq.defS^O}
(\Sc_N(t) f)(\vz_N)=f\big((\Hc_N(t))^{-1}\vz_N\big).
\end{equation}
\item\label{it.O5} \textbf{Collision.}
For each interaction between two particles, indexed by
$\alpha$, let
$$ J_\alpha:=\{i_\alpha,j_\alpha\}$$
be the corresponding particle pair. We define its \newterm{pairwise interaction interval} to be $I_\alpha=[t_\alpha^-,t_\alpha^+], \quad\! t_\alpha^-<t_\alpha^+,$ such that
\begin{equation*}
    \begin{aligned}
        &|x_{i_\alpha}(t_\alpha^-)-x_{j_\alpha}(t_\alpha^-)|=\varepsilon,
        \qquad
        (x_{i_\alpha}-x_{j_\alpha})
        \cdot(v_{i_\alpha}-v_{j_\alpha})
        \big|_{t=t_\alpha^-}\leq 0,\\
        &|x_{i_\alpha}(t_\alpha^+)-x_{j_\alpha}(t_\alpha^+)|=\varepsilon,
        \qquad
        (x_{i_\alpha}-x_{j_\alpha})
        \cdot(v_{i_\alpha}-v_{j_\alpha})
        \big|_{t=t_\alpha^+}\geq 0,
    \end{aligned}
\end{equation*}
and
$$|x_{i_\alpha}(t)-x_{j_\alpha}(t)|<\varepsilon, \qquad  \forall t\in(t_\alpha^-,t_\alpha^+).$$
Two pairwise interactions $\alpha$ and $\beta$ are said to be \textbf{connected} if
$$J_\alpha\cap J_\beta\neq\varnothing\qquad\text{and}\qquad I_\alpha\cap I_\beta\neq\varnothing.$$
Let $\sim$ be the smallest equivalence relation on the set of pairwise interaction satisfying the following two properties:
\begin{enumerate}
    \item If $\alpha$ and $\beta$ are connected, then $ \alpha\sim  \beta$; if $\alpha \sim \beta$ and $\beta \sim \gamma$, then $\alpha \sim \gamma$.
    \item Suppose that $ \beta\sim \gamma$ and there exists some interaction $\alpha$ satisfying
    $$J_\alpha\cap J_\beta \neq \emptyset,\quad \text{and} \quad\! J_\alpha \cap J_\gamma\neq\varnothing,$$
    and $ I_\alpha\subset [t_\beta^+,t_\gamma^-] \quad \text{or} \quad\![t_\gamma^+,t_\beta^-],$ then $ \alpha\sim  \beta$.\footnote{This condition is needed to account for configurations such as the following. Particles 1 and 2 initially interact. Before their interaction ends, particle 3 begins interacting with particle 2. Then particle 1 leaves, and particles 3 and 2 keep interacting. After that, particle 1 interacts with particle 4, and subsequently particle 1 goes back to interact with particle 2 or 3, before the interaction between particles 2 and 3 ends. In this case, we view the whole interaction as a whole collision.}
\end{enumerate}

A \newterm{collision} is a collection of all equivalent interactions under $\sim$. Its particle set consists of all particles involved in these interactions, and its collision interval is the union of every pairwise interaction interval. If a collision has two particles, then we call it \newterm{binary collision}; if it has three or more particles, we call it \newterm{multi-particle collision}.
\end{enumerate} 
\end{definition}
We list some basic properties of the short-range dynamics, which hold up to Lebesgue zero sets:
\begin{proposition}\label{prop.shortrange}  
Up to a set $Z$ of Lebesgue measure zero, for any $\vz_N^0 \in \Dc_N \setminus Z$, at each time there exists at most one pair of particles satisfying $|x_i - x_j|=\eps$ and $(x_i -x_j)\cdot(v_i-v_j) \leq 0$. The flow maps $\Hc_N(t)$ defined in Definition \ref{def.short_range} (\ref{it.O4}) are measure preserving homeomorphisms from $\Dc_N \setminus Z$ to itself, and satisfy the flow or semi-group property $\Hc_N(t+s)=\Hc_N(t)\Hc_N(s)$ for $t,s\geq 0$.  
\end{proposition}
\begin{proof} The first statement follows the same argument as Proposition 4.1.1 in \cite{GST14}. Volume preservation follows as an application of Liouville's theorem up to a zero measure set.
\end{proof}
 
\subsubsection{The grand canonical ensemble}\label{sec.grand_canon} In this paper, we study a particle system with short-range interactions and a random number of particles which have random initial states. We define a density function $W_0:\cup_N \Dc_N\to\Rb_{\geq 0}$ (or equivalently a sequence of density functions $W_{0,N}:\Dc_N\to\Rb_{\geq 0}$) that determines the law of the random initial data for the particle system with a short-range potential. We state the definition as follows:
\begin{definition}[The grand canonical ensemble] \label{def.grand_canon} Fix $0<\varepsilon\ll 1$ and a nonnegative function $f_0=f_0(z)$ with $\int f_0(z)\,\mathrm{d}z=1$. We define the \textbf{grand canonical ensemble} as follows.
\begin{enumerate}
\item\label{it.Can0} \emph{The grand canonical domain $\Dc$}. We define $\Dc:=\bigcup_{N}\Dc_N^0$ to be the \textbf{grand canonical domain}, where $\Dc_N^0 := \{\vz_N = (z_1,\cdots, z_N) \in \Rb^{2dN} : |x_i -x_j|\geq \eps \quad (\forall i\neq j)\} $. We also define the short-range dynamics $\vz(t)$ on $\Dc$, and the flow map $\Hc(t)$ on $\Dc$, by specifying them to be the $\vz_N(t)$ in Definition \ref{def.short_range} (\ref{it.O3}), and the $\Hc_N(t)$ in Definition \ref{def.short_range} (\ref{it.O4}), for $\vz^0=\vz_N^0\in\Dc_N^0$. Note that $\Hc(t)$ does not preserve $\Dc$ (particles may come within distance $\eps$), and it maps $\Dc$ into $\bigcup_N\Dc_N$.
\item \label{it.Can1} \emph{Random data and initial density}. Assume $\vz^0\in \Dc$ is a random variable, whose law is given by the \textbf{initial density functions} $(W_{0,N})$, in the sense that \begin{equation}\label{eq.inti_density}
\Pb\big(\vz^0=\vz_N^0\in A\subseteq\Dc_N\big)=\frac{1}{N!}\int_{A} W_{0,N}(\vz_N)\,\mathrm{d}\vz_N
\end{equation}
for any $N$ and any $A\subseteq\Dc_N$, where $(W_{0,N})$ is given by
\begin{equation}\label{eq.N_par_ensemble} 
    \frac{1}{N!}W_{0,N}(\vz_N):=\frac{1}{\Zc}\frac{\varepsilon^{-(d-1)N}}{N!}\prod_{j=1}^N f_0(z_j)\cdot \mathbbm{1}_{\Dc^0_N}(\vz_N).
\end{equation}
Here, $\mathbbm{1}_{\Dc_N^0}$ is the indicator function of the set $\Dc_N^0$, and the partition function $\Zc$ is defined as follows
\begin{equation}\label{eq.partition_func}
    \Zc:=1+\sum_{N=1}^\infty\frac{\varepsilon^{-(d-1)N}}{N!}\int_{\Rb^{2dN}}\prod_{j=1}^N f_0(z_j)\cdot \mathbbm{1}_{\Dc_N^0}(\vz_N)\,\mathrm{d}\vz_N.
\end{equation}
\item \label{it.Can2}\emph{Evolution of random data}. Let $\vz^0\in \Dc$ be the random variable defined in (\ref{it.Can1}) above, and let $\vz(t)=\Hc(t)\vz^0$ be the evolution of initial data $\vz^0$ by short-range dynamics, defined as in Definition \ref{def.short_range} (\ref{it.O3})--(\ref{it.O4}). We define the \textbf{density functions} $W_N(t,\vz_N)$ for the law of the random variable $\vz(t)$ by the relation
\begin{equation}\label{eq.inti_density_t}
\Pb\big(\vz(t)=\vz_N(t)\in A\subseteq\Dc_N\big)=\frac{1}{N!}\int_{A} W_{N}(t,\vz_N)\,\mathrm{d}\vz_N.
\end{equation} Then it is easy to see (using the volume preserving property in Proposition \ref{prop.shortrange}) that
\begin{equation}\label{eq.time_t_ensemble}W_N(t,\vz_N):=(\Sc_N(t)W_{0,N})(\vz_N).
\end{equation}
\item\label{it.Can3} \emph{The correlation function $f_s(t)$}. Given $s\in\Nb$, we define the $s$-particle (rescaled) \textbf{correlation function} $f_s(t,\vz_s)$ by the following formula
\begin{equation}\label{eq.s_par_cor}
f_s(t,\vz_s):=\varepsilon^{(d-1)s}\sum_{n=0}^\infty\frac{1}{n!}\int_{\Rb^{2dn}}W_{s+n}(t,\vz_{s+n})\,\mathrm{d}z_{s+1}\cdots \mathrm{d}z_{s+n},
\end{equation} 
where we abbreviate $\vz_s=(z_1,\cdots,z_s)$ and $\vz_{s+n}=(z_1,\cdots,z_{s+n})$ similar to $\vz_N$.
\item\label{it.Can4} \emph{The Boltzmann-Grad scaling law}. The formula (\ref{eq.N_par_ensemble}), particularly the choice of the power $\varepsilon^{-(d-1)N}$, implies that
\begin{equation}\label{eq.Bol-Grad}\Eb(N)\cdot \varepsilon^{d-1}\approx 1
\end{equation} up to error $O(\varepsilon)$, where $\vz^0\in\Dc$ is the random initial data defined in (\ref{it.Can1}) and $N$ is the random variable (i.e. number of particles) determined by $\vz^0=\vz_N^0\in\Dc_N^0$. In other words, we are considering the short-range potential dynamics with $\approx \varepsilon^{-(d-1)}$ particles (on average). This is referred to as the \textbf{Boltzmann-Grad scaling law}, and will be fixed throughout the paper.
\end{enumerate}
\end{definition}

\subsubsection{The Boltzmann equation} The kinetic theory predicts that the one-particle correlation function $f_1(t,z)=f_1(t,x,v)$ should solve the Boltzmann equation in the kinetic limit $\varepsilon\to 0$, which is defined as follows.
\begin{definition}[The Boltzmann equation]\label{def.boltzmann} We define the Cauchy problem for the \textbf{Boltzmann equation}, with initial data $f(0,z)=f_0(z)$ where $z=(x,v)\in\Rb^{2d}$, as follows:
\begin{equation}\label{eq.boltzmann}
\left\{
\begin{aligned}
(\partial_t+v\cdot\nabla_x)f&=\int_{\Rb^d}\int_{\Sb^{d-1}}\big((v-v_1) \cdot \nu \big)_+\cdot(f'f_1'-ff_1)\,\mathrm{d}\nu \mathrm{d}v_1;\\
f(0,x,v)&=f_0(x,v).
\end{aligned}
\right.
\end{equation} The right-hand side of (\ref{eq.boltzmann}) is referred to as the \textbf{collision operator}, where we denote
$$f=f(x,v),\quad f_1=f(x,v_1),\quad f'=f(x,v'),\quad f_1'=f(x,v_1')$$ and
\begin{equation}\label{eq.boltzmann2} 
    v'=v-\big((v-v_1)\cdot\omega\big)\omega,\quad v_1'=v_1+\big((v-v_1)\cdot\omega\big)\omega,
\end{equation} 
and $\big((v-v_1)\cdot\nu\big)_+$ is the positive part of $(v-v_1)\cdot\nu$. Here $\omega = \omega_{\Phi}(v-v_1,\nu)$ denotes a unit vector along the apsidal line associated with the two-body scattering generated by $\Phi$.
\end{definition}

The Boltzmann equation can be also rewritten in terms of $\d \omega$, by a change of variable $\nu \leftrightarrow \omega$:
\begin{equation}\label{eq.boltzmann_omega}
\left\{
\begin{aligned}
(\partial_t+v\cdot\nabla_x)f&=\int_{\Rb^d}\int_{\Sb^{d-1}}b(|v- v_1|,\omega)\cdot(f'f_1'-ff_1)\,\mathrm{d}\omega \mathrm{d}v_1;\\
f(0,x,v)&=f_0(x,v),
\end{aligned}
\right.
\end{equation}
where $b(|v- v_1|,\omega)$ is the scattering cross-section. The detailed derivation of $b(|v- v_1|,\omega)$ can be found in \cite{GST14}, and is summarized in Appendix \ref{app.aux}.
\subsubsection{The main result}\label{sec.main}
Throughout this paper, we assume the following properties for the potential $\Phi$:
\begin{assumption}\label{ass.phi1}
The potential $\Phi : \Rb^d \to \Rb$ is a radial, nonnegative, nonincreasing function supported in the unit ball of $\Rb^d$, of class $C^2$ in $\{r \in \Rb^d, 0 <|r| <1\}. $ Moreover, it is assumed that $\Phi$ is unbounded near zero, goes to zero at $|r| = 1$ with bounded derivatives, and that $|\nabla \Phi(r)| \neq 0$ on $|r| < 1$.   
\end{assumption}
Due to technical reasons, we impose the following assumptions on the potential. Assumption \ref{ass.phi2} is the same technical assumption as in \cite{GST14}.
\begin{assumption}\label{ass.phi2}
The potential $\Phi $ satisfies the following condition: \footnote{This is the monotonic assumption in \cite{GST14}.}
\begin{equation}
      r\Phi''(r) + 2 \Phi'(r) \geq 0 
\end{equation}
for all $r \in(0,1)$.
\end{assumption}
\begin{assumption}\label{ass.phi3} We further require that either $\Phi'(1^-) \leq -c_0 <0$ for some positive constant $c_0$, or there exists $0 <\eta_0 \ll 1$, such that for all $0< \eta \leq \eta_0$
    \begin{equation*}
        |\Phi'(1-\eta)| \gtrsim  \eta^{C_\Phi^0}
    \end{equation*} 
for some constant $C_\Phi^0 \geq d$.
\end{assumption}

We now state the main result of this paper. First we introduce the following norm for functions $f: \Rb^d \times \Rb^d \to \Cb$,
\begin{equation}\label{eq.boltzmann_decay_2}
\|f(x,v)\|_{\mathrm{Bol}^\beta}:=\sum_{k\in\Zb^d}\sup_{|x-k|\leq \sqrt{d},v\in\Rb^d}e^{\beta|v|^2}|f(x,v)|.
\end{equation} 
\begin{theorem}\label{th.main} Fix $d\geq 2, \beta>0$, and a nonnegative function $f_0=f_0(z)\geq 0$, with $\int f_0(z)\,\mathrm{d}z=1$. Suppose the solution $f(t,z)$ to the Boltzmann equation (\ref{eq.boltzmann}) exists on the time interval $[0,t_{\mathrm{fin}}]$, such that
\begin{equation}\label{eq.boltzmann_decay_1}
\big\|e^{2\beta  |v|^2}f(t)\big\|_{L_{x,v}^\infty}\leq A<\infty,\quad \forall t\in[0,t_{\mathrm{fin}}].
\end{equation}
Note that $t_{\mathrm{fin}}$ can be arbitrarily large. Suppose also that the initial data $f_0$ satisfies
\begin{equation}\label{eq.boltzmann_decay_3}\|f_0\|_{\mathrm{Bol}^{2\beta}}+\|\nabla_x f_0\|_{\mathrm{Bol}^{2\beta}}\leq B_0<\infty.
\end{equation}

Consider the $d$ dimensional short-range potential system (Definition \ref{def.short_range}), with random initial data $\vz^0$ given by the grand canonical ensemble (Definition \ref{def.grand_canon}), under the Boltzmann-Grad scaling law \eqref{eq.Bol-Grad}. Assume that the potential $\Phi$ satisfies Assumptions \ref{ass.phi1}--\ref{ass.phi3}. Let $\varepsilon$ be small enough depending on $(d,t_{\mathrm{fin}},\beta,A, B_0,\Phi)$. Then, uniformly in $t\in[0,t_{\mathrm{fin}}]$ and in $s\leq |\log\varepsilon|$, the $s$-particle correlation functions $f_s(t)$ defined as in \eqref{eq.N_par_ensemble}--\eqref{eq.s_par_cor}, satisfy
\begin{equation}\label{eq.maineqn}
\bigg\|f_s(t,\vz_s)-\prod_{j=1}^s f(t,z_j) \bigg\|_{L^1(\Rb^{2ds})}\leq \varepsilon^\theta,
\end{equation} 
where $\theta>0$ is a constant depending only on the dimension $d$ and $C_\Phi^0$. 
\end{theorem}
    
\begin{remark}\label{rem.theorem} This result is an analog of the main theorem in \cite{DHM24}. As in \cite{DHM24}, we allow $t_{\rm fin} \sim \log |\log \eps|$. 
\end{remark}
\begin{remark}[OpenAI's work \cite{OpenAI26}]\label{rem.oai}
Shortly before this manuscript was completed, the preprint \cite{OpenAI26} was released on October 6, 2026. It proves convergence of the marginals of any fixed number of particles (i.e. $s$ in Theorem \ref{th.main} does not depend on $\eps$), without a convergence rate, in dimension $d=3$, for stable radial potentials, a class that allows attractive parts. Both works follow the layered long-time strategy of Deng, Hani and Ma \cite{DHM24}, but the proofs are different: \cite{OpenAI26} uses an exact expansion over whole interaction components on a fine time mesh, together with energy and exposure estimates, and does not use molecules or a cutting algorithm; we work within the molecule framework of \cite{DHM24}, introduce boxes to handle unbounded numbers of recollisions, and adapt the cutting algorithm, which yields an explicit rate $\eps^\theta$, uniformly in $s\le|\log\eps|$, in every dimension $d\ge 2$. A first draft of this paper, containing the main results and the overall strategy, was almost complete before \cite{OpenAI26} appeared, and the later revisions did not make use of \cite{OpenAI26}. Our work is independent of \cite{OpenAI26}.
\end{remark} 
\subsection{Ideas of the proof.}\label{sec.idea_of_proof} The proof of Theorem \ref{th.main} follows the same overall strategy as the hard-sphere proof in \cite{DHM24}: we construct layered molecules (which correspond to the topological reduction of physical trajectories), derive the cumulant expansion, derive the combinatorial estimate on the number of molecules. We define elementary molecules and, classify them into \emph{good, normal, or bad} molecules depending on their estimates. Then we apply cutting algorithms to cut the whole molecule into elementary molecules, which produces sufficiently many good molecules while keeping the contribution of the bad molecules under control. In the present short-range potential setting, however, several new ingredients are needed. We summarize them below.
\begin{enumerate}
    \item {\bf Many collisions in time period $[\ell \tau, (\ell + 1) \tau]$.} For hard spheres, the Burago--Ferleger--Kononenko (BFK) theorem \cite{BFK98} provides a uniform bound for the number of collisions which depends only on the number of particles. Such a bound is not available for short-range potential: a fixed number of particles may have many collisions in some time period $[\ell \tau, (\ell + 1) \tau]$. If we apply the cutting algorithm directly to such a molecule, each collision cost a combinatorial factor, while the cutting produces only a few good molecules, which cannot absorb this loss. To overcome this, we treat such a pathological set, formed by $q\leq\mathsf{N}$ particles producing at least $G(q)$ collisions, as one single object, called a \emph{box} (Definitions \ref{def.pathological}--\ref{def.box}). Each box consists of a bottom part, a middle part and a top part. We keep the detailed structure only of the bottom and top parts, while the middle part, called an \emph{M-atom}, enters only through an indicator function. The key point is that the number of collisions inside an M-atom no longer enters the counting: after viewing each box as a single atom, at most $q$ particles carry at most $G(q)$ atoms. Moreover, the bottom and top parts have a bounded number of atoms and each contains a cycle, so a box cut with at most one fixed end gives a good molecule. In the algorithms below, we modify those in \cite{DHM24} to deal with molecules containing boxes.
    
    The intuitive explaintion of the definition of box is as follows: for a molecule if we can gain from part of the molecule, we do not need to further expand it, which leads to combinatorial loss without gaining.

    \item \textbf{Multi-particle collision.} In contrast to hard-sphere collisions, an interaction for a short-range potential lasts for a nonzero time interval, and during this interval additional particles may enter the interaction region. Hence an atom may contain more than two particles; we call such a collision a \emph{multi-particle collision}, and the corresponding atom a \emph{multi-particle atom}. To deal with this, the molecules, the cutting procedure and the counting argument allow multi-particle atoms (Definition \ref{def.molecule}), and the truncated dynamics rejects collisions with more than $d+1$ particles or with long interaction (Definition \ref{def.T_dynamics}). The key point is that a multi-particle atom is not only a new structure but also a source of gain, see \eqref{it.idea_6} below.

    \item \textbf{The Penrose argument with atom--atom or edge--atom connections.} In the hard-sphere case \cite{DHM24}, the non-overlap condition between two molecules can be reduced to the non-overlap of their edges:
    \begin{equation}
        \mathbbm{1}_{\Mb_j \not\sim \Mb_{j'}} = \prod_{e\in \Mb_j,\, e'\in \Mb_{j'}} \mathbbm{1}_{(\Mb_j,e) \not\sim(\Mb_{j'},e')} = \prod_{e\in \Mb_j,\, e'\in \Mb_{j'}} (1 - \mathbbm{1}_{(\Mb_j,e) \sim(\Mb_{j'},e')}) 
    \end{equation}
    where $e$ and $e'$ are edges of $\Mb_j$ and $\Mb_{j'}$ respectively. And the $\sim$ on the right-hand side above generates \emph{overlap} between two different edges of two different molecules by the Penrose argument. In the present setting, an overlap may involve either a bond or an atom (because the collision is no longer instantaneous). For short-range case, we have 
        \begin{equation}
        \mathbbm{1}_{\Mb_j \not\sim \Mb_{j'}} = \prod_{\af \in \Mb_j,\, \af'\in \Mb_{j'}} \mathbbm{1}_{(\Mb_j,\af) \not\sim(\Mb_{j'},\af')} = \prod_{\af \in \Mb_j,\, \af'\in \Mb_{j'}} (1 - \mathbbm{1}_{(\Mb_j,\af) \sim(\Mb_{j'},\af')}),
    \end{equation}
    where $\af$ and $\af'$ are {\bf edges or atoms} of $\Mb_j$ and $\Mb_j'$ respectively. Consequently, besides the usual bond--bond O-connections (O-atoms), the Penrose construction must also allow bond--atom and atom--atom O-connections, which we call CO-atoms or MO-atoms (Definitions \ref{def.molecule} and \ref{def.boxed_molecule}). A CO-atom has edge constituents and C-constituents (which are originally C-atoms); if it has an M-constituent (the middle part of a box), it is called an MO-atom. Note that a CO- or MO-atom may contain many constituents, hence many particles. The key point is that the Penrose argument of \cite{DHM24} still applies (Lemma \ref{lem.penrose_equality}): the only changes are that the overlaps may involve atom constituents, and that the O-connections at the same atom are ordered by its decoration (Definition \ref{def.molecule}).
     
    \item {\bf Partial cutting}: The O-connections create an additional difficulty for boxes. Each M-atom may be connected through O-connections to many edges. If one cuts such an atom naively, then every edge with O-connection attached to the M-atom would be cut. In the main algorithm of \cite{DHM24} (Section 13 there), the original molecule is reduced to a two-layer molecule $\Mb_{UD} = \Mb_U \cup \Mb_D$, and the gain comes from the bonds between $\Mb_U$ and $\Mb_D$ (UD connections); cutting an MO-atom with may UD-connections would break them without any gain. To overcome this, we introduce  
    \emph{partial cutting}: for a CO- or MO-atom, we only cut some of its constituents, and insert an auxiliary simple pair at each newly created O-connection. The key point is that this avoids breaking many bonds in a single cut, while the variables before and after the operation can still be matched as in Definition 8.12 of \cite{DHM24} (Proposition \ref{prop.cutting}).
    
    \item\label{it.idea_6} \textbf{Special atoms.} After partial cutting, new atoms appear, which are formed by a constituent of a CO- or MO-atom together with an auxiliary simple pair; we call them \emph{special atoms}. Their integrals do not satisfy the usual estimates, so the elementary molecules containing them need new arguments. The point is that for special $\{3_2\}$ and special $\{3_23_2\}$ molecules, we can only gain from the time variables (and not from the velocity or position variables), see Propositions \ref{prop.est_single} and \ref{prop.est_twoatom}. Hence, we introduce a new notion of \emph{short degeneracy} (Definition \ref{def.short_degeneracy}): a pair of ov-adjacent atoms in the same layer, at least one of which is a CO- or MO-atom, which are connected by an O-bond and whose times differ by at most $\eps^{\upsilon}$. If there are many disjoint short degenerate pairs, we gain from them directly (Proposition \ref{prop.short_degeneracy}); otherwise, we use specific arguments to obtain the gain from the special molecules.
    
    \item \textbf{Counting estimate and algorithm for multi-particle atoms and boxes.} Let $\Xf$ and $\Bf$ be the numbers of multi-particle atoms and of boxes (see Definition \ref{def.number_box_multi} for the precise definition of $\Xf$), and let $\Nf$ be the total number of edges at these atoms (each edge counted at most twice). Compared with \cite{DHM24}, the number of molecules has an extra factor depending on $\Xf$ and $\Bf$ (Proposition \ref{prop.layerrec3}):
    \begin{equation*}
        \# \Mb \leq C^{|\Mb| + \Nf} \cdot |\log \eps |^{C^*{(\rho + \Xf + \Bf)}}
    \end{equation*}
    where $\rho$ is the total number of cycles in $\Mb$ (not counting cycles inside the box). To absorb the new factor $|\log\eps|^{C^*(\Xf+\Bf)}$, we use that each multi-particle atom gives an extra factor $\eps^{\upsilon_0/2}$ (Proposition \ref{prop.est_single} \eqref{it.est_kl}), and that each box gives a good molecule. The key point for multi-particle atoms is the following dichotomy. If the collision duration is short, we gain from the time variables: for example, if particles $1,2,3$ collide, the time when the third particle enters the collision must be close to the time when the collision starts. If the collision duration is long, we gain from the velocities: by our assumptions on the potential, for a long collision, the velocities of the particles must be close (Proposition \ref{prop.upperbond.interaction}).

    \item\textbf{Modifications on algorithms in \cite{DHM24}.} The algorithms of \cite{DHM24} must be modified to deal with boxes, multi-particle atoms, CO- and MO-atoms and special atoms. First, we modify the algorithm \textbf{UP} of \cite{DHM24} (Definition \ref{def.alg_up}), and introduce its variants \textbf{UP1} and \textbf{UP2} (Definition \ref{def.alg_up1}), which view each box as a whole, cut one atom at a time, and apply the partial cutting \textbf{PM1} or \textbf{PM} to CO- and BO-atoms; their key properties are in Proposition \ref{prop.alg_up}. In particular, a multi-particle atom or a box cut with at most one fixed end gives a good molecule. Then, if the number of boxes or multi-particle atoms is large (say $\geq C\rho$), we apply \textbf{UP2} and gain from them directly; otherwise, we deal with the degenerate cases, the case of many cyclic clusters, and so on (Sections \ref{sec.strdeg}--\ref{sec.layer_select} and \ref{sec.maincr}). These algorithms follow \cite{DHM24}, but require many modifications, as the molecule structure is more complicated.
     
    \item \textbf{Estimate for the truncated dynamics.} Following \cite{DHM24}, we truncate the dynamics if a cluster has too many particles or too many recollisions. In addition, we truncate it if a collision involves more than $d+1$ particles, if a collision lasts too long (long interaction), or if a cluster contains too many boxes. The error term $f^{\rm err}_s$ is then controlled by the molecules $\Mb$ at which the truncation first occurs (Section \ref{sec.error}). As in Section 15 of \cite{DHM24}, we cut such an $\Mb$ into elementary molecules $\Mb_j$, so that the associated integral is bounded by the product of those of the $\Mb_j$. All of them are bounded by $|\log\eps|^{C^*}$, except one bad molecule which costs $\eps^{-(d-1)}$, and the goal is to find enough good molecules, each provides $\eps^{\sigma_j}$, such that $\sum_j \sigma_j  - (d-1) > 0 $; see \eqref{eq.extra_case_1}. A long interaction or an collision with at least $d+2$ particles is a single atom which already provides $\eps^{d-1+\upsilon}$ (Proposition \ref{prop.est.d+2}), and many boxes give many good molecules (Proposition \ref{prop.comb_est_extra2}). The main difficulty is the case of too many recollisions, $\Gamma<\rho(\Mb)\leq2\Gamma$ (here $\rho(\Mb)$ does not count recollisions inside boxes), for which we need the following new ingredients.
    
\emph{Boxes instead of the BFK bound.} In \cite{DHM24}, this case uses the BFK theorem to bound the number of collisions among a bounded number of particles. On the torus \cite{DHM25}, the BFK bound fails for long times but still holds for short times (in each period the theorem applies). Hence a pathological collision set there contains two collisions which are separated in time, which gives a long bond, and the gain comes from the estimate $\eps^{d-1}/|t_1-t_2|^d$ for a $\{33A\}$ molecule with collision times $t_1$ and $t_2$. Neither ingredient is available for short-range potentials: the BFK bound fails even for short times, see (1) above, and the above estimate for $\{33A\}$ molecules fails due to the singularity of the cross-section $b(|v_1-v_2|,\omega)$. To overcome this, we use the boxes as a substitute (Definitions \ref{def.box} \ref{def.boxed_molecule}): after viewing each box as a single atom, a connected set of atoms with at most $q$ particles has at most $G(q)$ atoms (Definitions \ref{def.pathological}--\ref{def.box}). 

\emph{Transversal sets, and the two-layer structure} As in Section 15 of \cite{DHM24}, we construct a sequence of connected transversal sets $A_j$ (Definition \ref{def.trans} and Proposition \ref{prop.trans}), and consider the connections between $A=A_j$ and $X(A)$, the set of atoms all of whose bonds go into $A$. If the process continues for $D$ steps, or if there are many such connections which is much larger than the number of cycles in $A$ (this is valid for hard-spehre case due to the BFK theorem), we argue as in \cite{DHM24}, using the algorithms \textbf{TRANSUP}, and \textbf{MAINTRUP} (the analogs algorithms there, dealing with multi-particle atoms and boxes) in the second case. 

In the short-range case, there may be few such connections. Then, by the bound on the number of atoms given by the boxes, some component $\Xc$ of $X(A)$ has at least $\mathsf{N}$ particles, which gives a two-layer structure $\Mb_U\cup\Mb_D$, where $\Mb_U=\Xc$ and $\Mb_D$ is the highest part of $A$, with at least $\mathsf{N}$ bonds between them.  In this case the above argument fails: \textbf{MAINTRUP} cuts all boxes in $A$, and a single box may carry up to $\mathsf{N}$ particle lines connected to $\Xc$, so these UD connections would be broken without any gain; however, as in the main algorithm \textbf{MAINUD} of \cite{DHM24}, we need to gain from them. Therefore, new algorithms are needed for this two-layer structure.

\emph{Canonical cycles.} In this two-layer structure, the algorithms \textbf{TRANSUP} and \textbf{MAINTRUP} gain roughly $\#_{\rm conn}-C\rho(A)$, where $\#_{\rm conn}$ is the number of bonds between $\Mb_U$ and $\Mb_D$; this fails here, since the number of cycles in $A$ may be much larger than $\#_{\rm conn}$. To overcome this, we always consider the highest part of $A$, i.e. the highest connected subset containing all particle lines of $A$; only this part is connected to $\Mb_U$ (Proposition \ref{prop_MUMD}). For this highest part, we introduce \emph{canonical cycles} (Definition \ref{def.canonical_cycle}): starting from a smallest cycle, we repeatedly merge the cycles which meet the current set, and we also regard each remaining box as a canonical cycle. The canonical cycles do not depend on the choices, and the molecule becomes a tree after viewing each canonical cycle as a single node (Proposition \ref{prop.canonical_cycle}). This gives a dichotomy. If $\Mb_D$ contains at least $K+1$ (where $K=(2d+1)/\upsilon$) canonical cycles, we cut this tree by \textbf{UP}, viewing each canonical cycle as a single node; then each node except the first is cut with exactly one fixed end, and each canonical cycle gives a good molecule (Lemma \ref{lem.case2_manycc}). Otherwise, $\Mb_D$ contains at most $K$ canonical cycles, which are treated below.  The key point, which is used at the end of the argument, is that passing to the highest part strictly reduces a lowest canonical cycle: its lowest atoms can be removed without losing a particle line, so it is never kept entirely (Lemma \ref{lem.case2_trim}). In particular, the highest part of a molecule is not a single canonical cycle, since otherwise a higher part could be chosen.

\emph{The induction on the highest parts.} We now have the highest part $\Mb^0_D$ of $A$, with at least $\mathsf{N}$ bonds to $\Mb_U$, and we may assume that it contains at most $K$ canonical cycles, where $K=(2d+1)/\upsilon$. The difficulty is that most of these bonds may belong to a single canonical cycle, and then we cannot gain from them directly. To overcome this, we construct a decreasing sequence $\Mb^0_D\supseteq\Mb^1_D\supseteq\cdots$ in at most $D:=8K^{2K+3}$ rounds, such that $\Mb^j_D$ has at least $N_{D-j}/(2(d+1))^{K+1}$ bonds to $\Mb_U$ (Definition \ref{def.construct_MD} and Proposition \ref{prop_MUMD}). In Round $j$, we stop and obtain the gain in each of the following cases:
\begin{itemize}
\item $\Mb^j_D$ contains at least $K+1$ canonical cycles (Lemma \ref{lem.case2_manycc});
\item $\Mb^j_D$ contains at least $2d^2$ multi-particle atoms outside the canonical cycles, each giving $\eps^{\upsilon_0/2}$;
\item for some earlier round $i$, there are many atoms all of whose children lie in $\Mb^{i-1}_D$, with at least one child in the removed part $\Mb^{i-1}_D\setminus\Mb^{i}_D$ and no such child in a canonical cycle (Lemma \ref{lem.case2_sigma});
\item each canonical cycle has fewer than $N_{D-j-1}$ bonds to $\Mb_U$; then the bonds to $\Mb_U$ are spread over many atoms, and we cut by the new algorithms \textbf{TRANSUP1} and \textbf{MAINTRUP2}.
\end{itemize}
Otherwise, some canonical cycle has at least $N_{D-j-1}$ bonds to $\Mb_U$, and we construct $\Mb^{j+1}_D$ as follows. Viewing each canonical cycle as a single node, $\Mb^j_D$ is a tree. We keep the smallest subtree containing all canonical cycles with bonds to $\Mb_U$, together with their ancestors, and then remove atoms from below: if a lowest node is a canonical cycle, we remove its lowest atoms as long as the set stays connected and keeps all particle lines (this is taking the highest part again); otherwise we remove a lowest atom and keep the component with the most bonds to $\Mb_U$. The key point is that each round loses at most a factor $(2(d+1))^{K+1}$ in the number of bonds to $\Mb_U$, while each threshold is much smaller than the previous one, $N_{D-j}=(60dN_{D-j-1})^{60d}$; so the construction can be repeated in all $D$ rounds.

\emph{Nested canonical cycles.} It remains to treat the case that we construct $\Mb^{j+1}_D$ in every round. A round in which a lowest node is a canonical cycle strictly reduces this canonical cycle (Lemma \ref{lem.case2_trim} \eqref{it.trim1}). The other rounds reduce the number of canonical cycles with bonds to $\Mb_U$, which is at most $K$, so at least $(D-1-K)/(K+1)$ rounds are of the first type (Lemma \ref{lem.case2_trim} \eqref{it.trim3}). Choosing each time the canonical cycle which contains the most canonical cycles reduced in later rounds, we lose a factor $K$ at each step, and since $D=8K^{2K+3}$, we obtain a nested sequence of canonical cycles $\Cc_1\supsetneq\Cc_2\supsetneq\cdots\supsetneq\Cc_{2K+2}$, where each $\Cc_{i+1}$ lies in the reduced part of $\Cc_i$ (Proposition \ref{prop.case2_chain}). Moreover, no atom outside $\Cc_i$ has both an ancestor and a descendant in $\Cc_i$ (Lemma \ref{lem.case2_convex}). To gain from this structure, we cut $\Cc_{2K+2}$ first, and then the layers $\Cc_i\setminus\Cc_{i+1}$ for $i=2K+1,\dots,1$, from the inside out, using the new algorithms \textbf{UP3} and \textbf{DOWN3}, which also allow M-atoms to be cut with fixed ends. The key point is that, since $\Cc_{i+1}$ lies in the reduced part of $\Cc_i$, each layer contains a component with no fully fixed atom which has either a cycle and a fixed end, or at least two fixed ends; this gives a good molecule (Proposition \ref{prop.alg_up} \eqref{it.up_proof_3}--\eqref{it.up_proof_4}). Hence we obtain at least $K$ good molecules and only one bad molecule, which is enough to close the proof.

\end{enumerate}
 
\subsection{Parameters and notations}\label{sec.notation}  By the assumptions (\ref{eq.boltzmann_decay_1}) and (\ref{eq.boltzmann_decay_3}) in Theorem \ref{th.main}, and using Proposition \ref{prop.pres_of_decay}, it follows that
\begin{equation}\label{eq.boltzmann_decay_4}
    \|f(t)\|_{\mathrm{Bol}^\beta}+\|\nabla_xf(t)\|_{\mathrm{Bol}^\beta}\leq B<\infty,\quad \forall t\in[0,t_{\mathrm{fin}}],
\end{equation} where $B$ is a constant depending only on $(d,t_{\mathrm{fin}},\beta,A,B_0)$. In this section we fix some parameters and define some notations, which will be used throughout the proof.
\begin{definition}[Parameters and notations]
\label{def.notation} We define the following parameters and notations.
\begin{enumerate}
\item\label{it.time_layer} \emph{Time layering}. Recall that $(d,t_{\mathrm{fin}},\beta,A,B_0,C_{\Phi}^0)$ has been fixed as in Theorem \ref{th.main}, and $B$ has been fixed as in (\ref{eq.boltzmann_decay_4}) depending only on these quantities. Throughout this paper we will assume $t=t_{\mathrm{fin}}$ (otherwise replace $t_{\mathrm{fin}}$ by $t$). Fix a large positive integer $\Lf$ depending on $(d,t_{\mathrm{fin}},\beta,A,B)$, and \textbf{subdivide the time interval} $[0,t_{\mathrm{fin}}]$ as follows:
\begin{equation}\label{eq.def_short_time}
\tau:=\frac{t_{\mathrm{fin}}}{\Lf};\qquad [0,t_{\mathrm{fin}}]=\bigcup_{\ell=1}^\Lf [(\ell-1)\tau,\ell\tau].
\end{equation}
\item\label{it.defC*} \emph{Convention for constants}. Throughout the proof, we will distinguish \textbf{two types of (large) constants}: those depending only on $(d,t_{\mathrm{fin}},\beta,A,B,C_\Phi^0)$ and not on $\Lf$ will be denoted by $C$, while those depending on $\Lf$ will be denoted by $C^*$. The notions of $\lesssim$ and $\lesssim_*$ etc. are defined accordingly. We assume $\varepsilon$ is small enough depending on $C^*$.

We also fix an \textbf{explicit sequence of large constants} $(C_1^*,C_2^*,\cdots, C_{15}^*)$ such that 
 \begin{equation}\label{eq.defC*}C_{j+1}^*\gg(\textrm{a suitable function of }C_j^*)\end{equation} all depending on $\Lf$. We understand that (i) for any large constant $C$ under the above convention, we always assume $\Lf^C\ll C_j^*$ for each $j$, and (ii) any large constant which we denote by $C^*$ will always be $\gg C_j^*$ for each $j$.
\item\label{it.defnot3} \emph{Other parameters}. \textbf{We fix} $\Gamma$ sufficiently large but depending only on $d$ and $C_\Phi^0$. \textbf{We fix} $\mathsf{N}$ sufficiently large but suitably smaller than $\Gamma$, and it also depends only on $d$ and $C_\Phi$. These constants will be specified in Section \ref{sec.error}. We also fix $\upsilon:=C_\Phi^{-2}\cdot 3^{-d-1}$ and $\upsilon_0 := 1/d - \upsilon$, where $C_\Phi = (3 C^0_\Phi +2)( 6 C^0_\Phi +d+3)$.

We fix a \textbf{sequence of parameters} $(A_\ell,\Lambda_\ell)\,(0\leq \ell\leq\Lf)$ by
\begin{equation}\label{eq.defLambdaseq}
A_{\Lf}=|\log\varepsilon|;\quad \Lambda_{\ell}=A_{\ell}^{10d},\quad A_{\ell-1}=\Lambda_{\ell}^{10d}.
\end{equation} It is obvious that \begin{equation}\label{eq.ineq_Al}A_\ell,\Lambda_\ell\leq |\log\varepsilon|^{C^*},\quad \Lambda_{\ell-1}\gg A_{\ell-1}\gg \Lambda_\ell^2A_\ell.\end{equation} We also fix \textbf{two sequences}
\begin{equation}
\label{eq.decayseq}\beta_\ell:=\bigg(\frac{9}{10}+\frac{\Lf-\ell}{10\Lf}\bigg)\frac{\beta}{2};\quad \theta_0:=C^{-2}_\Phi 3^{-d-8},\quad \theta_{\ell}:=\bigg(\frac{9}{10}+\frac{\Lf-\ell}{10\Lf}\bigg)\theta_0
\end{equation} for $1\leq\ell\leq\Lf$, and set $\theta=  \theta_0/2$ in Theorem \ref{th.main} (note that $\theta_0 <\upsilon/3$).
\item\label{it.vector} \emph{Set, vector and integral notations}. For any positive integers $a$ and $b$, define $$[a:b]=\{a,a+1,\cdots,b\},\quad[a]:=[1:a].$$ For any set $A$, we shall adopt the \textbf{vector notation} $\vz_A:=(z_j)_{j\in A}$, and $\vz_{a:b}:=\vz_{[a:b]}$ and $\vz_a:=\vz_{[a]}$, and $\mathrm{d}\vz_A=\prod_{j\in A}\mathrm{d}z_j$. The same notion applies for other variables, such as $\boldsymbol{t}_A=(t_j)_{j\in A}$ etc. The notion $\vz$ without any subscripts will be used to denote some $\vz_N$ where $N$ depends on the context.

\item \emph{Indicator functions and linear operators}. In the proof we will need to use a lot of \textbf{indicator functions}. These functions will be denoted by $\mathbbm{1}_{(\cdots)}$, which are defined by certain properties represented by $(\cdots)$ depending on the context, and take values in $\{0,1\}$.

We will view any function $F$ (including the above indicator functions) as a \textbf{linear operator} defined by multiplication. For any linear operators $\Sc_1$ and $\Sc_2$ we define $\Sc_1\circ\Sc_2$ as their composition. For any linear operator $\Sc$ and function $F$, the notion $\Sc\circ F$ means the composition of $\Sc$ with the multiplication operator defined by $F$. This is different from the result of applying $\Sc$ to the function $F$, which we denote by $\Sc (F)$.

\item\label{it.misc}\emph{Miscellaneous}. We use the notion $|A|$ or $\# A$ to denote the cardinality of finite sets $A$, and use $|A|$ to denote the Lebesgue measure of Borel sets $A\subseteq\Rb^p$ for some $p$. Define $X_+=\max(X,0)$ and $X_-=-\min(X,0)$.
\end{enumerate}
\end{definition}

\section{Molecules, truncated dynamics and the cumulant formula}\label{sec.summary_concepts} The framework of this section follows \cite{DHM24}, which treats hard spheres, and we keep its notation whenever possible. The short-range interaction requires several new notions: multi-particle atoms, O-connections, boxes and M-atoms. Molecules are introduced in Section \ref{sec.molecule}, and the topological reduction of trajectories to molecules in Section \ref{sec.topological_reduction}. The truncated dynamics is defined in Section \ref{sec.truncated_dynamic}. The goal of this section is the cumulant formula, Proposition \ref{prop.cumulant_formula}, which is stated and proved in Section \ref{sec.formula_cumulant_statement}. 
  
\subsection{Molecules}\label{sec.molecule} We start with describing the short-range dynamics into the molecule language. Molecules are topological reductions of actual particle trajectories. 
\begin{definition}[Molecules]\label{def.molecule}A \newterm{(layered) molecule} $\Mb=(\Mc,\Ec,\Pc,\Lc,\Df)$ is a diagram defined as follows. 
\begin{enumerate}
    \item \textit{Atoms, edges, bonds and ends.} $\Mc$ is a set of nodes referred to as \newterm{atoms}, and $\Ec$ is a set of \newterm{edges} consisting of \newterm{bonds} and \newterm{ends} to be specified below. For convenience, below we will abuse notation and use $\Mb$ instead of $\Mc$ to denote the set of atoms in $\Mb$. An atom represents a localized interaction (collision or overlap) event.  
    \begin{enumerate}
        \item\label{it.molecule_atom} An atom $\nf\in \Mb$ is either a \newterm{C-atom} (collision) or an \newterm{O-atom} (overlap), or a \newterm{CO}-atom (corresponding to interactions involving more than two particles, where some collide while others overlap.) denoted respectively by $\circ$, $\bullet$, and $\textcolor{gray}{\bullet}$. 
         Each atom is assigned with a unique {\bf layer} $\ell[\nf]$ which is an integer in $ [\underline{\ell}:\ell]$, where $\ell$ is the highest layer and $\underline{\ell}$ is the lowest layer of $\Mb$. The boundaries between layers are indicated by horizontal lines.
 
        \item\label{it.molecule_edge} Each edge in $\Ec$ is either a \newterm{bond} between two atoms $\nf$ and $\nf'$, or an \textbf{end} emanating from one atom $\nf$, or an \newterm{empty end} that does not involve any atom. Each edge is either \textbf{top} or \textbf{bottom} when viewed at each involved atom $\nf$; a top (resp. bottom) edge is drawn above (resp. below) the atom.
         Moreover, each end is also either \newterm{free} or \newterm{fixed}, with fixed ends indicated by a $\times$. 
          If two atoms are connected by more than one bond, these bonds are called {\bf parallel bonds}. In particular, a pair of parallel bonds connecting the same two atoms is called a {\bf double bond}. 
        \item\label{it.molecule_parent} For two atoms $\nf$ and $\nf'$, if a bond exists between $\nf$ and $\nf'$ that is bottom at $\nf$ and top at $\nf'$, we say $\nf$ is \textbf{parent} of $\nf'$ and $\nf'$ is \textbf{child} of $\nf$.
        \item\label{it.molecule_atom1} For each atom, it has at least two particles, the number of top edges equals the number of bottom edges equals the number of particles at this atom. If an atom has exactly two top edges and two bottom edges, we call it a {\bf binary atom}; otherwise we call it a {\bf multi-particle} atom. 
 
        \item Each O-atom must be a \emph{binary atom}. A CO-atom must be a \emph{multi-particle atom} consisting of a collection of C-atoms and edges (with at least one C-atom), referred to as its \textbf{constituents}, joined by {\bf O-connections}. These constituents are called \textbf{atom-constituents} (also called \textbf{C-constituents}) and \textbf{edge-constituents}, respectively. Each constituent is assigned with top and bottom edges, which forms edges of the CO-atom.
        
        Each O-connection corresponds to a pair of particles in $\nf$, with the two particles belonging to distinct constituents. In particular, \emph{no O-connection joins a constituent to itself}. Moreover, each edge constituent at a CO-atom must be O-connected to \emph{exactly one atom-constituent} (see also Proposition \ref{prop.mol_axiom} \eqref{it.axiom1_CO}).
        
        For an edge in  $\Ec$ connected to a CO-atom $\nf$, we call it a \textbf{C-edge} (resp. an \textbf{O-edge}) at $\nf$ if it is an edge of a C-constituent (resp. edge constituent) of $\nf$. 
        \item Each C-atom and C-constituent is marked as \newterm{long interaction}, or \newterm{short interaction}, which describes the corresponding  collision duration in physical trajectory, see Definition \ref{def.top_reduction} below. A C-atom is marked as \newterm{crossing atom} if the corresponding collision happens across different layers. We define bottom crossing and top crossing accordingly.
    \end{enumerate}
    \item\label{it.particle_line_same} \textit{Particle lines.} $\Pc$ is a collection of paths: each path starts from a bottom end at an atom, then keeps connecting to a parent by a bond, until reaching a top end at another atom. Each chain $\pb\in\Pc$ is called a \textbf{particle line}. We require that each edge belongs to a unique particle line. For a top edge and a bottom edge at an atom $\nf$, we say they are \textbf{serial} if they belong to the same particle line. Note that each particle line $\pb$ contains a unique bottom end and a unique top end.  If a molecule has no fixed end, the particle lines are identified with \emph{physical particles}. We further require that \emph{each O-edge at an atom $\nf$ must be serial with an O-edge at $\nf$.} 
    
    Recall $\ell$ defined in (\ref{it.molecule_atom}). We may define a \textbf{start layer} $\ell_1[\pb]\in[\underline{\ell}:\ell]$ and \textbf{finish layer} $\ell_2[\pb]\in[\underline{\ell}:\ell+1]$ for each particle line $\pb$. As $\pb$ contains a unique bottom (and unique top) end, we may also equivalently write $\ell_1[e]$ (resp. $\ell_2[e]$) for each bottom (resp. top) end $e$ instead of $\ell_1[\pb]$ and $\ell_2[\pb]$. We refer to the interval $[\ell_1[\pb]:\ell_2[\pb]]$ as the \textbf{lifespan} of the particle line $\pb$.

    \item\label{it.label_root} For a molecule $\Mb$, we define \textbf{set of all particle lines} in $\Mb$ to be $p(\Mb)$. Define a particle line to be \textbf{root particle line}, if its finish layer is $\ell_2[\pb] = \ell+1$. Define the set of root particle lines to be $r(\Mb)$. Moreover, we label the root particle lines by values in $[|r(\Mb)|]$. 

    \item\label{it.decoration} For each atom $\nf$, we define its \textbf{decoration} as follows. Denote by $\Pc_\nf$ the particles of $\nf$.  For a binary atom \(\nf\), the decoration is simply
    \begin{equation*}
        \Df(\nf) =\Cf (\nf) =\bigl((\pb_i,\pb_j)_C\bigr) \quad\text{or}\quad\! \Df(\nf)=\Of(\nf) = \bigl((\pb_i,\pb_j)_O\bigr),
    \end{equation*}
    where $\pb_i,\pb_j\in \Pc_\nf$, and the subscript indicates whether the interaction is a collision or an overlap.

    For a multi-particle C-atom $\nf$, its decoration records \emph{the order in which particles first become involved in the collision, together with their respective interaction partners}. More precisely, we define
    \[
    \Df(\nf)=\bigl((\pb_{i_1},\pb_{j_1})_C,\ldots,(\pb_{i_r},\pb_{j_r})_C\bigr).
    \]
    Here $\pb_{i_1},...,\pb_{i_r},\pb_{j_1},...,\pb_{j_r} \in \Pc_{\nf}$ and some of them may be the same. The collision sequence $\Cf(\nf)$ is ordered chronologically. Thus, for every $1\leq n<r$, the interaction represented by $(\pb_{i_n},\pb_{j_n})_C$ occurs before the interaction represented by $(\pb_{i_{n+1}},\pb_{j_{n+1}})_C$.

    Now let $\nf$ be a CO-atom with C-constituents $\nf_1,\ldots,\nf_m$ (and may also have edge constituents). We define its decoration by
    $$\Df(\nf)=\big(\Cf(\nf),\Of(\nf) \big)= \big( \Cf(\nf_1),\ldots,\Cf(\nf_m), (\pb_{i_1},\pb_{j_1})_O,\ldots,(\pb_{i_s},\pb_{j_s})_O \big).$$
    The decorations of different C-constituents and the O-connection are also ordered, while the ordering within each $\Cf(\nf_k)$ is preserved.
    \item\label{it.initial_link} \textit{Initial links.} Finally, $\Lc$ is a collection of pairs of bottom ends $(e,e')$ with $\ell_1[e]=\ell_1[e']=1$ (provided $\underline{\ell}=1$). Since each particle line has a unique bottom end by (\ref{it.particle_line_same}), we can equivalently view $\Lc$ as a collection of pairs of particle lines with start layer $1$. Each pair in $\Lc$ is called an \textbf{initial link}. 
    \item\label{it.ov_segment} \textit{Ov-segments.} For later use, we also define an \textbf{ov-segment} to be a path between two atoms (or ends) such that the two edges of the path at each internal atom are both O-edges at the atom, and this path is contained in a single particle line. We say an ov-segment is \textbf{maximal} if it is not contained in other ov-segments. 

    \item The objects $(\Mc, \Ec, \Pc,\Df)$ satisfy the following requirements. 
    \begin{enumerate}
        \item Starting from any atom $\nf$, it is impossible to get back to $\nf$ by iteratively taking parent or by iteratively taking child.
        \item If $\nf$ belongs to a particle line $\pb$ then we must have $\ell_1[\pb]\leq\ell [\nf]\leq  \ell_2[\pb]$. 
        \item If $\nf$ is a parent of $\nf'$, then $\ell[\nf'] \leq \ell[\nf]$.
    \end{enumerate}
\end{enumerate}
\end{definition} 

\begin{definition}[Clusters]\label{def.cluster} Given a molecule $\Mb$, we obtain a new molecule $\Mb'$ by deleting all O-connections from $\Mb$, leaving only C-atoms. The connected components $(\Mb_1, \ldots, \Mb_k)$ of $\Mb'$ are called the \textbf{clusters} of $\Mb$. We say that an O-connection $\of$ \textbf{joins} $\Mb_j$ and $\Mb_{j'}$ if the two particle lines at the O-connection connects the C-atoms in $\Mb_j$ and in $\Mb_{j'}$ respectively. We also define the auxiliary \textbf{cluster graph} to be the graph with nodes being clusters $\Mb_j$, and an edge is drawn between $\Mb_j$ and $\Mb_{j'}$ for each O-connection joining them. 
\end{definition}

We introduce some additional notions associated with molecules. 
\begin{definition}[More on molecules]\label{def.more_on_molecule} Let $\Mb = (\Mc,\Ec,\Pc,\Lc,\Df)$ be a molecule and $\nf,\nf'\in\Mb$ be atoms.
\begin{enumerate}
    \item\label{it.full_molecule} If a molecule $\Mb$ has no fixed end, we say $\Mb$ is a \newterm{full molecule}. 
    \item\label{it.fully_fixed} If the top or bottom edges of an atom are all fixed ends, we call it a \newterm{fully fixed atom}. If the top or bottom edges of an atom constituents are all fixed ends, we call it a \newterm{fully fixed atom-constituent}.
    \item \emph{Descendants and ancestors.} If $\nf$ is obtained from $\nf'$ by iteratively taking parents, then we say $\nf$ is an \newterm{ancestor} of $\nf'$, and $\nf'$ is a \newterm{descendant} of $\nf$. We use $S_\nf$ (resp. $Z_\nf$) to denote the set of descendants (resp. ancestors) of $\nf$. By convention, we understand that $\nf\in S_\nf$ (and $\nf\in Z_\nf$).
    \item\label{it.partial_order} \emph{The partial ordering of atoms and edges.} We also define a \newterm{partial ordering} between all atoms, by defining $\nf'\prec\nf$ if and only if $\nf'\neq\nf$ and is a descendant of $\nf$ (equivalently $\nf$ is ancestor of $\nf'$). In particular $\nf'\prec\nf$ if $\nf'$ is a child of $\nf$ (equivalently $\nf$ is parent of $\nf'$). 
    \item\label{it.minimax} Given any subset $A\subseteq\Mb$, we can define a \newterm{lowest} atom in $A$ to be any atom $\nf\in A$ that satisfies $\nf'\not\in A$ for each $\nf'\prec\nf$ (i.e. $\nf'\neq\nf$ is descendant of $\nf$). Define \newterm{highest} atoms in the same way. These are just minimal and maximal elements of $A$ under the partial ordering $\prec$, and note that they need not be unique for general $A$.
 
    \item\label{it.subsets_Ec} \emph{Sets $\Ec_*$ and $\Ec_{\mathrm{end}}$.} Define $\Ec_*$ to be the set of all edges that are either bonds or free ends (i.e. not counting fixed ends). Define also the sets $\Ec_{\mathrm{end}}$ (resp. $\Ec_{\mathrm{end}}^-$, $\Ec_{\mathrm{end}}^+$) to be the set of all ends (resp. bottom ends, top ends).  
    
    Note that the sets $\Ec_{\mathrm{end}}^-$ and $\Ec_{\mathrm{end}}^+$ are in bijection with the set $p(\Mb)$ of all particle lines, see Definition \ref{def.molecule} (\ref{it.particle_line_same}).
    \item\label{it.number_atom} \emph{Number of atoms and particle lines.} We will use $|\Mb|$, $|\Mb|_{\mathrm{O}}$, $|\Mb|_p$, $|\Mb|_r$ and $|\Mb|_{p\backslash r}$ to denote the number of atoms, O-connections, particle lines (i.e. $|p(\Mb)|$), root particle lines (i.e. $|r(\Mb)|$), and non-root particle lines (i.e. $|(p\backslash r)(\Mb)|$) of $\Mb$ respectively. 
    \item\label{def.layer_interval} Let $\Mb_{\ell'}$ be the layers of $\Mb$. Define the \textbf{layer intervals}
\[\Mb_{[\ell_1:\ell_2]}:=\bigcup_{\ell_1\leq\ell'\leq\ell_2}\Mb_{\ell'},\quad \Mb_{(\ell_1:\ell_2)}:=\bigcup_{\ell_1<\ell'<\ell_2}\Mb_{\ell'}\quad \Mb_{<\ell_1}=\bigcup_{\ell'<\ell_1}\Mb_{\ell'},\quad \textrm{etc.}\] After defining the refined layers $\Mb_{\zeta}^T$ (Definition \ref{def.layer_refine}) below, we will use the same notations $\Mb_{[\zeta_1:\zeta_2]}^T$ etc. for the intervals of refined layers.
    \item\label{it.submolecule} \emph{Sub-molecules.} Let $\Mb'\subseteq\Mb$ be a subset of atoms or C-constituents. We will regard $\Mb'$ as a molecule (referred to as a \newterm{sub-molecule} of $\Mb$),  by viewing each ov-segment between two atoms or constituents in $\Mb'$ as a bond, and view all other ov-segments between an atom in $\Mb'$ and an atom in $\Mb \backslash \Mb'$ as a free end at the atom in $\Mb'$. After that, we delete all O-connections between $\Mb'$ and $\Mb \backslash \Mb'$ (deleting all O-connections at an atom $\nf$ means we remove all the O-connections $\Of(\nf)$, but the constituents in $\nf$ are unchanged).
\end{enumerate}
\end{definition}

\begin{definition}\label{def.prsets} Let $\Mb$ be a molecule as in Definition \ref{def.molecule} with $\ell$ layers. For $\ell'\in[\underline{\ell}:\ell]$, define the sets
\begin{align}\label{eq.pset}
p(\Mb_{\ell'})&:=\{\pb:\ell_1[\pb]\leq\ell'\leq\ell_2[\pb]\},\\
\label{eq.rset}
r(\Mb_{\ell'})&:=\{\pb:\ell_1[\pb]\leq\ell'<\ell'+1\leq\ell_2[\pb]\}.
\end{align} 
We refer to the particle lines in $r(\Mb_{\ell'})$ as the \textbf{root particles of layer $\ell'$}. We refer to $p(\Mb_{\ell'})$ as the set of particle lines \textbf{occurring in layer $\ell'$}, and $r(\Mb_{\ell'})$ as the set of particle lines \textbf{crossing into layer $\ell'$}. 
 Define also $|\Mb_{\ell'}|_p=|p(\Mb_{\ell'})|$ and $|\Mb_{\ell'}|_r=|r(\Mb_{\ell'})|$ etc. similar to Definition \ref{def.more_on_molecule}.
\end{definition}
\begin{definition}\label{def.connectedvia} Let $\pb_1$ and $\pb_2$ be two particle lines, and $A\subseteq\Mb$ be a set. We say:
\begin{itemize}
\item $\pb_1$ \textbf{intersects} $A$, if there exists an atom in $A$ that contains $\pb_1$.
\item $\pb_1$ \textbf{belongs to} an atom $\nf$, if it intersects $\nf$.
\item $\pb_1$ and $\pb_2$ are \textbf{connected via $A$}, if $\pb_1$ and $\pb_2$ intersect the same component of $A$;
\item $\pb_1$ is \textbf{connected to a cycle within $A$}, if $\pb_1$ intersects a component of $A$ that contains a cycle;
\item $\pb_1$ \textbf{forms an initial link within $A$}, if there exists a particle line $\pb_1'$ which either equals $\pb_1$ or is connected to $\pb_1$ via $A$, such that the bottom end in $\pb_1'$ is involved in an initial link.
\end{itemize}
\end{definition}

\subsection{Topological reduction}\label{sec.topological_reduction} In this section, we rigorously formulate the notion of topological reduction. We start with the notion of admissible trajectories.

\begin{definition}[Admissible trajectories]\label{def.admissible}
Let $ \vz(t)=\vz_N(t)=(z_i(t))_{i}$ be any particle trajectory (a priori with no equation) in $[(\ell'-1) \tau, \ell'\tau]$ and $\Pf(t) := (\mathsf{P}_i(t))_{i \in[N]}$ for each particle $i$. We say  the pair $(\vz(t), \Pf(t))$ is \newterm{admissible} if 
\begin{enumerate}
    \item\label{it.admissible.1} $\Pf(t)$ satisfies the following properties
     \begin{enumerate}
        \item We always have $j \in \Pf_i(t)$ only if $|x_j(t) -x_i(t)|\leq \eps$;
       \item  we have that $i \in \Pf_i(t)$ for any $i \in [N]$; 
       \item  if $j \in \Pf_i(t)$ for some $j\neq i $, then $i \in \Pf_j(t)$. If $j \in\Pf_i(t)$ for some $t \in [(\ell'-1)\tau,\ell' \tau)$ then we also have $j \in \Pf_i(t_+)$ unless $|x_i(t_+) -x_j(t_+)| \geq \eps$. 
       \item  if $j\notin \Pf_i(t)$ at some time $t \in [(\ell'-1)\tau,\ell' \tau)$, then $j\notin \Pf_i(t_+)$ unless the pre-collision condition \eqref{it.pre-collision} holds at time $t$.
     \end{enumerate}
    \item The dynamics of $z_i (t)$ is given by 
     $$ \frac{\d x_i}{\d t}=v_i, \qquad \frac{\d v_i}{\d t} = \sum_{j\in\mathcal P_i(t)\setminus\{i\}} -\frac{1}{\eps}\nabla \Phi (\frac{x_i-x_j}{\eps}).$$
    In particular, if $\mathcal P_i(t)=\{i\}$, then 
    $$\frac{\d x_i}{\d t}=v_i,\qquad\frac{\d v_i}{\d t}=0.$$
    \item\label{it.pre-collision} At time $t \in [(\ell'-1)\tau, \ell'\tau]$, we say two particles satisfy the {\bf pre-collision condition} if: $ x_i(t) - x_j(t) = \eps \nu \ (\nu \in \Sb^{d-1})$, and $(v_i(t) - v_j(t)) \cdot \nu \leq 0.$ Moreover, we say a pairwise collision happens between $i$ and $j$ if $ j \in \Pf_i(t_+)$ (and also $i \in \Pf_j(t_+)$ by \eqref{it.admissible.1}).
    
    If $j \in \Pf_i((\ell'-1)\tau)$, then we also say a pairwise collision happens between $i$ and $j$.
    \item\label{it.overlap} An \newterm{overlap} between two particles $i\neq j$ is a time $t$ such that the pre-collision condition holds at time $t$ but $j \notin \Pf_i(t_+)$; or $|x_i((\ell'-1)\tau) -x_j((\ell'-1)\tau)| < \eps$ but $j \notin \Pf_i((\ell'-1)\tau)$. This means that the trajectory of particle $i$ is not influenced by $j$ at $t_+$. 
    \item\label{it.collision} We regard all pairwise collisions in $[(\ell'-1)\tau,\ell' \tau]$ as a single \newterm{collision} event as in Definition \ref{def.short_range} \eqref{it.O5}. 
\end{enumerate}

Moreover, we say a particle trajectory $\vz(t)$ is \textbf{admissible} in $[0,\ell \tau]$, if it is admissible in each interval $[(\ell'-1) \tau, \ell'\tau]$ for $\ell' \leq \ell$, and for each $2\leq\ell'\leq\ell$ the value of $\Pf_i$ at the end of the interval $[(\ell'-2)\tau,(\ell'-1)\tau]$ contains its value at the start of the interval $[(\ell'-1)\tau,\ell'\tau]$, for every $i$ (i.e. pairs may be removed from $\Pf$ at time $(\ell'-1)\tau$, but not added).
\end{definition}
 \begin{remark}
    In the admissible trajectory defined above,  if there exists a collision (as defined in Definition \ref{def.short_range} \eqref{it.O5}) such that $\ell' \tau$ belongs to its collision interval for some $\ell'$, then we will regard the interaction in $[(\ell'-1) \tau, \ell'\tau)$ and $[\ell' \tau, (\ell' +1)\tau]$ as \emph{two different interactions} (and correspond to \emph{two different atoms} in the molecule).
 \end{remark}
\begin{remark}\label{rem.conserved_energy}
    Define the corresponding energy of admissible trajectories:
    \begin{equation}\label{def.energy_admissible}
        E_{\Pf}(t) = \frac{1}{2} \sum_{i} |v_i(t)|^2 + \sum_{i< j, j \in \Pf_i(t)} \Phi(\frac{x_i(t) - x_j(t)}{\eps}).
    \end{equation}
    Then the admissible trajectory $\vz(t)$ has the conserved energy in each interval $[(\ell'-1)\tau,\ell' \tau]$, since $\Phi(|x|)=0$ for $|x|=1$, so adding or removing a pair from $\Pf$ at distance $\eps$ does not change $E_{\Pf}$.  At times $ \ell' \tau$, pairs can only be removed from $\Pf$, so $E_\Pf$ does not increase there. Hence the total energy $E_\Pf(t)$ is nonincreasing in time.
\end{remark}
\begin{definition}[Topological reduction]\label{def.top_reduction} Let $\vz(t)$ with $\Pf(t)$ be an admissible trajectory for $t\in[0,\ell\tau]$,
     and $\Mb$ be a molecule with $\ell$ layers (Definition \ref{def.molecule}). We say $\vz(t)$ \textbf{topologically reduces} to $\Mb = \cup_{1 \leq \ell' \leq \ell} \Mb_{\ell'}$, if there exists a bijection between maximal ov-segments of $\Mb$ and maximal linear pieces of the trajectory, which satisfies the following condition:
\begin{enumerate}
    \item If two (or more) maximal ov-segments intersect at a C-atom in $\Mb_{\ell'}$, then the corresponding maximal linear pieces of $\vz(t)$ have a binary (or multi-particle) collision in $[(\ell'-1)\tau,\ell' \tau]$. If the collision duration is more than $\eps^{\upsilon_0}$ (in the sense of Definition \ref{def.admissible} \eqref{it.collision}; recall that $\upsilon_0$ is defined in Section \ref{sec.notation} \eqref{it.defnot3}), or there exists a particle $j\in \Pf_i((\ell'-1)\tau)\backslash \{i\}$ (so they are in the collision at time $t=(\ell'-1)\tau$), such that $|v_i((\ell'-1)\tau)-v_j((\ell'-1)\tau)|\leq \eps^{1- \upsilon_0}$, then the corresponding C-atom is marked as long interaction, otherwise it is marked as short interaction. Moreover, for each C-atom, the decoration comes from the particle pairs that one of the particle gets into the interaction for the first time. Decorations are ordered according to the timing of these events.
    
    Moreover, if at time $(\ell'-1) \tau$ (or $\ell' \tau$), there exist particles $i,j$ in the atom or the constituent such that $j \in \Pf_i((\ell'-1)\tau)$ (or $ j \in \Pf_i (\ell'\tau)$) then the corresponding C-atom or C-constituent in layer $\ell'$ is marked as bottom crossing (or top crossing) . By the property of $\Pf$ in Definition \ref{def.admissible}, each bottom crossing atom in layer $\ell' > 1$ must be connected with a top crossing atom in layer $\ell'-1$. 
    
    If two maximal ov-segment intersects at an O-atom in $\Mb_{\ell'}$, then the corresponding maximal linear pieces have an overlap in $[(\ell'-1)\tau,\ell' \tau]$.  For each O-atom, the decoration is simply $(\pb_i,\pb_j)_\Of$, where $\pb_i,\pb_j$ are the corresponding particle lines.
    
     If one maximal ov-segment intersects at a C-constituent of a CO-atom with some other ov-segments, then the corresponding trajectories have a collision. Also, this C-constituent is marked as long interaction or short interaction as above. If one maximal ov-segment intersects at an edge-constituent of CO-atom in $[(\ell'-1)\tau,\ell' \tau]$ and has O-connection with some particle line in this atom, then the corresponding trajectories have an overlap in $[(\ell'-1)\tau,\ell' \tau]$. For each CO-atom, the decoration of each C-constituent is defined as above, while the decoration of each O-connection is determined by the overlaps between the corresponding pairs of particle lines. If two C-constituents are connected by an O-connection $(\pb_i,\pb_j)_\Of$, then the corresponding particles in each C-constituent collide as stated, and overlap happens between particles $\pb_i$ and $\pb_j$ while both collisions are ongoing (an overlap of a particle that has already left its collision is regarded as an overlap between an edge and an atom).
 
    \item The induced correspondence between atoms and collisions/overlaps preserves order (i.e. if $\nf$ is parent of $\nf'$ then the corresponding collision/overlap happens later in time), and \emph{is bijective on collisions} (it is injective on overlaps by time ordering, but might not be surjective).
    \item If two maximal ov-segments belong to the same particle line in $\Mb$, then the corresponding linear pieces belong to the trajectory of the same particle.
\end{enumerate}
\end{definition}
\begin{remark}\label{rem.top_reduction_uniqueness} 
The topological reduction $\Mb$ of a given trajectory $\vz(t)$ may not be unique, since we do not require a bijective correspondence between O-connections and the actual set of overlaps. The O-connection may represent only an arbitrary subset of all overlaps.  However there is always a \emph{unique} topological reduction that is C-molecule (i.e. contains only C-atoms). 
\end{remark}

\subsection{Truncated dynamics}\label{sec.truncated_dynamic} To introduce the truncated dynamics (Definition \ref{def.T_dynamics}), we first introduce the modified dynamics.

\begin{definition}[Modified dynamics]\label{def.modified_dynamics} We define the \textbf{modified dynamics}, governed by a given property $\mathtt{Pro}$, as follows.    
\begin{enumerate}
\item\label{it.modified_dynamics} \emph{Modified dynamics.} Let $\Xc_N$ be a set of pairs $(\vz,\Pf)$, where $\vz =\vz_{N} \in \Rb^{2dN}$ and $\Pf= (\Pf_i)_{i\in [N]}$ satisfies Definition \ref{def.admissible} \eqref{it.admissible.1} at a fixed time. Given an initial configuration $ (\vz^0, \Pf^0)\in \Xc_N$, we define the dynamics $\vz^{\mathtt{Pro}}(t)=\vz^{\mathtt{Pro}}(t,\vz^0,\Pf^0)$  as follows:
\begin{enumerate}
\item\label{it.modified_dynamics_0} We have $\vz^{\mathtt{Pro}}(0)=\vz^0$.
\item\label{it.modified_dynamics_1} Suppose $\vz^{\mathtt{Pro}}(t')$ is known for $t'\in[0,t]$ (so $\Pf_i(t')$ is known for every $i \in[N]$). Let $\boldsymbol{\gamma}$ be the admissible trajectory $\vz^{\mathtt{Pro}}(t')\,(t'\in[0,t])$. Let particles $\pb_i, \pb_j$ satisfy pre-collision condition at time $t$, and also $\mathtt{Pro}(\gamma,p_i,p_j) =\mathtt{True}$, then they collide as described in Definition \ref{def.admissible} (i.e. $j \in \Pf_i(t_+)$ and $i \in \Pf_j(t_+)$).  Otherwise, we have $j \notin \Pf_i(t_+)$.
\end{enumerate}
\item \emph{Flow map and flow operator.} We define the \textbf{flow map} $\Hc^{\mathtt{Pro}}(t)=\Hc_N^{\mathtt{Pro}}(t)$  by
\begin{equation}\label{eq.defH^M}
\Hc^{\mathtt{Pro}} (t):\Xc_N\to \Xc_N,\quad\Hc^{\mathtt{Pro}} (t)(\vz^0,\Pf^0)=(\vz^{\mathtt{Pro}}(t), \Pf(t)),
\end{equation} 
where $\vz^{\mathtt{Pro}}(t)=\vz_N^{\mathtt{Pro}}(t)$ is defined in (\ref{it.modified_dynamics}). The mapping $\Hc^{\mathtt{Pro}}(t)$ need not be injective; however, we can still define the \textbf{flow operator} $\Sc^{\mathtt{Pro}}(t)$ by
\begin{equation}\label{eq.defS^M}
    \Sc^{\mathtt{Pro}}(t)f(\vz,\Pf )=\sum_{(\vz^0,\Pf^0) :\Hc^{\mathtt{Pro}}(t)(\vz^0,\Pf^0)=(\vz,\Pf)}f(\vz^0,\Pf^0).
\end{equation}
\item \emph{Collision and overlap.} The notions of \textbf{collisions} and \textbf{overlaps} for any trajectory of modified dynamics, are defined as in Definition \ref{def.admissible}.
\item \emph{Label of particles.} In this definition, we assume that particles are labeled by $[N]=\{1,\cdots,N\}$ (see the convention in Definition \ref{def.notation} (\ref{it.vector})). For a general labeling set $\lambda$, the formulation remains the same, except that the subscript $N$ in expressions like $\vz_N(t)$, $\Hc_N(t)$, $\Sc_N(t)$ etc. is replaced by $\lambda$.
\end{enumerate}
\end{definition}

\begin{definition}[The extended dynamics]\label{def.extended_dynamic}
Given $(\vz^0,\Pf^0)$, the extended dynamics $\bigl(\vz^E(t),\mathcal \Pf^E(t)\bigr)$ is the modified dynamics (Definition \ref{def.modified_dynamics}) with $\mathtt{Pro}(\boldsymbol{\gamma},\pb_i,\pb_j) \equiv \mathtt{True}$. Define $\mathcal H^E_N(t)$ and $\mathcal S^E_N(t)$ accordingly.
\end{definition}

We now introduce the definition of box, this is introduced to deal with the case where a finite number of particles collide too many times.
\begin{definition}\label{def.pathological}
   Recall the definition of $\mathsf{N}$ in Definition \ref{def.notation} \eqref{it.defnot3}. Define the function 
   $$G(q) = (2 q+1)\cdot \binom{q}{2} \cdot (32 q^{3/2})^{q^2}.$$
   Fix an integer $3 \leq q \leq \mathsf{N}$. We define \emph{pathological trajectories} as follows: we say that the particles $\pb_1,\dots,\pb_q$ undergo pathological trajectories, if the corresponding C-molecule $\Mb$ has a connected sub-molecule $\Mb'$ such that (i) its set of particle lines is exactly $\{\pb_1,\dots,\pb_q\}$; (ii) it has at least $G(q)$ atoms; (iii) there exists no atom in $\Mb \backslash \Mb'$ that has both descendants and ancestors in $\Mb'$.\footnote{The function $G$ is only used as a threshold, chosen such that the bottom and top parts of a box have a bounded number of atoms (see Part 1 of the proof of Proposition \ref{prop.layerrec3}); since no bound of Burago--Ferleger--Kononenko type \cite{BFK98} is available, the number of collisions in the middle part of a box is not bounded.}
\end{definition}

\begin{definition}\label{def.box}  Let $\vz(t)$ be an admissible trajectory in $[(\ell-1)\tau,\ell\tau]$, and $\Mb$ the corresponding C-molecule. We reduce its pathological collisions by regarding each maximal pathological-collision sub-molecule as a single effective interaction, called a \emph{box}. The resulting molecule is called a \emph{boxed C-molecule}. The reduction is performed as follows. 
    \begin{enumerate}
        \item\label{it.explore_box} We start with $q=3$. Among all connected \emph{maximal} sub-molecules $\Mb'$ containing $q$ particles that undergo pathological trajectories as in Definition \ref{def.pathological}, we choose the one whose earliest collision occurs first, and regard it as a single box (such choice is unique up to a zero measure set by Proposition \ref{prop.shortrange}.)
        
        Each box is divided into three parts: the bottom part, the middle part, and the top part. The bottom part is constructed as follows. Consider the minimal sub-molecules $A$ in $\Mb'$ such that each particle line in $\Mb'$ is connected within $A$, and \emph{$A$ has no descendants in $\Mb'\backslash A$}. If there is more than one such $A$, we choose the one whose latest collision is the earliest; if this still does not determine $A$, we compare the second latest collisions, and so on. If $A$ already contains a cycle, then we identify $A$ as the bottom part of the box; if not, we choose the lowest atom $\pf$ in $\Mb'\backslash A$ which corresponds to the earliest collision in atoms $\Mb'\backslash A$. Then the bottom part of the box is defined to be $A \cup\{\pf\}$. It is clear that for each admissible trajectory, this choice of the bottom part is unique. The top part can be defined similarly, with time reversed and descendants replaced by ancestors. The remaining part is called the middle part. We treat it as a whole without distinguishing its internal collision structure. Two boxes are regarded as the same if they have the same bottom and top parts, regardless of the internal structures of their middle parts.

        \item We then view the box constructed in \eqref{it.explore_box} \emph{as a single atom}, and repeat the procedure until there exists no such sub-molecule with $q$. Then we go to \eqref{it.explore_box} with $q+1$, or we stop exploring if $q = \mathsf{N}$ and call the resultant graph a boxed molecule.
    \end{enumerate}
\end{definition} 
   \begin{remark}\label{rem.uniqueness_boxed}
    By Definition \ref{def.box}, we see that for each admissible trajectory $\vz(t)$, it reduces to a unique boxed molecule. Indeed, up to a zero measure set, the starting times of different pairwise interactions are different (Proposition \ref{prop.shortrange}), so different atoms correspond to different collision times. Moreover, each step decreases the number of atoms (each box being counted as one atom), so the procedure terminates, and its result is unique.
 \end{remark}
\begin{remark}
    We now explain the motivation of defining the box. Pathological trajectories contain a large number of recollisions. In the combinatorial estimates of Section \ref{sec.elem_int}, each recollision costs a factor of $|\log\eps|^{C^*}$, and we need comparable good molecules to compensate this loss. However, if we count the quantity $\sigma$ as in Proposition \ref{prop.alg_up}, most atoms involved in pathological collisions are fully fixed when they are cut, whatever cutting algorithm we choose. Therefore, for those pathological trajectories, we only expand its bottom and top parts, which already give enough gain. 
\end{remark}
\begin{definition}[boxed molecule]\label{def.boxed_molecule}
We define a boxed molecule $\Mb_B=(\Mc_B, \Ec,\Pc,\Df)$ similarly to Definition \ref{def.molecule}, but with boxes and the following adaptations: 
\begin{enumerate}
    \item We view the middle part of a box as a whole (i.e. as a single object), called an {\bf M-atom}.  We call an atom {\bf MO-atom} if it is formed by joining O-connections among some edges or C-atoms and contains at least one M-atom. The corresponding constituent is called M-constituent (which is also regarded as an atom-constituent, but not as a C-constituent). The decorations of an MO-atom are inherited from the decorations of all atoms contained within it in $\Mb_B$. Each atom in the same box is contained in the same layer. An MO-atom $\mf$ may have C-atoms, M-atoms, and edges as constituents. Moreover, if an MO-atom or a CO-atom has at least two atom-constituents, then each atom constituent must be O-connected with some other atom constituent. Each M-atom or M-constituent is not marked as long interaction or short interaction.
    
    \item For each molecule in a single layer, it is formed by joining O-connections among clusters. Each box originally belongs to one cluster in a single layer. In certain situations, we \emph{treat a box as a whole}, by viewing each C-atom or M-atom forming the box before joining any O-connections as a single object, and we call it a {\bf B-atom} (which is a sub-molecule of $\Mb_B$). We can also define {\bf BO-atom} if the B-atom has O-connection with other constituents, the decoration of a BO-atom is inherited from $\Mb_B$.  
    
    By \emph{viewing the atoms in a box individually}, we mean treating each of its constituent atoms separately.
    \item\label{it.recollision_number_box} \emph{Recollision number.} Given a subset $A\subseteq\Mb_B$, we define its recollision number by \emph{viewing boxes as a whole}, $\rho(A)=\#_{\mathrm{bond}(A)}-|A|+\#_{\mathrm{comp}(A)}$, where $\#_{\mathrm{bond}(A)}$ and $\#_{\mathrm{comp}(A)}$ are respectively the number of bonds between atoms in $A$, and the number of connected components of $A$. We say a set $A$ is \textbf{recollisional} (or equivalently \textbf{cyclic}) if $\rho(A)>0$. Note that $\rho(A)$ \textbf{does not count cycles inside boxes.}
    \item\label{it.number_box} For a boxed C-molecule $\Mb$, we denote by $\#_{\mathrm{box}}(\Mb)$ the number of maximal boxes in $\Mb$. Namely, if a box is contained in another Box, then only the latter is counted. 
    
    For a general boxed molecule $\Mb$, we denote by $\#_{\mathrm{box}}(\Mb)$ the total number of boxes in all clusters and all layers.
   \item\label{it.box_connect} In the top/bottom part of a box, each particle line is connected, there exists at least one cycle. 
   \item All definitions on molecules in Section \ref{sec.molecule} also apply to the boxed molecule. We also use $|\Mb_B|$ to define the number of atoms in $\Mb_B$, but in the sense that \emph{viewing each box as one atom}.
\end{enumerate}
\end{definition} 
\begin{remark}
    We also write a boxed molecule as $\Mb$ and say it is a molecule when there is no confusion.
\end{remark}
To define the truncated dynamics, we also introduce an auxiliary trajectory.
\begin{definition}\label{def.aux_time}
    Suppose $\vz(t')$ is known for $t'\in [0,t]$ for some $t>0$ (so $\Pf(t')$ is known). Let $\gamma$ be the corresponding admissible trajectory and $t$ be a pre-collision time of particles $\pb_i$ and $\pb_j$. We define an auxiliary trajectory $\boldsymbol{\gamma}'$ obtained by setting $\mathtt{Pro}(\gamma,\pb_i,\pb_j) = \mathtt{True}$, but $\mathtt{Pro}(\gamma,\pb_i,\pb_j) = \mathtt{False}$ for every pre-collision time $t' > t$. Define the boxed molecule reduced from $\boldsymbol{\gamma}'$ by $\Mb'_B$. Define the duration of the collision containing this pairwise collision as in Definition \ref{def.admissible} \eqref{it.collision}; it determines whether the corresponding C-atom is marked as long or short interaction (Definition \ref{def.top_reduction}).
\end{definition}
\begin{definition}[Well-prepared initial state]\label{def.well-prepared}
    Given an initial state of particles $\vz^0 \in \Rb^{2d N}$ and initial interaction set $\Pf^0_i$ for each $i $. Define 
$$\widetilde{\Pf}_i^0 = \{j|\textrm{ there exists } \{ j_0,\cdots,j_r\}\textrm{ with } j_0 = j\textrm{ and } j_r =i,\textrm{ such that } \ \ j_{k+1} \in \Pf_{j_k}^0 \textrm{ for every } 0\leq k < r\}.$$
 Then an initial data $(\vz^0,\Pf^0)$ is said to be well-prepared if $|\widetilde{\Pf}_i^0| \leq d+1$ for every $i \in [N]$, and for every $j \in \Pf^0_i\backslash \{i\}$ we have $|v_i-v_j|\geq \eps^{1-\upsilon_0}$.
\end{definition}
\begin{definition}[The truncated dynamics]\label{def.T_dynamics} Let $(\Lambda,\Gamma,\mathsf{N},\upsilon_0)$ be fixed. 

    We first define the truncated dynamics for well-prepared $(\vz^0,\Pf^0)$, we define the \textbf{truncated dynamics} (or \textbf{T-dynamics}), denoted by $\vz^{T}(t)$, to be the modified dynamics in Definition \ref{def.modified_dynamics}. At each pre-collision time $t$ of particles $\pb_i$ and $\pb_j$, we define $\Mb'_B$ as in Definition \ref{def.aux_time}. Then $\mathtt{Pro}(\gamma,\pb_i,\pb_j) = \mathtt{True}$, if $\Mb'_B$ has no C-atom with more than $d+1$ particles, and no atom is marked as \emph{long interaction}, and one of the following conditions holds:
\begin{enumerate}\label{eq.trunc_dym}
    \item\label{it.truncated_dynamic1} If we have $\Mb_{B,i} = \Mb_{B,j}$, and $\#_{\mathrm{box}}(\Mb'_{B,i}) \leq 2 \Gamma$ and $\rho(\Mb_{B,i}')\leq\Gamma$; 
    \item\label{it.truncated_dynamic2} If we have $\Mb_{B,i}\neq \Mb_{B,j}$ and $|\Mb_{B,i}|_p + |\Mb_{B,j}|_p \leq\Lambda$ and $\rho(\Mb_{B,i})+\rho(\Mb_{B,j})\leq\Gamma$  and $\#_{\mathrm{box}}(\Mb_{B,i}) + \#_{\mathrm{box}}(\Mb_{B,j}) \leq 2 \Gamma$, 
\end{enumerate}
where $\Mb_{B,i}$ is the connected components of $\Mb_{B}$ in which the particle line $\pb_i$ occurs, and $\Mb_{B}$ is the unique boxed C-molecule which is reduced from $\boldsymbol{\gamma}$. Moreover, $\Mb'_{B}$ is the unique boxed C-molecule reduced from $\boldsymbol{\gamma}'$, and $\boldsymbol{\gamma}'$ is the admissible trajectory that $\mathtt{Pro}(\gamma,\pb_i,\pb_j) = \mathtt{True}$ holds, and $\Mb'_{B,i}$ is the connected component that contains $\pb_i$ in $\Mb'$. Here $\#_{\mathrm{box}}$ is the number of maximal boxes, as in Definition \ref{def.boxed_molecule} \eqref{it.number_box} for boxed C-molecules.  

For $(\vz^0,\Pf^0)$ that is not well-prepared, we define the truncated dynamics by $\mathtt{Pro}( \gamma, \pb_i,\pb_j) = \mathtt{False}$ for all pre-collision time.
\end{definition}
\begin{remark}\label{rem.T_dyn} In the truncated dynamics, we force the property that each cluster has size $\leq\Lambda$ and recollision number $\leq\Gamma$, and the number of boxes $\leq 2 \Gamma$. Each C-atom or C-constituent involves at most $d+1$ particles, and no atom is marked as long interaction, by neglecting any collisions that may violate this property. In case \eqref{it.truncated_dynamic2} of Definition \ref{def.T_dynamics}, we have $\rho(\Mb'_{B,i})=\rho(\Mb_{B,i})+\rho(\Mb_{B,j})$ and $\#_{\mathrm{box}}(\Mb'_{B,i})=\#_{\mathrm{box}}(\Mb_{B,i})+\#_{\mathrm{box}}(\Mb_{B,j})$. Indeed, the new collision $\nf$ joins two different components, so it does not belong to any cycle and does not create any new box. 
    
\end{remark}
\begin{definition}\label{def.projection_map}
    Given $\vz \in \Rb^{2d N}$, define $\Pf_{\mathrm{dist}} = (\Pf_{i})_{i}$, such that $j \in \Pf_{i}$ if and only if $|x_i -x_j| < \eps$. For a function $f(z)$, we define $(\Lb f) (z,\Pf) = f(z) \cdot \mathbbm{1}_{\Pf = \Pf_{\mathrm{dist}}} (\Pf)$, and for a function $f(\vz,\Pf)$, we define $(\Pi f)(\vz) = \sum_{\Pf} f(\vz,\Pf)$. Therefore, we have $\Pi \circ \Lb f(\vz) = f(\vz)$ for any $f= f(\vz)$.
\end{definition}
\begin{proposition}\label{prop.dynamics_E_T_property} The T-dynamics satisfy the following properties.
    \begin{enumerate}
        \item\label{it.dynamics_E_T_property_1} Up to some Lebesgue zero set of $\vz^0$, the E- and T- dynamics in Definitions \ref{def.extended_dynamic} and \ref{def.T_dynamics} are well-defined and produce admissible trajectories in Definition \ref{def.admissible}. 
        \item\label{it.dynamics_E_T_property_2} Define $(\Hc^E,\Sc^E)$ and $(\Hc^T,\Sc^T)$ as special cases of (\ref{eq.defH^M}) and (\ref{eq.defS^M}). Then, when the initial data $\vz^0 \in \Dc^0_N$, the E-dynamics coincides with the O-dynamics. Moreover, $\Pi \Sc^E_N(t) (\Lb f)  = S_N (t) f$ if $f$ is supported in $\Dc^0_N$.
        \item\label{it.dynamics_E_T_property_3}  The $\Hc^E_N$ and $\Sc^E_N$ satisfy the following semi-group property for $t,s\geq 0$:
        \begin{equation}
            \Hc^E_N(t+s) = \Hc^E_N(t)\Hc^E_N(s),\quad \Sc^E_N(t+s) = \Sc^E_N (t) \Sc^E_N(s).
        \end{equation}
        In general, the $\Hc^T_N(t)$ and $\mathcal{S}^T_N(t)$ might not satisfy the flow property for $t,s\geq 0$:
        \begin{equation}\label{eq.semigroup}
            \mathcal{H}_N^T(t+s)\neq \mathcal{H}_N^T(t)\mathcal{H}_N^T(s),\qquad\mathcal{S}_N^T(t+s)\neq \mathcal{S}_N^T(t)\mathcal{S}_N^T(s).
        \end{equation}
        \item\label{it.dynamics_E_T_property_4} Both $\Sc^E_N(t)$ and $\Sc^T_N(t)$ preserve $L^1$ integral on the phase space of $(\vz,\Pf)$:
        \begin{equation}\label{eq.L1pres}
            \int \sum_{\Pf}(\Sc^E_N(t) f) \d x \d v = \int \sum_{\Pf} (\Sc^T_N(t) f) \d x \d v =  \int \sum_{\Pf} f \d x \d v.
        \end{equation}
        \end{enumerate}
\end{proposition}
\begin{proof} The proof is the same as Appendix A of \cite{DHM24}, with the following changes. The main difference is the state is now the set of pairs $(\vz,\Pf)$, so all flow operators act on functions of $(\vz,\Pf)$ and the $L^1$ integral in \eqref{eq.L1pres} is $\int\sum_{\Pf}$. 
    
Note that the force may be discontinuous at distance $\eps$ (see Assumption \ref{ass.phi1}), so the E- and T-dynamics are defined up to a measure zero set of grazing configurations, as in Proposition \ref{prop.shortrange}. For fixed $\Pf$, on each time interval where $\Pf$ does not change, the flow is a smooth Hamiltonian flow and preserves the Lebesgue measure. $\Pf$ only changes when a pair is at distance $\eps$, where the state $\vz$ is continuous. Hence $\Hc^E_N(t)$ and $\Hc^T_N(t)$ are piecewise compositions of measure preserving maps, and they are injective on each piece; summing over the pieces as in \eqref{eq.defS^M} gives \eqref{eq.L1pres}.
    
For \eqref{it.dynamics_E_T_property_2}, if $\vz^0\in\Dc_N^0$ then $\Pf^0=(\{i\})_i$ is trivial, then every pre-collision configuration leads to a collision in the E-dynamics, which is exactly the O-dynamics.
\end{proof}
We now introduce the domain $\mathcal D^{T}_N$, which characterizes the initial data such that the truncated dynamics coincide with the O-dynamics (original dynamics).
\begin{definition}[The truncated and error domain $\mathcal D^{T}_N$ and $\Ec_N^{T}$]\label{def.truncated_domain} Given $(\Lambda,\Gamma,\upsilon_0,\mathsf{N})$, we define the \textbf{truncated domain} $\mathcal D^{T}_N$ as follows:
\begin{equation}\label{eq.truncated_domain} 
       \begin{aligned}
         \mathcal D^{T}_N =  \bigg\{(\vz^0,\Pf^0) : \textrm{ $(\vz^0,\Pf^0)$ is well-prepared, }  \Hc^{T}(t)(\vz^0,\Pf^0) =  \Hc^{E}(t)(\vz^0,\Pf^0)  \ \textrm{for all}\  t \in [0,\tau]     \bigg\}, 
       \end{aligned}
\end{equation} where $\vz^E(t)$ is the dynamics under the E-dynamics. We also define the \textbf{error domain} $\Ec_N^{T}$ as 
\begin{equation}\label{eq.error_domain}
    \Ec_N^{T}:=\Xc_N\backslash\Dc_N^{T}.
\end{equation}
Their indicator functions are denoted by $\mathbbm{1}_{\Dc_N^{T}}$ and $\mathbbm{1}_{\Ec_N^{T}}$, and we have $\mathbbm{1}_{\Ec_N^{T}} = 1-\mathbbm{1}_{\Dc_N^{T}}$. For a truncated dynamics with $\Lambda = \Lambda_\ell$ (defined in Section \ref{sec.notation}), we define the corresponding $T_\ell$-dynamics, $\Hc^{T_\ell}$, $\Sc^{T_\ell}$,  $\Dc^{T_\ell}_N $, and $ \Ec^{T_\ell}_N$ (or simply $\Dc^{\ell}_N $ and $ \Ec^{\ell}_N$), etc.
\end{definition}
 
     \begin{proposition}[Truncated correlation function] Recall the maps $\Pi$ and $\Lb$ in Definition \ref{def.projection_map}. We define the \textbf{truncated correlation function} $\widetilde{f}_s$, and the \textbf{error correlation function} $f_s^{\mathrm{err}}$, by
    \begin{align}\label{eq.truncated_correlation}
        \widetilde{f}_s(\ell\tau, \boldsymbol{z}_s) &= \varepsilon^{s(d-1)}\sum_{n = 0}^{\infty} \frac{1}{n!} \int_{\Rb^{2dn}}\Pi \widetilde{W}_{s+n}(\ell\tau,\boldsymbol{z}_{s+n})\,\mathrm{d}\vz_{[s+1:s+n]},\\
        \label{eq.truncated_error}
        f_s^{\mathrm{err}}(t_{\mathrm{fin}}, \boldsymbol{z}_s) &= \varepsilon^{s(d-1)}\sum_{n = 0}^{\infty} \frac{1}{n!} \int_{\Rb^{2dn}} \Pi W^{\mathrm{err}}_{s+n}(t_{\mathrm{fin}},\boldsymbol{z}_{s+n})  \,\mathrm{d}\vz_{[s+1:s+n]},
    \end{align} 
   where  
    \begin{equation}\label{eq.W_widetilde_l}
        \widetilde{W}_N(\ell\tau,\vz_{s+n},\Pf_{s+n}) =   \bigg(\prod_{\ell'=1}^{\ell}\mathcal{S}^{\ell'}_N(\tau)\bigg) ( \Lb  W_{0,N}),
    \end{equation}
    and
    \begin{equation}\label{eq.W_err}
        W^{\mathrm{err}}_N(t_{\mathrm{fin}}) = \sum_{\ell = 1}^{\mathfrak{L}}\bigg((\Sc_N^E(\tau))^{\Lf-\ell}\circ\mathcal{S}^{\mathrm{err}}_{\ell,N}(\tau)\bigg) \circ  \widetilde{W}_N((\ell-1)\tau)
    \end{equation}
     \label{def.ftrunc} 
    with $\Sc^\ell_N:= \Sc^{T_\ell}_N$, and
    \begin{equation}
        \Sc^{\mathrm{err}}_{\ell,N}:= (\Sc_N^E(\tau)-\Sc^{\ell}_N(\tau))\circ \mathbbm{1}_{\Ec^\ell_N}.
    \end{equation}

    Then we have the following decomposition 
    $$f_s (t_{\mathrm{fin}},\vz_s) =  \widetilde{f}_s(t_{\mathrm{fin}}, \boldsymbol{z}_s) +f_s^{\mathrm{err}}(t_{\mathrm{fin}}, \boldsymbol{z}_s) .$$
     \end{proposition}

     \begin{proof}
        As $\Hc^E = \Hc^\ell$ on the support of $\Dc^\ell_N$, we have $\Sc^E(\mathbbm{1}_{\Dc^\ell}  f) = \Sc^\ell(\mathbbm{1}_{\Dc^\ell} f)$, hence $\Sc^E = \Sc^\ell + \Sc^{\rm err}_{\ell , N}$. Moreover, by Proposition \ref{prop.dynamics_E_T_property} \eqref{it.dynamics_E_T_property_2} and similar calculation as in Section 4.3 in \cite{DHM24}, we have the following decomposition
        \begin{equation}
            W_N(t_{\mathrm{fin}}) = \Sc_N(t_{\mathrm{fin}}) W_{0,N} = \Pi \Sc^E_N(t_{\mathrm{fin}})(\Lb W_{0,N})= \Pi \widetilde{W}_N (t_{\mathrm{fin}}) + \Pi W^{\mathrm{err}}_N(t_{\mathrm{fin}}) . 
        \end{equation} 
     \end{proof}
\subsection{Statement of the cumulant formula and proof of the cumulant formula}\label{sec.formula_cumulant_statement} In this section, we provide a rigorous statement of the cumulant formulas, see Proposition \ref{prop.cumulant_formula}, and the proof of Proposition \ref{prop.cumulant_formula}.
We first define the following sets of molecules, which are needed in Proposition \ref{prop.cumulant_formula}.
\begin{definition}[Set of the boxed molecules $\Tc$ and $\Fc$]\label{def.set_T_F} Given $\ell\leq\Lf$, $\vLambda_{\ell}:=(\Lambda_{\ell'})_{\ell'\leq\ell}$, $\Gamma$ and $\mathsf{N}$, we define the following sets of molecules. Recall the notions in Definition \ref{def.molecule} and the sets $p(\Mb_{\ell'})$ and $r(\Mb_{\ell'})$ in Definition \ref{def.prsets}.
\begin{enumerate}
    \item\label{it.set_F} \emph{The set $\mathcal{F}_{\boldsymbol{\Lambda}_\ell}$.} Define $\mathcal{F}_{\boldsymbol{\Lambda}_\ell}$ to be the boxed molecules $\Mb$ that satisfy the following properties:
    \begin{enumerate}
        \item\label{it.set_F_l_1} $\Mb$ has the lowest layer $\underline{\ell}=1$ and the highest layer $\ell$, with $|r(\Mb_{\ell'})|\leq A_{\ell'}$ for each $\ell'\in [1:\ell]$. Each component of $\Mb_{\ell'}$ contains $<\Lambda_{\ell'}$ clusters (see Definition \ref{def.cluster}), and for each cluster $\Mb_j$ we have $|\Mb_j|_p\le \Lambda_{\ell'}$ and $\rho(\Mb_j)\le \Gamma$ (Definition \ref{def.more_on_molecule} and Definition \ref{def.boxed_molecule} \eqref{it.recollision_number_box}). Each C-atom and C-constituent contains at most $d+1$ particles, and each box contains at most $\mathsf{N}$ particles. Each C-atom or C-constituent is marked as short interaction, and each cluster contains at most $2\Gamma$ boxes.
        \item\label{it.set_F_l_crossing} $\Mb_1$ has no bottom crossing atom. For $\ell'\geq2$, each bottom crossing atom of $\Mb_{\ell'}$ contains two particle lines of $r(\Mb_{\ell'-1})$ that belong to the same top crossing atom of $\Mb_{\ell'-1}$, i.e. it is the continuation of this collision (Definition \ref{def.top_reduction}); a top crossing atom of $\Mb_{\ell'-1}$ need not be continued in $\Mb_{\ell'}$.
        \item\label{it.set_F_l_2} For each $\ell'\in[1:\ell]$, each particle line in $p(\Mb_{\ell'})$ is either connected to a particle line in $r(\Mb_{\ell'})$ via $\Mb_{\ell'}$, or belongs to $r(\Mb_{\ell'})$. Each component of $\Mb_{\ell'}$ contains a particle line in $r(\Mb_{\ell'})$. 
        \item\label{it.set_F_l_3} The cluster graph (Definition \ref{def.cluster}) of each $\Mb_{\ell'}$ is a forest. The initial link graph, formed by initial link conditions between particle lines $\pb\in p(\Mb_1)$, is also a forest.
        \item\label{it.set_F_l_4} For $\ell'\in[2:\ell]$, each particle line in $r(\Mb_{\ell'})$ is either connected to another particle line in $r(\Mb_{\ell'})$ via $\Mb_{\ell'}$, or connected to a particle line in $r(\Mb_{\ell'-1})$ via $\Mb_{\ell'}$ (in the sense of Definition \ref{def.connectedvia}), or belongs to $r(\Mb_{\ell'-1})$. For $\ell'=1$, the set $r(\Mb_{0})$ in the above definition should be replaced by the set of particle lines that are involved in an initial link.
    \end{enumerate}
        \item\label{it.set_F_err} \emph{The set $\mathcal{F}_{\boldsymbol{\Lambda}_\ell}^{\mathrm{err}}$.} $\mathcal{F}_{\boldsymbol{\Lambda}_\ell}^{\mathrm{err}}$ is the set of boxed molecule $\Mb$ as in (\ref{it.set_F}) that satisfy conditions \eqref{it.set_F_l_1}--\eqref{it.set_F_l_4}, except that \eqref{it.set_F_l_4} is assumed only for $\ell'\leq\ell-1$. Moreover, for $\ell'=\ell$, instead of  \eqref{it.set_F_l_1} we assume that each component of $\Mb_{\ell}$ contains $\leq\Lambda_{\ell}$ clusters and there exists a unique component that contains \emph{exactly} $\Lambda_{\ell}$ clusters.
    \item\label{it.set_T} \emph{The set $\mathcal{T}_{\Lambda_\ell}$ and $\mathcal{T}_{\Lambda_\ell}^{\mathrm{err}}$.} Given $\Lambda_\ell$, define $\mathcal{T}_{\Lambda_\ell}$ to be the set of boxed molecules $\Mb$  such that it has only one layer $\ell$ (i.e. $\underline{\ell}=\ell$) and no initial link, has only one root particle labeled by $1$, and satisfies \eqref{it.set_F_l_1}--\eqref{it.set_F_l_3}. Moreover, we require that $\Mb$ contains no bottom crossing atom. (The root particle line may end at a top crossing atom; the other particles of such an atom are non-root particle lines.) 
    
    Define also $\mathcal{T}_{\Lambda_\ell}^{\mathrm{err}}$ to be the set of $\Mb$ such that  $\Mb$ is obtained from a molecule in $\Tc_{\Lambda_{\ell}}$ by joining two edges or particles of a C-atom or C-constituent $\af_1$ and $\af_2$ with an extra O-connection.
    \item\label{it.set_F_err_2} \emph{The set $\mathcal{F}_{\boldsymbol{\Lambda}_\ell}^{\mathrm{trc.err}}$.} Define $\mathcal{F}_{\boldsymbol{\Lambda}_\ell}^{\mathrm{trc.err}}$ to be the set of boxed molecules $\Mb$ such that it has $\ell$ layers and satisfies (\ref{it.set_F_l_1})--(\ref{it.set_F_l_4}) for $\ell'\leq \ell-1$. Moreover, the layer $\Mb_{\ell}$ is connected (apart from possible empty ends which are root particle lines) and has no O-connection, and satisfies that 
    \begin{equation}\label{eq.trc_err}
        |\Mb_{\ell}|_p\leq 2 \Lambda_{\ell},\ \mathrm{and}\ \rho(\Mb_{\ell})\leq 2\Gamma, \ \mathrm{and}\ \#_{\mathrm{box}}(\Mb_\ell) \leq 4 \Gamma, 
    \end{equation}
     moreover, either $|\Mb_{\ell}|_p>\Lambda_{\ell}$  or $\rho(\Mb_{\ell})>\Gamma $  or $\#_{\mathrm{box}}(\Mb_\ell) > 2\Gamma $ or there exists a C-atom with more than $d+1$ particle lines, or a C-atom marked as long interaction. In any case, each C-atom of $\Mb_\ell$ has at most $2(d+1)$ particle lines, and at most one C-atom of $\Mb_\ell$ has more than $d+1$ particle lines or is marked as long interaction.
\end{enumerate}
\end{definition}

 We  define the following objects.

\begin{definition}[Partial dynamics]\label{def.partial_dynamics} Recall $\vz=\vz_N=(z_1,\cdots,z_N)$ as in Definition \ref{def.notation} (\ref{it.vector}). Let $\lambda\subseteq[N]$ and $\vz_\lambda=(z_j)_{j\in\lambda}$. Given initial configuration $\vz(0)=\vz^0$ and $\Pf^0$, let $\vz(t)$ be a modified dynamics, then we can define the following notions.
\begin{enumerate}
    \item\label{it.dynamics_partial} \emph{The partial dynamics $\vz_\lambda(t)$.} Consider $\lambda\subseteq[N]$ as a set of particles. We define the \textbf{partial dynamics} $\vz_\lambda(t)$ to be the dynamics following the same evolution rule but only taking into account particles in $\lambda$, and with initial configuration $\vz_\lambda^0=(z_j^0)_{j\in\lambda}$, with $\Pf_{\lambda,i}^0 = \{j\in \lambda  \mid j \in \Pf^0_i\}$, and $\Pf_\lambda^0 := (\Pf_{\lambda,i}^0)_{i \in \lambda}$. When $\lambda = [N]$, we refer to it as the \textbf{full dynamics}.
    \item\label{it.partial_flow_map} \emph{The partial flow maps $\Hc_\lambda(t)$ and $\Sc_\lambda(t)$.} We define the \textbf{partial flow maps} $\Hc_\lambda(t)$ by $\Hc_\lambda(t)(\vz_\lambda^0)=\vz_\lambda(t)$. It can be extended to $(\vz^0,\Pf^0)$ with $\vz^0 = \vz^0_N$ by $\Hc_\lambda(t)(\vz^0,\Pf^0) = \Hc_\lambda(t)(\vz_\lambda^0,\vz^0_{\lambda^c},\Pf_\lambda(t), \Pf^0 |_{\lambda^c})$, where $\Pf_{\lambda^c}:= (\Pf_i \cap \lambda^c)_{i \in \lambda^c}$. We define the \textbf{partial flow operators} $\Sc_\lambda(t)$ as follows, 
    \begin{equation}\label{eq.transport_op_partial}
    \begin{gathered}
        (\Sc_\lambda(t)f)(\vz_\lambda,\Pf_\lambda) = \sum_{(\vz_\lambda^0, \Pf_{\lambda}^0):\Hc_\lambda(t)(\vz_\lambda^0,\Pf^0_\lambda)=(\vz_\lambda,\Pf_\lambda)}f(\vz_\lambda^0,\Pf_\lambda^0),
        \\
        (\Sc_\lambda(t)f)(\vz,\Pf)  =\sum_{(\vz^0,\Pf^0):\Hc_\lambda(t)(\vz^0,\Pf^0)=(\vz,\Pf)}f(\vz^0,\Pf_\lambda^0).
    \end{gathered}
    \end{equation}
    Finally, when $\lambda=[s]$, we will denote $\vz_\lambda=\vz_s$, and similarly $\Hc_\lambda=\Hc_s$ and $\Sc_\lambda=\Sc_s$.
    \item \emph{Composition of partial dynamics.} Given disjoint subsets $(\lambda_1, \cdots, \lambda_k)$, the \textbf{composition} of the partial dynamics of these subsets is the  dynamics specified by the flow map $\prod_{j=1}^k\Hc_{\lambda_j}(t)$.
\end{enumerate}
\end{definition}

\begin{definition}
    Let $\Mb_B$ and $\Mb'_B$ be boxed C-molecules with only one layer $\ell$, and $\Mb$, $\Mb'$ be C-molecules before they reduce to boxed molecules with only one layer $\ell$. Let $\vz^T(t)$ be the truncated dynamics on $[(\ell-1)\tau,\ell\tau]$ and $\vz^0 = \vz^T((\ell-1)\tau)$ be the initial configuration. 
        \begin{enumerate}
            \item  We define the indicator function $\mathbbm{1}_{\Mb}^T (\vz^0,\Pf^0)$ to be 1 if and only if the trajectory $\vz_{p(\Mb)}^T(t)$ topologically reduces to $\Mb$. We define the indicator function $\mathbbm{1}_{\Mb_B}^T (\vz^0,\Pf^0)$ to be
            \begin{equation}\label{eqn.1_Mb_B_T}
                \mathbbm{1}_{\Mb_B}^T (\vz^0,\Pf^0) = \sum_{\Mb \textrm{ that can be reduced into } \Mb_B}\mathbbm{1}_{\Mb}^T (\vz^0,\Pf^0).
            \end{equation}  
            \item We write $\Mb_B \nsim \Mb_B'$ if and only if there exists no $j \in p(\Mb_B)$ and $i \in p(\Mb_B')$, such that $j \in \Pf^0_i$, and in the trajectory $\vz^T_{p(\Mb_B) \cup p(\Mb_B')}(t)$, there exists no collisions between any particles in $p(\Mb_B)$ and $p(\Mb_B')$. Define $\mathbbm{1}_{\Mb_B\nsim \Mb'_B}(\vz^0,\Pf^0)$ to be 1 if and only if $\Mb_B \nsim \Mb_B'$. 
         
            \item Denote $\af$ and $\af'$ to be edges or particles of some atoms in $\Mb_B$ and $\Mb_B'$, then for each trajectory that is reduced to $\Mb_B$ or $\Mb_B'$, $\af$ and $\af'$ correspond a unique part of the trajectory (either a linear segment, or a curved trajectory which is involved in a collision or a trajectory of some particle in an M-atom). 
            
            We write $(\Mb_B,\af)\nsim (\Mb'_B,\af')$, if (i) when $\af,\af'$ contain the states of particles $i,j$ at time $(\ell-1)\tau$, we have $i\notin \Pf_j((\ell-1)\tau)$; and (ii) under the composition dynamics $\mathcal{H}_{p(\Mb_B)}^T \circ \mathcal{H}_{p(\Mb'_B)}^T$ either (a) the trajectories corresponding to $\af$ and $\af'$ never enter a pre-collisional configuration at any time $t_*$; or (b) $\af$ and $\af'$ do enter a pre-collisional configuration under the dynamics $\mathcal{H}_{p(\Mb_B)}^T \circ \mathcal{H}_{p(\Mb'_B)}^T$ at time $t_*$, but $\mathtt{Pro}(\gamma,\pb,\pb') = \mathtt{False}$. Here $\pb,\pb'$ are the particles containing $\af$ and $\af'$ respectively, $\gamma$ is the trajectory of the $\mathcal{H}_{p(\Mb_B)}^T \circ \mathcal{H}_{p(\Mb'_B)}^T$ dynamics up to time $t_*$, and $\mathtt{Pro}$ is the condition in the truncated dynamics (Definition \ref{def.T_dynamics}) with particle set $p(\Mb_B) \cup p(\Mb'_B)$.
            
            We define the indicator function $\mathbbm{1}_{(\Mb_B,\af)\nsim (\Mb'_B,\af')}(\vz^0)$ to be $1$ if $(\Mb_B,\af)\nsim (\Mb'_B,\af')$, and $0$ otherwise. Note that these objects are well-defined only when the trajectory of particles in $p(\Mb_B)$ and $p(\Mb_B')$ can be reduced to $\Mb_B$ and $\Mb_B'$.
        \end{enumerate}
\end{definition}
Below, we define the following dynamics prescribed by boxed molecules. This definition differs from the one in \cite{DHM24}, due to the possible complicated trajectory in one atom.
\begin{definition}[Boxed Molecule prescribed dynamics]\label{def.molecule_truncated_dynamics} We define the boxed molecule described dynamics as follows, depending on the type of molecules.
     \begin{enumerate}
    \item\label{it.molecule_prescribed_molecule1} Given a boxed C-molecule $\Mb_B$ with only one layer $\ell$, which is reduced from some truncated trajectory $\vz^T_{p(\Mb_B)(t)}$, with $\Pf_{p(\Mb_B)(t)}$. If $\{\Mb_{B,1},\cdots,\Mb_{B,k}\}$ are connected components of $\Mb_B$, the dynamics, with particle set being $p(\Mb)$, \textbf{prescribed by $\Mb_B$} is defined to be 
        $$(\vz_{\Mb_B},\Pf_{p(\Mb_B)}) (t+(\ell-1)\tau)= \Hc_{\Mb_B}(t)(\vz^0,\Pf^0)= 
         \prod_{j=1}^{k} \Hc^T_{p(\Mb_{B,j})}(t)(\vz^0,\Pf^0),  \quad t \in[0,\tau], $$
         where $(\vz^0,\Pf^0)$ denotes the state of particles at time $(\ell - 1)\tau$.
        
    If $\Mb_B$ contains O-connections, and each cluster $\Mb_{B,j}$ in $\Mb_B$ is reduced from some trajectory $\vz^T_{\Mb_{B,j}}(t)$, the dynamics prescribed by $\Mb_B$ is defined to be 
     \begin{equation}
       (\vz_{\Mb_B},\Pf_{p(\Mb_B)})(t+(\ell-1)\tau)= \Hc_{\Mb}(t)(\vz^0,\Pf^0) = \prod_{j=1}^{k} \Hc^T_{p(\Mb_{B,j})}(t)(\vz^0,\Pf^0), \quad t \in[0,\tau] ,
    \end{equation}
     where $(\vz^0,\Pf^0)$ denotes the state of particles at time $(\ell - 1)\tau$.
    
    If $\Mb_B$ is a boxed molecule containing multiple layers $[\underline{\ell},\ell]$, and for every $\ell'\in [\underline{\ell},\ell] $, we have $\Mb_{\ell'}$ satisfies the above assumption, then the dynamics prescribed by $\Mb_B$ with initial state $(\vz^0,\Pf^0)$ at time $(\underline{l}-1)\tau$, is
      \begin{equation}\label{eq.multi_layer_flow}
        \Hc_{\Mb_B}((\ell - \underline{\ell}+1)\tau) := \Tc_2 \circ  \Hc_{\Mb_{\ell}}(\tau)\circ \cdots \circ \Hc_{\Mb_{\underline{\ell}}}(\tau) \circ \Tc_1,  
      \end{equation}
    and similarly for intermediate times. Here $\Hc_{\Mb_{\ell'}}$ acts on $(\vz_{p(\Mb_{\ell'})},\Pf|_{p(\Mb_{\ell'})})$ (Definition \ref{def.partial_dynamics}) and sets $\Pf_\pb=\{\pb\}$ for $\pb\notin p(\Mb_{\ell'})$; $\Tc_1$ moves each $\pb$ linearly from $(\underline{\ell}-1)\tau$ to $(\ell_1[\pb]-1)\tau$, and $\Tc_2$ moves each $\pb$ linearly from $\ell_2[\pb]\tau$ to $\ell\tau$. Both $\Tc_1$ and $\Tc_2$ only act on trivial interaction sets: $\Pf^0$ is trivial, and for $\ell_2[\pb]<\ell$, $\Pf_\pb$ is made trivial by $\Hc_{\Mb_{\ell_2[\pb]+1}}$.

    Thus the layers are connected by restriction: at time $\ell'\tau$ the pairs containing a particle line that finishes at layer $\ell'$ are removed (e.g. a root and a non-root at a top crossing atom), and no pair is added. So the trajectory is admissible (Definition \ref{def.admissible}), and a bottom crossing atom of $\Mb_{\ell'+1}$ is the continuation of a top crossing atom of $\Mb_{\ell'}$.
    
       \item\label{it.molecule_prescribed_molecule2}  Let $\Mb_B\in \mathscr{F}^{\mathrm{trc.err}}_{\Lambda_\ell}$, where $\mathscr{F}^{\mathrm{trc.err}}_{\Lambda_\ell}$ is the collection of all boxed molecules $\Mb$ of a single layer $\ell$ that is connected apart from possible empty ends (which are root particle lines), has only C-atoms, and satisfies the conditions of $\Mb_\ell$ in Definition \ref{def.set_T_F} \eqref{it.set_F_err_2}. We define the dynamics \textbf{prescribed by $\Mb_B$} to be
        $$ \Hc_{\Mb_B}(t)(\vz^0,\Pf^0) \widetilde{\Hc}^T_{p(\Mb_{B })}(t)(\vz^0,\Pf^0) \quad t \in[0,\tau],  $$
         where $\vz^0$ denotes the state of particles at time $(\ell - 1)\tau$ and $\widetilde{\Hc}^T_{p(\Mb_{B })}$ is defined as follows:
         \begin{enumerate}
            \item $(\vz^0,\Pc^0)$ is not well-prepared, then $\Pc_i^0$ is defined as in Definition \ref{def.well-prepared}, and $\mathtt{Pro} (\gamma,\pb_1,\pb_2)= \mathtt{False}$ at every pre-collision time.
            \item $\vz^0$ is well-prepared, $\widetilde{\Hc}^T_{p(\Mb_{B })}$ is defined to be the same as ${\Hc}^T_{p(\Mb_{B })}$, but with $\mathtt{Pro} (\gamma,\pb_1,\pb_2)= \mathtt{True}$ at the first time that $\mathtt{Pro}(\gamma,\pb_1,\pb_2)= \mathtt{False}$ in $\vz^T(t)$, and all $\mathtt{Pro}(\gamma,\pb_1,\pb_2) = \mathtt{False}$ for every pre-collision time after that.
         \end{enumerate} 
\end{enumerate}
\end{definition}
    With the above definitions, we can obtain the following lemma, which corresponds to Lemma 5.5 in \cite{DHM24}.
   \begin{lemma}\label{lem.factorization} We have the following decompositions.
      \begin{enumerate}
         \item\label{it.factorization1} 
         We have
         \begin{equation}\label{eqn.sum_MbB}
             1 = \sum_{\Mb_B} \mathbbm{1}_{\Mb_B}^T
         \end{equation}
         where the summation is taken over all boxed C-molecules $\Mb_B$ with $p(\Mb_B) = [N]$ that satisfy the truncated property in Definition \ref{def.T_dynamics}.
         \item\label{it.factorization2} If $\{\Mb_{B,1},\dots,\Mb_{B,k}\}$ are the connected components of $\Mb_B$, then we have
         \begin{equation}\label{eqn.factorization_1M} 
                 \mathbbm{1}^T_{\Mb_B} = \big(\prod_{j=1}^k \mathbbm{1}^T_{\Mb_{B,j}} \big)\cdot\big(\prod_{1\leq j<j' \leq k} \mathbbm{1}_{\Mb_{B,j} \nsim \Mb_{B,j'}}\big),
                \end{equation}
         In particular, on $\supp \mathbbm{1}^T_{\Mb_B}$, every $j \in \Pf^0_i$ implies particles $i,j$ belong to the same component. Moreover, if $f$ is supported on $\mathbbm{1}_{\Mb_B}^T $,
                \begin{equation}\label{eqn.factorization_S}
                 \Sc^T_N(t)f   = \prod^k_{j=1} \Sc_{\Mb_{B,j}}(t)f .
         \end{equation}
         
         \item\label{it.factorization3} In the set $(\supp \mathbbm{1}^{T}_{\Mb_B})\cap (\supp \mathbbm{1}^T_{\Mb'_B})$, we have 
         $$\mathbbm{1}_{\Mb_B\nsim \Mb_B'} = \prod_{\af \in \Mb_B,\af' \in \Mb_B'}\mathbbm{1}_{(\Mb_B,\af)\nsim (\Mb_B',\af')} =  \prod_{\af \in \Mb_B,\af' \in \Mb_B'}(1 - \mathbbm{1}_{(\Mb_B,\af)\sim (\Mb_B',\af')}),$$
         where $\af,\af'$ are edges in $\Mb_B$ and $\Mb_B'$, or particle lines contained in some atoms in $\Mb_B$ and $\Mb_B'$. 
        \item\label{it.factorization4} By putting together \eqref{it.factorization1}--\eqref{it.factorization3}, we get
        \begin{equation}
            \Sc^T_N(t) = \sum_{k\leq N} \sum_{\Mb_{B,1},\dots,\Mb_{B,k}} \big(\prod^k_{j=1} \Sc_{\Mb_{B,j}}(t)\circ \mathbbm{1}^T_{\Mb_{B,j}}\big)\circ\big( \prod_{1\leq j < j' \leq k}\prod_{\af \in \Mb_B,\af' \in \Mb_B'}\mathbbm{1}_{(\Mb_{B,j},\af)\nsim (\Mb_{B,j'},\af')} \big).
        \end{equation}
     \end{enumerate}
   \end{lemma}
     \begin{proof}[Proof of Lemma \ref{lem.factorization}]
         The equation \eqref{eqn.sum_MbB} holds as each admissible trajectory reduces to a unique boxed C-molecule, see Remarks \ref{rem.top_reduction_uniqueness} and \ref{rem.uniqueness_boxed}. The proof of  \eqref{it.factorization3} follows by the same argument as in Lemma 5.5 \cite{DHM24} (see equation (5.9) there), with one difference: a pre-collision configuration may occur while a particle is inside an atom, so the product runs over pairs $(\af,\af')$ of edges and particle pieces of atoms; and \eqref{it.factorization4} follows from putting \eqref{it.factorization1}--\eqref{it.factorization3} together. Note that \eqref{eqn.factorization_S} directly follows by our definition of the boxed molecule prescribed dynamics in Definition \ref{def.molecule_truncated_dynamics}.
        
         Now we prove \eqref{eqn.factorization_1M}. We first note that if $j \in \Pf^0_i$ but particles $i,j$ belong to different components, then both side vanishes, so we may assume that $\Pf^0 = (\Pf |_{p(\Mb_B,j)})_j$.
         If $\mathbbm{1}^T_{\Mb_B}(\vz^0,\Pf^0) = 1$, then $\vz^T_{p(\Mb_B)}(t)$ must topologically reduce to some molecule $\Mb=(\Mb_1,\dots,\Mb_k)$ which can be further reduced into $\Mb_B$. Hence, $\mathbbm{1}^T_{\Mb}(\vz^0,\Pf^0) = 1$. By the same proof of Lemma 5.6 in \cite{DHM24}, we have 
         $$\mathbbm{1}^T_{\Mb}(\vz^0,\Pf^0) = \big(\prod_{j=1}^k \mathbbm{1}^T_{\Mb_{j}} (\vz^0,\Pf^0) \big)\cdot\big(\prod_{1\leq j<j' \leq k} \mathbbm{1}_{\Mb_{j} \nsim \Mb_{j'}} (\vz^0,\Pf^0)\big).$$
          Note that $\mathbbm{1}_{\Mb_{B,j} \nsim \Mb_{B,j'}}(\vz^0,\Pf^0) =  \mathbbm{1}_{\Mb_j \nsim \Mb_{j'}}(\vz^0,\Pf^0) $, where $\Mb_j$ and $\Mb_{j'}$ are the topological reduction of $\vz^T_{p(\Mb_j)}(t)$ and $\vz^T_{p(\Mb_{j'})}(t)$ with initial configuration $\vz^0$. We have
        \begin{equation}
             \begin{aligned}
                1 = &\mathbbm{1}^T_{\Mb}(\vz^0,\Pf^0) =   \big(\prod_{j=1}^k \mathbbm{1}^T_{\Mb_{j}}(\vz^0,\Pf^0)\big)\cdot\big(\prod_{1\leq j<j' \leq k} \mathbbm{1}_{\Mb_{j} \nsim \Mb_{j'}}(\vz^0,\Pf^0)\big)\\
                \leq & \big( \prod_{j=1}^k \mathbbm{1}^T_{\Mb_{B,j} }(\vz^0,\Pf^0)\big)\cdot\big(\prod_{1\leq j<j' \leq k} \mathbbm{1}_{\Mb_{B,j} \nsim \Mb_{B,j'}}(\vz^0,\Pf^0)\big).
             \end{aligned}
         \end{equation}
        This implies that (as both sides only take values $0$ and $1$)
        \begin{equation}\label{eqn.proof.fact1}
             \mathbbm{1}^T_{\Mb_B} \leq \big(\prod_{j=1}^k \mathbbm{1}^T_{\Mb_{B,j}}\big)\cdot\big(\prod_{1\leq j<j' \leq k} \mathbbm{1}_{\Mb_{B,j} \nsim \Mb_{B,j'}}\big).
        \end{equation}
        On the other hand, if the right-hand side of \eqref{eqn.factorization_1M} equals $1$ at $(\vz^0,\Pf^0)$, then we must have $\mathbbm{1}^T_{\Mb_{j}}(\vz^0,\Pf^0) = 1$ where $\Mb_j$ can be reduced into $\Mb_{B,j}$ and $\big(\prod_{1\leq j<j' \leq k} \mathbbm{1}_{\Mb_{j} \nsim \Mb_{j'}}(\vz^0,\Pf^0)\big) = 1$. Hence, we have $\mathbbm{1}^T_{\Mb}(\vz^0,\Pf^0) =1$ for $\Mb = (\Mb_1,\cdots,\Mb_k)$. Moreover, $\Mb$ further reduces to $\Mb_B$ as $\Mb_j$ reduces to $\Mb_{B,j}$ for each $j$. Hence, we have $\mathbbm{1}^T_{\Mb_B}(\vz^0,\Pf^0)=1$ and 
        \begin{equation}\label{eqn.proof.fact2}
             \mathbbm{1}^T_{\Mb_B} \geq \big(\prod_{j=1}^k \mathbbm{1}^T_{\Mb_{B,j}}\big)\cdot\big(\prod_{1\leq j<j' \leq k} \mathbbm{1}_{\Mb_{B,j} \nsim \Mb_{B,j'}}\big).
        \end{equation}
        Combining \eqref{eqn.proof.fact1} and \eqref{eqn.proof.fact2}, we obtain \eqref{eqn.factorization_1M}.

    \end{proof}
   \begin{remark} We have the following remarks on Lemma \ref{lem.factorization}.
     \begin{enumerate}
         \item In Lemma 5.5 of \cite{DHM24}, the authors used the decomposition
         \begin{equation}\label{eqn.sum_Mb}
             1 = \sum_{\Mb} \mathbbm{1}_{\Mb}^{\Lambda,\Gamma},
         \end{equation} 
         where the summation is taken over all C-molecules $\Mb$ with $p(\Mb) = N$ that satisfy the truncated property in \cite{DHM24}. But in this paper, we use the decomposition in \eqref{eqn.sum_MbB}. This means that we do not further perform the decomposition
         $$\mathbbm{1}_{\Mb_B}^T = \sum_{\Mb \text{ that can be reduced to } \Mb_B } \mathbbm{1}_{\Mb}^T.$$
         This choice leads to fewer possibilities in the summation and thus reduces the combinatorial loss in Proposition \ref{prop.layerrec3} below.
         \item Note that in (5.5) \cite{DHM24}, it shows that in the set $(\supp \mathbbm{1}^T_{\Mb})\cap (\supp \mathbbm{1}^T_{\Mb'})$
         $$\mathbbm{1}_{\Mb\nsim \Mb'} = \prod_{e \in \Mb,e' \in \Mb'}\mathbbm{1}_{(\Mb,e)\nsim (\Mb',e')} = \prod_{e \in \Mb,e' \in \Mb'}(1 - \mathbbm{1}_{(\Mb,e)\sim (\Mb',e')}),$$
         for edges $e$ and $e'$. But here, due to the possibility of multi-particle collisions, in the set $(\supp \mathbbm{1}^{T}_{\Mb_B})\cap (\supp \mathbbm{1}^T_{\Mb'_B})$, we have 
         $$\mathbbm{1}_{\Mb_B\nsim \Mb_B'} = \prod_{\af \in \Mb_B,\af' \in \Mb_B'}\mathbbm{1}_{(\Mb_B,\af)\nsim (\Mb_B',\af')} =  \prod_{\af \in \Mb_B,\af' \in \Mb_B'}(1 - \mathbbm{1}_{(\Mb_B,\af)\sim (\Mb_B',\af')}),$$
         where $\af,\af'$ are edges in $\Mb_B$ and $\Mb_B'$, or particles of some atoms in $\Mb_B$ and $\Mb_B'$.  Therefore, after the Penrose argument as in Section 4.2 \cite{DHM24}, it is possible that two or more atoms are joined by several O-connections. Therefore, it is possible that a CO-atom may have many particle lines, even each C-constituent has at most $d+1$ particles by the truncated dynamics.
     \end{enumerate}
   \end{remark}
  
\begin{definition}[Associated transport operators]\label{def.associated_op_nonlocal} Let $\Mb_B$ be a boxed molecule. We define the following objects:
\begin{enumerate}
    \item\label{it.prescribed_flow} \emph{Prescribed dynamics and flow operators.} We use $\vz_{\Mb_B}(t)$ to denote the $\Mb_B$-prescribed dynamics in Definition \ref{def.molecule_truncated_dynamics}, with $\Hc_{\Mb_B}(t)$ and $\Sc_{\Mb_B}(t)$ etc. defined accordingly as in Definition \ref{def.modified_dynamics}. Note that the subscript $N$ is omitted, as $N=|\Mb|_p$ can be inferred from $\Mb$.
    \item\label{it.associated_ind} \emph{Indicator functions $\mathbbm{1}_{\Mb_B}$ and $\mathbbm{1}_{\Lc}^\eps$.}  We define the indicator function $\mathbbm{1}_{\Mb_B}(\vz^0,\Pf^0)$ such that
    \begin{equation}\label{eq.ind_1MB}
        \begin{aligned}
            \mathbbm{1}_{\Mb_B}(\vz^0,\Pf^0)=1  \Leftrightarrow \textrm{the trajectory $\vz_{\Mb_B}(t)$, with well-prepared initial configuration $(\vz^0,\Pf^0)$, reduces to $\Mb_B$.}
        \end{aligned}
    \end{equation} 

    Define also the indicator function $\mathbbm{1}_{\Lc}^\varepsilon$, associated with the initial link structure $\Lc$, by 
    \begin{equation}\label{eq.1_L}
        \mathbbm{1}_\Lc^{\varepsilon} = \prod_{(e_1, e_2)\in \Lc} \mathbbm{1}_{e_1\sim_{\mathrm{in}} e_2}^\varepsilon\qquad\textrm{and}\qquad \mathbbm{1}_{e_1\sim_{\mathrm{in}} e_2}^\varepsilon = \left\{\begin{aligned}
            &1 \quad &&|x_{e_1}-x_{e_2}|\le \varepsilon,
            \\
            &\varepsilon^\upsilon \quad &&|x_{e_1}-x_{e_2}|\ge \varepsilon,
        \end{aligned}\right.
    \end{equation} 
    where $\upsilon$ is defined in Section \ref{sec.notation}.
    \item \emph{The operator $(\Sc\circ\mathbbm{1})_{\Mb_B}$ and $|(\Sc\circ\mathbbm{1})_{\Mb_B}|$.} Define the \textbf{associated transport operators} $(\Sc\circ\mathbbm{1})_\Mb$ and $|(\Sc\circ\mathbbm{1})_\Mb|$ by 
    \begin{equation}\label{eq.associated_op}
        (\Sc\circ\mathbbm{1})_{\Mb_B} = (-1)^{|\Mb_B|_{\mathrm{O}}}\cdot\Sc_{\Mb_B}((\ell-\underline{\ell}+1)\tau)\circ\mathbbm{1}_{\Mb_B}\quad\text{and}\quad|(\Sc_{\Mb_B}\circ\mathbbm{1})_{\Mb_B}| = \Sc_{\Mb_B}((\ell-\underline{\ell}+1)\tau)\circ \mathbbm{1}_{\Mb_B},
    \end{equation} where $\ell$ is the highest layer and $\underline{\ell}$ is the lowest layer of $\Mb_B$, and recall $|\Mb_B|_{\mathrm{O}}$ is the number of O-connections in $\Mb_B$, see Definition \ref{def.more_on_molecule} \eqref{it.number_atom}.
\end{enumerate}
\end{definition}

\begin{definition}[Associated integrals $\Ic\Nc_\Mb$ and $|\Ic\Nc_\Mb|$]\label{def.associated_int} Let $\Mb_B$ be a boxed molecule, let the lowest layer of $\Mb_B$ be $\underline{\ell}$. Define the \textbf{associated integrals} $\Ic\Nc_{\Mb_B}$ and $|\Ic\Nc_{\Mb_B}|$ as follows:
    \begin{align}\label{eq.associated_integral_molecule}
        \Ic\Nc_{\Mb_B}(\vz_{r({\Mb_B})}) &\coloneqq \varepsilon^{-(d-1)|{\Mb_B}|_{p\backslash r}}\int_{\Rb^{2d|{\Mb_B}|_{p\backslash r}}} \Pi (\Sc\circ\mathbbm{1})_{{\Mb_B}}\, (Q\cdot \mathbbm{1}_{\{\Pf \ \rm trivial\}}) \,\mathrm{d}\vz_{(p\backslash r)({\Mb_B})},\\\label{eq.associated_integral_molecule_abs}
        |\Ic\Nc_{\Mb_B}|(\vz_{r({\Mb_B})}) &\coloneqq \varepsilon^{-(d-1)|{\Mb_B}|_{p\backslash r}}\int_{\Rb^{2d|{\Mb_B}|_{p\backslash r}}} \Pi |(\Sc\circ\mathbbm{1})_{{\Mb_B}}|\, (|Q| \cdot \mathbbm{1}_{\{\Pf \ \rm trivial\}})\,\mathrm{d}\vz_{(p\backslash r)({\Mb_B})},
    \end{align}
    where  $\vz_{(p\backslash r)({\Mb_B})}:=(\vz_\pb)_{\pb\in (p\backslash r)({\Mb_B})}$, $\Pi$ is defined in Definition \ref{def.projection_map}, $\mathbbm{1}_{\Pf\ \rm trivial}= 1$ if and only if $\Pf_i =\{i\}$ for each $i$, and
    \begin{equation}\label{eq.extraQ}
    Q:=\mathbbm{1}^\eps_{\Lc} \prod_{\pb\in p({\Mb_B})} f^\Ac((\ell_1[\pb]-1)\tau,x_\pb'+(\ell_1[\pb]-\underline{\ell})\tau\cdot v_\pb',v_\pb'),
    \end{equation}
where $\vz_\pb= (x_\pb,v_\pb)$ and $\vz_{\pb}' = (x_\pb',v_\pb')$ denotes the states of particle $\pb$ at time $\ell \tau$ and $(\underline{\ell} -1)\tau$ respectively, in the $\Mb_B$ prescribed dynamics. For an interaction set $\Pf_{r(\Mb_B)}$ on $r(\Mb_B)$, we also define
\begin{equation}\label{eq.associated_integral_molecule_P}
|\Ic\Nc_{\Mb_B}|(\vz_{r(\Mb_B)},\Pf_{r(\Mb_B)}):=\varepsilon^{-(d-1)|{\Mb_B}|_{p\backslash r}}\int\sum_{\Pf:\ \Pf|_{r(\Mb_B)}=\Pf_{r(\Mb_B)}}|(\Sc\circ\mathbbm{1})_{{\Mb_B}}|\,(|Q| \cdot \mathbbm{1}_{\{\Pf \ \rm trivial\}})(\vz,\Pf)\,\mathrm{d}\vz_{(p\backslash r)({\Mb_B})},
\end{equation}
where the sum runs over the interaction sets $\Pf$ on $p(\Mb_B)$ at time $\ell\tau$, and $\Pf|_{r(\Mb_B)}$ is as in Definition \ref{def.partial_dynamics}; $\Ic\Nc_{\Mb_B}(\vz_{r(\Mb_B)},\Pf_{r(\Mb_B)})$ is defined in the same way. Then $\sum_{\Pf_{r(\Mb_B)}}|\Ic\Nc_{\Mb_B}|(\vz_{r(\Mb_B)},\Pf_{r(\Mb_B)})=|\Ic\Nc_{\Mb_B}|(\vz_{r(\Mb_B)})$.
 \end{definition}

Now, we are finally ready to state the main cumulant expansion formula in this section, in Proposition \ref{prop.cumulant_formula} below.
\begin{proposition}\label{prop.cumulant_formula} Recall $\Lambda_\ell$ and $A_\ell$ defined in \eqref{eq.defLambdaseq}, and $t_{\mathrm{fin}}=\Lf\tau$. Assume the initial density function $W_{0,N}(\vz_N)$ is given by \eqref{eq.N_par_ensemble}. Then we can define suitable quantities $\widetilde{f}_s(\ell\tau,\vz_s)$ (for all $\ell\in[\Lf]$ and $s\leq A_\ell$) and $f_s^{\mathrm{err}}(t_{\mathrm{fin}},\vz_s)$ (for $s\leq A_\Lf$), see (\ref{eq.truncated_correlation})--(\ref{eq.truncated_error}), such that we have the decomposition
\begin{equation}\label{eq.f_s_dec}
    f_s(t_{\mathrm{fin}}, \boldsymbol{z}_s) = \widetilde{f}_s(t_{\mathrm{fin}}, \boldsymbol{z}_s) + f_s^{\mathrm{err}}(t_{\mathrm{fin}}, \boldsymbol{z}_s),
\end{equation}
and for each $\ell\in[0:\Lf]$ and $s\leq A_\ell$, we have
\begin{equation}\label{eq.cumulant_expansion}
    \widetilde{f}_{s}(\ell\tau,\vz_s)=\sum_{H\subseteq [s]}(f^\Ac(\ell\tau))^{\otimes([s]\backslash H)}\cdot E_H(\ell\tau,\vz_H)+\mathrm{Err}(\ell\tau, \vz_s).
\end{equation} 
In (\ref{eq.f_s_dec})--\eqref{eq.cumulant_expansion} we have the following estimates:
\begin{enumerate}
    \item\label{it.cumulant_formula_1} $(f^\Ac(\ell\tau))^{\otimes([s]\backslash H)}$ is a function of $\vz_{[s]\backslash H}$, which is the tensor product of the single variable function $f^\Ac(\ell\tau, z)$. This function is defined inductively in $\ell$ by $f^\Ac(0)=f_0$ and
    \begin{equation}\label{eq.fAterm}
        f^{\mathcal{A}}(\ell\tau,z) = \sum_{\Mb\in \mathcal{T}_{\Lambda_\ell}}\Ic\Nc_\Mb(z) + f^{\mathcal{A},\mathrm{err}}(\ell\tau,z),\quad |f^{\mathcal{A},\mathrm{err}}(\ell\tau,z)| \le|\log\varepsilon|^{C^*} \sum_{\Mb\in \mathcal{T}_{\Lambda_\ell}^{\mathrm{err}}} |\Ic\Nc_\Mb|(z).
    \end{equation}
    \item The cumulant $E_H(\ell \tau)$ is bounded by 
    \begin{equation}\label{eq.cumulant_formula}
        |E_H(\ell\tau,\vz_H)| \le 
        \sum_{\Mb\in \mathcal{F}_{\boldsymbol{\Lambda}_\ell},r(\Mb) = H} |\Ic\Nc_\Mb|(\vz_H).
    \end{equation}
    \item We have $\mathrm{Err}=\mathrm{Err}^1+\mathrm{Err}^2$ (where $\ell\geq 1$, for $\ell=0$ we have $\mathrm{Err}(0)=0$), and
    \begin{align}\label{eq.cumulant_formula_err}
        |\mathrm{Err}^1(\ell\tau, \vz_s)|&\leq \sum_{\Mb\in \mathcal{F}_{\boldsymbol{\Lambda}_\ell}^{\mathrm{err}},r(\Mb) = [s]}|\Ic\Nc_\Mb|(\vz_s),\\
\label{eq.cumulant_formula_err2}
        \|\mathrm{Err}^2(\ell\tau)\|_{L^1} &\le \varepsilon^{-2(d-1)\Lambda_{\ell}^2A_{\ell}}\cdot \|\mathrm{Err}((\ell-1)\tau)\|_{L^1}.
    \end{align}
    \item The truncation error $f^{\mathrm{err}}(t_{\mathrm{fin}}, \vz_s)$ in (\ref{eq.f_s_dec}) satisfies
    \begin{equation}\label{eq.cumulant_formula_err_trunc}
        \|f^{\mathrm{err}}(t_{\mathrm{fin}}, \vz_s)\|_{L^1}\leq \sum_{\ell=1}^\Lf\bigg(\sum_{\Mb\in \mathcal{F}_{\boldsymbol{\Lambda}_\ell}^{\mathrm{trc.err}},r(\Mb) = [s]}\big\||\Ic\Nc_\Mb|(\vz_s)\big\|_{L^1}+\varepsilon^{-8(d-1)\Lambda_\ell}\cdot\|\mathrm{Err}((\ell-1)\tau)\|_{L^1}\bigg).
    \end{equation}
\end{enumerate}
\end{proposition}
 
The proof of Proposition \ref{prop.cumulant_formula} uses the Penrose equality (Lemma \ref{lem.penrose_equality}) and the recurrence formula for $\widetilde f_s$ (Proposition \ref{prop.recurrence_f}), which are the analogs of Propositions 5.10 and 5.12 of \cite{DHM24}.

As in Definition 5.11 of \cite{DHM24}, we denote by $\Fc_{\Lambda_\ell}$ (resp. $\Fc^{\rm err}_{\Lambda_\ell}$) the set of molecules $\Mb$ with only one layer $\ell$ and $r(\Mb)=[s]$ that satisfy Definition \ref{def.set_T_F} \eqref{it.set_F} (resp. \eqref{it.set_F_err}) in the layer $\ell$.
\begin{lemma}[Penrose equality]\label{lem.penrose_equality}
For $N\geq s$ we have
\begin{equation}\label{eq.penrose_equality}
\Sc^{\ell}_N(\tau)=\sum_{\substack{\Mb\in\Fc_{\Lambda_\ell}\\ p(\Mb)\subseteq[N]}}(\Sc\circ\mathbbm1)^{\rm Pen}_\Mb\circ\Sc^{\ell}_{[N]\setminus p(\Mb)}(\tau)+\sum_{\substack{\Mb\in\Fc^{\rm err}_{\Lambda_\ell}\\ p(\Mb)\subseteq[N]}}O_1\Big(|(\Sc\circ\mathbbm1)_\Mb|\circ\Sc^{\ell}_{[N]\setminus p(\Mb)}(\tau)\Big),
\end{equation}
where $|O_1(A)|\leq A$, and the operators $(\Sc\circ\mathbbm1)^{\rm Pen}_\Mb$ satisfy properties (2a)--(2c) of Definition 5.7 of \cite{DHM24}. In particular $|(\Sc\circ\mathbbm1)^{\rm Pen}_\Mb g|\leq|(\Sc\circ\mathbbm1)_\Mb|(|g|)$, and $(\Sc\circ\mathbbm1)^{\rm Pen}_\Mb$ factorizes over the components of $\Mb$.
\end{lemma}
\begin{remark}\label{rem.factor_P}
Let $\Mb_B$ be a boxed C-molecule in \eqref{eqn.sum_MbB}, and let $i,j$ be particles in different components of $\Mb_B$. By Definition \ref{def.admissible}, $j\in\Pf_i(t)$ means that $i$ and $j$ are colliding at time $t$, so $j\notin\Pf_i(t)$ for all $t\in[(\ell-1)\tau,\ell\tau]$. Hence, if $\lambda$ is a union of components of $\Mb_B$, then $\Pf(t)=(\Pf(t)|_\lambda,\Pf(t)|_{\lambda^c})$, where $\Pf|_\lambda:=(\Pf_i\cap\lambda)_{i\in\lambda}$ as in Definition \ref{def.partial_dynamics}, and the particles in $\lambda$ and in $\lambda^c$ follow two partial dynamics. In Lemma \ref{lem.penrose_equality} we apply this with $\lambda=p(\Mb)$: the two operators act on $(\vz_{p(\Mb)},\Pf|_{p(\Mb)})$ and on $(\vz_{[N]\setminus p(\Mb)},\Pf|_{[N]\setminus p(\Mb)})$, and they do not depend on whether $j\in\Pf_i$ for $i\in p(\Mb)$ and $j\notin p(\Mb)$.
\end{remark}
\begin{proof}[Proof of \ref{lem.penrose_equality}]
The proof of Proposition 5.10 of \cite{DHM24} applies, starting from Lemma \ref{lem.factorization} \eqref{it.factorization4} instead of (5.6) of \cite{DHM24}. We point out the differences.
\begin{enumerate}
\item The set $E$ consists of pairs $(\af,\af')$ of edges and particle pieces of atoms in different clusters, as in Lemma \ref{lem.factorization} \eqref{it.factorization3}. The linear ordering is extended by ordering atom pieces along each particle line in the same way as edges. An O-connection at an atom produces a CO- or MO-atom. The indicator $\mathbbm{1}_\prec$ orders the O-connections along the same ov-segment, and also the O-connections at the same atom, according to the decoration $\Of(\nf)$ in Definition \ref{def.molecule} \eqref{it.decoration}. The correspondence $(\Mb,F,\prec)\leftrightarrow\Mb_F$ and the construction of $\pi$ and $FB(F)$ are unchanged.
\item The interaction sets $\Pf$ factorize over the components, see Remark \ref{rem.factor_P}.
\end{enumerate}
\end{proof}

\begin{proposition}[Recurrence formula for $\widetilde f_s$]\label{prop.recurrence_f}
For $\ell\in[0:\Lf]$ and $\Pf_s=(\Pf_i)_{i\in[s]}$, define
\begin{equation}\label{eq.f_tilde_link}
\widetilde f_s(\ell\tau,\vz_s,\Pf_s):=\eps^{(d-1)s}\sum_{n\geq0}\frac1{n!}\int\sum_{\Pf:\ \Pf|_{[s]}=\Pf_s}\widetilde W_{s+n}(\ell\tau,\vz_{s+n},\Pf)\,\d\vz_{[s+1:s+n]},
\end{equation}
where $\Pf |_{[s]}$ is as in Definition \ref{def.partial_dynamics}, and the sum is over all $\Pf$ satisfying Definition \ref{def.admissible} \eqref{it.admissible.1}; for a general set of labels $Q$ in place of $[s]$ we define $\widetilde f_s(\ell\tau,\vz_Q,\Pf_Q)$ in the same way. Then $\widetilde f_s(\ell\tau,\vz_s)=\Pi\widetilde f_s(\ell\tau,\vz_s,\cdot)$. For $\ell\geq1$ and $\Mb$ with one layer $\ell$ and $r(\Mb)=[s]$, define
\begin{equation}\label{eq.f_M_link}
\widetilde f_{\Mb}(\ell\tau,\vz_s,\Pf_s):=\eps^{-(d-1)|\Mb|_{p\setminus r}}\int\sum_{\Pf_{p(\Mb)}|_{[s]}=\Pf_s}(\Sc\circ\mathbbm1)^{\rm Pen}_\Mb\,\big( \widetilde f_{|\Mb|_p}((\ell-1)\tau, \vz'_{p(\Mb)},\Pf'_{p(\Mb)}) \big)\,\d\vz_{(p\setminus r)(\Mb)},
\end{equation}
where $(\vz'_{p(\Mb)},\Pf_{p(\Mb)}')$ denotes the state of $p(\Mb)$ at time $(\ell-1)\tau$, and define $\widetilde f_{\Mb,\rm upp}$ in the same way with $|(\Sc\circ\mathbbm1)_\Mb|$ in place of $(\Sc\circ\mathbbm1)^{\rm Pen}_\Mb$. Then
\begin{equation}\label{eq.recurrence_f}
\widetilde f_s(\ell\tau,\vz_s,\Pf_s)=\sum_{[\Mb]\in\Fc_{\Lambda_\ell}}\widetilde f_{\Mb}(\ell\tau,\vz_s,\Pf_s)+\mathrm{Err},\qquad |\mathrm{Err}|\leq\sum_{[\Mb]\in\Fc^{\rm err}_{\Lambda_\ell}}\widetilde f_{\Mb,\rm upp}(\ell\tau,\vz_s,\Pf_s),
\end{equation}
where $[\Mb]$ is the equivalence class under relabeling of non-root particle lines, and the sum over $[\Mb]$ is interpreted as in (5.34) of \cite{DHM24}.
\end{proposition}
\begin{proof}
The proof of Proposition 5.12 of \cite{DHM24} applies, with (5.17) there replaced by \eqref{eq.penrose_equality}. The only difference is in the step from (5.38) to (5.39) of \cite{DHM24}, where the particles outside $p(\Mb)$ are integrated out. By Remark \ref{rem.factor_P}, at time $\ell\tau$ we have $\Pf=(\Pf|_{p(\Mb)},\Pf|_{[N]\setminus p(\Mb)})$, so the sum over $\Pf$ with $\Pf|_{[s]}=\Pf_s$ becomes the sum over $\Pf|_{p(\Mb)}$ in \eqref{eq.f_M_link} and a free sum over $\Pf|_{[N]\setminus p(\Mb)}$. By the $L^1$ invariance \eqref{eq.L1pres} of $\Sc^\ell_{[N]\setminus p(\Mb)}(\tau)$, the free sum and the integral over the particles outside $p(\Mb)$ give the same sum and integral at time $(\ell-1)\tau$.
\end{proof}

\begin{proposition}[Recurrence formula for cumulants]\label{prop.recurrence_E}
Fix $\ell\geq1$ and $(A_\ell,\Lambda_\ell)$. Assume that for each $Q_{\ell-1}$ with $s_+=|Q_{\ell-1}|\le\Lambda_\ell^2A_\ell$ we have, at $t=(\ell-1)\tau$,
\begin{equation}\label{eq.ansatz_lminus1}
\begin{aligned}
    \widetilde f_{s_+}\big((\ell-1)\tau,\vz_{Q_{\ell-1}},\mathsf P_{Q_{\ell-1}}\big) =\sum_{H_{\ell-1}\subseteq Q_{\ell-1}}(f^{\mathcal A}((\ell-1)\tau))^{\otimes(Q_{\ell-1}\setminus H_{\ell-1})} \mathbbm{1}_{\{\mathsf P_{Q_{\ell-1}\setminus H_{\ell-1}}\ \rm trivial\}}\\
    \times E_{H_{\ell-1}}\big((\ell-1)\tau,\vz_{H_{\ell-1}},\mathsf P_{H_{\ell-1}}\big)  +\mathrm{Err}\big((\ell-1)\tau,\vz_{Q_\ell},\Pf_{Q_\ell}\big),
\end{aligned}
\end{equation}
where $ \widetilde f_{s_+}$ is as in \eqref{eq.f_tilde_link}, $E_\emptyset=1$, $E_H$ is symmetric under permutations acting on $(\vz_H,\mathsf P_H)$, and $f^{\mathcal A}$ is a function of a single $z=(x,v)$ (hence does not depend on $\Pf$). For a single-layer boxed molecule $\Mb$ and
$H\subseteq p(\Mb)$ define
\[
\begin{aligned}
    \Ic\Nc^{\rm Pen}_{\Mb,H}(\vz_{r(\Mb)},\mathsf P_{r(\Mb)})
    :=\varepsilon^{-(d-1)|\Mb|_{p\setminus r}}\int\sum_{\Pf=\Pf_{p(\Mb)}:\ \mathsf P|_{r(\Mb)}=\mathsf P_{r(\Mb)}}
    (\mathcal S\circ\mathbbm{1})^{\rm Pen}_{\Mb}\Big[(f^{\mathcal A}((\ell-1)\tau))^{\otimes(p(\Mb)\setminus H)}
    \mathbbm{1}_{\{\mathsf P'_{p(\Mb)\setminus H}\ \rm trivial\}}\\
    \times E_H((\ell-1)\tau,\vz'_H,\mathsf P'_H)\Big] \d \vz_{(p\setminus r)(\Mb)},
\end{aligned}
\]
where the sum runs over the interaction sets $\mathsf P$ on $p(\Mb)$ at time $\ell\tau$, $\Pf_{Q \ \rm trivial} $ means $\Pf_i = \{i\}$ for all $i \in Q$, and $(\vz',\Pf')$ is the state at time $(\ell-1)\tau$ , and define
$|\Ic\Nc_{\Mb,H}|$ in the same way with $|(\mathcal S\circ\mathbbm{1})_{\Mb}|$, $|f^{\mathcal A}|$
and $|E_H|$. Then for each $Q_\ell$ with $s=|Q_\ell|\le A_\ell$,
\[
\widetilde f_s(\ell\tau,\vz_{Q_\ell},\mathsf P_{Q_\ell})=\sum_{H_\ell\subseteq Q_\ell}
(f^{\mathcal A}(\ell\tau))^{\otimes(Q_\ell\setminus H_\ell)}\,\mathbbm{1}_{\{\mathsf P_{Q_\ell\setminus H_\ell}\ \rm trivial\}}\,
E_{H_\ell}(\ell\tau,\vz_{H_\ell},\mathsf P_{H_\ell})+\mathrm{Err}(\ell\tau,\vz_{Q_{\ell}}, \Pf_{Q_{\ell}}),
\]
where
\begin{enumerate}
\item\label{it.recE_fA} $f^{\mathcal A}(\ell\tau,z)=\sum_{[\Mb]\in\mathcal T_{\Lambda_\ell}}\Ic \Nc^{\rm Pen}_{\Mb,\emptyset}(z)$ (note that $r(\Mb)=\{i\}$ consists of a single particle, so $(\{i\})$ is the only interaction set on $r(\Mb)$, and we write $\Ic\Nc^{\rm Pen}_{\Mb,\emptyset}(z):=\Ic\Nc^{\rm Pen}_{\Mb,\emptyset}(z,(\{i\}))$);
\item $E_{H_\ell}$ is symmetric under permutations acting on $(\vz_{H_\ell},\mathsf P_{H_\ell})$, and 
\begin{equation}\label{eq.EH}
    |E_{H_\ell}(\ell\tau,\vz_{H_\ell},\mathsf P_{H_\ell})|\leq \sum_{\substack{[\Mb_\ell] \in \Fc_{\Lambda_\ell}, r(\Mb_{\ell}) = H_{\ell}\\
     H_{\ell-1} \subseteq p(\Mb_{\ell}), r(\Mb_\ell)\searrow H_{\ell-1}}} |\Ic\Nc_{\Mb_\ell,H_{\ell-1}}|(\vz_{H_\ell},\mathsf P_{H_\ell}). 
\end{equation}
 
\item $\mathrm{Err}=\mathrm{Err}_1+\mathrm{Err}_2$ with
$$|\mathrm{Err}_1(\ell\tau,\vz_{Q_{\ell}}, \Pf_{Q_{\ell}})|\le\sum_{[\Mb_\ell]\in\mathcal F^{\rm err}_{\Lambda_\ell},\,r(\Mb_\ell)=Q_\ell,\,H_{\ell-1}\subseteq p(\Mb_\ell)} |\Ic\Nc_{\Mb_\ell,H_{\ell-1}}|(\vz_{Q_\ell},\Pf_{Q_{\ell}}),$$
and
$$\|\mathrm{Err}_2(\ell\tau,\vz_{Q_\ell},\Pf_{Q_\ell})\|_{L^1}\le\varepsilon^{-2(d-1)\Lambda_\ell^2A_\ell}\|\mathrm{Err}((\ell-1)\tau,\vz_{Q_{\ell-1}},\Pf_{Q_{\ell-1}})\|_{L^1},$$
where $\|\cdot\|_{L^1}$ means $\int \d \vz_Q \sum_{\mathsf P_Q}$.
\end{enumerate}
\end{proposition}
\begin{proof}
The proof of Proposition 5.16 of \cite{DHM24} applies, with Proposition 5.12 there replaced by Proposition \ref{prop.recurrence_f}: for $[\Mb]\in\Fc_{\Lambda_\ell}$ we have $|p(\Mb)|\leq\Lambda_\ell^2A_\ell$, so we may insert \eqref{eq.ansatz_lminus1} with $Q_{\ell-1}=p(\Mb)$ into \eqref{eq.f_M_link}. The term $\mathrm{Err}((\ell-1)\tau)$ gives $\mathrm{Err}_2$, $\mathrm{Err}_1$ is obtained as in \cite{DHM24}, and the main terms are $\Ic\Nc^{\rm Pen}_{\Mb,H_{\ell-1}}$. We point out the differences, which come from the interaction sets.
\begin{enumerate}
\item\label{it.non-bottom} Every particle of a bottom crossing atom of $\Mb$ with $\Pf'_i((l-1)\tau)\neq \{i\}$ belongs to $H_{\ell-1}$, and every bottom crossing atom contains such a particle. Indeed, a bottom crossing atom requires $j\in\Pf'_i$ at time $(\ell-1)\tau$ for some $i\neq j$ in this atom. By the factor $\mathbbm{1}_{\{\Pf'_{p(\Mb)\setminus H_{\ell-1}}\ \rm trivial\}}$ in $\Ic\Nc^{\rm Pen}_{\Mb,H_{\ell-1}}$ and the symmetry of $\Pf'$, both $i$ and $j$ lie in $H_{\ell-1}$. In the multi-layer molecule, $H_{\ell-1}=r(\Mb_{\ell-1})$, so the particles which are still interacting at time $(\ell-1)\tau$ are roots of layer $\ell-1$, where they belong to the top crossing atom of the same collision. For $\ell=1$ there is no bottom crossing atom, since $\Lb W_{0,N}$ is supported on $\Pf=(\{i\})_i$.

\item Outside $H_{\ell-1}$ the interaction sets $\Pf'$ at time $(\ell-1)\tau$ are trivial, and $(\Sc\circ\mathbbm{1})^{\rm Pen}_\Mb$ factorizes over the components of $\Mb$ (Lemma \ref{lem.penrose_equality}). As in Part 2 of the proof of Proposition 5.16 of \cite{DHM24}, we say that a component of $\Mb$ is \emph{single} if it contains exactly one root and no particle of $H_{\ell-1}$. A single component gives a factor $\Ic\Nc^{\rm Pen}_{\Mb_j,\emptyset}$, and these factors sum to $f^\Ac(\ell\tau)$. The other components form $\Mb_\ell$, with $r(\Mb_\ell)=H_\ell$ and $r(\Mb_\ell)\searrow H_{\ell-1}$ (i.e. each root of $\Mb_\ell$ lies in a component containing another root or a particle line of $H_{\ell-1}$, Definition 5.14 of \cite{DHM24}), and give $E_{H_\ell}$, which is symmetric since the sum is over all labelings.
\item The root $i$ of a single component either interacts with no other particle at time $\ell\tau$, or it is at a top crossing atom whose other particles are non-root particle lines, which are integrated out at time $\ell\tau$. In both cases $\Pf_i|_{Q_\ell}=\{i\}$, which gives the factor $\mathbbm{1}_{\{\Pf_{Q_\ell\setminus H_\ell}\ \rm trivial\}}$.
\end{enumerate}
\end{proof}
\begin{remark}
    Here we remark that bottom crossing atom is connected to a top crossing atom, but top fixed end is not necessarily connected to a bottom of the higher layer. In particular, the total energy is nonincreasing. 
\end{remark}

\begin{proof}[Proof of Proposition \ref{prop.cumulant_formula}] \begin{enumerate}
\item\label{it.phasespace} \emph{Phase space.} All transport operators act on functions of $(\vz,\Pf)$. Since $W_{0,N}$ is supported in $\mathcal D_N^0$, $\Lb W_{0,N}$ is supported on $\Pf=(\{i\})_i$, which is well-prepared, and by Proposition \ref{prop.dynamics_E_T_property} \eqref{it.dynamics_E_T_property_2}--\eqref{it.dynamics_E_T_property_3}, $W_N(t_{\rm fin})=\Pi \circ(\Sc_N^E(\tau))^\Lf\circ\Lb W_{0,N}$. The $L^1$ invariance (Proposition \ref{prop.dynamics_E_T_property} \eqref{it.dynamics_E_T_property_4}) is $\int\sum_{\Pf}$, which is the form used in (5.38) and (5.98) of \cite{DHM24}.

\item \emph{Initial cumulant expansion.} We can directly argue as Proposition 5.18 of \cite{DHM24}, as $W_{0,N}$ contains only the positional exclusion $\mathbbm{1}_{\mathcal D_N^0}$ and $\Pf$ is trivial at $t=0$.

As in Proposition 5.18 of \cite{DHM24}, for each $Q_0$ with $|Q_0|\leq A_0$ we have
\begin{equation}\label{eq.initial_cumulant}
\widetilde f_{|Q_0|}(0,\vz_{Q_0},\Pf_{Q_0})=\sum_{H_0\subseteq Q_0}(f^\Ac(0))^{\otimes(Q_0\setminus H_0)}\,\mathbbm{1}_{\{\Pf_{Q_0\setminus H_0}\ \rm trivial\}}\cdot E_{H_0}(0,\vz_{H_0},\Pf_{H_0}),\qquad f^\Ac(0)=f_0,
\end{equation}
with $\mathrm{Err}(0)=0$ and
\begin{equation}\label{eq.initial_cumulant_bound}
E_{H_0}(0,\vz_{H_0},\Pf_{H_0})=E_{H_0}(0,\vz_{H_0})\,\mathbbm{1}_{\{\Pf_{H_0}\ \rm trivial\}},\qquad |E_{H_0}(0,\vz_{H_0})|\leq(f^\Ac(0))^{\otimes H_0}\sum_{\Lc}\mathbbm{1}^\eps_\Lc,
\end{equation}
where $\Lc$ is a set of pairs of particle lines in $H_0$ which forms a forest with no isolated vertex over the vertex set $H_0$, and $\mathbbm{1}^\eps_\Lc$ is as in \eqref{eq.1_L}. Indeed, $\Lb W_{0,N}$ is supported on $\Pf=(\{i\})_i$ (see \eqref{it.phasespace}), so by \eqref{eq.truncated_correlation}, \eqref{eq.W_widetilde_l} and \eqref{eq.f_tilde_link}, we have $\widetilde f_{|Q_0|}(0,\vz_{Q_0},\Pf_{Q_0})=f_{|Q_0|}(0,\vz_{Q_0})\,\mathbbm{1}_{\{\Pf_{Q_0}\ \rm trivial\}}$. So the expansion of $f_{|Q_0|}(0,\vz_{Q_0})$ and the bound for $E_{H_0}(0,\vz_{H_0})$ are (5.70)--(5.71) of \cite{DHM24}.

\item \emph{The term $f^{\Ac,\mathrm{err}}$.} By Proposition \ref{prop.recurrence_E} \eqref{it.recE_fA}, $f^\Ac(\ell\tau)=\sum_{\Tc_{\Lambda_\ell}}\Ic\Nc^{\rm Pen}_{\Mb,\emptyset}$, where $(\Sc\circ\mathbbm{1})^{\rm Pen}_{\Mb}$ contains, in addition to $(\Sc\circ\mathbbm{1})_{\Mb}$, the indicator functions that no two edges or particle pieces of atoms of $\Mb$ overlap in a way which is not recorded in $\Mb$ (a single-root molecule has no O-connection with another component). Exactly as in (5.83)--(5.85) of \cite{DHM24}, expanding $\mathbbm{1}_{\not\sim}=1-\mathbbm{1}_{\sim}$ for these pairs gives
\[
\Ic\Nc^{\rm Pen}_{\Mb,\emptyset}(z)=\Ic\Nc_\Mb(z)+\sum_{(\af_1,\af_2)}O_1\big(|\Ic\Nc_{\widetilde\Mb}|(z)\big),\qquad |O_1(A)|\leq A,
\]
where $\widetilde\Mb=\Mb_{\af_1,\af_2}\in\Tc^{\mathrm{err}}_{\Lambda_\ell}$ is obtained from $\Mb$ by adding one O-connection joining $\af_1$ and $\af_2$ (edges, or particles of a C-atom or C-constituent). Since the number of $\Mb$ for fixed $\widetilde\Mb$ is at most $|\log\eps|^{C^*}$, this gives \eqref{eq.fAterm}.

\item \emph{Multi-layer molecules.} For $\ell'\in[1:\ell]$, a single-layer molecule $\Mb_{\ell'}$ with $r(\Mb_{\ell'})=H_{\ell'}$ and $H_{\ell'-1}\subseteq p(\Mb_{\ell'})$, and a function $g(\vz_{H_{\ell'-1}},\Pf_{H_{\ell'-1}})\geq0$, define
\begin{equation}\label{eq.proof_multi_op}
\begin{aligned}
\mathrm{Op}_{\Mb_{\ell'}}[g](\vz_{H_{\ell'}},\Pf_{H_{\ell'}}):=\ &\eps^{-(d-1)|\Mb_{\ell'}|_{p\setminus r}}\int\sum_{\Pf:\ \Pf|_{H_{\ell'}}=\Pf_{H_{\ell'}}}|(\Sc\circ\mathbbm1)_{\Mb_{\ell'}}|\Big[|f^\Ac((\ell'-1)\tau)|^{\otimes(p(\Mb_{\ell'})\setminus H_{\ell'-1})}\\
&\times\mathbbm1_{\{\Pf'_{p(\Mb_{\ell'})\setminus H_{\ell'-1}}\ \rm trivial\}}\,g(\vz'_{H_{\ell'-1}},\Pf'_{H_{\ell'-1}})\Big]\,\d\vz_{(p\setminus r)(\Mb_{\ell'})},
\end{aligned}
\end{equation}
where the sum runs over the interaction sets $\Pf = \Pf_{p(\Mb_{\ell'})}$ at time $\ell'\tau$, and $(\vz',\Pf')$ is the state of $p(\Mb_{\ell'})$ (so $\Pf' = \Pf'_{p(\Mb_{\ell'})}$ at time $(\ell'-1)\tau$. Then \eqref{eq.EH} reads $|E_{H_\ell}|\leq\sum_{[\Mb_\ell]}\mathrm{Op}_{\Mb_\ell}[|E_{H_{\ell-1}}|]$, with the sum as in \eqref{eq.EH}. Iterating this and using \eqref{eq.initial_cumulant_bound} gives, as in (5.86)--(5.87) of \cite{DHM24},
\begin{equation}\label{eq.proof_cumulant_multi_3_2}
|E_{H_\ell}(\ell\tau,\vz_{H_\ell},\Pf_{H_\ell})|\leq\sum_{\substack{([\Mb_{\ell'}])_{\ell'\in[1:\ell]}: \,[\Mb_{\ell'}]\in \mathcal{F}_{\Lambda_{\ell'}}\\ r(\Mb_{\ell'}) = H_{\ell'},\,H_{\ell'-1}\subseteq p(\Mb_{\ell'})\\r(\Mb_{\ell'})\searrow H_{\ell'-1}}}\sum_{\Lc}\big(\mathrm{Op}_{\Mb_\ell}\circ\cdots\circ\mathrm{Op}_{\Mb_1}\big)\Big[|f^\Ac(0)|^{\otimes H_0}\,\mathbbm1^\eps_\Lc\,\mathbbm1_{\{\Pf_{H_0}\ \rm trivial\}}\Big],
\end{equation}
where $\Lc$ runs over the forests on $H_0$ with no isolated vertex, as in \eqref{eq.initial_cumulant_bound}.

\item \emph{Constructing $\Mb$ from $\Mb_{\ell'}$.} As in (5.88)--(5.91) of \cite{DHM24}, we merge the $\Mb_{\ell'}$ into one molecule $\Mb$ with $r(\Mb)=H_\ell$, by joining the top end in $\Mb_{\ell'}$ and the bottom end in $\Mb_{\ell'+1}$ of each particle line in $H_{\ell'}$. The sum over $([\Mb_{\ell'}])_{\ell'}$ and $\Lc$ becomes the sum over $\Mb$, the measures combine into $\eps^{-(d-1)|\Mb|_{p\setminus r}}\d\vz_{(p\setminus r)(\Mb)}$, and
\begin{equation}\label{eq.proof_multi_p}
p(\Mb)=\bigsqcup_{\ell'=1}^\ell\big(p(\Mb_{\ell'})\setminus H_{\ell'-1}\big)\sqcup H_0,\qquad \pb\in p(\Mb_{\ell'})\setminus H_{\ell'-1}\iff\ell_1[\pb]=\ell'.
\end{equation}
 The properties of $\Fc_{\vLambda_\ell}$ in Definition \ref{def.set_T_F} \eqref{it.set_F} are checked as in \cite{DHM24}, except the condition on bottom crossing atoms in \eqref{it.set_F_l_1}, which holds for every nonzero term in \eqref{eq.proof_cumulant_multi_3_2}. Indeed, a bottom crossing atom of $\Mb_{\ell'}$ contains $i\neq j$ with $j\in\Pf'_i$ at time $(\ell'-1)\tau$. By the indicator in \eqref{eq.proof_multi_op} (or $\mathbbm1_{\{\Pf_{H_0}\ \rm trivial\}}$ if $\ell'=1$), $\ell'\geq2$ and $i,j\in H_{\ell'-1}$. Since $\Pf'_{H_{\ell'-1}}$ is the restriction of the output of $\mathrm{Op}_{\Mb_{\ell'-1}}$, $i$ and $j$ are still colliding at time $(\ell'-1)\tau$ in layer $\ell'-1$, so they belong to the same top crossing atom of $\Mb_{\ell'-1}$. Hence \eqref{eq.proof_cumulant_multi_3_2} is bounded by
\begin{equation}\label{eq.proof_multi_sumM}
|E_{H_\ell}(\ell\tau,\vz_{H_\ell},\Pf_{H_\ell})|\leq\sum_{\Mb\in\Fc_{\vLambda_\ell},\,r(\Mb)=H_\ell}\big(\mathrm{Op}_{\Mb_\ell}\circ\cdots\circ\mathrm{Op}_{\Mb_1}\big)\Big[|f^\Ac(0)|^{\otimes H_0}\,\mathbbm1^\eps_\Lc\,\mathbbm1_{\{\Pf_{H_0}\ \rm trivial\}}\Big](\vz_{H_\ell},\Pf_{H_\ell}),
\end{equation}
where $\Mb_{\ell'}$ is the layer $\ell'$ of $\Mb$ and $\Lc$ is the set of initial links of $\Mb$. 

\item \emph{The flow map.}  As in (5.92)--(5.94) of \cite{DHM24}, we show that the summand in \eqref{eq.proof_multi_sumM} equals $|\Ic\Nc_\Mb|(\vz_{H_\ell},\Pf_{H_\ell})$ in \eqref{eq.associated_integral_molecule_P}. Let
\[
F_{\ell'}(\vz',\Pf'):=\mathbbm1_{\Mb_{\ell'}}(\vz',\Pf')\,|f^\Ac((\ell'-1)\tau)|^{\otimes(p(\Mb_{\ell'})\setminus H_{\ell'-1})}(\vz')\,\mathbbm1_{\{\Pf'_{p(\Mb_{\ell'})\setminus H_{\ell'-1}}\ \rm trivial\}}.
\]
Since $\Hc_{\Mb_{\ell'}}(\tau)$ is a map, the sum over $\Pf$ and the sum over preimages in \eqref{eq.proof_multi_op} form one sum over the states $(\vz',\Pf')$ of $p(\Mb_{\ell'})$ at time $(\ell'-1)\tau$:
\[
\begin{aligned}
&\mathrm{Op}_{\Mb_{\ell'}}[g](\vz_{H_{\ell'}},\Pf_{H_{\ell'}})\\
&=\eps^{-(d-1)|\Mb_{\ell'}|_{p\setminus r}}\int\sum_{\substack{(\vz',\Pf'):\ \Hc_{\Mb_{\ell'}}(\tau)(\vz',\Pf')=(\vz,\Pf),\\ \Pf|_{H_{\ell'}}=\Pf_{H_{\ell'}}}}F_{\ell'}(\vz',\Pf')\,g(\vz'_{H_{\ell'-1}},\Pf'_{H_{\ell'-1}})\,\d\vz_{(p\setminus r)(\Mb_{\ell'})}.
\end{aligned}
\]
In $\mathrm{Op}_{\Mb_{\ell'+1}}\circ\mathrm{Op}_{\Mb_{\ell'}}$, the input of layer $\ell'+1$ is the restriction to $H_{\ell'}$ of the output of layer $\ell'$, with $\Pf'_\pb=\{\pb\}$ for the new particles by $F_{\ell'+1}$; by Definition \ref{def.molecule_truncated_dynamics} this is how $\Hc_{\Mb_{\ell'+1}}(\tau)$ acts, so the two sums form one sum over the preimages under $\Hc_{\Mb_{\ell'+1}}(\tau)\circ\Hc_{\Mb_{\ell'}}(\tau)$ (the pairs at time $\ell'\tau$ containing a particle of $(p\setminus r)(\Mb_{\ell'})$ are made trivial by $\Hc_{\Mb_{\ell'+1}}(\tau)$, so the sum over them is part of it), and the integral over $\vz_{(p\setminus r)(\Mb_{\ell'})}$ can be taken outside, since the later layers do not move these particles. Repeating this, we obtain
\begin{equation}\label{eq.proof_flow_all}
\begin{aligned}
&\big(\mathrm{Op}_{\Mb_\ell}\circ\cdots\circ\mathrm{Op}_{\Mb_1}\big)[G_0](\vz_{H_\ell},\Pf_{H_\ell})\\
&=\eps^{-(d-1)|\Mb|_{p\setminus r}}\int\sum_{\substack{(\vz^0,\Pf^0):\ \Hc_{\Mb_\ell}(\tau)\circ\cdots\circ\Hc_{\Mb_1}(\tau)\circ\Tc_1(\vz^0,\Pf^0)=(\vz,\Pf),\\ \Pf|_{H_\ell}=\Pf_{H_\ell}}}\mathbbm1_{\{\Pf^0\ \rm trivial\}}\,G_0\prod_{\ell'=1}^\ell F_{\ell'}\,\d\vz_{(p\setminus r)(\Mb)},
\end{aligned}
\end{equation}
where $G_0:=|f^\Ac(0)|^{\otimes H_0}\mathbbm1^\eps_\Lc\mathbbm1_{\{\Pf_{H_0}\ \rm trivial\}}$ is the function in \eqref{eq.proof_multi_sumM}, and $F_{\ell'}$ is evaluated at the state of $p(\Mb_{\ell'})$ at time $(\ell'-1)\tau$.

Now $\mathbbm1_\Mb=\prod_{\ell'}\mathbbm1_{\Mb_{\ell'}}$, since the topological reduction is defined layer by layer; the states at times $(\ell'-1)\tau$ are well-prepared, since by Definition \ref{def.set_T_F} \eqref{it.set_F_l_1} each bottom crossing atom has at most $d+1$ particles and is marked as short interaction. With \eqref{eq.extraQ} and \eqref{eq.proof_multi_p}, the integrand in \eqref{eq.proof_flow_all} is $(\mathbbm1_\Mb\,|Q|\,\mathbbm1_{\Pf\ \rm trivial})(\vz^0,\Pf^0)$. Finally, $\Hc_\Mb(\ell\tau)=\Tc_2\circ\Hc_{\Mb_\ell}(\tau)\circ\cdots\circ\Hc_{\Mb_1}(\tau)\circ\Tc_1$ by \eqref{eq.multi_layer_flow}, and $\Tc_2$ does not change $\Pf$ and has Jacobian $1$. Changing variables by $\Tc_2$ and splitting the sum back into the sum over $\Pf$ and the sum over preimages, \eqref{eq.proof_flow_all} equals
\[
\begin{aligned}
&\eps^{-(d-1)|\Mb|_{p\setminus r}}\int\sum_{\Pf:\ \Pf|_{H_\ell}=\Pf_{H_\ell}}\Sc_\Mb(\ell\tau)\big(\mathbbm1_\Mb\,|Q|\,\mathbbm1_{\Pf\ \rm trivial}\big)(\vz,\Pf)\,\d\vz_{(p\setminus r)(\Mb)}\\
&=|\Ic\Nc_\Mb|(\vz_{H_\ell},\Pf_{H_\ell}).
\end{aligned}
\]

\emph{Conclusion.} By \eqref{eq.proof_multi_sumM},
\[
|E_{H_\ell}(\ell\tau,\vz_{H_\ell},\Pf_{H_\ell})|\leq\sum_{\Mb\in\Fc_{\vLambda_\ell},\,r(\Mb)=H_\ell}|\Ic\Nc_\Mb|(\vz_{H_\ell},\Pf_{H_\ell}),
\]
and summing over $\Pf_{H_\ell}$ gives \eqref{eq.cumulant_formula}. The same argument, with $\Mb_\ell\in\Fc^{\rm err}_{\Lambda_\ell}$ as in Part 4 of the proof of Proposition 3.25 of \cite{DHM24}, gives $|\mathrm{Err}^1(\ell\tau,\vz_s,\Pf_s)|\leq\sum_{\Mb\in\Fc^{\rm err}_{\vLambda_\ell},\,r(\Mb)=[s]}|\Ic\Nc_\Mb|(\vz_s,\Pf_s)$, hence \eqref{eq.cumulant_formula_err}, and it applies in the same way to $\Fc^{\rm trc.err}_{\vLambda_\ell}$ in \eqref{eq.cumulant_formula_err_trunc}.

\item By induction, we have
 $$\widetilde{f}_s(\ell\tau, \vz_s,\Pf_s) = \sum_{H \subseteq [s]} f^\Ac(\ell \tau)^{\otimes ([s]\backslash H)} \mathbbm{1}_{\{\Pf_{[s]\backslash H} \rm trivial\}}\cdot E_H(\ell \tau,\vz_H ,\Pf_H) + \mathrm{Err}(\ell \tau, \vz_s,\Pf_s).$$
 
We now sum over $\Pf_s$, using $\widetilde{f}_s(\ell\tau,\vz_s) = \sum_{\Pf_s} \widetilde{f}_s(\ell\tau,\vz_s,\Pf_s)$ (Proposition \ref{prop.recurrence_f}). For fixed $\vz_s$ there are finitely many $\Pf_s$, so all sums are finite. By the symmetry of $\Pf_s$, the factor $\mathbbm{1}_{\{\Pf_{[s]\setminus H}\ \rm trivial\}}$ means that $\Pf_s$ consists of the trivial sets $\{i\}$ for $i\in[s]\setminus H$ and of an interaction set $\Pf_H$ on $H$. Hence
\[
\begin{aligned}
    \sum_{\Pf_s}(f^\Ac(\ell\tau))^{\otimes([s]\setminus H)}\mathbbm{1}_{\{\Pf_{[s]\setminus H}\ \rm trivial\}}\,E_H(\ell\tau,\vz_H,\Pf_H) =&(f^\Ac(\ell\tau))^{\otimes([s]\setminus H)}\,E_H(\ell\tau,\vz_H),\\
     E_H(\ell\tau,\vz_H):=&\sum_{\Pf_H}E_H(\ell\tau,\vz_H,\Pf_H),
\end{aligned}
\]
and we set $\mathrm{Err}(\ell\tau,\vz_s):=\sum_{\Pf_s}\mathrm{Err}(\ell\tau,\vz_s,\Pf_s)$. This gives \eqref{eq.cumulant_expansion}. The estimates follow by summing over $\Pf_H$ (resp. $\Pf_s$): for \eqref{eq.cumulant_formula} and \eqref{eq.cumulant_formula_err}, we use that $\sum_{\Pf_H}|\Ic\Nc_\Mb|(\vz_H,\Pf_H)=|\Ic\Nc_\Mb|(\vz_H)$ (Definition \ref{def.associated_int}); for \eqref{eq.cumulant_formula_err2}, the norm $\|\cdot\|_{L^1}$ is $\int\d\vz\sum_{\Pf}$ as in Proposition \ref{prop.recurrence_E}, and $\|\sum_{\Pf}g(\cdot,\Pf)\|_{L^1(\d\vz)}\leq\int\sum_{\Pf}|g|\,\d\vz$.

\item \emph{Error term.} $(\vz^0,\Pf^0)\in\mathcal E^T_N$ if and only if either $(\vz^0,\Pf^0)$ is not well-prepared (a single atom with $|\widetilde{\Pf^0}|>d+1$, or a long interaction), or there is a first rejected collision at some time $t_*$. In the second case the rejection is due to $|\Mb|_p$, $\rho$, $\#_{\rm box}$, an atom with $>d+1$ particles, or a long interaction, and the molecule at $t_*^+$ is read off from $\gamma'$, which coincides with $\widetilde{\mathcal H}^T$ in Definition \ref{def.molecule_truncated_dynamics}. Thus $\mathbbm{1}_{\Mb_B}(\vz^0)=1$, which gives \eqref{eq.cumulant_error} as in (5.96) of \cite{DHM24}.
\end{enumerate}
\end{proof}

At the end of this section, we introduce the following proposition.
\begin{proposition}
\label{prop.mol_axiom} Let $\Mb\in\Fc_{\vLambda_\ell}$ as in Definition \ref{def.set_T_F} (\ref{it.set_F}) and $r(\Mb)=H$. Then we have the following:
\begin{enumerate}
\item\label{it.axiom1-} The total number of atoms and particle lines in $\Mb$ is bounded by $|\log\varepsilon|^{C^*}$.
\item\label{it.axiom1} If two atoms are connected by two different ov-segments, then at least one of them is a multi-particle atom; or they are in two different layers $\ell'-1$ and $\ell'$, and they correspond to the same collision in the sense of Definition \ref{def.short_range} \eqref{it.O5}, but the moment $(\ell'-1)\tau$ belongs to the corresponding collision duration; see Definitions \ref{def.admissible} and \ref{def.top_reduction}.
\item\label{it.axiom1_CO} Each CO-atom may have several C-, M- or edge constituents. Each edge constituent must be O-connected exactly one atom constituent. If there are at least two atom constituents, each atom constituent must be O-connected with other atom constituents. 
\end{enumerate}
Properties (\ref{it.axiom2})--(\ref{it.axiom3}) below correspond to the clusters (Definition \ref{def.cluster}) in each layer $\Mb_{\ell'}\,(1\leq \ell'\leq \ell)$.
\begin{enumerate}
    \setcounter{enumi}{3}
\item\label{it.axiom2} Any cycle in a single layer $\Mb_{\ell'}$ must be formed by ov-segments connecting C-atoms or CO-atoms in the same cluster. In particular, each connected cluster graph is a tree. For any cycle in a single layer $\Mb_{\ell'}$ and any O-bond at an atom in this cycle, if the O-bond at this atom belong to this cycle, then the other O-bond which is serial to this O-bond at the same atom, must also belong to the cycle. This means for a cycle in a single layer, there are no turns at those O-bonds.
\item\label{it.axiom3} For any cluster $\Qb$ in each layer $\Mb_{\ell'}$, we have that $\rho(\Qb)\leq\Gamma$ (Definition \ref{it.recollision_number_box}). We have that the number of boxes in each cluster is at most $2 \Gamma$. Recall that a cluster $\Qb$ or component is cyclic (or recollisional) if $\rho(\Qb)>0$, i.e. if it contains a cycle.
\end{enumerate}
Properties (\ref{it.axiom4})--(\ref{it.axiom5}) below correspond to the particle line sets $p(\Mb_{\ell'})$ and $r(\Mb_{\ell'})$ (Definition \ref{def.prsets}).
\begin{enumerate}
    \setcounter{enumi}{5}
\item\label{it.axiom4} Each particle line $\pb$ must belong to some $p(\Mb_{\ell'})$ for some $\ell'$. If $\pb$ intersects $\Mb_{\ell'}$ (Definition \ref{def.connectedvia}), then $\pb\in p(\Mb_{\ell'})$. If $\pb$ contains an atom in both $\Mb_{\ell_1}$ and $\Mb_{\ell_2}$, then $\pb\in p(\Mb_{\ell'})$ for all $\ell'\in[\ell_1:\ell_2]$.

\item \label{it.axiom4+} We have $r(\Mb_\ell)=H$, and $r(\Mb_{\ell'})=p(\Mb_{\ell'})\cap p(\Mb_{\ell'+1})$ for $1\leq \ell'\leq\ell-1$. Define also $r(\Mb_0)=H_0$.
\item\label{it.axiom5} For $1\leq\ell'\leq\ell$, each component of $\Mb_{\ell'}$ must intersect a particle line in $r(\Mb_{\ell'})$. Each particle line in $p(\Mb_{\ell'})$ is either connected to a particle line in $r(\Mb_{\ell'})$ via $\Mb_{\ell'}$, or belongs to $r(\Mb_{\ell'})$. 
\end{enumerate}
Property (\ref{it.axiom7}) below corresponds to the boxes.
\begin{enumerate}
        \setcounter{enumi}{8}
\item\label{it.axiom7} The Box has three parts. The bottom part and the top part satisfy  that: all particles are connected, and contains at least one cycle. The middle part is regarded as one atom, where we refer to M-atom.
\end{enumerate}
Property (\ref{it.axiom6}) below ensures the connectedness of particle lines in $r(\Mb_{\ell'})$. Recall Definition \ref{def.connectedvia} and Definition \ref{def.more_on_molecule} \eqref{def.layer_interval}.
\begin{enumerate}
        \setcounter{enumi}{9}
\item\label{it.axiom6} For each $1\leq\ell'\leq \ell$, and any particle line $\pb\in r(\Mb_{\ell'})$, this $\pb$ is either connected to another particle line $\pb'\in r(\Mb_{\ell'})$ via $\Mb_{[1:\ell']}$, or connected to a cycle within $\Mb_{[1:\ell']}$.
\end{enumerate}

\textbf{The $\Fc_{\vLambda_\ell}^{\mathrm{err}}$ case.} Finally, suppose $\Mb\in\Fc_{\vLambda_\ell}^{\mathrm{err}}$ and $r(\Mb)=H$. Then $\Mb$ satisfies all of the above properties (\ref{it.axiom1-})--(\ref{it.axiom6}), except that (\ref{it.axiom6}) is true only for $\ell'\leq\ell-1$.

\end{proposition}

\begin{proof}
 \eqref{it.axiom1-} follows from Part 1 of the proof of Proposition \ref{prop.layerrec3} below. Items \eqref{it.axiom2}--\eqref{it.axiom5} and \eqref{it.axiom6} follow from Definition \ref{def.set_T_F} and the Penrose construction in exactly the same way as the corresponding properties of molecules in \cite{DHM24} (stated after Definition 3.20 there); boxes are viewed as single atoms, which does not affect these arguments. \eqref{it.axiom7} follows from the construction of the bottom and top parts in Definition \ref{def.box}, see also Definition \ref{def.boxed_molecule} \eqref{it.box_connect}. We now prove  \eqref{it.axiom1} follows from the definition of topological reduction in Definition \ref{def.top_reduction}.
 
 \emph{Proof of \eqref{it.axiom1_CO}.} In the Penrose argument (Lemma \ref{lem.penrose_equality}), an O-connection between an edge $e$ and an atom constituent $\af$ creates a CO- or MO-atom, which contains $e$ as an edge constituent O-connected to $\af$ (Definition 3.10 of \cite{DHM24}). If $e$ were O-connected to two different atom constituents, the two O-connections would create two atoms instead of one. Hence each edge constituent is O-connected to exactly one atom constituent. The last statement of \eqref{it.axiom1_CO} is part of Definition \ref{def.boxed_molecule} (1).
\end{proof}

\section{Statements of cumulant and $f^\Ac$ estimates}\label{sec.reduction}

    In this section we present the proof of Theorem \ref{th.main}. The proof requires the following three main Propositions. 
 
\begin{proposition}
\label{prop.est_fa} Let $f^\Ac$ be defined as in Proposition \ref{prop.cumulant_formula} \eqref{it.cumulant_formula_1}, inductively via (\ref{eq.fAterm}), with $f^\Ac(0,z)=f_0(z)$; also recall $\beta_\ell$ and $\theta_\ell$ defined in (\ref{eq.decayseq}) and $\|\cdot\|_{\mathrm{Bol}^{\beta_\ell}}$ defined by \eqref{eq.boltzmann_decay_2}. Then we have
\begin{equation}\label{eq.approxfA_1}
    f^\Ac(\ell\tau,z)=f(\ell\tau,z)+(f^\Ac)_{\mathrm{rem}}(\ell\tau,z),
\end{equation}
 where $f(t,z)$ is the solution to the Boltzmann equation as in Theorem \ref{th.main}, and $(f^\Ac)_{\mathrm{rem}}$ satisfies that
\begin{equation}\label{eq.approxfA_2}
    \|(f^\Ac)_{\mathrm{rem}}(\ell\tau)\|_{\mathrm{Bol}^{\beta_\ell}}\leq\varepsilon^{\theta_\ell}.
\end{equation}
\end{proposition}
\begin{proof} The proof is given in Section \ref{sec.f_Ac}.
\end{proof}

The second proposition, provides upper bounds for $E_H$ and $\mathrm{Err}_{\ell}^1$. This is the crucial proposition, and will occupy most parts of the rest of proof. 
\begin{proposition}
\label{prop.cumulant_est} Recall the $E_H(\ell\tau,\vz_H)$ and $\mathrm{Err}^1(\ell\tau,\vz_s)$ introduced in Proposition \ref{prop.cumulant_formula}, and also the notion $C_j^*$ in (\ref{eq.defC*}) in Definition \ref{def.notation} (\ref{it.defC*}). Assume that the estimates (\ref{eq.approxfA_1})--(\ref{eq.approxfA_2}) hold with $\ell$ replaced by any $\ell'\leq\ell-1$. Then we have
\begin{equation}\label{eq.cumulant_est_1}
\|E_H(\ell\tau,\vz_H)\|_{L^1}\leq\varepsilon^{\upsilon/3+(C_{14}^*)^{-1}\cdot |H|}
\end{equation} for any $H\neq\varnothing$, and
\begin{equation}\label{eq.cumulant_est_2}
\|\mathrm{Err}^1(\ell\tau,\vz_s)\|_{L^1}\leq \tau^{\Lambda_{\ell}/10}.
\end{equation}
\end{proposition}
\begin{proof}
    The proof will occupy Sections \ref{sec.local_associated}--\ref{sec.treat_integral}, and will be finished in Section \ref{sec.summary}, pending a crucial combinatorial statement, Proposition \ref{prop.comb_est}. Then, the proof of Proposition \ref{prop.comb_est} and will occupy Sections \ref{sec.cuttingI}--\ref{sec.maincr}.
\end{proof}

In the last proposition we prove upper bounds for $f_s^{\mathrm{err}}(t_{\mathrm{fin}})$.
\begin{proposition}
\label{prop.cumulant_error} Recall the $f_s^{\mathrm{err}}(t_{\mathrm{fin}},\vz_s)$ introduced in Proposition \ref{prop.cumulant_formula}. Then we have
\begin{equation}\label{eq.cumulant_error}
    \|f_s^{\mathrm{err}}(t_{\mathrm{fin}},\vz_s)\|_{L^1}\leq \varepsilon^{\upsilon/3}.
\end{equation}
\end{proposition}
\begin{proof} The proof will occupy Section \ref{sec.error}.
\end{proof}

Now we prove Theorem \ref{th.main} assuming Propositions \ref{prop.est_fa}--\ref{prop.cumulant_error}.
\begin{proof}[Proof of Theorem \ref{th.main} assuming Propositions \ref{prop.est_fa}--\ref{prop.cumulant_error}] With Propositions \ref{prop.est_fa}--\ref{prop.cumulant_error}, Theorem \ref{th.main} follows by the same argument as in Section 6 \cite{DHM24}. 
\end{proof}

\section{Preparations for $E_H$ estimates} \label{sec.local_associated}
 In the rest of this paper, we \textbf{only consider the boxed molecule}, instead of the original molecule. Therefore, for notational simplicity, we also denote a boxed molecule by $\Mb$, and \textbf{also call it a molecule}.

In Sections \ref{sec.local_associated}--\ref{sec.maincr}, we prove Proposition \ref{prop.cumulant_est}. We begin with the upper bounds on $E_H$ and $\mathrm{Err}^1$ in (\ref{eq.cumulant_formula})--(\ref{eq.cumulant_formula_err}) in Proposition \ref{prop.cumulant_formula}, which have the form $\sum_{\Mb}|\Ic\Nc_\Mb|$. To bound these terms, it suffices to (i) control the number of terms $\Mb$ in the summation, and (ii) control the size of $|\Ic\Nc_\Mb|$ for each individual $\Mb$.

\subsection{Upper bounds of combinatorial factors}\label{sec.int_dia_exp} In this subsection we prove an upper bound on the number of molecules $\Mb$ in Proposition \ref{prop.layerrec3}. For this we need to introduce the following parameter $\rho$ for molecules $\Mb$:
\begin{definition}[The parameter $\rho$]\label{def.parameter_rho_old} Let molecule $\Mb$ whose lowest and highest layer are $1$ and $\ell$, respectively, and denote its $\ell'-$th layer by $\Mb_{\ell'}$. Define the parameters $s_{\ell'}$ and $\Rf_{\ell'}$ by $s_{\ell'} \coloneqq |r(\Mb_{\ell'})|$ and $\Rf_{\ell'} \coloneqq \rho(\Mb_{\ell'})$ for $\ell'\in[1:\ell]$; also define $s_0 = |H_0|$ where $H_0$ is the initial cumulant set. Now define the parameter $\rho$ by
\begin{equation}\label{eq.def_rho_old}
    \rho = \sum_{\ell'=0}^\ell s_{\ell'}+\sum_{\ell'=1}^\ell\Rf_{\ell'}.
\end{equation}
We interpret the quantity $\rho$ as the sum of $\Rf_{\ell'}$ (the total number of recollisions within each single layer), together with the sum of $s_{\ell'}$ (the total number of particles lines crossing different layers).
\end{definition}
Each C-atom or C-constituent of CO-atoms involves at most $d+1$ particles, and each box involves at most $\mathsf{N}$ particles. However, O-connections may produce CO-atoms or MO-atoms involving more than $\mathsf{N}$ particles. This motivates the following definition.
\begin{definition}\label{def.number_box_multi}
    Let $\Mb$ be a molecule defined as above. We denote by $\Bf = \#_{\mathrm{box}}(\Mb)$ the total number of boxes (after deleting all O-connections) defined in Definition \ref{def.boxed_molecule} \eqref{it.number_box}. We denote by $\Xf$ the total number of multi-particle C-atoms, plus  the number of C-constituents of CO-atoms or MO-atoms in $\Mb$, \emph{excluding those contained within any box}\footnote{This means that those multi-particle C-atom, or a C-constituents in a box are not counted. If a C-constituent does not belong to any box, and is O-connected with an atom/edge in a box, then it is counted}. 
    For each multi-particle atom, M-atom or MO-atom $\nf$, we denote the number of its particle lines by $\Nf_\nf$, and set 
    $$\Nf := \sum_{\nf} \Nf_\nf,$$ 
    The corresponding parameter at layer $\ell'$ are denoted by $\Bf_{\ell'}$, $\Xf_{\ell'}$ and $\Nf_{\ell'}$, respectively.
\end{definition}
We now state and prove the main combinatorial upper bound in this section.
\begin{proposition}\label{prop.layerrec3} Let $\Mb$ be a molecule and $\rho$ be as in Definition \ref{def.parameter_rho_old}. Let $\Bf$ and $\Xf$ be as in Definition \ref{def.number_box_multi}. Assume $H$ is a non-empty set and $|H|\le |\log\varepsilon|^{C^*}$. Consider all molecules $\Mb\in \mathcal{F}_{\boldsymbol{\Lambda}_\ell}$ (recall that we are summing over non-labeled molecules $\Mb$) that satisfy $r(\Mb)=H$, with the size $|\Mb|$ (the number of all atoms, including those in boxes) of $\Mb$ fixed, and the $\rho, \Bf,\Xf,\Nf$ defined above also fixed. Then the number of choices for $\Mb$, together with the labelings of root particles, is bounded by
\begin{equation}\label{eq.layercomb}
    \sum_{\Mb} 1\leq C^{|\Mb| + \Nf}\cdot |\log\varepsilon|^{C^*(\rho+\Xf + \Bf)}.
\end{equation}
Here $C$ depends only on $d$ and $C_\Phi$. The same result is true if we replace $\mathcal{F}_{\boldsymbol{\Lambda}_\ell}$ by $\mathcal{F}_{\boldsymbol{\Lambda}_\ell}^{\mathrm{err}}$ or $\Tc_{\Lambda_\ell}$ or $\Fc_{\vLambda_\ell}^{\mathrm{trc.err}}$.
\end{proposition}
\begin{proof}
     We only consider the case $\Mb\in \mathcal{F}_{\boldsymbol{\Lambda}_\ell}$, as the other cases are similar. Note that $\Mb$ is constructed from its layers $\Mb_{\ell'}$ by the following process:
\begin{itemize}
\item Label the root particles, and draw a molecule $\Mb_{\ell'}$ with a single layer $\ell'$, for each $\ell'\in[1:\ell]$;
\item For each $\ell'\in[1:\ell-1]$, identify a subset of $p(\Mb_{\ell'})$ with a subset of $p(\Mb_{\ell'+1})$, call this identified subset $H_{\ell'}$ and $|H_{\ell'}|=s_{\ell'}$;
\item For each particle line in $H_{\ell'}$, merge the corresponding top end in $\Mb_{\ell'}$ and bottom end in $\Mb_{\ell'+1}$ to form a bond;
\item Finally choose a set $H_0\subseteq p(\Mb_1)$, and fix the initial links $\Lc$ between particle lines in $H_0$.
\end{itemize}

We now estimate the number of choices for $\Mb$, by considering each of the above steps.

\textbf{Proof part 1}. In this part we prove that $|\Mb|\leq|\log\varepsilon|^{C^*}$ and $|\Mb|_p\leq |\log\varepsilon|^{C^*}$. 
We first show that each box (with no O-connection) contains at most $C=C(\mathsf N)$ atoms. By the structure of a box, it suffices to prove that its bottom part contains at most $C$ atoms. This can be proved by induction on the number of particles in the box.
 
First, suppose that the box has $3$ particle lines. By our construction of box, we stop exploring its bottom part once \emph{all particles are connected in the box, and at least one cycle has been formed}. Hence, in this case, the bottom part contains at most $3$ atoms. 
Now suppose that the bottom part of every box containing at most $q-1$ particles has a bounded number of atoms. Consider the box containing $q$ particles. Then its bottom part has at most $G(q-1)$ atoms, in the sense that viewing each previously constructed box as one single atom. Otherwise, before all $q$ particles become connected, there would exist a connected component containing $q'<q$ particles and more than $G(q')$ internal atoms, as $G$ grows super-exponentially. However, according to our construction, these $q'$ particles and their internal atoms would already have been reduced into a box. This contradicts the construction of the box. Moreover, each box in the bottom part which is viewed as a single atom, also has an upper bound on the number of atoms by induction. Since $q\leq \mathsf N$, the number of atoms in the bottom part is bounded by a constant depending only on $\mathsf N$. The top part can be treated in the same way, and therefore each box contains at most $C=C(\mathsf N)$ atoms. As a result, the number of atoms in each box has an absolute upper bound $C=C(\mathsf{N})$.

By Definition \ref{def.set_T_F} (\ref{it.set_F_l_1})--(\ref{it.set_F_l_2}) we know that each cluster of $\Mb_{\ell'}$ has at most $\Lambda_{\ell'}$ particle lines and at most $\Lambda_{\ell'}+ C(\mathsf{N}) \cdot\Gamma $ atoms (using the upper bound of $\Gamma$ for the recollision number, and each cluster has at most $2 \Gamma$ boxes). Then each component has at most $\Lambda_{\ell'}^2$ particle lines and at most $2\Lambda_{\ell'}^2$ collisions (as each component has at most $\Lambda_\ell$ clusters). Then each $\Mb_{\ell'}$ has at most $A_{\ell'}\Lambda_{\ell'}^2$ particle lines and $2A_{\ell'}\Lambda_{\ell'}^2$ collisions (as each component must contain a root particle line, and there are $s_{\ell'}\leq A_{\ell'}$ of them). All these values are bounded as $\leq|\log\varepsilon|^{C^*}$ by Definition \ref{def.notation} (\ref{it.defnot3}), hence the result.

\textbf{Proof part 2.} In this part we consider the first step of labeling and constructing individual $\Mb_{\ell'}$. Note that labelings lead to $|H|!\leq|\log\varepsilon|^{C^*\rho}$ choices as $|H|=|r(\Mb)|=|r(\Mb_\ell)|=s_\ell\leq\rho$. Next we may fix the values of $(|\Mb_{\ell'}|)_{\ell'\in[1:\ell]}$. As these numbers add up to the fixed $|\Mb|$, we have at most $\binom{|\Mb|+\Lf}{\Lf}$ choices which is acceptable in (\ref{eq.layercomb}). Similarly, we can also fix the values of $(s_{\ell'})$, $(\rho(\Mb_{\ell'}))$, $\Bf_{\ell'}$ and $\Xf_{\ell'} $ at an acceptable price.  

Now we construct $\Mb_{\ell'}$. The molecule $\Mb_{\ell'}$ can be constructed as follows: we start by constructing trees, then forest, and then adding cycles within those forests. For the tree contain every atom in $\Mb_{\ell'}$, the number of multi-particle C-atoms, M-atoms plus the number of C- or M- constituents is no more than $\Xf_{\ell'} + C(\mathsf{N}) \Bf_{\ell'}$. By Lemma \ref{lem.forest_count}, the choice of forests is no more than $C^{|\Mb_{\ell'}| + \Nf_{\ell'}|}  |\log \eps|^{C^*(\Xf_{\ell'}+\Bf_{\ell'})}$. 

Moreover, in $\Mb_{\ell'}$, we know the total number of cycle is no more than $\rho(\Mb_{\ell'}) + C \Bf_{\ell'}$. Hence, once the tree parts are fixed, we choose at most $C \Bf$ atoms from those atoms or constituents with more than $3$ particle lines to be M-atoms or M-constituents. This leads to at most $|\log \eps|^{C^* \Bf}$ loss, which is also acceptable. We then connect $\leq \rho(\Mb_{\ell'}) + C \Bf_{\ell'} $ distinct edges pairs one by one; each edge pair has at most $|\Mb_{\ell'} + C \Bf_{\ell'}|^2 |\log \eps|^{2 C^*} \lesssim|\log \eps|^{C^*} $ choices (first choose two atoms, then choose an edge at each atom. Note that each atom has at most $|\log \eps|^{ C^*}$ edges), which leads to $|\log \eps|^{ C^*(\rho(\Mb_{\ell'}) + \Bf_{\ell'})} $ in total. All these factors are acceptable.

\textbf{Proof part 3.} In this part we consider the steps of matching and merging particle lines/ends in $\Mb_{\ell'}$ and $\Mb_{\ell'+1}$. Recall that $|H_{\ell'}|=s_{\ell'}$; selecting a set of this cardinality in $\Mb_{\ell'}$ and $\Mb_{\ell'+1}$ leads to at most \[\binom{|\Mb_{\ell'}|_p}{s_{\ell'}}\binom{|\Mb_{\ell'+1}|_p}{s_{\ell'}}\leq |\log\varepsilon|^{C^*s_{\ell'}}\] choices. Putting different $\ell'$ together, we again get a factor of $|\log\varepsilon|^{C^*\rho}$ which is acceptable.

\textbf{Proof part 4.} In this part we consider the last step of constructing $H_0$ and $\Lc$. Note that $|H_0|=s_0$, and $|\Lc|$ (i.e. the number of initial links in $\Lc$) is bounded by $2|H_0|$ as $\Lc$ is a forest. The number of choices of $H_0$ and $\Lc$ is then trivially bounded by $|\log\varepsilon|^{C^*s_0}$, which is also acceptable.

By putting Parts 1--4 above together that accounts for all steps in the construction of $\Mb$, we have proved (\ref{eq.layercomb}) and Proposition \ref{prop.layerrec3}.
\end{proof}

\begin{lemma}\label{lem.forest_count}
    Let $\Mb$ be a molecule, containing no cycle. Denote the number of multi-particle C-atoms, M-atoms plus the number of C- or M- constituents by $\Xf'$ (this is different from $\Xf'$ in Definition \ref{def.number_box_multi}). Let $\Nf_\nf$ be the number of particle lines at each multi-particle atom, and $\Nf = \sum_\nf \Nf_\nf$. Let the number of atoms be $|\Mb|$. Then the number of such molecules is at most $C^{|\Mb| + \Nf} |\log \eps|^{C^* (\Xf'+1)}$, where $C$ depends on $\Nf$.
\end{lemma}
\begin{proof}
\textbf{Proof part 1.} In this part we obtain the upper bound of number of molecules $\Mb$ assume that $\Mb$ is a tree.  Define the total number of the required molecule is $B_{|\Mb|, \Xf',\Nf}$. If $\Mb$ has no atoms with more than three particles, then by the same arguments in Proposition 7.2 \cite{DHM24}, we have $B_{\Mb,0,0}\leq C^{|\Mb|}$, which is acceptable. We also denote $B_{|\Mb|, \Xf',\Nf}=0$  if $|\Mb| =0$ and $(\Xf',\Nf)\neq (0,0)$; moreover we define $B_{0,0,0} = 1$. Note that we always have $\Nf \geq 2 \Xf'$ by definition, and we assume $B_{|\Mb|,\Xf',\Nf} =0$ if $\Nf < 2 \Xf'$.

Now assume $\Mb$ has at least one multi-particle atom $\nf$. Assume it has $l$ particle lines (hence $2l$ edges). Then from $\nf$, we can construct $\Mb$ as follows. First, the decoration of $\nf$ is encoded by a plane rooted tree $T_\nf$ whose vertices represent particle lines in $\nf$ (there are $l$ of them). The tree has three types of edges: M-edges, C-edges and O-edges. Furthermore, every subtree $A$ consisting only C-edges has at most $\mathsf{N}$ vertices, and all C-edges within $A$ are ordered. Then the number of such tree is no more than 
\begin{equation*}
    3^l \times 4^l \times \textrm{the number of plane rooted trees} \times (\mathsf{N}!)^{l/2}.
\end{equation*}
Here, $3^l$ determines the choice of every edge,  the number of plane rooted trees determines the shape of the tree, and the factor $(\mathsf{N}!)^{l/2}$ is because in each subtree $A$, we order the C-edges which leads to a factor $\mathsf{N}!$, and the number of such subtree is no more than $l/2$. Moreover, the number of plane rooted trees is bounded by the Catalan number which is no more than $C^l$. The factor $4^l$ accounts for the time order of the O-connections that the plane structure does not record. At a vertex $p$ with $k_p$ O-children, the O-children are already in time order, and if the edge from $p$ to its parent is an O-edge, its position among them is one of $k_p+1\le\deg p$ possibilities. Since $\sum_p\deg p=2(l-1)$, this gives at most $\prod_p\deg p\le\prod_p2^{\deg p}=4^{l-1}$ choices.
Therefore, the decoration of $\nf$ contributes at most $C^l$ choices. After $\nf$ being determined, for each pair of serial edges $e_j$ and $e_j'$ at $\nf$, we add two sub-molecules $\Mb_j$ and $\Mb_j'$ (which may be empty set) with $\Mb_j$ connected to the top edge and the other connected to the bottom edge. Denote the number of such two paired sub-molecules to be $Q_{m,x,n}$, then we have 
     \begin{equation}\label{eqn.Q}
        Q_{m,x,n} \leq  \sum_{i=0}^{m} \sum_{r=0}^{x} \sum_{k=0}^{n} (5i+\mathsf{N} \cdot r + 1)\cdot (5(m-i) + \mathsf{N}\cdot (x-r) + 1)\cdot B_{i,r,k} B_{m-i,x-r,n-k}. 
     \end{equation}
Now we explain the factor $(5i+\mathsf{N}\cdot r + 1)\cdot (5(m-i) + \mathsf{N}\cdot (x-r) + 1)$ in $Q_{m,x,n}$. For each molecule $\Mb_j$ contributes to $B_{i,r,k}$  with $i\geq 1$ atoms, we first choose either a binary atom, or a C- or an M-constituent of a multi-particle atom, or an O-edge of a multi-particle atom in $\Mb_j$, which would be connected to $\nf$ by $e_j$. If we choose a binary atom, then we have at most $2i$ choices. If it is a C- or an M- constituent, then we have at most $\mathsf{N} \cdot r$ choices. If it is an O-edge at a multi-particle atom, and it is serial with a O-free end at the atom, then choose different such O-edge at the same atom leads to the same molecule $\Mb$. Hence, the choice of such edges is at most $i$. If the bottom O-edge is serial with a bond, then the total such bottom particle lines is no more than $\#_{\mathrm{bond}} = i - 1 $ (since $\#_{\mathrm{bond}} = |\Mb_j| - \#_{\mathrm{comp}}$). If it is a C- or an M-constituent, then the choice of bottom edges is no more than $\mathsf{N}r$. Therefore, the total contribution of choosing such bottom edges is no more than $5i+\mathsf{N}r + 1 $.

   Assume the number of atom constituents in $\nf$ is $x_0$, hence $1 \leq x_0 \leq [l/2]$. We also note that $3 \leq l \leq \Nf$. As a result, we have
    \begin{equation}\label{eqn.B}
       \begin{aligned}
         B_{|\Mb|, \Xf',\Nf} \leq  \sum_{l =3}^{\Nf} \sum_{x_0=1}^{[l/2]} \sum_{\substack{m_1 + \cdots + m_{l} = |\Mb|-1 \\ 
             x_1 + \cdots + x_{l} = \Xf'-x_0 \\
        n_1+\cdots + n_{l} = \Nf -l}} 
             C^l  \prod_{j=1}^l Q_{m_j, x_j, n_j},
       \end{aligned}
    \end{equation} 
    where $m_j \geq 0, x_j \geq 0,n_j\geq 0$. Now we apply Lemma \ref{lem.tree-count}, to obtain that $B_{|\Mb|, \Xf',\Nf} \leq C^{|\Mb|+\Nf} |\log \eps|^{C^* \Xf'}$.
   
\textbf{Proof part 2.} Define the total number of the required molecule that is a forest is $F_{\Mb, \Xf',\Nf}$. Then
 \begin{equation}\label{eq.FMX}
    F_{|\Mb|, \Xf',\Nf} \leq \sum_{k=1}^{|\Mb|} \sum_{\substack{m_1 + \cdots + m_k =|\Mb|, m_k \geq 1 ,\\ x_1 + \cdots + x_k = \Xf', x_k \geq 0,\\
             n_1+\cdots + n_k = \Nf}}  \prod_{j=1}^k  B_{m_j, x_j,n_j}.
 \end{equation} 
For fixed $k\leq |\Mb|$, $\#\{m_1 + \cdots + m_k = |\Mb|,\ m_i \ge 1\} = \binom{|\Mb|-1}{k-1} \leq 2^{|\Mb|},$ and $\#\{x_1 + \cdots + x_k = \Xf',\ x_i \ge 0\} = \binom{\Xf' + k - 1}{k - 1} \leq 2^{|\Mb|+\Xf'}$. Similarly, $\#\{n_1 + \cdots + n_l = \Nf,n_i \ge 0\} \leq C^{\Nf}$. Moreover, the number of empty ends (these empty ends are undistinguishable) in $\Mb_{\ell'}$ is no more than $|\log \eps|^{C^*}$. Put the estimate of $B_{m_j,x_j,n_j}$ into \eqref{eq.FMX} we obtain the desired estimate.
\end{proof}

\subsection{Reducing to the integral $\Ic_\Mb(Q_\Mb)$}\label{sec.int_reduce}

In this subsection, we rewrite the quantity $|\Ic\Nc_\Mb|$ in Proposition \ref{prop.cumulant_formula} in the form $\Ic_\Mb(Q_\Mb)$, using Dirac delta distributions to localize each individual collision and overlap, see Definition \ref{def.associated_op} and Proposition \ref{prop.local_int}. As we will see in Section \ref{sec.oper}, this representation behaves well under the cutting operation, and allows us to obtain efficient estimates.
 
\begin{definition}
    Consider a molecule with a single atom with $n$ particles $\pb_1,\ldots,\pb_n$ for some integer $n\geq 2$. Denote the initial state of these particles by $\vz_n = (z_1,z_2,...z_n)$. We define 
    \begin{equation}
     (Z_1(t),Z_2(t),...Z_n(t))
    \end{equation} 
    where $Z_j(t) = (X_j(t),V_j(t)) = (X_j(t;z_1,...z_n),V_j(t;z_1,...z_n))$ to be the molecule described dynamics. We omit the dependence on $(z_1,...z_n)$ when there is no confusion.
\end{definition}

\begin{definition}[Associated operator of a molecule]\label{def.associated_op} Consider the molecule $\Mb=(\Mc,\Ec,\Pc,\Lc,\Df)$, we define the following objects. 
\begin{enumerate}
    \item\label{it.associated_vars} \emph{Associated variables $\vz_\Ec$ and $\vt_\Mb$.} We associate each atom $\nf\in \Mb$ with a time variable $t_\nf$ and each edge $e\in \Ec$ with a position-velocity vector $z_e=(x_e,v_e)$. We refer to them as \newterm{associated variables} and denote the collections of these variables by $\vz_\Ec=(z_e:e\in\Ec)$ and $\vt_\Mb=(\vt_\nf:\nf\in\Mb)$, where $\vt_\nf = (t_{ij}\in \Rb:{(i,j)\in \Df(\nf)})$, where $\Df(\nf)$ is defined in Definition \ref{def.molecule} \eqref{it.decoration}.
    \item\label{it.associated_dist} \emph{Associated distributions $\boldsymbol{\Delta}_\nf$.} Given an atom $\nf\in\Mb$, let $(e_1,...,e_n)$ and $(e_1',...,e_n')$ be the bottom edges and top edges at $\nf$ respectively, such that $e_j$ and $e_j' $ are serial for $(j=1,..,n)$. We define the \newterm{associated distribution} $\Dirac_\nf =\boldsymbol{\Delta}(z_{e_1},...,z_{e_n},z_{e_1'},...,z_{e_n'},\vt_\nf)$ as follows. 
    \begin{enumerate}
        \item If $\nf$ is a C-atom.
        \begin{enumerate}
        \item $\nf$ is a binary atom, then
        \begin{equation}\label{eq.associated_dist_C}
        \begin{aligned}
            \boldsymbol{\Delta}_\nf&:=\boldsymbol{\delta}\big(x_{e_1'} + t^+_\nf  v_{e_1'} -X_{e_1}(t^+_\nf)\big)\cdot \boldsymbol{\delta}\big(x_{e_2'} +t^+_\nf v_{e_2'} -X_{e_2}(t^+_\nf)\big) \\
            &\times \boldsymbol{\delta}\big(|x_{e_1}-x_{e_2}+t_{12}(v_{e_1}-v_{e_2})|-\varepsilon\big)\cdot \big[(v_{e_1}-v_{e_2})\cdot \nu\big]_-\\
           &\times\boldsymbol{\delta}(v_{e_1'}-V_{e_1}(t^+_\nf) ) \cdot \boldsymbol{\delta}\big(v_{e_2'}-V_{e_2}(t^+_\nf)\big).
        \end{aligned}
        \end{equation} 
        Here  $\nu_{12}:=\varepsilon^{-1}(x_{e_1}-x_{e_2}+t_{12}(v_{e_1}-v_{e_2}))$ is a unit vector, $z_-=-\min(z,0)$, and $t^+_\nf$ denotes the time that the collision ends. Moreover, we have $|t_{12}- t^+_\nf| \leq \eps^{\upsilon_0}$ if $\nf$ marked as short interaction. Specifically, for a binary collision happens inside the interval $((\ell[\nf]-1)\tau, \ell[\nf]\tau)$, $z_{e_1'}$ and $z_{e_2'}$ are given by \eqref{eq.relation_collision} below.
        \item $\nf$ is a multi-particle atom. Then
                \begin{equation}\label{eq.associated_dist_C_multi} 
        \begin{aligned}
            \boldsymbol{\Delta}_\nf  := &\prod_{i=1}^{n}\Big(\boldsymbol{\delta}\big(x_{e_i'} + t^+_\nf  v_{e_i'} -X_{e_i}(t^+_\nf)\big)\boldsymbol{\delta}\big(v_{e_i'}-V_{e_i}(t^+_\nf)\big)\Big)\\
              &\prod_{(i,j) \in \Cf(\nf)} \Big( \boldsymbol{\delta}\big(|X_{e_i}(t_{ij}) -X_{e_j}(t_{ij})|-\varepsilon\big)\cdot \big[(V_{e_i}(t_{ij})-V_{e_j}(t_{ij})) \cdot \nu_{ij}\big]_- \chi_{ij}(t_{ij})\cdot \mathbbm{1}_{|t_{ij} - t^-_\nf| \leq \eps^{\upsilon_0}} \Big).
        \end{aligned}
        \end{equation} 
        Here  \begin{equation}\label{eqn.chi_ij}
            \chi_{ij}(t_{ij}):= \mathbbm{1}_{|X_{e_i}(t) - X_{e_j}(t)| > \eps  \text{ for all } t \in [(\ell[\nf]-1)\tau, t_{ij})}. 
        \end{equation}
        $\Cf(\nf) $ (and $\Df(\nf)$, $\Of(\nf)$ below) is defined in Definition \ref{def.molecule} \eqref{it.decoration}, $t^-_\nf = \min_{(i,j)\in \Cf(\nf)} ( t_{ij}) $ and  $t^+_\nf$ is the moment when the collision starts and ends, in the sense of Definition \ref{def.admissible}. 
        In particular, $X_{e_i}(t_{ij}) = x_{e_i} + t_{ij} v_{e_i}$ if $\pb_i$ has not entered the collision before time $t_{ij}$. 
        \end{enumerate}
        \item If $\nf$ is an O-atom, then
        \begin{equation}\label{eq.associated_dist_O}
        \begin{aligned}
            \boldsymbol{\Delta}(z_{e_1},z_{e_2},z_{e_1'},z_{e_2'},t_{12})&:=\boldsymbol{\delta}(x_{e_1'}-x_{e_1})\cdot \boldsymbol{\delta}(x_{e_2'}-x_{e_2}) \\ 
            &\times \boldsymbol{\delta}\big(|x_{e_1}-x_{e_2}+t_{12}^-(v_{e_1}-v_{e_2})|-\varepsilon\big)\cdot \big[(v_{e_1}-v_{e_2})\cdot \nu_{12}\big]_-\\&\times\boldsymbol{\delta}\big(v_{e_1'}-v_{e_1}\big)\cdot\boldsymbol{\delta}\big(v_{e_2'}-v_{e_2}\big).
        \end{aligned}
        \end{equation} Note that (\ref{eq.associated_dist_O}) forces $z_{e_1}=z_{e_1'}$ and $z_{e_2}=z_{e_2'}$.
        \item If $\nf$ is a CO-atom. The CO-atom is formed by joining O-connection between several C-atoms $\nf_1,\dots,\nf_m$ or edges $e_{k_1},...e_{k_l}$. Then 
        \begin{equation}\label{eq.associated_dist_CO_multi} 
        \begin{aligned}
            \boldsymbol{\Delta}_\nf:=&\prod_{i=1}^{m}\boldsymbol{\Delta}_{\nf_i}\cdot\prod_{1\leq j\leq l}\big(\dirac(x_{e_{k_j}'} - x_{e_{k_j}}) \cdot\dirac(v_{e_{k_j}'} - v_{e_{k_j}})\big) \\
            &\times \prod_{(i,j) \in \Of(\nf)}\Big( \boldsymbol{\delta}\big(|X_{e_i}(t_{ij}) -X_{e_j}(t_{ij})|-\varepsilon\big)\cdot \big[(V_{e_i}(t_{ij})-V_{e_j}(t_{ij})) \cdot \nu_{ij}\big]_- \chi_{ij}(t_{ij}) \mathbbm{1}_{|t_{ij} - t^-_{ij}| \leq \eps^{\upsilon_0}} \Big) .
        \end{aligned}
        \end{equation}
       Moreover, $|t_{ij} - t^-_{ij}|\leq \eps^{\upsilon_0}$, where $t^-_{ij} = \min (t_{\nf'}^-,t^-_{\nf''})$, $\nf'$ and $\nf''$ are atoms containing particles $i$ and $j$, respectively. 
       In particular, if ${e_i}$ is an edge constituent at $\nf$, then $V_{e_i}(t) = v_{e_i} $, and $X_{e_i}(t) = x_{e_i} + t v_{e_i}$.  
 
    \item  If $\nf$ is an M-atom in layer $\ell' $, then,
  \begin{equation}\label{eq.associated_dist_M}
    \Dirac_\nf: = \prod_{1 \leq j \leq n}\dirac (x_{e_j'} + t^+ v_{e_j'} - X_{e_j}(t^+) )\dirac (v_{e_j'} - V_{e_j}(t^+))\cdot  \mathbbm{1}_{\Mf_\nf}(\vz_n) 
  \end{equation}
where $\Mf_\nf$ denotes the set of all possible boxed molecules $\Mb$ that can appear as the middle part of a box in Definition \ref{def.box}. More precisely, $\mathbbm{1}_{\Mf_\nf}(\vz_n)=1$ if and only if the trajectory of the $n$ particles of the M-atom, with incoming states $\vz_n$ (given by its bottom edges) and evolving by the truncated dynamics of these $n$ particles, reduces to a middle part of a box (Definition \ref{def.box}) whose bottom and top parts are the ones given in $\Mb$; and $(X_{e_j}(t^+),V_{e_j}(t^+))$ is the state of this trajectory at the time $t^+$ when the middle part ends.
    \item  If $\nf$ is a MO-atom, and is obtained by joining O-connections among several C-atoms or M-atoms $\nf_1,\cdots, \nf_m$, or bonds $e_{k_1},...e_{k_l}$. Then
\begin{equation}\label{eq.associated_dist_MO}
      \begin{aligned}
        \Dirac_\nf: = &\prod_{i=1}^m \Dirac_{\nf_i} \cdot \prod_{1\leq j\leq l}\big(\dirac(x_{e_{k_j}}' - x_{e_{k_j}}) \cdot\dirac(v_{e_{k_j}}' - v_{e_{k_j}})\big)  \\
        &\prod_{(i,j) \in \Of(\nf)} \Big(\boldsymbol{\delta}\big(|X_{e_i}(t_{ij}) -X_{e_j}(t_{ij})|-\varepsilon\big)\cdot \big[(V_{e_i}(t_{ij})-V_{e_j}(t_{ij})) \cdot \nu_{ij}\big]_- \chi_{ij}(t_{ij})\cdot\mathbbm{1}_{|t_{ij} - t^-_{ij}| \leq \eps^{\upsilon_0}} \Big) .
      \end{aligned}
\end{equation}
Moreover, we have $|t_{ij} - t^-_{ij}| \leq \eps^{\upsilon_0}$, where $t^-_{ij} = \min (t_{\nf'}^-,t^-_{\nf''})$, and $\nf'$ and $\nf''$ are C-constituents containing particles $i$ and $j$, respectively.

    \end{enumerate}
    \item\label{it.associated_domain} \emph{The associated domain $\mathcal{D}_\Mb$.} The associated domain $\mathcal{D}_\Mb$ is defined by
    \begin{equation}\label{eq.associated_domain}
        \begin{aligned}
            \Dc_\Mb:=\big\{\vt_\Mb=(\vt_\nf: \nf\in\Mb): (\ell'-1)\tau - \eps^{\upsilon_0}<t^-_\nf<\ell'\tau\mathrm{\ if\ }\nf\in\Mb_{\ell'};\\
             t_{\nf}^-<t_{\nf^+}^- +    \eps^{\upsilon_0} \mathrm{\ if\ }\nf^+\mathrm{\ is\ parent\ of\ }\nf \big\}.
        \end{aligned}
    \end{equation} 
    \item\label{it.associated_int} \emph{The associated operator $\Ic_\Mb$.}  With the associated variables $\vz_\Ec=(z_e:e\in\Ec)$ and $\vt_\Mb$, we define the \newterm{associated integral operator} $\Ic_\Mb$ by 
    \begin{equation}\label{eq.associated_int_op}
        \begin{aligned}
            \Ic_\Mb(Q):=\varepsilon^{-(d-1)(|\Ec_*|-\sum_{l \geq 2} l|\Mb_l|)}\int_{\mathbb{R}^{2d|\Ec_*|} \times \Dc^\Mb} 
            \prod_{\nf\in\Mb}\boldsymbol{\Delta}_\nf \ Q(\vz_\Ec,\vt_\Mb) \,\mathrm{d}\vz_{\Ec_*} \,\d \vt_{\Mb} . 
        \end{aligned}
    \end{equation}   
    
    Here $\Ec_*$ is the set of edges of $\Mb$ that are \emph{bonds or free ends}, and $|\Mb_l|$ denotes the number of atoms (including M- and MO-atoms) that contains $l$ particles. Note that in (\ref{eq.associated_int_op}) we are only integrating in the variables $\vz_{\Ec_*}$ for bonds and free ends, thus $\Ic_\Mb(Q)$ is a function (or distribution) of \emph{all fixed ends of $\Mb$}.
\end{enumerate}
\end{definition} 
We have the following Proposition regarding binary collision, see Section 8.2 of \cite{GST14}.
\begin{proposition}\label{prop.scattering_map}
   In \eqref{eq.associated_dist_C}, we have
    \begin{equation}\label{eq.relation_collision}
        \begin{aligned}
        x_{e_1'} + (t_{12}+ \Delta t)  v_{e_1'} &= x_{e_2} + t_{12} v_{e_2}  + \eps(\omega_{12}\cdot \nu_{12}) \omega_{12} + \frac{\Delta t}{2}(v_{e_1}+v_{e_2}), \\
        x_{e_2'} +(t_{12}+\Delta t) v_{e_2'} &= x_{e_1} + t_{12} v_{e_1} -\eps(\omega_{12}\cdot \nu_{12}) \omega_{12} + \frac{\Delta t}{2}(v_{e_1}+v_{e_2}), \\
        v_{e_1'} &:= v_{e_1}- \omega_{12}\cdot (v_{e_1}-v_{e_2})\omega_{12}, \\
        v_{e_2'} &:= v_{e_2} + \omega_{12}\cdot (v_{e_1}-v_{e_2})\omega_{12},
        \end{aligned} 
    \end{equation}
    where $\nu_{12} = \eps^{-1} \big( x_{e_1} -x_{e_2} + t_{12} (v_{e_1} - v_{e_2})\big)$, and $\omega_{12}= \omega_{\Phi}(v_{e_1} -v_{e_2}, \nu)$ is the line of apsides of the collision. Moreover, $\Delta t$ is the collision duration, as in Definition \ref{def.short_range} \eqref{it.O5}, which implicitly depends on $v_{e_1} -v_{e_2}$ and $\nu$. We refer the reader to Lemma 8.1.2 in \cite{GST14} for more details.
\end{proposition}

We now state the main result of this subsection.
\begin{proposition}\label{prop.local_int} Let $\Mb$ be a full molecule (so $\Ec_*=\Ec$) with $r(\Mb) = H$. Let $|\Ic\Nc_\Mb|$ be the associated integral of $\Mb$ defined in \eqref{eq.associated_integral_molecule_abs}. Then, we have
\begin{equation}\label{eq.local_int}
    \|\, |\Ic\Nc_\Mb|(\vz_{H})\|_{L^1} \leq  \varepsilon^{(d-1)|H|}\cdot\Ic_\Mb(Q_\Mb).
\end{equation}
where $Q_\Mb$ is the function given by 
\begin{equation}\label{eq.Q_M}
    Q_\Mb = \mathbbm{1}_{\Dc_\Mb}(\vt_\Mb) \cdot \prod_{e\in \Ec_{\mathrm{end}}^-}\big|f^\Ac\big((\ell_1[e]-1)\tau,x_e-(\ell_1[e]-1)\tau\cdot v_e,v_e\big)\big|.
\end{equation}
Here in \eqref{eq.Q_M}, $\Ec_{\mathrm{end}}^-$ is the set of bottom ends of $\Mb$ (see Definition \ref{def.more_on_molecule} \eqref{it.subsets_Ec}), and $\ell_1[e]$ is the initial layer for bottom end $e$ in Definition \ref{def.molecule} \eqref{it.particle_line_same}. The indicator functions $\mathbbm{1}_{\Dc_\Mb}$ is defined in \eqref{eq.associated_domain}.
\end{proposition}

\begin{proof} First we prove $|\Mb|_p = |\Ec| - \sum_{l \geq 2} l |\Mb_l|$. This is because each atom in $\Mb_l$ has exactly $2 l $ edges; by adding them together and noticing that each bond (between two atoms) is counted twice and each end is counted once, we obtain $\sum_{l \geq 2} 2l|\Mb_l| = 2|\Ec| - |\Ec_{\mathrm{end}}^-| - |\Ec_{\mathrm{end}}^+|$. Using also that $|\Ec_{\mathrm{end}}^-| = |\Ec_{\mathrm{end}}^+|=|\Mb|_p$ (as each particle line has a unique top end and a unique bottom end, see Definition \ref{def.molecule} (\ref{it.particle_line_same})), we conclude that $|\Mb|_p = |\Ec| - \sum_{l \geq 2} l |\Mb_l|$, as desired.

By $|\Mb|_p = |\Ec| - \sum_{l \geq 2} l|\Mb_l|$, and note that $|H|=|\Mb|_r$, we see that the powers of $\varepsilon$ on both sides of (\ref{eq.local_int}) match. 
Therefore, in the definition of $|\Ic\Nc_\Mb|$ in (\ref{eq.associated_integral_molecule_abs}), by integrating in the $\vz_{r(\Mb)}$ variables, using the definition $|(\Sc\circ\mathbbm{1})_\Mb|$ in \eqref{eq.associated_op} and the $L^1$ preserving property of $\Sc_\Mb$, we see that
\begin{equation}\label{eq.proof_local_int_*}
\begin{aligned}
    \big\||\Ic\Nc_\Mb|\big\|_{L^1}= \varepsilon^{-(d-1)|\Mb|_{p\backslash r}}\int_{\Rb^{2d|\Mb|_p}}\mathbbm{1}_\Mb\cdot (|Q(\vz_{p(\Mb)}')|)\,\mathrm{d}\vz_{p(\Mb)}',
\end{aligned}
\end{equation} where $Q$ is as in (\ref{eq.extraQ}). We may identify a particle line $\pb\in p(\Mb)$ with its bottom end $e\in\Ec_{\mathrm{end}}^-$; as such, we see that (\ref{eq.local_int}) would follow if we can prove 
\begin{equation}\label{eq.proof_local_int_1_*}
\int_{\Rb^{2d|\Ec\backslash \Ec_{\mathrm{end}}^-|} }\mathbbm{1}_{{\Dc_\Mb}}(\vt_\Mb)\cdot \prod_{\nf\in\Mb}\boldsymbol{\Delta}(z_{e_1},z_{e_2},z_{e_1'},z_{e_2'},t_\nf)\,\mathrm{d}\vz_{\Ec\backslash\Ec_{\mathrm{end}}^-}  \geq \mathbbm{1}_{\Mb}\big(\vz_{\Ec_{\mathrm{end}}^-}\big),\quad a.e.,
\end{equation}
where we identify $\pb\in p(\Mb)$ with the bottom end $e\in\Ec_{\mathrm{end}}^-$ and identify $\vz_\pb'$ in (\ref{eq.proof_local_int_*}) with $z_e$. This can be proved as in Proposition 7.5 in \cite{DHM24}. For binary atoms the argument is exactly the same. For the other atoms we argue as follows, by induction starting from the lowest atoms.
For a multi-particle C-atom $\nf$, with its incoming variables fixed, integrating \eqref{eq.associated_dist_C_multi} in the outgoing variables $(x_{e_i'},v_{e_i'})$ evaluates the corresponding Dirac functions; for each $(i,j)\in\Cf(\nf)$, thanks to the factor $\chi_{ij}$, the integral in $t_{ij}$ of $\dirac(|X_{e_i}(t_{ij})-X_{e_j}(t_{ij})|-\eps)[(V_{e_i}-V_{e_j})\cdot\nu_{ij}]_-\chi_{ij}$ only counts the first time when $\pb_i,\pb_j$ reach distance $\eps$, so it equals $1$ (outside a null set) if this pairwise collision happens; finally, in the support of $\mathbbm{1}_\Mb$ the atom is marked as short interaction, so the time indicators $\mathbbm{1}_{|t_{ij}-t^-_\nf|\leq\eps^{\upsilon_0}}$ are equal to $1$. Hence $\int\Dirac_\nf\prod_i\d x_{e_i'}\d v_{e_i'}\d\vt_\nf\geq\mathbbm{1}_\nf$, where $\mathbbm{1}_\nf$ is the indicator that the collision of $\nf$ happens as described by its decoration. 
 For an M-atom, integrating \eqref{eq.associated_dist_M} in the outgoing variables gives $\mathbbm{1}_{\Mf_\nf} = \mathbbm{1}_\nf$. The CO- and MO-atoms are treated in the same way, applying the above to their constituents and to the factors of the O-connections (whose time variables are integrated in the same way). Applying this inductively proves \eqref{eq.proof_local_int_1_*}. 

\end{proof}
 
\section{Operations on molecules}

\subsection{Cutting, deleting and splitting} \label{sec.oper}
In this subsection, we introduce the operations of cutting, deletion, and splitting. These definitions basically follow those in \cite{DHM24}, with suitable modifications to incorporate partial cutting. 

\begin{definition}[Deleting an O-connection]\label{def.delete} Let $\Mb$ be a molecule and \(\mathfrak o\) be an \(O\)-connection in \(\Mb\), joining two objects \(\af\) and \(\af'\), where each of \(\af,\af'\) may be an edge, a \(C\)-atom, or an M-atom. By deleting $\of$, we remove the \(O\)-connection \(\mathfrak o\) from \(\Mb\), while keeping \(\af\) and \(\af'\) unchanged. In particular, deleting an \(O\)-connection does not create or remove any free end. 
\end{definition}
 
\begin{definition}\label{def.reg} We say $\Mb$ is \textbf{regular}, if for each atom $\nf \in \Mb$,\begin{enumerate}
    \item\label{it.regular1} any O-fixed end at $\nf$ must be serial with an O-free end $\nf$ (we say they form a \textbf{simple pair}), and
    \item\label{it.regular2}  the C-fixed ends are on the same side of $\nf$, i.e., are all top or bottom.
\end{enumerate}
This property will be assumed for all the molecules studied below, and will be preserved throughout the reduction process.
\end{definition}
We remark that property \eqref{it.regular1} is preserved by Lemma \ref{lem.regularity} below, and property \eqref{it.regular2} will be preserved by the choice of algorithms in Sections \ref{sec.cuttingI}--\ref{sec.error}.
\begin{definition}[Cutting]\label{def.cutting} Let $\Mb$ be a regular molecule and $A\subseteq\Mb$ be a set of atoms. The operation of \newterm{cutting} $A$ \newterm{as free} (resp. $\Mb\backslash A$ \newterm{as fixed}) turns $\Mb$ into a new molecule $\Mb'$ which is the disjoint union of two molecules, i.e. $\Mb'=\Mb_1\sqcup\Mb_2$. For convenience, we may also say that $\Mb$ is ``cut into two smaller molecules $\Mb_1$ and $\Mb_2$". We describe the procedure as follows. 
\begin{enumerate}
\item\label{it.cutting_1} The atom sets of $\Mb_1$ and $\Mb_2$ are $A$ and $\Mb\backslash A$ respectively. If $\Mb$ contains no O-connections (so each maximal ov-segment is just a bond), then we simply break each bond $e$ between an atom in $A$ and an atom in $\Mb\backslash A$, by removing it and drawing a free end $e_1$ (resp. fixed end $e_2$) at the corresponding atom in $A$ (resp. $\Mb\backslash A$). We form $\Mb_1$ (resp. $\Mb_2$) by atoms in $A$ (resp. $\Mb\backslash A$) and these edges.
\item\label{it.cutting_2} Now, suppose $\Mb$ contains O-connections. For each maximal ov-segment $\sigma$ that contains an atom in $A$ and an atom in $\Mb\backslash A$ (other maximal ov-segments are not affected), let the highest and lowest atoms of $\sigma$ be $\nf^+$ and $\nf^-$, and let $(\pf_1,\cdots,\pf_s)$ be the atoms in $A\cap\sigma$ (ordered high to low). Then, we take the following steps, which extend the steps in (\ref{it.cutting_1}) from bonds to ov-segments:
\begin{enumerate}
\item\label{it.cutting_2a} First remove every bond in $\sigma$. Then, for $1\leq j\leq s-1$, draw a bond between $\pf_j$ and $\pf_{j+1}$. Then, if $\pf_1\neq\nf^+$, draw a top free end at $\pf_1$; if $\pf_s\neq\nf^-$, draw a bottom free end at $\pf_s$. Let $\Mb_1$ be formed by atoms in $A$ and these edges (constructed from different $\sigma$).
\item\label{it.cutting_2b} Next, if $\nf^+\in\sigma\backslash A$, draw a bottom fixed end at $\nf^+$; if $\nf^-\in\sigma\backslash A$, draw a top fixed end at $\nf^-$. Then, for any $\nf^\pm\neq \nf\in\sigma\backslash A$, draw a simple pair at $\nf$ consisting of a free end and a (serial) fixed end. Let $\Mb_2$ be formed by atoms in $\Mb\backslash A$ and these edges (constructed from different $\sigma$).
\end{enumerate}
\end{enumerate}
\end{definition}

\begin{definition}[Partial cutting]\label{def.partial_cutting}
Let $\Mb$ be a regular molecule. Let \(\nf\in \Mb\) be a  CO-atom or an  MO-atom formed by joining objects
\(\af_1,\ldots,\af_m\) through O-connections. 
\begin{enumerate}
\item\label{it.partial_cutting1} Suppose $\af_i$ is a C-constituent or an M-constituent. We \emph{partially cut $\af_i$} by deleting every O-connection $\of=(\af_i,\af_k)$ incident to $\af_i$. For each deleted $\of$, we draw a simple pair $\sigma_{\of}$ at $\af_k$ and add the O-connection $(\af_k,\sigma_{\of})$. Then we cut $\af_i$ as free from $\Mb$. This step may break one atom $\nf$ into several atoms. 
\item\label{it.partial_cutting2} Suppose that \(\af_i\) and \(\af_j\) are \(O\)-connected. We \emph{partially cut of \(\af_i\) and \(\af_j\) together} by deleting only the O-connections
\[
    \of=(\af_l,\af_k),\qquad l\in\{i,j\},\quad k\notin\{i,j\}.
\]
For each such deleted connection, introduce an auxiliary simple pair
\(\sigma_{\of}\) at \(\af_k\) and add \((\af_k,\sigma_{\of})\).  The O-connection \((\af_i,\af_j)\) is kept unchanged. Let the new formed atom that contains $\af_i,\af_j$ be $\nf'$. We cut $\nf'$ as free.
\item\label{it.partial_cutting3} Let $A$ be a set of atoms and $\nf\in A$ be a CO-atom or an MO-atom. Assume that there exists exactly one atom constituent $\af_i$ that still connects to the other atoms of $A$, while no other constituent of $\nf$ is. We say cut $A$ as free but partially at $\nf$ with $\af_i$, means we first delete every O-connections $\of=(\af_i,\af_k)$ and add simple pairs as in \eqref{it.partial_cutting1}. Then cut each atom $(A\setminus \nf)\cup \af_i$ as in Definition \ref{def.cutting}
\item\label{it.partial_cutting_box} Let $A$ be a collection of atoms, C-constituents and M-constituents (for example, $A$ can be a cluster). We say (partially) cut $A$ as free, means the following operation: first delete all O-connections joining a C- or an M-constituents in $A$ to a constituent outside $A$; the affected constituents in $A$ then become new independent atoms in $A$. For each deleted O-connection, add a simple pair to the corresponding constituent outside $A$. Finally, cut the resulting atoms associated with $A$ as free. 
\end{enumerate}
\end{definition}
 
\begin{proposition}\label{prop.cutting} We define suitable mapping (as in Definition 8.12 \cite{DHM24}) to match variables before and after the operation ((partial) cutting and deleting). In particular, for a special atom $\nf$ formed by partially cutting $\af_i$ consists of a constituent $\af_j$ and an auxiliary simple pair $(e,e')$, we define its distribution by
    \begin{equation}
        \Dirac =\Dirac_{\af_j} \cdot \dirac(v_{e'} - v_e) \dirac(x_{e'} - x_e) \dirac(|X_{e_i}(t_{ij}) - X_{e_j}(t_{ij})| - \eps) [(V_{e_i} -V_{e_j})\cdot \nu_{ij}]_- \chi(t_{ij})\mathbbm{1}_{|t_{ij} - t^-| \leq \eps^{\upsilon_0}} ,
    \end{equation}
    where $\nu_{ij} : = \eps^{-1} (X_{e_i} (t_{ij} - X_{e_j}(t_{ij})))$. Here $X_{e_i(\cdot)}$ denotes the (curved) trajectory of some constituent cut before $\nf$, which is a function of the variables of the molecule $\Mb_1$ containing $\af_i$, and is fixed in $\Ic_{\Mb_2}$. The indicator function $\mathbbm{1}_{(\cdots)}$ is omitted if $\af_i$ is an M-constituent. 
    
    With the above definitions, we have the followings.
     \begin{enumerate}
    \item Let $\Mb'$ be formed by deleting one O-connection. Then we have
       \begin{equation}
           \Ic_{\Mb}(Q) \leq \Ic_{\Mb'}(Q), \quad \forall Q \geq 0.
       \end{equation}
    \item Assume that a molecule $\Mb$ is cut into $\Mb_1$ and $\Mb_2$, where $\Mb_1$ and $\Mb_2$ are cut as free and fixed respectively. Then we have
       \begin{equation}\label{eq.cutting}
           \Ic_{\Mb}(Q) = \Ic_{\Mb_1}\circ \Ic_{\Mb_2} (Q).
       \end{equation}
     \end{enumerate}
\end{proposition}
 \begin{proof}
    Recall the formula \eqref{eq.associated_int_op}. Note that the quantity $|\Ec_*| - \sum_{l \geq 2} l |\Mb_l|$ denotes the number of particle lines that both start and end at free ends. None of the above operations changes this quantity. 
    
    Moreover, by defining the map as in Definition 8.12 \cite{DHM24}, $\Ic_{\Mb}$ is the same as $\Ic_{\Mb'}$ except with an extra indicator function $\mathbbm{1}_{\mathrm{precol}} (Z_{e_i}(t_{ij}),Z_{e_j}(t_{ij}),t_{ij})$, which is $\leq 1$. Hence, $\Ic_{\Mb}(Q) \leq \Ic_{\Mb'}(Q)$. Moreover, cutting and partial cutting merely change the order of integration, so \eqref{eq.cutting} follows directly from Fubini's theorem.
 \end{proof}

\begin{definition}[Ordering components after cutting]\label{def.ordered_disjoint} Suppose a molecule $\Mb$ becomes a new molecule $\Mb'$ after a sequence of cutting operations; we refer to the connected components of $\Mb'$ as final components. We define a total ordering $\prec_{\mathrm{cut}}$ between these final components inductively, as follows:
\begin{enumerate}
    \item\label{it.disjoint_1} In the first cutting operation, suppose $\Mb$ is cut into $\Mb_1$ (as free) and $\Mb_2$ (as fixed). Then for any final component $\Xc,\Yc$ of $\Mb'$, we require that $\Xc\prec_{\mathrm{cut}}\Yc$ if $\Xc\subseteq\Mb_1$ and $\Yc\subseteq\Mb_2$.
    \item For the final components that are contained in $\Mb_1$ (resp. in $\Mb_2$), order them inductively using $\prec_{\mathrm{cut}}$ defined for the subsequent sequence of cutting operations applied to $\Mb_1$ (resp. to $\Mb_2$).
    \end{enumerate}
    
    If $\Mb'$ is the result of $\Mb$ after (deleting and) cutting, and $\Mb'$ is the disjoint union of its final components, i.e. $\Mb' = \Mb_1 \sqcup \cdots \sqcup \Mb_k$ with $\Mb_1 \prec_{\mathrm{cut}} \cdots \prec_{\mathrm{cut}} \Mb_k$, then we define the \textbf{associated integral operator} by
    \begin{equation}\label{eq.associated_op_disjoint}
        \Ic_{\Mb'} = \Ic_{\Mb_1}\circ \Ic_{\Mb_2}\circ\cdots\circ \Ic_{\Mb_k}.
    \end{equation}
\end{definition}

\begin{definition}[Ov-notions]\label{def.ov_connect} Recall molecules $\Mb$ and ov-segments in Definition \ref{def.molecule}. 
\begin{enumerate}
\item \emph{Ov-adjacency.} We define two atoms $\nf_1$ and $\nf_2$ to be \textbf{ov-adjacent}, if they are connected by an ov-segment; if $\nf_1$ is an ancestor (resp. descendant) of $\nf_2$ along this ov-segment, we say it is an \textbf{ov-parent} (resp. \textbf{ov-child}) of $\nf_2$.
\item \emph{Ov-connectedness.} For a subset $A\subseteq\Mb$, we say it is \textbf{ov-connected} if for any two atoms $\nf,\nf'\in A$ there exist atoms $\nf_j\in A \,(0\leq j\leq q)$ such that $(\nf_0,\nf_q)=(\nf,\nf')$ and $\nf_j$ is ov-adjacent with $\nf_{j+1}$. The maximal ov-connected subsets of $A$ are called \textbf{ov-components} of $A$.
\item\emph{Ov-distance.} Finally, define the \textbf{ov-distance} of two subsets $B,B'\subseteq \Mb$ to be the minimal $q$ such that there exist atoms $\nf_j\in \Mb\,(0\leq j\leq q)$ such that $(\nf_0,\nf_q)\in B\times B'$ and $\nf_j$ is ov-adjacent with $\nf_{j+1}$. Note that the ov-distance is $0$ if $B\cap B'\neq\varnothing$, and is $\infty$ if they are contained in different connected components\footnote{The ov-components of the \emph{whole molecule} $\Mb$ are the same as the components. However the ov-components of \emph{subsets} $A\subseteq\Mb$ might not be the same as its components.} of $\Mb$.
\end{enumerate}  
\end{definition}

\begin{lemma}[]\label{lem.regularity}
    The property in Definition \ref{def.reg} \eqref{it.regular1} is preserved after deleting an O-connection or cutting a subset $A$ as fixed or free, or partially cut an atom. 
\end{lemma}
\begin{proof} Simply note that, the forbidden configuration, where a O-fixed end is serial with another fixed end or a bond, cannot be constructed in the deleting or cutting process. This can be directly verified for each step in Definitions \ref{def.delete}, \ref{def.cutting} and \ref{def.partial_cutting}.
\end{proof}

\begin{definition}\label{def.splitting} In \textbf{splitting}, we do not change $\Mb$ but split the support of $Q$ in (\ref{eq.associated_int_op}) as follows: decompose $1$ into finitely many indicator functions $1=\sum_{j}\mathbbm{1}_{\Sc_j}$, and decompose $Q=\sum_j Q\mathbbm{1}_{\Sc_j}$. Each of these terms corresponds to a \textbf{sub-case} where we make extra \textbf{restrictions} on the support of $Q$ provided by the conditions defining $\Sc_j$.
\end{definition}

\begin{definition}[Operation sequence]\label{def.cutting_algorithm} Given a molecule $\Mb$ and $Q$, define an \textbf{operation sequence} to be a sequence of operations acting on $\Mb$ or $Q$, such that each one is either \textbf{(partial) cutting} a subset of $\Mb$ as in Definitions \ref{def.cutting} and \ref{def.partial_cutting}, or \textbf{deleting} an O-connection in $\Mb$ as in Definition \ref{def.delete}, or \textbf{splitting} into finitely many sub-cases (by decomposing $1$ into indicator functions) as in Definition \ref{def.splitting}. If an operation sequence does not contain splitting, then we call it a \textbf{cutting sequence}.
\end{definition}

Finally, we define the elementary molecules, which are the result of the cutting sequences we will construct in the algorithms below, and the main object of study in Section \ref{sec.treat_integral}:
\begin{definition}[Elementary molecules]\label{def.elementary} We define a molecule $\Mb$ to be elementary, if it is regular, and satisfies one of the properties:
\begin{enumerate}
\item   {\bf \{$k_l$\} ($2 \leq l \leq k \leq 2l$)  molecules}: $\Mb$ contains a single atom with $k$ free ends and $l$ particles. If the single atom is an M-atom, then we must have $k=l$ (this is guaranteed in Sections \ref{sec.cuttingI}--\ref{sec.maincr} by Proposition \ref{prop.alg_up} \eqref{it.up_proof_1}: each M-atom is cut either as fully fixed, or with no fixed end together with a simple pair; in Section \ref{sec.error}, the algorithms \textbf{UP3} and \textbf{DOWN3} may also cut an M-atom with $k>l$, or together with a $\{3_2\}$ atom, and such molecules are also regarded as elementary there, see Remark \ref{rem.up3}); if the single atom is an MO-atom, then it must be formed by joining an O-connection between an M-atom and a simple pair, and  the M-constituent has no fixed end (so $k=2l-1$). \textbf{$\{2_2\}, \{3_2\}, \{4_2\}$ molecules} are specific examples. We also call an atom \(\nf\) in an arbitrary molecule \(\Mb\) a \(\{k_l\}\) atom if the total number of its free ends and bonds is \(k\), and it contains \(l\) particles. Moreover, we define {\bf special $\{3_2\}$ atom} to be those $\{3_2\}$ atoms originally from a multi-particle atom, and formed by adding an auxiliary simple pair during the partial cutting process. Each such special $\{3_2\}$ atom is regarded as an O-atom (and we do not distinguish its top/bottom edges). A molecule containing a single special $\{3_2\}$ atom is called a {\bf special $\{3_2\}$ molecule}. We also recall that an atom such that $k = l$ is a fully fixed atom.
 \item {\bf $\{3_23_2\}$ molecules}: $\Mb$ contains two atoms connected by a bond, which are both $\{3_2\}$. We say $\Mb$ is \textbf{$\{3_23_2A\}$ molecule} if we can cut one atom as free such that the other atom becomes a $\{2_2\}$ molecule as defined above. Otherwise we say $\Mb$ is \textbf{$\{3_23_2B\}$ molecule}. Moreover, if one or two atoms are special $\{3_2\}$ atoms, then we say it is a special $\{3_2 3_2\}$ molecule. Note that each special atom $\{3_2 3_2\}$ molecule must be a $\{3_2 3_2A \}$ molecule.
\item\label{it.elementary_32kl}  {\bf $\{3_2 k_l\}$ $(l \geq 2)$ molecules}: $\Mb$ contains two atoms that are connected (by one or two bonds). One atom is a $\{3_2\}$ atom while the other one is a $\{k_l\}$ atom. Moreover, if we can cut the $\{3_2\}$ atom as free, then the other atom is still regular, and becomes a fully fixed atom. In particular, if $l=2$ then we must have $k=4$, and there exists a double bond between these two atoms. Moreover, we have $k=l+1$ if the $\{3_2\}$ atom has exactly one bond to the $\{k_l\}$ atom, and $k= l+2$ if $\{3_2\}$ atom has a double bond to the $\{k_l\}$ atom.  We also call it a special molecule if the $\{3_2\}$ atom is special.
\item\label{it.elementary_422l2}{\bf $\{4_2(2l)_l\}$ molecules} $\Mb$ contains two connected atoms. One is a $\{4_2\}$ atom and the other is a $\{(2l)_l\}$ atom. No atom is a special atom. In particular, if there exists a double bond between these two atoms, then $l \neq 2$. 
\end{enumerate}
\end{definition}

\begin{remark}\label{rem.layer_interval} Let $A$ be a layer interval $\Mb_{[\ell_1:\ell_2]}$ (Definition \ref{def.more_on_molecule} \eqref{def.layer_interval}), or any union of connected components of a layer interval, or the same for thin layers (defined in Definition \ref{def.layer_refine} below). It is easy to see that: for any atoms $\nf,\nf'\in A$ that belong to the same particle line, and any atom $\mf$ in the same particle line between $\nf$ and $\nf'$, we must have $\mf\in A$. This implies that (i) any particle line $\pb$ can intersect at most one component of $A$, and (ii) if two particle lines $\pb_1,\pb_2$ are both connected to a particle line $\pb$ via $A$ (Definition \ref{def.connectedvia}), then $\pb_1$ and $\pb_2$ are also connected via $A$.
\end{remark}
\begin{remark}\label{rem.full_cut} The following fact will be convenient to note: suppose we cut $A\subseteq\Mb$ \emph{as free}, then this may divide the remaining set $\Mb\backslash A$ into connected components; however, any \emph{full} component of $\Mb\backslash A$ must also be a full component of $\Mb$ before cutting $A$ as free (in other words, any \emph{newly created} connected component must have at least one fixed end). This is because, when any bond or ov-segment is broken in cutting $A$ as free, we get a fixed end at \emph{every atom} in this ov-segment that belongs to $\Mb\backslash A$. 
\end{remark}

\section{Estimating the integral $\Ic_\Mb(Q_\Mb)$}\label{sec.treat_integral}

In this section, we complete the proof of Proposition \ref{prop.cumulant_est}, pending the proof of the pure combinatorial Proposition \ref{prop.comb_est} (which is in Sections \ref{sec.cuttingI}--\ref{sec.maincr}).  The operation sequence is constructed in Proposition \ref{prop.comb_est} using the algorithms in Sections \ref{sec.cuttingI}--\ref{sec.maincr}.
\subsection{Integrals for elementary molecules}\label{sec.elem_int} In this subsection, we prove upper bounds for $\Ic_{\Mb}$ for an elementary molecule $\Mb$. 
\begin{proposition}[Upper bounds for one-atom $\Ic_\Mb$]\label{prop.est_single} Let $\Ic_\Mb$ be defined as in \eqref{eq.associated_int_op}, where $\Mb$ is elementary molecule defined in Definition \ref{def.elementary}, and contains only one atom $\nf$ marked as short interaction, with $2l$ ($l \geq 2$) edges $e_1,\dots,e_{2l}$ and $l\leq k\leq 2l$ free ends. We assume that $e_j$ ($1\leq j \leq l$) are bottom edges of $\nf$, and each $e_j$ is serial with $e_{l+j}$. Denote by $\pb_i\ (1\leq i\leq l)$ the particles corresponding to $e_i$. Let $Q$ be a function supported in $|v_j|\leq|\log\varepsilon|^{C^*}$ and $|x_j-x_j^0|\leq |\log\varepsilon|^{C^*}$, where $(x_j,v_j)=z_{e_j}$ and $x_j^0$ are fixed vectors, and also $\vt_\Mb\in\Dc_\Mb$ with $\Dc_\Mb$ as in (\ref{eq.associated_domain}).
\begin{enumerate}
\item\label{it.est_full_fixed}
Suppose $\Mb$ is a \{$l_l$\} molecule as in Definition \ref{def.elementary}, then 
\begin{equation}\label{eq.int_elem_est_deg2}
    \|\Ic_\Mb(Q)\|_{L^\infty}\leq \|Q\|_{L^\infty}.
\end{equation}
\item\label{it.est_3} Suppose $\Mb$ is a \{$3_2$\} molecule (not special) as in Definition \ref{def.elementary}, by symmetry we may assume $e_1$ is fixed and $(e_2,e_3,e_4)$ are free. Then
\begin{equation}\label{eq.intmini_3}
\Ic_\Mb(Q)  = \int_{\Dc_\Mb\times\Sb^{d-1}\times\Rb^d}\big[(v_1-v_2)\cdot\nu_{12}\big]_-\cdot |Q|\,\mathrm{d}t_{12}\mathrm{d}\nu_{12}\mathrm{d}v_2.
\end{equation} 
Here in \eqref{eq.intmini_3}, $(z_1,z_2,z_3,z_4)$ in the input variables of $Q$ are functions of the integral variables $(t_{12}, \nu_{12}, v_2)$ given by \eqref{eq.relation_collision}, with $(x_1,v_1)$ fixed.

In addition, assume $Q$ is supported in some set depending on some external parameters $(x^*,t^*)$, and the support satisfies one of the following:
\begin{enumerate}[(a)]
\item\label{it.est31} $|t_{12}-t^*|\leq\eps^{\upsilon}$;
\item\label{it.est32} for some $i\neq j$, or for some $j\neq1$, respectively,
\begin{equation}\label{eq.restriction.3}
    \min\big(|x_i - x_j|, |v_i - v_j|\big) \leq \eps^{C_\Phi\upsilon}
    \quad\text{or}\quad
    \min\big(|v_j - v^*|, |x_j - x^*|\big) \leq \eps^{C_\Phi \upsilon};
\end{equation}
\item\label{it.est33} $[t^-_\nf,t^+_\nf] \subset ((\ell[\nf]-1)\tau,\ell[\nf] \tau)$, and
    \begin{equation}
        \frac{|(v_1-v_2)\wedge \nu_{12}|}{|v_1-v_2|} \leq \eps^\upsilon,
    \end{equation}
     where $\nu_{12} = \eps^{-1} (x_{1} -x_{2} + t_{\nf}^-(v_{1} -v_{2}))$ (this ratio is equal to $|(v_3-v_4)\wedge \nu_{34}|/|v_3-v_4|$ with $\nu_{34} =  \eps^{-1} (x_{3} -x_{4} + t^+_\nf(v_{3} -v_{4}))$, see the proof).
\end{enumerate}
 Then in each case, we have
$$ \|\Ic_\Mb(Q)\|_{L^\infty} \leq \eps^{\upsilon} \cdot |\log\varepsilon|^{C^*} \cdot \|Q\|_{L^\infty}.$$
\item\label{it.est_special_3} Suppose $\Mb$ is a special $\{3_2\}$ atom, then 
\begin{equation}\label{eq.intmini_3_special}
    \Ic_\Mb(Q)  = \int_{\Dc_\Mb\times\Sb^{d-1}\times\Rb^d}\big[(V_1(t_{12})-v_2)\cdot\nu_{12}\big]_-\cdot |Q| \mathrm{d}t_{12}\mathrm{d}\nu_{12}\mathrm{d}v_2,
\end{equation}
where $\nu_{12} = \eps^{-1}(X_1(t_{12}) - (x_2 +   t_{12} v_2))$. Moreover, assume $Q$ is supported in some set depending on external parameters $t^*$, and the support satisfies
\begin{equation}
    |t_{12} - t^*| \leq \eps^{\upsilon}.
\end{equation}
Then $$ \|\Ic_\Mb(Q)\|_{L^\infty} \leq \eps^{\upsilon} \cdot |\log\varepsilon|^{C^*} \cdot \|Q\|_{L^\infty}.$$
\item\label{it.est_kl} Let $\Mb$ be a (special) $\{k_l\}$ molecule (not an M- or an MO-atom) with $l\ge3$ and $l < k\leq 2l$, whose atom $\nf$ is marked as short interaction. Then $$\|\Ic_\Mb(Q)\|_{L^\infty}\le|\log\eps|^{C^*}\|Q\|_{L^\infty}$$
if $l<k<2l$, and $$\|\Ic_\Mb(Q)\|_{L^\infty}\le\eps^{-(d-1)+(l-2)\upsilon_0}|\log\eps|^{C^*}\|Q\|_{L^\infty}$$
if $k=2l$. Moreover, if $l<k<2l$ and one of the following holds:
\begin{enumerate}[(a)]
\item\label{it.kgeql+2} $k\ge l+2$; or
\item\label{it.k=l+1} $k=l+1$, $\Mb$ is not a special molecule and $Q$ is supported in $|v_{e_i}-v_{e_j}|\ge\eps^{\upsilon_0/2}$ for every pair of fixed ends $e_i\neq e_j$ at $\nf$,
\end{enumerate}
then
\begin{equation}\label{eq.int_elem_est_general}
\|\Ic_\Mb(Q)\|_{L^\infty}\le\eps^{\upsilon_0/2}\,|\log\eps|^{C^*}\,\|Q\|_{L^\infty}.
\end{equation} 
In particular, if $k \geq l+2$, we have $$\|\Ic_\Mb(Q)\|_{L^\infty}\le\eps^{\upsilon_0}\,|\log\eps|^{C^*}\,\|Q\|_{L^\infty}.$$
\item\label{it.est_kl2} Suppose $\Mb$ is a $\{(l+1)_l\} \ (l\geq 3)$ molecule with a C-atom, and $Q$ is supported in  $ |v_{e_i} - v_{e_j}| \geq \eps^{1 - \upsilon}$ for every pair of fixed ends $e_i \neq e_j$ at $\nf$, then
\begin{equation}\label{eq.int_elem_est_general_degeneracy}
    \|\Ic_\Mb(Q)\|_{L^\infty}\leq \varepsilon^{\upsilon}|\log\varepsilon|^{C^*}\cdot\|Q\|_{L^\infty},
\end{equation}

\item\label{it.est_3_extra} In \eqref{it.est_3} above, suppose we do not assume (a)--(c). Instead, we assume that $[t^-_\nf,t^+_\nf] \subset ((\ell[\nf]-1)\tau,\ell[\nf] \tau)$ and $Q$ is supported in the set  (depending on some external parameters $(x^*,v^*)$)  
\begin{equation}\label{eq.est_3_extra}
    \inf_{|t| \leq |\log \eps|^{C^*}} |x_j - x^* - t(v_j-v^*)| \leq \mu
\end{equation}
for some $\eps \lesssim \mu \lesssim 1$ and $j \neq 1$, and 
\begin{equation}\label{eq.est_3_extra_1}
     \max(|x_1 - x^*|, |v_1 - v^*|) \gtrsim \mu',  
\end{equation}for some $\eps \lesssim \mu \lesssim \mu' \lesssim 1$. Then $\|\Ic(Q)\|_{L^\infty}\leq (\frac{\mu}{\mu'})^{1/C_\Phi}  \cdot |\log \eps|^{C^*}\cdot \|Q\|_{L^\infty}$.
\end{enumerate}
\end{proposition}

\begin{proof} We may assume $\|Q\|_{L^\infty}=1$. The proof of \eqref{it.est_full_fixed} directly follows from \eqref{eq.associated_int_op}. As for the proof of \eqref{it.est_3} and \eqref{it.est_special_3}, we can derive \eqref{eq.intmini_3} and \eqref{eq.intmini_3_special} by \eqref{eq.associated_dist_C} and making a change of variable $x_2  \leftrightarrow \nu_{12} $, where $\nu_{12} = \eps^{-1}(x_1 - x_2 + t_{12}(v_1 -v_2))$.
    
    We now prove that, if $[t^-_\nf,t^+_\nf] \subset ((\ell[\nf]-1)\tau,\ell[\nf] \tau)$, we have
    \begin{equation}
        \frac{|(v_1-v_2)\wedge \nu_{12}|}{|v_1-v_2|} =\frac{|(v_3-v_4)\wedge \nu_{34}|}{|v_3-v_4|}.
    \end{equation}
    Indeed, by Proposition \ref{prop.scattering_map}, we have $v_3 - v_4 = v_1 -v_2 - 2(v_1 -v_2)\cdot \omega_{12} \omega_{12}$, and $\nu_{34} = - \nu_{12} +2 (\omega_{12} \cdot \nu_{12})\omega_{12}$. Therefore, we have 
    \begin{equation}
      \begin{aligned}
         |(v_3-v_4)\wedge \nu_{34}|^2  = |v_3 - v_4|^2 - |(v_3-v_4) \cdot \nu_{34}|^2 = |v_1 - v_2|^2 - |(v_1-v_2) \cdot \nu_{12}|^2 = |(v_1-v_2)\wedge \nu_{12}|^2,
      \end{aligned}
    \end{equation}
    Here $|a \wedge b|^2 = |a|^2 \cdot |b|^2 - (a\cdot b)^2$, and we used $|v_3 - v_4| = |v_1 -v_2|$, and $(v_3 -v_4)\cdot\nu_{34} = - (v_1 -v_2)\cdot \nu_{12}$. Hence the result.

Moreover, from Proposition \ref{prop.scattering_map}, we know that 
\begin{equation}\label{eq.dis_13}
     |x_3 - x_1 + t_{12}(v_3 -v_1)| = |- \eps \nu_{12} +\eps (\omega_{12} \cdot \nu_{12})\omega + \frac{\Delta t_{12}}{2}(v_1 + v_2 - 2 v_3)|   \lesssim \eps, 
\end{equation}
where we used the fact that $|v_1 + v_2 - 2 v_3| = |v_4 -v_3| = |v_1 -v_2|$, and $|v_1 -v_2|\cdot \Delta t \lesssim \eps$ by Proposition \ref{prop.upperbond.interaction}. Similarly,
\begin{equation}\label{eq.dis_14}
     |x_4 - x_2 + t_{12}(v_4 -v_2)| \lesssim \eps.
\end{equation} Then by a change of variable $x_{2} \leftrightarrow \nu_{12}$, we can argue the same way as in Propositions 9.1--9.2 of \cite{DHM24}. More precisely, the integral is taken in $(t_{12},\nu_{12},v_2)$ as in \eqref{eq.intmini_3}, with weight $[(v_1-v_2)\cdot\nu_{12}]_-\leq|\log\eps|^{C^*}$. In case \eqref{it.est31}, the gain $\eps^\upsilon$ comes from the integration in $t_{12}$. In case \eqref{it.est33}, the set of $\nu_{12}$ with $|(v_1-v_2)\wedge\nu_{12}|\leq\eps^\upsilon|v_1-v_2|$ has measure $\lesssim\eps^{(d-1)\upsilon}$, which gives the gain. 

In case \eqref{it.est32}, denote by $\lambda:=\eps^{C_\Phi\upsilon}$. If the restriction involves only $x_1,v_1$ and the incoming variables $x_2,v_2$, then it restricts $(t_{12},v_2)$ to a set of measure $\lesssim\lambda|\log\eps|^{C^*}$, which gives the gain $\lambda\leq\eps^\upsilon$, as in \cite{DHM24}.
If it involves an outgoing variable $x_3,v_3,x_4$ or $v_4$, then by Proposition \ref{prop.scattering_map}, \eqref{eq.dis_13} and \eqref{eq.dis_14}, for fixed $t_{12}$ it restricts $(v_2,\omega_{12})$ to a set $E$ with $|E|\lesssim\lambda|\log\eps|^{C^*}$ (for example, $|v_1 -v_3| \leq \lambda$ means $|(v_1 -v_2)\cdot \omega_{12}| \leq \lambda$ and for fixed $\omega$, $v_2$ belong to a set with measure $\leq \lambda \cdot |\log\eps|^{C^*}$). Then by the relation $[(v_1-v_2)\cdot\nu_{12}]_- \d\nu_{12}=b(|v_1-v_2|,\omega_{12})\,\d\omega_{12}$ (Proposition \ref{prop.coll.kernel}) and using Proposition \ref{prop.est.b} \eqref{it.est.b1}, we obtain the gain $|\log\eps|^{C^*}\lambda^{1/C_\Phi}=|\log\eps|^{C^*}\eps^{\upsilon}$. 

{\bf Proof of \eqref{it.est_kl}.} By time reversal, we may assume that the fixed ends are the bottom ends of $\pb_i$, $i\in A$, with $|A|=2l-k$. The top ends are determined by the exit factors of $\Dirac_\nf$, so integrating them only changes the arguments of $Q$:
\begin{equation*}
\begin{aligned}
\Ic_\Mb(Q)\le\eps^{-(k-l)(d-1)}\int\prod_{(i,j)\in\Df(\nf)}&\Big(\dirac\big(|X_i(t_{ij})-X_j(t_{ij})|-\eps\big)\big[(V_i(t_{ij})-V_j(t_{ij}))\cdot\nu_{ij}\big]_-\\
&\times\chi(t_{ij})\,\mathbbm 1_{|t_{ij}-t^-_\nf|\le\eps^{\upsilon_0}}\Big)\,|Q|\prod_{j\notin A}dx_j\,dv_j\,d\vt_\nf .
\end{aligned}
\end{equation*}
We treat the factors of $\Df(\nf)$ in chronological order. The first one, say $(\pb_1,\pb_2)_C$, has $t_{12}=t^-_\nf$. If $e_1$ and $e_2$ are both fixed, integrating in $t_{12}$ gives at most $1$. If $e_1$ is free, the change of variables $x_1=x_2+t_{12}(v_2-v_1)+\eps\nu_{12}$ gives $\eps^{d-1}$. After it, $X_1(t)=x_2+\mathbf a(t,\cdot)$ and $X_2(t)=x_2+\mathbf b(t,\cdot)$, so these particles depend on $x_2$ only through a translation. In general, for each factor $(i,j) \in \Df(\nf)$:
\begin{itemize}
\item if $e_j$ is free, or a particle $\pb_m$ already in the collision with $\pb_i$ has a free bottom end (correspond to $(i,m)$), and the position has not been used yet, we change this position to $\nu_{ij}=\eps^{-1}\big(X_i(t_{ij})-X_j(t_{ij})\big)$. Here $X_j(t_{ij})=x_j+t_{ij}v_j$. In the first case,
\begin{equation}\label{eqn.x2x3change}
x_j=X_i(t_{ij})-t_{ij}v_j-\eps\nu_{ij}.
\end{equation}
In the second case, the particles already in the collision depend on $x_m$ only through a translation, $X_i(t)=x_m+\mathbf c(t,\cdot)$, and
\[
x_m=x_j+t_{ij}v_j+\eps\nu_{ij}-\mathbf c(t_{ij},\cdot).
\]
For fixed $t_{ij}$, this change of variables is a translation, and $\dirac(|w|-\eps)\,dw=\eps^{d-1}d\nu_{ij}$. Since $[(V_i-V_j)\cdot\nu_{ij}]_-\le|\log\eps|^{C^*}$ on $\supp Q$, we obtain a factor $\eps^{d-1}|\log\eps|^{C^*}$, and $t_{ij}$ becomes a free variable. If $(i,j)$ is not the first factor, the restriction $|t_{ij}-t^-_\nf|\le\eps^{\upsilon_0}$ gives an extra factor $\eps^{\upsilon_0}$;

\item otherwise we integrate in $t_{ij}$; thanks to $\chi$, the result is at most $1$.
\end{itemize}
Every particle with a free bottom end eventually gets into the collision together with a particle with a fixed end, so each free position is used exactly once. Hence, if there is a fixed end, the $k-l$ factors $\eps^{d-1}$ cancel the prefactor, and $\Ic_\Mb(Q)\le|\log\eps|^{C^*}$, since all remaining variables are bounded on $\supp Q$. If $k=2l$, the last free position is integrated over a ball of radius $|\log\eps|^{C^*}$, which leaves $\eps^{-(d-1)}$, and the $l-2$ non-first factors give $\eps^{(l-2)\upsilon_0}$.

\emph{Case $k\ge l+2$.} At most one free position is used at the first factor, so at least one is used later. This gives $\eps^{\upsilon_0}$.

\emph{Case $k=l+1$.} Let $\pb_p$ be the only particle with a free bottom end. If the first factor is between two particles with fixed ends, then $t^-_\nf$ is determined by the fixed variables and $x_p$ is used later, which gives $\eps^{\upsilon_0}$. Otherwise the first factor is $(\pb_i,\pb_p)$ with $i\in A$; then $x_p$ is used there and $t^-_\nf$ is free. Since $l\ge3$, there is a second factor, at time $t'\in[t^-_\nf,t^-_\nf+\eps^{\upsilon_0}]$, which is the first interaction of some $\pb_j$ with $j\in A\setminus\{i\}$.
\begin{itemize}
\item If the partner of $\pb_j$ is a particle $\pb_{j'}$ with a fixed end that is not yet in the collision, then $t'$ is the time at which the lines $x_j+tv_j$ and $x_{j'}+tv_{j'}$ reach distance $\eps$, so $t'$ is determined by the fixed variables. Hence $t^-_\nf\in[t'-\eps^{\upsilon_0},t']$, which gives $\eps^{\upsilon_0}$.
\item Otherwise $\pb_j$ gets into the collision of $\pb_i$ and $\pb_p$ at time $t'$. Since all velocities are at most $|\log\eps|^{C^*}$ on $\supp Q$, for $t\in[t^-_\nf,t^-_\nf+\eps^{\upsilon_0}]$ we have $ |X_i(t)-x_i-tv_i|\lesssim |\log\eps|^{C^*}\eps^{\upsilon_0}$. Similarly, since $|X_p(t^-_\nf)-X_i(t^-_\nf)|=\eps$,  we have $|X_p(t)-x_i-tv_i| \lesssim |\log\eps|^{C^*}\eps^{\upsilon_0}.$
At time $t'$, the particle $\pb_j$ is still on its line $x_j+tv_j$, therefore 
\[
\big|(x_i-x_j)+t'(v_i-v_j)\big|\lesssim  |\log\eps|^{C^*}\eps^{\upsilon_0}.
\]
The vector $(x_i-x_j)+t(v_i-v_j)$ depends only on the fixed variables. Note that  $|v_i-v_j|\ge\eps^{\upsilon_0/2}$. So the set of $t'$ satisfying this inequality is contained in an interval of length
\[
\lesssim \frac{|\log\eps|^{C^*}\eps^{\upsilon_0} }{|v_i-v_j|}\lesssim|\log\eps|^{C^*}\eps^{\upsilon_0/2}.
\]
Since $t^-_\nf\in[t'-\eps^{\upsilon_0},t']$, the time $t^-_\nf$ lies in an interval of length $\lesssim|\log\eps|^{C^*}\eps^{\upsilon_0/2}$, and integrating in $t^-_\nf$ gives the factor $|\log\eps|^{C^*}\eps^{\upsilon_0/2}$.
\end{itemize}
In both cases, all later factors are integrated in time and are bounded by $1$. This proves \eqref{eq.int_elem_est_general}.

For a CO-atom, the O-connections in $\Of(\nf)$ have factors of the same form, with $|t_{ij}-t^-_{ij}|\le\eps^{\upsilon_0}$, where $t^-_{ij}$ is determined by earlier variables; the same argument applies, with $l$ the total number of particles. 

\textbf{Proof of \eqref{it.est_kl2}} We may assume the fixed ends are the bottom ends of $\pb_1,\dots,\pb_{l-1}$, and the bottom end of $\pb_l$ is free. As in the proof of \eqref{it.est_kl}, all collision factors except the one where $\pb_l$ enters are at most $1$. For the entry of $\pb_l$, say by colliding with $\pb_m$, the change of variables $x_l\leftrightarrow(t_{ml},\nu_{ml})$ gives $\eps^{d-1}$, which cancels the normalization. Hence
\[
\Ic_\Mb(Q)\le\|Q\|_{L^\infty}\int\big[(V_m(t_{ml})-v_l)\cdot\nu_{ml}\big]_-\,\mathbbm 1_{\nf}\ \d\nu_{ml}\,\d v_l\,\d t_{ml},
\]
where $\mathbbm 1_\nf$ is the indicator function that the remaining collisions of $\nf$ take place.

Consider the first collision of $\nf$ that does not involve $\pb_l$ or its first partner. If this collision is between two particles that have not interacted before, its time is determined by the fixed ends. Since all first-entry times of $\nf$ lie in an interval of length $\eps^{\upsilon_0}$, $t_{ml}$ lies in a fixed interval of length $2\eps^{\upsilon_0}$, and we gain $\eps^{\upsilon_0}$.
 
Otherwise, $\pb_2$ must enter before the interaction of $\pb_1,\pb_l$ ends, since otherwise this interaction would be a binary collision by itself. Hence $|t-t'|\lesssim\eps|v_1-v_l|^{-1}$ by Proposition \ref{prop.upperbond.interaction}. Since $|V_1-v_1|,|V_l-v_l|\le|v_1-v_l|$ by conservation of momentum and energy and $\Phi\ge0$, both $\pb_1$ and $\pb_l$ stay within $C\eps$ of their free lines during the interaction. Therefore
\[
|x_1-x_l+t(v_1-v_l)|\le\eps,\qquad |x_2-x_l+t'(v_2-v_l)|\lesssim\eps .
\]
Together with $|t-t'||v_1-v_l|\lesssim\eps$, we get $|x_1-x_2+t'(v_1-v_2)|\lesssim\eps$. Hence $t'$ lies in an interval of length $\lesssim\eps|v_1-v_2|^{-1}$, and $t$ in a set of measure $\lesssim\eps|v_1-v_2|^{-1}+\eps|v_1-v_l|^{-1}$.  
\[
\|\Ic_\Mb(Q)\|_{L^\infty}\lesssim\big(\eps|v_1-v_2|^{-1}+\eps\big)|\log\eps|^{C^*}\le\eps^{\upsilon}|\log\eps|^{C^*},
\]
by the assumption $|v_1-v_2|\ge\eps^{1-\upsilon}$.

{\bf Proof of \eqref{it.est_3_extra}.} We know $j \in \{2,3,4\}$ and $Q$ is supported in the set where \eqref{eq.est_3_extra} holds. From \eqref{eq.dis_13} and \eqref{eq.dis_14}, we know that $|(x_j + t_{12} v_j) - (x_1 + t_{12} v_1)| \leq O(\eps)$. From \eqref{eq.est_3_extra}, we get (note that $\varepsilon \lesssim \mu$)
\begin{equation}\label{proof.est_3_extra1}
    |(x_j - x^*) \times (v_j - v^*)| \leq \mu |\log\varepsilon|^{C^*} |v_j - v^*| \implies |\big(x_1 - x^* + t_{12}(v_1 - v_j)\big) \times (v_j - v^*)| \leq \mu |\log\varepsilon|^{C^*} |v_j - v^*|. 
\end{equation}
Let $v_j - v_1 := u_j$, then \eqref{proof.est_3_extra1} implies that
\begin{equation}\label{proof.est_3_extra2}
    |w \times( u_j + (v_1 -v^*))| \leq \mu |\log\varepsilon|^{C^*} |v_j - v^*| \leq \mu |\log\varepsilon|^{C^*},  
\end{equation}
where $w := x_1 - x^* + t_{12}(v_1 - v^*)$. By a dyadic decomposition we may assume $|w| \sim \mu''$ for some dyadic $\mu'' \in [\mu, \mu']$ (or $|w| \lesssim \mu$ if $\mu'' = \mu$, or $|w| \gtrsim \mu'$ if $\mu'' = \mu'$). Note that for fixed $w$, the inequality \eqref{proof.est_3_extra2} restricts $u_j$ to a tube in $\mathbb{R}^d$ with size $\leq |\log\varepsilon|^{C^*}$ in one direction and size $\mu(\mu'')^{-1} |\log\varepsilon|^{C^*}$ in all other $(d-1)$ directions; denote this tube by $\mathcal{X}$. Denote also $v_2 - v_1 := u$, then $u_j \in \{u,\, (u\cdot \omega_{12})\omega_{12},\, u - (u\cdot \omega_{12})\omega_{12}\}$ (corresponding to $j \in \{2,3,4\}$ respectively). Note that here $\omega_{12}=\omega_{12}(\nu_{12}, |v_1- v_2|)$.
 
We claim that: for fixed $t_{12}$ (hence fixed $w$), the condition $u_j \in \mathcal{X}$ restricts $(v_2,\omega_{12})$ to a set $E$ with measure $\leq (\frac{\mu}{\mu''})|\log\varepsilon|^{C^*}$, so the corresponding integral is $\lesssim (\frac{\mu}{\mu''})^{\frac{1}{C_\Phi}}|\log\varepsilon|^{C^*}$ by Proposition \ref{prop.est.b} \eqref{it.est.b1}. In fact, if $j=2$ this is trivial using the volume of $\mathcal{X}$. Now consider $j=3$. We apply the change of variable $u \to a \cdot \omega_{12} +\eta$ with $a = (u \cdot \omega_{12}) $ and $\eta \in \Pi^\perp_{\omega_{12}}$. Then $\d u = \d a \d \eta$. 

Define $\Omega_a:=\{\omega_{12}\in\Sb^{d-1}: a \omega_{12} \in \Xc \}$ and we denote the middle line of $\Xc$ by $L$ and its distance to the origin is $x_0$. Since $|\eta| \leq|\log\eps|^{C^*}$, we have $|E|\lesssim|\log\eps|^{C^*}\int|\Omega_a|\d a$. If $a \leq x_0 - \mu (\mu'')^{-1} $, then $\Omega_{a} = \emptyset$ and does not contribute to the integral. If $x_0 - \mu (\mu'')^{-1} \leq a \leq x_0 + \mu (\mu'')^{-1}$, we can gain volume from $\d a$, which is $O(\mu(\mu'')^{-1})$. Now we consider the case $a \geq x_0 + \mu (\mu'')^{-1}$. By a direct geometric estimate, we see that $\Omega_a$ has two symmetric regions, each of them is given by the intersection of $\Sb^{d-1}$ with a cone with angle
$$\arcsin{\frac{x_0 + \mu (\mu'')^{-1}}{a} } - \arcsin{\frac{x_0 - \mu (\mu'')^{-1}}{a}}.$$
Hence, 
\begin{equation}
    |E| \lesssim |\log \eps|^{C^*} \int |\arcsin{\frac{x_0 + \mu (\mu'')^{-1}}{a} } - \arcsin{\frac{x_0 - \mu (\mu'')^{-1}}{a} }|^{1/C_\Phi} \d a \lesssim |\log\eps|^{C^*}(\mu /\mu'')^{1/C_\Phi}.
\end{equation}

Finally, if $j=4$ and dimension $d=2$, then the same argument applies with $\omega$ replaced by $\omega^\perp\in \Sb^{d-1}$ which is orthogonal to   $\omega$; if $j=4$ and dimension $d\geq 3$, we may first fix $\omega$ and $u_\omega$, then $u^\perp \in \Pi_\omega^\perp \cap \Xc$, and it is easy to see that the measure of the set $\Pi_\omega^\perp \cap \Xc \subset \Xc$ is bounded by $\mu(\mu'')^{-1}|\log\varepsilon|^{C^*}$. This proves the first claim in either case.

We next claim that: for fixed $\mu''$, the condition $|w| \sim \mu''$ restricts $t_{12}$ to a subset of $\mathbb{R}$ with measure $\leq (\mu''(\mu')^{-1})^{1/C_\Phi}|\log\varepsilon|^{C^*}$. In fact, recall $w = x_1 - x^* + t_{12}(v_1 - v^*)$. If $|v_1 - v^*| \geq \mu'|\log\varepsilon|^{-1}$, then the desired bound follows by solving a linear inequality in $t_{12}$; if $|v_1 - v^*| \leq \mu'|\log\varepsilon|^{-1}$, then by the assumption $\max(|x_1 - x^*|,|v_1 - v^*|) \gtrsim \mu'$, we know that $|x_1 - x^*| \gtrsim \mu'$, and hence $\mu'' \sim |w| \gtrsim \mu'$, in which case the desired bound becomes obvious (the set of $t_{12}$ is trivially bounded by $|\log\varepsilon|^{C^*}$). This proves the second claim.

By putting the above two claims together (and summing over $\mu''$ which loses at most a logarithm and can be absorbed by $|\log\varepsilon|^{C^*}$), we then conclude the desired result.
\end{proof}

\begin{proposition}[Better upper bounds for two-atom $\Ic_\Mb$]\label{prop.est_twoatom} Let $\Ic_\Mb$ be defined as in \eqref{eq.associated_int_op}, where $\Mb$ is a regular molecule containing only two atoms $\nf_1, \nf_2$ and both are marked as short interaction. We also denote by $e$ the edge between $\nf_1$ and $\nf_2$, which equals some $e_j$, and $t_j = t^-_{\nf_j}$ for $j=1,2$. Let $Q$ be a function supported in $|v_j|\leq|\log\varepsilon|^{C^*}$ and $|x_j-x_j^0|\leq |\log\varepsilon|^{C^*}$ for each edge, as in Proposition \ref{prop.est_single}. Consider the following cases:
\begin{enumerate}[(a)]
\item\label{it.est_33A}  Suppose $\Mb$ is a \{$3_2 3_2A$\} molecule with atoms $(\nf_1,\nf_2)$, and if we cut $\nf_1$ as free then $\nf_2$ would become a fully fixed atom. Let $e$ be the bond connecting them and $(e_1,e_7)$ be the two fixed ends at $\nf_1$ and $\nf_2$ respectively. Assume moreover that either (i) $\ell[\nf_1]\neq\ell[\nf_2]$ for the two atoms $\nf_1$ and $\nf_2$, or (ii) the support of $Q$ is a subset of
\begin{equation}\label{eq.intmini_8-}
\max(|x_1-x_7|,|v_1-v_7|)\geq \eps^{1- C_\Phi \upsilon} \quad \mathrm{or}\quad |t_1-t_2|\geq \eps^{\upsilon}
\end{equation} 
where $(x_j,v_j) = z_{e_j}$ as above.
\item\label{it.est_special_33A} Suppose $\Mb$ is a special \{$3_2 3_2$A\} molecule or a special $\{3_2 k_l\}$ molecule.
Assume either $|t_1 - t_2| \geq  \eps^{ \upsilon}$, or $\ell[\nf_1] \neq \ell[\nf_2]$.
\item\label{it.est_3k} Suppose $\Mb$ is a $\{3_2 k_l\}$ molecule with atoms $(\nf_1,\nf_2)$, where $\nf_1$ is a binary atom and $\nf_2$ is a multi-particle atom (note there may be double bond between $\nf_1$ and $\nf_2$). Moreover, if we cut $\nf_1$ then $\nf_2$ becomes a fully fixed atom. Let $e_1$ being the fixed end at $\nf_1$, and $e_i' \ (1 \leq i \leq l-1)$ being the fixed ends at $\nf_2$. Assume moreover that either (i) $\ell[\nf_1]\neq\ell[\nf_2]$ for the two atoms $\nf_1$ and $\nf_2$, or (ii) the support of $Q$ is a subset of
\begin{equation} 
\max_{1 \leq i \leq l-1} (|x_1-x'_i|,|v_1-v'_i|) \sim \mu, 
\end{equation} 
 for some  $\mu \gtrsim \eps^{1 - C_\Phi \upsilon}$.

\item\label{it.est_33B} Suppose $\Mb$ is a \{$3_23_2$B\} molecule with two C-atoms. Assume the support of $Q$ is a subset of 
\begin{equation}\label{eq.intmini_9-}
|t_1-t_2|\geq \eps^{\upsilon} \quad\text{and}\quad |v_i-v_j|\geq \eps^{C_\Phi \upsilon}\quad\textrm{for\ any\ edge\ }(e_i,e_j)\textrm{\ at\ the\ same\ atom},
\end{equation}
as well as 
\begin{equation}\label{eq.33b_except_crossing}
 |t_j - \ell'\tau| \geq \eps^{\upsilon_0} \ (j=1,2) \textrm{ for } \ell'= \ell[\nf] \textrm{ or } \ell[\nf]+1.
\end{equation} 
Here \eqref{eq.33b_except_crossing} implies no atom is a crossing atom, as all atoms are marked as short interaction.
\item\label{it.est_44} Suppose $\Mb$ is a \{$4_24_2$\} molecule with atoms $(\nf_1,\nf_2)$ and $(e_1,e_7)$ being two free ends at $\nf_1$ and $\nf_2$ respectively, such that $\Mb$ becomes a \{$3_2 3_2A$\} after turning these two free ends into fixed ends. Moreover assume the support of $Q$ is a subset of 
\begin{equation}\label{eq.intmini_10-}
\max(|x_1-x_7|,|v_1-v_7|)\leq \eps^{1-C_\Phi\upsilon}.
\end{equation} 
\item\label{it.est_4k} Suppose $\Mb$ is a \{$4_2(2l)_l$\} molecule with atoms $(\nf_1,\nf_2)$ and $e_1$ and $e'_i (1 \leq i \leq l-1)$ being free ends at $\nf_1$ and $\nf_2$ respectively, such that $\Mb$ becomes a \{$3_2k_l$\} molecule after turning these free ends into fixed ends. Moreover assume the support of $Q$ is a subset of 
\begin{equation} 
\max (|x_1-x'_i|,|v_1-v'_i|) \sim \mu 
\end{equation}
for some $i$ and some $ \mu  \lesssim \eps^{1-C_\Phi \upsilon}$.
\end{enumerate}

Then, in any of above cases, we have 
\begin{equation}\label{eq.intmini_8}
    \|\Ic_\Mb(Q) \|_{L^\infty}\leq \eps^{\upsilon} |\log \eps |^{C^*} \cdot \| Q\|_{L^\infty}.
\end{equation}
\end{proposition}
\begin{proof}
The proof is the same way as Proposition 3.3 in \cite{DHM25} and Proposition 9.3 in \cite{DHM24}. We only outline the ideas and the differences.

For the proof of \eqref{it.est_33A}, we first consider the case that $[t^-_{\nf_1},t^+_{\nf_1}] \notin ((\ell[\nf]-1)\tau,\ell[\nf] \tau)$, then we simply cut the $\nf_1$ as free and we can gain a factor $\eps^{\upsilon}$ from the integration in $t_1$. If $[t^-_{\nf_1},t^+_{\nf_1}] \in ((\ell[\nf]-1)\tau,\ell[\nf] \tau)$, then the proof is the same as the corresponding part in the Proof of Proposition 3.3 in \cite{DHM25}, by using Proposition \ref{prop.est_single} \eqref{it.est_3_extra}. 

{\bf Proof of \eqref{it.est_special_33A}.} We first consider the special $\{3_2 3_2A\}$ molecule case. We assume $ \nf_1$ is a C-atom, the other case is easier. By definition of the special molecule, we have
\begin{equation*}
    \Ic_\Mb(Q) = \int_{\Rb \times \Sb^{d-1} \times \Rb^d} [(V_1(t_1)-v_2)\cdot \nu_{12}]_- \cdot \mathbbm{1}_{\mathrm{col}}(z_j,z_7) \cdot |Q| \d t_1 \d \nu_{12} \d v_2.
\end{equation*}
where $\mathbbm{1}_{\mathrm{col}}(z_j,z_7)\leq1$ is the indicator function that the line $x_j+tv_j$ reaches distance $\eps$ from the given trajectory $X_7$ in the pre-collision configuration at some time $t_2$. If $\ell[\nf_1]\neq\ell[\nf_2]$ and $|t_1-t_2|<\eps^{\upsilon}$, then $t_1$ lies within $\eps^{\upsilon}$ of a layer boundary, and integrating in $t_1$ gives the factor $\eps^{\upsilon}$.

Now assume $|t_1-t_2|\ge\eps^{\upsilon}$. We have
\begin{equation*}
|x_j + t_1 v_j -X_1(t_1)|  \lesssim \eps,\qquad
|x_j + t_2 v_j -X_7(t_2)| \lesssim \eps,
\end{equation*}
where $X_1(t)$ and $X_7(t)$ are given trajectories of $\pb_1$ and $\pb_7$ (which may be curves). Let $\Yc_{t_1}(s):=\frac{X_1(t_1)-X_7(s)}{t_1-s}$. Subtracting the two inequalities gives $|v_j-\Yc_{t_1}(t_2)|\lesssim\eps/|t_1-t_2|$. Hence, for fixed $t_1$, $v_j$ belongs to the set
\begin{equation*}
T_{t_1}:=\big\{v:\ |v-\Yc_{t_1}(s)|\lesssim\eps/|t_1-s| \text{ for some } s \text{ with } |t_1-s|\ge\eps^{\upsilon}\big\}.
\end{equation*}
Take a grid $\{s_m\}$ of spacing $\eps$. The set $T_{t_1}$ is covered by the balls $B\big(\Yc_{t_1}(s_m),C\eps/|t_1-s_m|\big)$ with $|t_1-s_m|\ge\eps^\upsilon/2$. Since the $k$-th such grid point on each side of $t_1$ satisfies $|t_1-s_m|\ge\eps^\upsilon/2+k\eps$, we have
\[
|T_{t_1}|\lesssim\sum_{m}\Big(\frac{\eps}{|t_1-s_m|}\Big)^d\lesssim\sum_{k\ge0}\Big(\frac{\eps}{\eps^{\upsilon}+k\eps}\Big)^d\lesssim\eps^{(d-1)(1-\upsilon)}.
\]
 Let $u=V_1(t_1)-v_2$, and write $u=a\omega_{12}+\eta$ with $a=u\cdot\omega_{12}$ and $\eta\perp\omega_{12}$. Then $v_3=V_1(t_1)-a\omega_{12}$, and $\d v_2=\d a\,\d\eta$ for fixed $\omega_{12}$. Fix $\rho>0$ and $v_2$. If $|a|<\rho$, then $\omega_{12}$ lies in a  domain of measure $\lesssim\min(1,\rho/|u|)$, and by Proposition \ref{prop.est.b} this part contributes $\lesssim\rho^{1/C_\Phi}$. If $|a|\ge\rho$, then we have $\d v_3=|a|^{d-1}\d a\,\d\omega_{12}$. For fixed $v_2$, let $A(v_2):=\{\omega_{12}:\ |a|\ge\rho,\ v_3\in T_{t_1}\}$. Since $\d v_2=\d a\,\d\eta$ for fixed $\omega_{12}$, the condition $v_3\in T_{t_1}$ does not depend on $\eta$, and $\eta$ is bounded, Fubini and the polar coordinates $\d v_3=|a|^{d-1}\d a\,\d\omega_{12}$ give
\[
\int|A(v_2)|\,\d v_2\lesssim\int_{|a|\ge\rho}\mathbbm 1_{v_3\in T_{t_1}}\d a\,\d\omega_{12}=\int_{T_{t_1}}\frac{\d v_3}{|a|^{d-1}}\le\rho^{1-d}|T_{t_1}|.
\]
By Proposition \ref{prop.est.b} applied for each fixed $v_2$, and then H\"older in $v_2$, this part contributes
\[
\int\Big(\int_{A(v_2)}b(|u|,\omega_{12})\,\d\omega_{12}\Big)\d v_2\lesssim\int|A(v_2)|^{\frac1{C_\Phi}}\d v_2\lesssim\Big(\int|A(v_2)|\,\d v_2\Big)^{\frac1{C_\Phi}}\lesssim\big(\rho^{1-d}|T_{t_1}|\big)^{\frac1{C_\Phi}}.
\]
Combining the cases $|a|\geq \rho$ and $|a|\leq \rho$, we have 
\begin{equation*}
\int[(V_1(t_1)-v_2)\cdot\nu_{12}]_-\mathbbm 1_{v_3\in T_{t_1}}\d\nu_{12}\d v_2\lesssim\rho^{\frac1{C_\Phi}}+\big(\rho^{1-d}|T_{t_1}|\big)^{\frac1{C_\Phi}}\lesssim|T_{t_1}|^{\frac1{dC_\Phi}}\lesssim\eps^{\frac{(d-1)(1-\upsilon)}{dC_\Phi}}\le\eps^{\upsilon}.
\end{equation*}

If $v_j =v_4$, the same argument applies with $v_3$ replaced by $v_4=V_1(t_1)-\eta$.
For a special $\{3_2k_l\}$ molecule, the same argument applies with $X_7$ replaced by fixed ends of $\nf_2$.
 
{\bf Proof of \eqref{it.est_3k}.} Case (i) is treated as in \eqref{it.est_33A}, and now we consider case (ii). Let $e=e_j$ be the bond between $\nf_1$ and $\nf_2$.

If the first collision of $\nf_2$ is between two fixed ends, then $t_2$ depends only on the fixed data. If $|t_1-t_2|<\eps^{\upsilon}$, integrating in $t_1$ gives $\eps^{\upsilon}$. If $|t_1-t_2|\ge\eps^{\upsilon}$, we can gain $\eps^{\upsilon}$ as in the proof of \eqref{it.est_special_33A}.  

If $\pb_j$ is in the first collision of $\nf_2$, say with $\pb'_i$, then $\pb'_i$ has not interacted before $t_2$, so $|x_j-x'_i+t_2(v_j-v'_i)|\leq \eps$. By Proposition \ref{prop.est_single} \eqref{it.est_3_extra}, we get $(\eps/\mu)^{1/C_\Phi}\lesssim \eps^{\upsilon}$.

{\bf Proof of \eqref{it.est_33B}.}  By assumption each collision happens inside the interval $(\ell[\nf_i] - 1)\tau,\ell[\nf_i]\tau)$ for $i=1,2$.

\textbf{Case 1.} This is the case when both fixed ends are serial with the common edge $e$. Here by Proposition \ref{prop.scattering_map},
\begin{align*}
\Ic_\Mb(Q)&=\varepsilon^{-(d-1)}\int_{\Rb^{10d}\times\Dc_\Mb}\dirac\big(x_3+ (t_1+ \Delta t_1)v_3 -x_2 - t_1 v_2 - \eps(\omega_{12} \cdot \nu_{12}) \omega_{12} - \frac{\Delta t_1}{2 }(v_1 +v_2))\big)\\
&\times\dirac\big(x_4+ (t_1+ \Delta t_1)v_4 - x_1 - t_1 v_1 + \eps(\omega_{12} \cdot \nu_{12}) \omega_{12} - \frac{\Delta t_1}{2 }(v_1 +v_2)\big)\\
& \times \dirac\big(v_3-v_1+[(v_1-v_2)\cdot\omega_{12}]\omega_{12}\big)\cdot \dirac\big(v_4-v_2-[(v_1-v_2)\cdot\omega_{12}]\omega_{12}\big)\\
&\times \dirac\big(|x_1-x_2+t_1(v_1-v_2)|-\varepsilon\big)\cdot\big[(v_1-v_2)\cdot\nu_{12}\big]_-\\
&\times \dirac\big(x_6+ (t_2+ \Delta t_2) v_6 - x_3 - t_2 v_3 - \eps(\omega_{53} \cdot \nu_{53}) \omega_{53} - \frac{\Delta t_2}{2 }(v_3 +v_5) \big)\\   
&\times\dirac\big(x_7 + (t_2+ \Delta t_2)v_7 - x_5 - t_2 v_5 + \eps(\omega_{53} \cdot \nu_{53}) \omega_{53} - \frac{\Delta t_2}{2 }(v_3 +v_5)\big)\\
&\times \dirac\big(v_6-v_5+[(v_5-v_3)\cdot\omega_{53}]\omega_{53}\big)\dirac\big(v_7-v_3-[(v_5-v_3)\cdot\omega_{53}]\omega_{53}\big)\\
&\times\dirac\big(|x_5-x_3+t_2(v_5-v_3)|-\varepsilon\big)\cdot\big[(v_5-v_3)\cdot\nu_{53}\big]_- \cdot Q(x_1,v_1,\cdots,x_7,v_7,t_1,t_2)\prod_{j=2}^6\mathrm{d}x_j\mathrm{d}v_j\prod_{j=1}^2\mathrm{d}t_j,
\end{align*}
 where $\nu_{12}:=\varepsilon^{-1}(x_1-x_2+t_1(v_1-v_2))$ and $\nu_{53}:=\varepsilon^{-1}(x_5-x_3+t_2(v_5-v_3))$ are two unit vectors. Here we first integrate in $(x_4,v_4,x_6,v_6)$, then make the change of variables $x_5\leftrightarrow\nu_{53}$ to reduce to
\begin{equation}
\label{eq.intmini_proof_11}
\begin{aligned}
\Ic_\Mb(Q)&=\int_{\Rb^{4d}\times\Dc_\Mb}\dirac\big(x_3+ (t_1+ \Delta t_1)v_3 -x_2 - t_1 v_2 - \eps(\omega_{12} \cdot \nu_{12}) \omega_{12} - \frac{\Delta t_1}{2 }(v_1 +v_2)\big)\\
&\times\dirac\big(v_3-v_1+[(v_1-v_2)\cdot\omega_{12}]\omega\big) \cdot \dirac\big(|x_1-x_2+t_1(v_1-v_2)|-\varepsilon\big)\\
&\times \dirac\big(x_7 + (t_2+ \Delta t_2)v_7 - x_5 - t_2 v_5 + \eps(\omega_{53} \cdot \nu_{53}) \omega_{53} - \frac{\Delta t_2}{2 }(v_3 +v_5)\big)\cdot\big[(v_1-v_2)\cdot\nu_{12}\big]_- \prod_{j=2}^3\mathrm{d}x_j\mathrm{d}v_j\prod_{j=1}^2\mathrm{d}t_j\\
&\times\int_{\Sb^{d-1}\times\Rb^d}\dirac\big(v_7-v_3-[(v_5-v_3)\cdot\omega_{53}]\omega_{53}\big) \cdot \big[(v_5-v_3)\cdot\nu_{53}\big]_-\cdot Q(\cdots)\,\mathrm{d}\nu_{53}\mathrm{d}v_5,
\end{aligned}
\end{equation} 
where the exact expressions of the input variables of $Q$ are omitted for simplicity. We first apply a change of variable $\nu_{53} \leftrightarrow \omega_{53}$, and $[(v_5-v_3)\cdot \nu_{53}]_- \d \nu_{53} = b(|v_5- v_3|,\omega_{53})\d \omega_{53}$.
Note that in \eqref{eq.intmini_proof_11}, the $\dirac$ functions imply that $\omega_{53}=(v_7-v_3)/|v_7-v_3|$ and $(v_5-v_3)\cdot\omega_{53}=(v_7-v_3)\cdot\omega_{53}$, we may decompose $v_5-v_3$ into $[(v_5-v_3)\cdot\omega_{53}]\omega_{53}$, and a vector $w\in\Pi_{\omega_{53}}^\perp$ which is orthogonal to $\omega_{53}$, i.e., $v_5-v_3 = a \omega_{53} + w$ with $a =(v_5-v_3)\cdot\omega_{53}$. Then, with fixed $v_3$, we have $\d v_5 = \d a \d w$. Then, by defining $u' = a \cdot \omega$, we have $\d v_5 \d \omega_{53} = a^{-(d-1)} \d u' \d w$. This leads to
\begin{equation}
\label{eq.intmini_proof_12}
\begin{aligned}
\Ic_\Mb(Q)&=\int_{\Rb^{4d}\times\Rb^2}\dirac\big(x_3+ (t_1+ \Delta t_1)v_3 -x_2 - t_1 v_2 - \eps(\omega_{12} \cdot \nu_{12}) \omega_{12} - \frac{\Delta t_1}{2 }(v_1 +v_2)\big)\\
&\times\dirac\big(v_3-v_1+[(v_1-v_2)\cdot\omega_{12}]\omega\big)  \cdot\dirac\big(|x_1-x_2+t_1(v_1-v_2)|-\varepsilon\big)\\
&\times \dirac\big(x_7 + (t_2+ \Delta t_2)v_7 - x_5 - t_2 v_5 + \eps(\omega_{12} \cdot \nu_{12}) \omega_{12} - \frac{\Delta t_2}{2 }(v_3 +v_5)\big)\cdot\big[(v_1-v_2)\cdot\nu_{12}\big]_- \\
&\times |v_7-v_3|^{-(d-1)}\prod_{j=2}^3\mathrm{d}x_j\mathrm{d}v_j\prod_{j=1}^2\mathrm{d}t_j\cdot\int_{\Pi_{\omega_{53}}^\perp} b(|v_5-v_3|,\omega_{53})Q(\cdots)\,\mathrm{d}w\mathrm,
\end{aligned}
\end{equation} where $\omega_{53}=(v_7-v_3)/|v_7-v_3|$ as above, $\Pi_{\omega_{53}}^\perp$ is the orthogonal complement of $\omega_{53}$, and the input variables of $Q$ are again omitted. The factor $|v_7-v_3|^{-(d-1)}$ comes from the Jacobian coming from viewing $[(v_5-v_3)\cdot\eta]\eta$ as the polar coordinates of $v_7-v_3$. Now for fixed $(x_3,v_3)$, we make the change of variables $x_2\leftrightarrow \nu_{12}$, and argue similarly as above, to get
\begin{equation}
\label{eq.intmini_proof_13}
\begin{aligned}
\Ic_\Mb(Q)&=\varepsilon^{d-1}\int_{\Rb^{2d}\times\Rb^2}\dirac\big(x_3+ (t_1+ \Delta t_1)v_3 -x_2 - t_1 v_2 - \eps(\omega_{12} \cdot \nu_{12}) \omega_{12} - \frac{\Delta t_1}{2 }(v_1 +v_2)\big)\\
&\dirac\big(x_7 + (t_2+ \Delta t_2)v_7 - x_5 - t_2 v_5 + \eps(\omega_{53} \cdot \nu_{53}) \omega_{53} - \frac{\Delta t_2}{2 }(v_3 +v_5)\big)\\
&\times|v_7-v_3|^{-(d-1)}|v_1-v_3|^{-(d-1)}\,\mathrm{d}x_3\mathrm{d}v_3\mathrm{d}t_1\mathrm{d}t_2\cdot\int_{\Pi_{\omega_{53}}^\perp\times\Pi_{\omega_{12}}^\perp} b(|v_5-v_3|,\omega_{53})b(|v_1-v_2|,\omega_{12}) Q(\cdots)\,\mathrm{d}w\mathrm{d}\theta,
\end{aligned}
\end{equation} where $\omega_{12}=-(v_3-v_1)/|v_3-v_1|$ and $\Pi_{\omega_{12}}^\perp$ is the orthogonal complement of $\omega_{12}$. Finally, by integrating the $\dirac$ functions in \eqref{eq.intmini_proof_13}, first in $x_3$ and then in $v_3$ and then in $(t_1,t_2)$, using that $x_5-x_3+ t_2(v_5-v_3) = \eps\nu_{53}$, we get
\begin{equation}
\label{eq.intmini_proof_14}
\begin{aligned}
    \Ic_\Mb(Q)\leq&\varepsilon^{d-1}\int_{\Rb^2}|t_1-t_2|^{-d}|v_7-v_3|^{-(d-1)}|v_1-v_3|^{-(d-1)}\,\mathrm{d}t_1\mathrm{d}t_2\\
    & \times \int_{\Pi_{\omega_{53}}^\perp\times\Pi_{\omega_{12}}^\perp} b(|v_5- v_3|,\omega_{53}) \cdot b(|v_1- v_2|,\omega_{12})Q(\cdots)\,\mathrm{d}w\mathrm{d}\theta,
\end{aligned}
\end{equation} 
where each of $(v_3,\eta,\omega)$ is a function of $(t_1,t_2)$ and the fixed parameters. The weight $|t_1-t_2|^{-d}|v_7-v_3|^{-(d-1)}|v_1-v_3|^{-(d-1)}$ in $\Ic_\Mb(Q)$ is bounded above by $\varepsilon^{-1/4}$ due to \eqref{eq.intmini_9-} and $0<\upsilon<C_\Phi^{-2} 3^{-d-1}$, and can be absorbed by $\varepsilon^{d-1}$. For the cross-sections, note that $|(v_1- v_2) \cdot \omega_{12}| = |v_1 -v_3| \geq \eps^{C_\Phi \upsilon}$ and $|(v_5-v_3)\cdot \omega_{53}| \geq \eps^{C_\Phi \upsilon}$. Hence, by  Proposition \ref{prop.est.b} \eqref{it.est.b2}, $b(|v_i-v_j|,\omega)$ are bounded by $\leq \eps^{-1/4}$, and can also be absorbed by our choice of $\upsilon$, hence the result.

\textbf{Case 2.} This is the case when the bottom fixed edge is serial with $e$, but the top fixed edge is not. Here
\begin{align*}
\Ic_\Mb(Q)&=\varepsilon^{-(d-1)}\int_{\Rb^{10d}\times\Dc_\Mb}\dirac\big(x_3+ (t_1+ \Delta t_1)v_3 -x_2 - t_1 v_2 - \eps(\omega_{12} \cdot \nu_{12}) \omega_{12} - \frac{\Delta t_1}{2 }(v_1 +v_2)\big)\\
&\times\dirac\big(x_4+ (t_1+ \Delta t_1)v_4 - x_1 - t_1 v_1 + \eps(\omega_{12} \cdot \nu_{12}) \omega_{12} - \frac{\Delta t_1}{2 }(v_1 +v_2)\big)\\
& \times \dirac\big(v_3-v_1+[(v_1-v_2)\cdot\omega_{12}]\omega_{12}\big)\cdot \dirac\big(v_4-v_2-[(v_1-v_2)\cdot\omega_{12}]\omega_{12}\big)\\
&\times \dirac\big(|x_1-x_2+t_1(v_1-v_2)|-\varepsilon\big)\cdot\big[(v_1-v_2)\cdot\nu_{12}\big]_-\\
&\times \dirac\big(x_7+ (t_2+ \Delta t_2) v_7 - x_3 - t_2 v_3 - \eps(\omega_{53} \cdot \nu_{53}) \omega_{53} - \frac{\Delta t_2}{2 }(v_3 +v_5) \big)\\   
&\times\dirac\big(x_6 + (t_2+ \Delta t_2)v_6 - x_5 - t_2 v_5 + \eps(\omega_{12} \cdot \nu_{12}) \omega_{12} - \frac{\Delta t_2}{2 }(v_3 +v_5)\big)\\
&\times \dirac\big(v_7-v_5+[(v_5-v_3)\cdot\omega_{53}]\omega_{53}\big)\dirac\big(v_6-v_3-[(v_5-v_3)\cdot\omega_{53}]\omega_{53}\big)\\
&\times\dirac\big(|x_5-x_3+t_2(v_5-v_3)|-\varepsilon\big)\cdot\big[(v_5-v_3)\cdot\nu_{53}\big]_- \cdot Q(x_1,v_1,\cdots,x_7,v_7,t_1,t_2)\prod_{j=2}^6\mathrm{d}x_j\mathrm{d}v_j\prod_{j=1}^2\mathrm{d}t_j,
\end{align*}
where $\nu_{12} $ and $\nu_{53} $ are as above. Here we first integrate in $(x_4,v_4)$, make the change of variables $x_2\leftrightarrow \nu_{12}$ and argue as in Case 1 above (i.e. decomposing $v_1-v_2$ into $[(v_1-v_2)\cdot\omega_{12}]\omega_{12}$ and a vector $\theta \in\Pi_{\omega_{12}}^\perp$ orthogonal to $\omega_{12}$), to reduce
\begin{equation}
\label{eq.intmini_proof_16}
\begin{aligned}
\Ic_\Mb(Q)&=\int_{\Rb^{6d}\times\Rb^2}\dirac\big(x_7+ (t_2+ \Delta t_2) v_7 - x_3 - t_2 v_3 - \eps(\omega_{53} \cdot \nu_{53}) \omega_{53} - \frac{\Delta t_2}{2 }(v_3 +v_5) \big) \\ 
&\times\dirac\big(x_6 + (t_2+ \Delta t_2)v_6 - x_5 - t_2 v_5 + \eps(\omega_{53} \cdot \nu_{53}) \omega_{12} - \frac{\Delta t_2}{2 }(v_3 +v_5)\big)\\
&\times \dirac\big(v_7-v_5+[(v_5-v_3)\cdot\omega_{53}]\omega_{53}\big)\dirac\big(v_6-v_3-[(v_5-v_3)\cdot\omega_{53}]\omega_{53}\big)\\
&\times\dirac\big(|x_7-x_6+(t_2+ \Delta t)(v_6-v_7)|-\varepsilon\big)\cdot\big[(v_5-v_3)\cdot\nu_{53}\big]_-  \\
 &\times \dirac\big(x_3+ (t_1+ \Delta t_1)v_3 -x_2 - t_1 v_2 - \eps(\omega_{12} \cdot \nu_{12}) \omega_{12} - \frac{\Delta t_1}{2 }(v_1 +v_2)\big)\\
 &\times |v_1-v_3|^{-(d-1)}\prod_{j\in\{3,5,6\}}\mathrm{d}x_j\mathrm{d}v_j\prod_{j=1}^2\mathrm{d}t_j\cdot\int_{\Pi_{\omega_{12}}^\perp} b(|v_1-v_2|,\omega_{12})Q(\cdots)\,\mathrm{d}\theta,
\end{aligned}
\end{equation} where $\omega_{12}= - (v_3-v_1)/|v_3-v_1|$ and $\Pi_{\omega_{12}}^\perp$ is the orthogonal complement of $\omega_{12}$. Now in \eqref{eq.intmini_proof_16}, we make the change of variables $x_6\leftrightarrow \nu_{67} = \eps^{-1}(x_6 - x_7 + (t_2 +\Delta t_{2})(v_6-v_7) )= \nu_{53} - 2(\omega_{53} \nu_{53}) \omega_{53}$. Noting that 
$$(v_5-v_3)\cdot \nu_{53}=(v_6-v_7)\cdot \nu_{67} \quad \text{as well as}\quad\! v_3 - v_7 = v_6 - v_7 -\big((v_6-v_7)\cdot \omega_{35} \big) \omega_{35},$$ 
we have
\begin{equation}
\label{eq.intmini_proof_17}
\begin{aligned}
    \Ic_\Mb(Q)&=\varepsilon^{d-1}\int_{\Rb^{d}\times\Sb^{d-1}\times\Rb^2}\dirac\big(x_7+ (t_2+ \Delta t_2) v_7 - x_3 - t_2 v_3 - \eps(\omega_{53} \cdot \nu_{53}) \omega_{53} - \frac{\Delta t_2}{2 }(v_3 +v_5)\big)\\
     &\times \big[(v_6-v_7)\cdot\nu_{67}\big]_- 
    \,\mathrm{d}v_6\mathrm{d}\nu_{67}\mathrm{d}t_1\mathrm{d}t_2 \cdot |v_1-v_3|^{-(d-1)}\cdot\int_{\Pi_{\omega_{12}}^\perp} b(|v_1-v_2|,\omega_{12}) Q(\cdots)\,\mathrm{d}\theta,
\end{aligned}
\end{equation}
where $(x_3,v_3)$ is determined by $(v_6,\eta,t_1,t_2)$ and the fixed parameters. Note that with \emph{fixed} $(t_1,t_2,\eta)$, we have 
\[x_7+ (t_2+ \Delta t_2) v_7 - x_3 - t_2 v_3 - \eps(\omega_{53} \cdot \nu_{53}) \omega_{53} - \frac{\Delta t_2}{2 }(v_3 +v_5)\\
= (t_1 -t_2 -\frac{\Delta t_2}{2})\sigma +\gamma \]
 where $\gamma$ depends only on $(t_1,t_2,\eta)$ and the fixed parameters, and $\sigma$ is the projection of $v_6-v_7$ to the orthogonal complement $\Pi_{\omega_{35}}^\perp$ of $\omega_{35}$. Then the $\dirac$ function in \eqref{eq.intmini_proof_17} can be written as $\dirac_\Rb(\gamma \cdot \omega_{53})\dirac_{\Pi_{ \omega_{53}}^{\perp}}(\gamma+(t_1-t_2 - \frac{\Delta t_2}{2})  \sigma)$ where $\Pi_{\omega_{53}}^{\perp}$ is the orthogonal complement of $ \omega_{53}$. The latter $\dirac_{\Pi_{\omega_{53}}^{\perp}}$ yields $|t_1-t_2 - \frac{\Delta t_2}{2}|^{-(d-1)}$ upon integrating in $\sigma$, and leaves us with the remaining component $\alpha:=(v_6-v_7) \cdot {\omega_{53}}$ from the integral in $v_6$. Note also that 
\begin{align*}
\gamma\cdot {\omega_{53}}&={\omega_{53}}\cdot \big[x_7+ (t_2+ \Delta t_2) v_7 - x_3 - t_2 v_3 - \eps(\omega_{53} \cdot \nu_{53}) \omega_{53} - \frac{\Delta t_2}{2 }(v_3 +v_5)\big]\\
&={\omega_{53}}\cdot \big[x_7 -x_1+ t_1 (v_7-v_1) +\frac{\Delta t_1}{2} (v_3 -v_4) + \eps \nu_{12}\big]  - \frac{\Delta t_2}{2}\alpha + (\nu_{53}\cdot \omega_{53})  \\
\end{align*} 
Denote by $\zeta:=x_7 -x_1+ t_1 (v_7-v_1) +\frac{\Delta t_1}{2} (v_3 -v_4) + \eps \nu_{12}$. We also have $|\zeta|\gtrsim \varepsilon^{(1+C_\Phi)\upsilon}$ since $|\frac{\Delta t_1}{2}(v_3 -v_4)| \lesssim \eps$ by Proposition \ref{prop.upperbond.interaction}, and 
\begin{multline}
x_7 -x_1+ t_1 (v_7-v_1)=(x_7+t_2v_7)-(x_1+t_1v_1)+(t_1-t_2)v_7\\
=(x_3+t_2v_3)-(x_3+t_1v_3) +(t_1-t_2)v_7  + \eps \mathbf{a} =(t_1-t_2)(v_7-v_3) +\eps  \mathbf{a} 
\end{multline}
due to \eqref{eq.intmini_9-}, here $\mathbf{a}$ denotes some $d$ dimension vector with length $O(1)$. As a result, we now have that
\begin{equation}
\label{eq.intmini_proof_18}
\begin{aligned}
    \Ic_\Mb(Q)&=\varepsilon^{d-1}\int_{\Rb \times\Sb^{d-1}\times\Rb^2}[\alpha]_-\cdot|t_1-t_2 -\frac{\Delta t_2}{2}|^{-(d-1)} |v_1-v_3|^{-(d-1)} \\
    &\times\dirac({\omega_{53}}\cdot \zeta- \frac{\Delta t_{2}}{2}\alpha + (\nu_{53}\cdot \omega_{53})\,\mathrm{d}\alpha\mathrm{d}{\omega_{53}}\mathrm{d}t_1\mathrm{d}t_2\cdot\int_{\Pi_\omega^\perp} b(|v_1-v_2|,\omega_{12} Q(\cdots)\,\mathrm{d}\theta,
\end{aligned}
\end{equation}
 where $v_3$ and $\omega_{12}$ are determined by $(\alpha,{\omega_{53}},t_1,t_2)$ and the fixed parameters. Note that the value $\nu_{53} \cdot \omega_{53}$ does not depend on $\omega_{53}$ by rotation invariant. The integral in ${\omega_{53}}$ can be performed by decomposing ${\omega_{53}}$ into a vector ${\omega_{53}}^\circ$ parallel to $\zeta$ and a vector ${\omega_{53}}^\perp$ orthogonal to $\zeta$. The part of integrand involving ${\omega_{53}}^\circ$ is 
 $$|\zeta|^{-1}\dirac_\Rb({\omega_{53}}^\circ-\frac{ \zeta}{|\zeta|^2}\cdot \big(\frac{\Delta t_2 \alpha}{2} + (\nu_{53}\cdot \omega_{53})\big) ),$$
and the integral in ${\omega_{53}}^\perp$ can be reparametrized into an integral on the unit sphere $\Sb^{d-2}$ in the orthogonal complement of $\zeta$ up to a multiplicative factor of $(1-\varepsilon |\zeta|^{-1})^{d-2}\sim 1$. Thus
\begin{equation}
\label{eq.intmini_proof_18+}
\begin{aligned}
\Ic_\Mb(Q)\leq\varepsilon^{d-1}\int_{\Rb^\times\Sb^{d-2}\times\Rb^2}[\alpha]_-|t_1-t_2|^{-(d-1)}|v_1-v_3|^{-(d-1)}\frac{|\zeta|^{-1}}{(1-\varepsilon |\zeta|^{-1})^{d-2}}\,\mathrm{d}\alpha\mathrm{d}\eta^\perp \mathrm{d}t_1\mathrm{d}t_2\\
\times \int_{\Pi_\omega^\perp} b(|v_1-v_2|,\omega_{12})Q(\cdots)\,\mathrm{d}\theta,
\end{aligned}
\end{equation} 
where ${\omega_{53}}$ is determined by $(t_1, t_2, {\omega_{53}}^\perp)$ and the fixed parameters, and $v_3$ and $\omega$ are determined by $(\alpha,{\omega_{53}},t_1,t_2)$ and the fixed parameters. The weight $|t_1-t_2|^{-(d-1)}|v_1-v_3|^{-(d-1)}|\zeta|^{-1} \cdot b(|v_1-v_2|,\omega_{12})$ in $\Ic_\Mb(Q)$ are bounded by $\varepsilon^{-2/3}$ due to \eqref{eq.intmini_9-} and Proposition \ref{prop.est.b} \eqref{it.est.b2}, and can be absorbed by $\varepsilon^{d-1}$, and hence \eqref{eq.intmini_8} follows.

\textbf{Case 3.} This is the case where neither fixed end is serial with $e$. The proof is similar to the corresponding part in the proof of Proposition 9.3 \cite{DHM24}, using Proposition \ref{prop.est.b} to bound the cross-section as in Case 1.

{\bf Proof of \eqref{it.est_44} and \eqref{it.est_4k}.} These are easy and follows from \eqref{it.est_33A} and \eqref{it.est_3k}. Suppose that $\Mb$ is a $\{4_2  (2l)_l\}$ molecule for $l\geq 2$. We first cut $\nf_1$ as free, and then $\nf_2$ becomes a $(2l-1)_l$ molecule (or a $(2l-2)_l$ molecule if there exists a double bond between $\nf_1$ and $\nf_2$, but in this case we must have $l\geq 3$). Then we have for $l=2$, $|v_7 - v_1|\leq \eps^{1-C_\Phi\upsilon}$. Hence the integral with respect to $v_7$ gives a factor $\eps^{d(1 - C_\Phi \upsilon)}$, which can absorb the factor $\eps^{-(d-1)}$ in $\Ic_{\nf-1}$ and provides an extra $\eps^\upsilon$ by our choice of $\upsilon$. 

Now if $l > 2$, we can also gain at least $\mu^d \leq \eps^{d(1- C_\Phi \upsilon)}$ in integrating velocity variables in $\nf_2$, which is also enough to obtain the desired result.

\end{proof}

\subsection{Classification of components and Definition of excess.}
Now we define the normal, good, and bad components, based on Propositions \ref{prop.est_single}--\ref{prop.est_twoatom}. The cutting algorithm aims to produce good components while avoiding bad ones. These properties depend not only on the molecule $\Mb$, but also on the restriction on the support of $Q$, which comes from splitting (Definition \ref{def.splitting}).

\begin{definition}[Normal, bad and good molecules]
\label{def.good_normal} Suppose $\Mb'$ is formed from a full molecule $\Mb$ by an operation sequence (Definition \ref{def.cutting_algorithm}) and picking one sub-case at each splitting, and $\Mb'$ contains only elementary molecules (Definition \ref{def.elementary}). The initial link structure $\Lc$ between the bottom ends of $\Mb'$ are inherited from that of $\Mb$ (as the operations do not affect any bottom end). Note also that in the splitting steps, by selecting a sub-case each time, we may have made certain extra restrictions to the support of $Q$ in (\ref{eq.associated_int_op}).

Based on all this information, we now define the notions of good, normal and bad molecules. Below, for any  molecule with one binary atom $\{\nf_1\}$, we denote the 4 edges by $e_1,\dots,e_4$, where $e_1$ is serial with $e_3$, and $e_2$ is serial with $e_4$. For two-atom molecule $\{\nf_1,\nf_2\}$, we denote the 7 edges by $e_1,\cdots,e_7$ with $e_1$ and $e_7$ being the fixed ends in the case of \{$3_23_2$\} molecules.
Recall also the constant $\upsilon $ is defined in Section \ref{sec.notation}.
\begin{enumerate}
    \item The following molecules are \newterm{normal}:
    \begin{enumerate}
        \item Any $\{l_l\}$ molecule.
        \item Any $\{3_2\}$ molecule, except those described in \eqref{it.good_1} and \eqref{it.good_1_1} below.
        \item\label{it.normal_33A} Any $\{3_23_2A\}$ molecule, except those described in \eqref{it.good_2} below. 
    \end{enumerate}
    
    \item The following molecules are \newterm{bad}:
    \begin{enumerate}
        \item Any \{$(2l)_l$\} molecule.
        \item Any empty end.
    \end{enumerate}
    
    \item The following molecules are \newterm{good}:
    \begin{enumerate}
    \item\label{it.good_1} Any non-special \{$3_2$\} molecule $\Xc=\{\nf\}$ with fixed end $e_1$, where either the extra restriction includes one of the following conditions:
    \begin{equation}\label{eq.good_normal_1}
    |t_\nf-t_{\nf'}|\leq \varepsilon^{\upsilon}\quad\mathrm{or}\quad|v_{e_i}-v_{e_j}|\leq\varepsilon^{C_\Phi \upsilon}\quad\mathrm{or}\quad|x_e-x_{e'}|\leq\varepsilon^{1- C_\Phi \upsilon}\quad\mathrm{or}\quad|x_{e_i}-x_{e_j}|\leq\varepsilon^{1- C_\Phi\upsilon}
    \end{equation}
    or $(e_i,e_j)\in\Lc$, or $(e,e')\in\Lc$ (i.e. the relevant pair of ends are initial linked as in Definition \ref{def.molecule}; see also the $\mathbbm{1}_{\Lc}^{\varepsilon}$ factor in Definition \ref{def.associated_op}).

     Here in \eqref{eq.good_normal_1} we assume that (i) $\nf'$ is an atom in a molecule $\Yc$ with $\Yc\prec_{\mathrm{cut}}\Xc$ (Definition \ref{def.ordered_disjoint}), and $e'$ is a free end or bond at $\nf'$; (ii) $e$ is a free end at $\nf$ and $e_i\neq e_j$ are two ends at $\nf$; (iii) if $\nf$ is an O-atom, then $e'$ is not serial with $e_1$, and $e_i$ not serial with $e_j$. In this case, $t_{\nf'}$ and $x_{e'}$ are fixed variables, so the extra restrictions make $\Xc$ good.
     \item\label{it.good_1_1} Any special $\{3_2\}$ molecule with the restriction $|t^-_\nf - t^*| \leq \eps^{\upsilon}$. Any special $\{3_2\}$ that is O-connected with a short C-constituent.
    \item\label{it.good_multiple}  Any (special) $\{k_l\}$ molecule with $l \geq 3$ and $l+2 \leq k\leq 2l$. Any non-special $\{(l+1)_l\}$ atom such that $|v_i -v_j| \geq \eps^{\upsilon_0/2}$ for every pair of fixed ends.
    \item\label{it.good_2} Any \{$3_2 3_2A$\} molecule $\Xc=\{\nf_1,\nf_2\}$ with fixed ends $(e_1,e_7)$, where either the extra restriction includes one of the following conditions:
    \begin{equation}\label{eq.good_normal_2}|x_{e_1}-x_{e_7}|\geq\varepsilon^{1- C_\Phi\upsilon}\quad \mathrm{or}\quad |v_{e_1}-v_{e_7}|\geq\varepsilon^{1-C_\Phi \upsilon}\quad\mathrm{or}\quad|t_{\nf_1}-t_{\nf_2}|\geq\varepsilon^{\upsilon},
    \end{equation}
    or $\ell[\nf_1] \neq \ell[\nf_2]$. In this case $\Xc$ is good as a \{$3_2 3_2A$\} molecule without strong degeneracy.

    \item Any special $\{3_2 3_2 A\}$ molecule $\Xc = \{\nf_1,\nf_2\} $ such that $\ell[\nf_1] \neq \ell[\nf_2]$, or $|t_{\nf_1} - t_{\nf_2}| \geq \eps^\upsilon$.
    \item Any $\{3_2 (l+1)_l\}$,  $\{3_2 (l+2)_l\}$ molecule stated in Proposition \ref{prop.est_twoatom} \eqref{it.est_3k}.
    \item Any $\{4_2 (2l)_l\}$ molecule stated in Proposition \ref{prop.est_twoatom} \eqref{it.est_4k}.
    \item Any \{$3_23_2$B\} molecule $\Xc=\{\nf_1,\nf_2\}$ of two C-atoms with fixed ends $(e_1,e_7)$, where the extra restriction includes both of the following conditions: 
    \begin{equation}\label{eq.good_normal_4}|v_{e_i}-v_{e_j}|\geq \varepsilon^{C_\Phi\upsilon}\quad\textrm{and}\quad |t_{\nf_1}-t_{\nf_2}|\geq \varepsilon^{ \upsilon}.
    \end{equation} 
Here $(e_i,e_j)$ are two different ends at the same atom ($\nf_1$ or $\nf_2$). In this case $\Xc$ is good as a \{$3_23_2$B\} molecule without weak degeneracy.
    \item\label{it.good_44} Any \{$4_24_2$\} molecule $\Xc=\{\nf_1,\nf_2\}$, where the extra restriction includes both of the following conditions: 
    \begin{equation}\label{eq.good_normal_3}
        |x_{e_1}-x_{e_7}|\leq \varepsilon^{1-C_\Phi\upsilon}\quad\mathrm{and}\quad |v_{e_1}-v_{e_7}|\leq \varepsilon^{1-C_\Phi\upsilon}.
    \end{equation} 
    Here $(e_1,e_7)$ are two free ends at $\nf_1$ and $\nf_2$ respectively, such that $\Mb$ becomes a \{$3_2 3_2A$\} molecule after turning these two free ends into fixed ends. In this case the strong extra restrictions (\ref{eq.good_normal_3}) makes $\Xc$ good.
    \item\label{it.good_3} 
Any $\{{4}_2\}$ molecule $\Xc=\{\nf\}$, where either the extra restriction includes one of the following conditions:
    \begin{equation}\label{eq.good_normal_4+}
   |v_{e_i}-v_{e_j}|\leq\varepsilon^{C_\Phi \upsilon}\quad \mathrm{or}\quad|x_e-x_{e'}|\leq\varepsilon^{1-C_\Phi\upsilon}\quad\mathrm{or}\quad |x_{e_i}-x_{e_j}|\leq\varepsilon^{1-C_\Phi \upsilon},
    \end{equation}
    or $(e_i,e_j)\in\Lc$, or $(e,e')\in\Lc$.

     Here in \eqref{eq.good_normal_4+} we assume that (i) $\nf'$ is an atom in a molecule $\Yc$ with $\Yc\prec_{\mathrm{cut}}\Xc$ (Definition \ref{def.ordered_disjoint}), and $e'$ is a free end or bond at $\nf'$; (ii) $e$ is a free end at $\nf$; (iii) $(e_i,e_j)$ are two different ends at $\nf$, and are non-serial when $\nf$ is O-atom. In this case, $x_{e'}$ is a fixed variable, so the extra restrictions make $\Xc$ good. 

    \item  Any MO-atom, formed by one M-constituent and a simple pair, such that there exist no fixed end at the M-constituent.\footnote{The estimate of such atom is proved in Proposition \ref{prop.alg_up} below.}
    \item Finally, by convention, we also view an empty end $e$ as good, if $(e,e')\in\Lc$ for some $e'$, i.e. $e$ is involved in some initial link. In this case the initial link condition (see the $\mathbbm{1}_{\Lc}^{\varepsilon}$ factor in Definition \ref{def.associated_op}) will ensure the gain and thus we view $e$ as good.
    \end{enumerate}
\end{enumerate}
\end{definition}

\begin{definition}\label{def.excess}
Let $\Mb$ be any elementary molecule as in Definition \ref{def.elementary}; note that here we are allowed to add various restrictions to the support of $Q$ as in Propositions \ref{prop.est_single}-\ref{prop.est_twoatom}.
 
For dyadic number $\sigma \lesssim 1$, we say $\Mb$ has excess $\sigma$, if $\Ic_\Mb(Q)$ satisfies
\[ \|\Ic_{\Mb}(Q)\|_{L^\infty} \leq \sigma \cdot |\log \eps|^{C^*}, \]
or  \[\|\Ic_{\Mb}(Q)\|_{L^\infty} \leq \eps^{-(d -1)}\sigma \cdot |\log \eps|^{C^*}, \]
if $\Mb$ is a bad molecule.
\end{definition}
 \begin{remark}
    By Proposition \ref{prop.est_single}, each good $\{k_l\} \ (l\geq 3)$ molecule has excess $\eps^{\upsilon_0/2}$. Moreover, a $\{k_l\} \ (k\geq l+2, l\geq 3)$ molecule must be good with excess $\eps^{\upsilon_0}$.
 \end{remark}
\subsection{Proof of Proposition \ref{prop.cumulant_est}}\label{sec.summary} In this subsection we finish the proof of Proposition \ref{prop.cumulant_est}, assuming Proposition  \ref{prop.comb_est} below, which states that there exists an operation sequence that satisfies the nice property (\ref{eq.overall_alg}).  
\begin{proposition}
\label{prop.comb_est} Let $\Mb=(\Mc,\Ec,\Pc,\Lc,\Df)$ be a full (unlabeled) molecule as in Definition \ref{def.set_T_F}, which belongs to either $ \mathcal{F}_{\boldsymbol{\Lambda}_\ell}$ or $\mathcal{F}_{\boldsymbol{\Lambda}_\ell}^{\mathrm{err}}$. Let $r(\Mb)=H$ and $Q=Q_\Mb$ as in (\ref{eq.Q_M}). Let also $\rho$ be defined as in (\ref{eq.def_rho_old}) in Definition \ref{def.parameter_rho_old}.

Then, there exist an operation sequence (Definition \ref{def.cutting_algorithm}) which produces at most $C^{|\Mb|}\cdot |\log\varepsilon|^{C^*\rho}$ total sub-cases, such that in each sub-case the $\Mb'$ contains only good, normal and bad elementary molecules (Definition \ref{def.good_normal}). Moreover, if $\Mb\in \mathcal{F}_{\boldsymbol{\Lambda}_\ell}$, we have the following estimate:
\begin{equation}\label{eq.overall_alg}(\upsilon/2)\cdot\#_{\mathrm{good}}+(d-1)(|H|-\#_{\mathrm{bad}})\geq (C_{13}^*)^{-1}\cdot (\rho + \Bf + \Xf),
\end{equation} where $\#_{\mathrm{good}}$ and $\#_{\mathrm{bad}}$ are the number of good and bad molecules respectively. If $\Mb\in \mathcal{F}_{\boldsymbol{\Lambda}_\ell}^{\mathrm{err}}$, then (\ref{eq.overall_alg}) holds with an extra term $-C^* |H|$ on the right hand side.
\end{proposition}
\begin{proof} The proof will occupy Sections \ref{sec.cuttingI}--\ref{sec.maincr}.
\end{proof}
 
\begin{proof}[Proof of Proposition \ref{prop.cumulant_est} assuming Proposition \ref{prop.comb_est}] The proof is similar to those in \cite{DHM24}. 

    \textbf{Proof of Part 1: properties of final components and decomposition of the support of $Q$}. Since we allow the combinatorial loss $|\log \eps |^{C^* (\rho+\Bf +\Xf)}$. We also have     \begin{equation}\label{eq.not_3_comp_upp_3}
        \#_{\neg\mathrm{N}\{3\}}\le C(\#_{\mathrm{good}}+\rho + \Bf + \Xf),
    \end{equation} where $\#_{\neg\mathrm{N}\{3\}}$ is the number of final components in $\Mb'$ that are either not normal or not \{$3_2$\} molecules, and $\#_{\mathrm{good}}$ and $\#_{\mathrm{bad}}$ are the number of good and bad final components in $\Mb'$ respectively.

Indeed, define $\sigma:=\#_{\mathrm{bo/fr}}- \sum_{l \geq 2} (2 l -1) |\Mb_l|$ and consider its increment under cutting operations, where $\#_{\mathrm{bo/fr}}$ is the number of bonds plus the number of free ends (not counting fixed ends), and $|\Mb_l|$ is the number of (remaining) $l$-atoms in $\Mb$. It is straightforward to check that cutting a $k_l$ atom increases the value of $\sigma$ by $ (2l -1) - k$. Cutting a $3_2 (l+1)_l$ increases the value by $l - 1$. We also claim that partial cut a $k_l$ atom as free from a ${k_0}_{l_0}$ atom increases the value of $\sigma$ by $ (2l -1) - k$. Indeed, let the partial cutting generates ${k_j}_{l_j} \ (j \geq 1)$ atoms in $\Mb$, we know that
\begin{equation} 
        \sum_{j \geq 1} (k_j - 1) + k = k_0, \qquad  \sum_{j} (l_j - 1) + l = l_0 .
\end{equation}
Hence,
\begin{equation}
    \Delta \sigma = \sum_{j \geq 1} (k_j - (2 l_j -1)) - \big(k_0 -(2 l_0 -1)\big) = (2l -1) - k.
\end{equation}
In view of the initial and final values of $\sigma$, we then get  
\begin{equation}\label{eq.incr1}
\begin{aligned}
       &\sum_{2\leq l\leq k\leq 2l} \#_{\{k_l\}}(2 l - k -1) + \sum_{l \geq 2} (l -1) ( \#_{\{3_2 (l+1)_l \} +\#_{3_2 (l+2)_l}})\\
       =& - \sigma_{\mathrm{init}} = -\sum_{l \geq 2} |\Mb_l| + \#_{\mathrm{bond}} + \#_{\mathrm{fixed}} =-1 +  \#_{\mathrm{cycle}} + \#_{\mathrm{fixed}},
\end{aligned}
\end{equation}
 where $\sigma_{\mathrm{init}}$ is the initial value of $\sigma$ (with final value $\sigma_{\mathrm{fin}}=0$), $\#_{\mathrm{cycle}}$ is the number of cycles (initially) in $\Mb$ and $\#_{\mathrm{fixed}}$ is the number of fixed ends (initially) at $\Mb$.

Since at first the number of cycles is no more than $\rho + C \Bf$, $\#_{\neg\mathrm{N}\{3\}}$ is no more than $C( \#_{\mathrm{good}}+  \#_{\mathrm{bad}} + \rho + \Bf + \Xf)$. We conclude the result by the assumption that $\#_{\mathrm{bad}} \leq C(\#_{\mathrm{good}} + H)$ and $|H| \leq \rho$.

Moreover, we can decompose $$Q_\Mb  \leq \sum_{\alpha}w_\alpha Q_\alpha,$$
 where $\alpha$ is the summation index, $w_\alpha\geq 0$ is a scalar, and $Q_\alpha$ are nonnegative functions that have localized support in the $x_e$ and $v_e$ variables, specified as below:
\begin{enumerate}
    \item\label{it.proof_int_est_from_alg_step1_1} Each $Q_\alpha$ has the following form for some $(x_e^*)$ with $x_e^*\in\Zb^d$: 
    \begin{equation}\label{eq.proof_int_est_from_alg_step1_*2}
    \begin{aligned}
        Q_\alpha= \mathbbm{1}_\Dc\cdot \mathbbm{1}_\Lc^{\varepsilon}\cdot\prod_{e\in \Ec_{\mathrm{end}}^-} \big(e^{-\beta|v_e|^2/3}\cdot \oneb_{|x_e-x_e^*|\leq |\log\varepsilon|^{C^*},\,|v_e|\leq |\log\varepsilon|^{C^*}}\big).
    \end{aligned}
    \end{equation}
\item\label{it.proof_int_est_from_alg_step1_2} The $w_\alpha$ satisfies
    \begin{equation}\label{eq.proof_int_est_from_alg_step1_*1}
        \sum_\alpha w_\alpha\le C^{|\Mb|}.
    \end{equation}
    \item\label{it.proof_int_est_from_alg_step1_3} In the support of $Q_\alpha$ and $\prod_{\nf\in\Mb}\boldsymbol{\Delta}_\nf$ in $\Ic_{\Mb}$ (see \eqref{eq.associated_int_op}), for each edge $e\in \Ec$, we have $|x_{e}-x_{e'}^*|\leq|\log\varepsilon|^{C^*}$ and $|v_e|\leq |\log\varepsilon|^{C^*}$, where $e'\in \Ec_{\mathrm{end}}^-$ is the bottom end of the particle line containing $e$. 
\end{enumerate}

\textbf{Proof part 2.} Auxiliary results: Proposition 9.8 in \cite{DHM24}. This proposition is used in (9.67) to obtain the upper bound for $B_{N\{3_2\}}$. Here $\{3_2\}$ molecules can also be produced through partial cutting.

Let $\Mb'$ be obtained from $\Mb$ by the operation sequence, and let $\Sf$ be the set of the following $\{3_2\}$ molecules $\Mb'_j$ of $\Mb'$:
\begin{enumerate}
\item normal $\{3_2\}$ molecules whose atom is not a crossing atom;
\item special $\{3_2\}$ molecules formed from an M-constituent.
\end{enumerate}
Special $\{3_2\}$ molecules formed from a C-constituent, and crossing $\{3_2\}$ molecules, are excluded and discussed in part 4. Let $\nf_j$ be the unique atom of $\Mb'_j$. We relate these molecules to $\Mb$ \emph{before cutting}, which are binary atoms in $\Mb$, or an edge constituent O-connected to an M-constituent.
\begin{enumerate}
\item Let $W$ be the set of binary atoms and binary C-constituent of $\Mb$ that are not crossing atoms and becomes a normal $\{3_2\}$ molecule in $\Sf$ after cutting.
\item Let $\mathcal O_M$ be the set of O-connections $\of=(\pb,\pb')_O$ of $\Mb$ between a particle $\pb$ of an M-constituent and a particle $\pb'$ of
  an edge-constituent, with overlap time $t_\of$, that become a special $\{3_2\}$ atom in $\Sf$ after cutting.
\item Let $A$ be the set of the remaining non-crossing binary atoms and binary C-constituents, and the remaining O-connections between a particle of an M-constituent and a particle of an edge constituent, i.e. those not in $W \cup \Oc_\Mb$.
\end{enumerate}
Each $\Mb'_j\in\Sf$ comes either from an atom $\nf_j\in W$ or, by partial cutting, from some $\of\in\mathcal O_M$. In the second case, by
\eqref{eq.intmini_3_special}, its weight is $((V_\pb(t_\of)-v_{\pb'})\cdot\nu)_-$, where $V_\pb(t)$ is the velocity of $\pb$ inside the box and $v_{\pb'}$ is the velocity of $\pb'$. Hence
\begin{equation}
B_{N\{3_2\}}\leq\prod_{\Mb'_j\in\Sf}|v_{e_1(\nf_j)}-v_{e_2(\nf_j)}| \leq \prod_{\nf\in W}\big(1+|v_{e_1(\nf)}-v_{e_2(\nf)}|\big) \prod_{\of\in\mathcal O_M}\big(1+|V_\pb(t_\of)-v_{\pb'}|\big),
\end{equation}
where, for $\Mb'_j$ coming from $\of$, we write $v_{e_1(\nf_j)}:=V_\pb(t_\of)$
and $v_{e_2(\nf_j)}:=v_{\pb'}$.

\emph{Time variables.} By Definition \ref{def.associated_op} \eqref{it.associated_vars}, an atom $\nf$ has the time variables $t_\nf=(t_{ij})_{(i,j)\in \Df(\nf)}$, and only $t^-_\nf=\min_{(i,j)}(t_{ij})$ is constrained by \eqref{eq.associated_domain}. For a binary atom, $t^-_\nf=t_{12}$. Below, each atom in $\Sf$ carries the single time variable $t^-_\nf$, and each $\of\in\mathcal O_M$ carries $t_\of$. 

We modify Proposition 9.8 as follows.  
\begin{proposition}
    Recall that initially there exists no potential energy as the initial data is supported in $\Dc^0_N$ in Definition \ref{def.grand_canon}, and that the energy $E_\Pf$ (Remark \ref{rem.conserved_energy}) is nonincreasing along admissible trajectories. Assume $|v_e|\leq |\log \eps|^{C^*}$ for every edge $e$, and fix $\gamma$ with $|\log \eps|^{-1} \leq \gamma \leq 1$. Then there exists $\mathscr{B} \ge1$, depending only on $\Mb$, such that
    \begin{equation}\label{eq.weight_est_1}
    \prod_{\nf\in W}(1+|v_{e_1(\nf)}-v_{e_2(\nf)}|) \prod_{\of\in\mathcal O_M}(1+|V_\pb(t_\of)-v_{\pb'}|) \leq\mathscr{B} \cdot\exp\Big(\gamma\sum_{e\in\Ec_{\rm end}^-}|v_e|^2\Big),
    \end{equation}
    and
    \begin{equation}\label{eq.weight_est_2}
    \mathscr{B} \cdot\Big(\int_{\widetilde\Dc}1\d \vt\Big)\le \big(C\gamma^{-1/2}\tau\big)^{|W|+|\mathcal O_M|}  |\log\eps|^{C^*(\rho+\Bf+\Xf+|A|)}.
    \end{equation}
    Here the integral is over $(t_\nf)_{\nf\in W}$ and $(t_\of)_{\of\in\mathcal O_M}$, and $\widetilde\Dc$ is the projection of the set $\Dc_\Mb$ (defined in \eqref{eq.associated_domain}) onto these time variables.
\end{proposition} 
The proof is the same as \cite{DHM24}, with the following changes.
\begin{enumerate}
\item\label{it.energy} \emph{Energy.} 
\begin{enumerate}
\item The conservation of energy $|v_{e_1}|^2 + |v_{e_2}|^2 = |v_{e_1'}|^2 + |v_{e_2'}|^2$ holds at every atom that is not a crossing atom;
\item Since there is no initial potential energy, at any time, the kinetic energy of any set of distinct particles is at most $\sum_{\Ec^-_{\rm end}}|v_e|^2$. The total energy is non-increasing.
\end{enumerate}
The first fact gives (9.87) and (9.100) of \cite{DHM24} as there. 

\item\label{item.timedomain} \emph{Time domain.} As in Proof part 1 of \cite{DHM24}, consider first a binary tree $T$ with no crossing atoms. Each atom $\nf$ of $T$ has the single time variable $t^-_\nf$ ($=t_\of$ if $\nf=\nf_\of$). We have $t^-_\nf<t^-_{\nf'}+\eps^{\upsilon_0}$ if $\nf$ is a child of $\nf'$:
\begin{enumerate}
\item for atoms other than the $\nf_\of$, because all these interactions are short;
\item along a chain, because $t_{\of_{k-1}}<t_{\of_k}$ by the order in the decoration;
\item along the line of $\pb'$, because the overlap happens between leaving the lower atom and entering the upper atom of that line.
\end{enumerate}
Put $s_\nf=t^-_\nf-\eps^{\upsilon_0}\,\mathrm{dep}(\nf)$, where
$\mathrm{dep}(\nf)$ is the number of ancestors of $\nf$ in $T$. Then the
$s_\nf$ are strictly ordered and lie in an interval of length
$\tau+\eps^{\upsilon_0}(|T|+1)$. Hence (9.86) becomes
\begin{equation}
\int_{\Dc_T}1\d \vt_T\le\big(\tau+\eps^{\upsilon_0}(|T|+1)\big)^{|T|}
\Cf_T^{-1}\le e\,\tau^{|T|}\,\Cf_T^{-1},
\end{equation}
since $\eps^{\upsilon_0}|T|^2\le\tau$, where $\Cf_T=\prod_{\nf\in T}\sigma(\nf)$.
In the same way, (9.91) holds with an extra factor $e$. 

\item \emph{Forest case.} Consider a forest without box and multi-particle atoms. As in Proof part 2 of \cite{DHM24}, we divide it into binary trees and upside down binary trees, and apply \eqref{item.timedomain} to each of them (for upside down trees, by time reversal). This gives (9.96) and (9.97), with an extra factor $e$ for each tree, which is at most $C^{|W|+|\Oc_M|+|A|}$.

\item \emph{Single layer case.} In each cluster $\mathcal C$ of $\Mb_{\ell'}$ (with its crossing atoms removed, see \eqref{it.multi-layer}) we remove:
\begin{enumerate}
\item at most $\rho(\mathcal C)\le\Gamma$ C-atoms or C-constituents, so that the remaining cluster has no cycle;
\item all atoms belonging to boxes (not containing those edge-constituents O-connected to them). There are at most $C\Gamma$, since each cluster has
  at most $2\Gamma$ boxes;
\item the crossing atoms;
\item all constituents that contributes to $\Xf$, and edge constituents that is O-connected to any C-constituents.
\end{enumerate}
The total number of removed atoms is at most
$C(\Rf_{\ell'}+\Bf_{\ell'} + \Xf_{\ell'} + s_{\ell'}+s_{\ell'-1})$. 
After that, we turn those any bond connected to deleted atoms and other atoms $\nf$ into a free end at $\nf$. If there exists a O-connection between an edge and a removed box, then we add free end at those special atoms. For those  $\of=(\pb,\pb')_O\in\mathcal O_M$ O-connected to the same particle at an M-constituent, we draw a auxiliary chain among them, keeping the order in decoration $t_{\of_1}<\dots<t_{\of_m}$. 

After that, each $\nf_\of$ is thus a binary atom, exactly like an O-atom in \cite{DHM24}. The only difference is that the velocity $V_\pf(t)$ changes at different atom $\nf_\of$ in the same chain. They are estimated as
\[
|V_\pb(t)|^2\le\sum_{e \text{ is bottom edge of the M-atom}}|v_e|^2
=\sum_{e \text{ is bottom edge of the M-atom}}|v_e|^2 .
\]
So (9.87) holds for $\nf_\of$:
\[
(1+|V_\pb(t_\of)-v_{\pb'}|)^2\le2+4\Big(\sum_{e \text{ is bottom edge of the M-atom} }|v_e|^2+|v_{\pb'}|^2\Big)\leq 2+4\sum_{(\mathfrak l,e)}|v_e|^2.
\]
where the sum is taken over all pairs $(\mathfrak{l},e)$ of lowest atoms $\lf$ in $\Cc$ and bottom edges at $\lf$. Since we order the O-connections inside the same M-atom, we see that $\sigma(\nf)$ is also strictly increases, and (9.90) still valid, and the estimate for the forest case holds.

Now, applying the forest case, we get
\begin{multline}
\prod_{\nf\in\Mb_{\ell'}}(1+|v_{e_1}-v_{e_2}|)\le\Cf_{\ell'}^{5/6} C^{|\Mb_\ell'|}|\log\eps|^{C^*(\Rf_{\ell'}+\Bf_{\ell'} +\Xf_{\ell'})}\\
\times\exp\Big(C\sum_{e\in\Ec_{\rm end}(\Mb_{\ell'})}|v_e|^{4/3}\Big) \exp\Big(C\sum_{e\in\Ec^*_{\ell'}}|v_e|^{4/3}\Big),
\end{multline}
\begin{equation}
\int_{\widetilde\Dc}1\d t_{\Mb_{\ell'}\setminus A}\le (C\tau)^{|W_{\ell'}|+|\Oc_{\Mb_\ell'}|}\,\Cf_{\ell'}^{-1}\,|\log\eps|^{C^*|A\cap\Mb_{\ell'}|}.
\end{equation}
Here $\Mb_{\ell'}$ includes the atoms $\nf_\of$ of layer $\ell'$.
$\Ec^*_{\ell'}$ is the set of edges at removed atoms that are not crossing
atoms; it contains the edges of the boxes, which are also the ends of the
chains in (2). Each cluster contains at most $C\Gamma$ removed atoms of this
kind, and each edge of a box is used by at most $\mathsf N+1$ of these sums.
So, as in (9.100)--(9.101), the last factor is absorbed by the one before it.

\item\label{it.multi-layer} \emph{General multi-layer case.} We follow Proof part 4 of \cite{DHM24}.  As $W$ and $\mathcal O_M$ contain no crossing atoms, we remove the crossing atoms from each layer and apply (4) to each layer. Multiplying over $\ell'$ gives, as in (9.102),
\[
\text{LHS of \eqref{eq.weight_est_1}}\le\Cf^{5/6}\,C^{|W|+|\mathcal O_M|+|A|}\,|\log\eps|^{C^*(\rho+\Bf+\Xf)}\exp\Big(C\sum_{\ell'}\sum_{e\in\overline\Ec_{\rm end}(\Mb_{\ell'})}|v_e|^{4/3}\Big),
\]
where $\Cf=\prod_{\ell'}\Cf_{\ell'}$ and $\overline\Ec_{\rm end}(\Mb_{\ell'})$ is the set of ends of $\Mb_{\ell'}$ after removing the crossing atoms. Each $e\in\overline\Ec_{\rm end}(\Mb_{\ell'})$ is an end of $\Mb$ (case 1), a bond connecting two layers (case 2), or a bond to a crossing atom (case 2 if the crossing atom is connected to another layer, and case 1 otherwise).  Case 1 is the same as (9.104) in \cite{DHM24}, with the energy identity replaced by the inequality
$$\sum_{\rm case 1} |v_e|^2 \leq 2 \sum_{\Ec^-_{\rm end}}|v_e|^2.$$
as the energy is conserved inside each layer and the initial potential energy is $0$. For case 2, each such edge meets a particle of some $r(\Mb_{\ell'})$, and a crossing atom has at most $d+1$ particles, so by (9.105) there are at most $C\rho$ of them. For fixed $\ell'$, the edges of case 2 at the bottom (resp. top) of $\Mb_{\ell'}$ are free velocities or velocities of disjoint crossing atoms at one boundary time, so by (1)(b) their total is at most $\sum_{\Ec^-_{\rm end}}|v_e|^2$. This replaces energy conservation in (9.107), and (9.106)--(9.108) follow as in \cite{DHM24}, at the cost $(C^*)^{C\rho}$.
\end{enumerate}


\textbf{Proof part 3.} Auxiliary results: Proposition 9.9 in \cite{DHM24}. This proposition is used in (9.73) to bound the volume of the variables $(v_{e_2(\nf_j)})$ of the molecules in $\Sf$.

For $\Mb'_j\in\Sf$, \eqref{eq.intmini_3} and \eqref{eq.intmini_3_special} are integrals in $(t_{\nf_j},\nu_j,v_{e_2(\nf_j)})$. With the convention of part 2, ($v_{e_1(\nf_j)}:=V_\pb(t_\of)$ and $v_{e_2(\nf_j)}:=v_{\pb'}$ if $\Mb'_j$ comes from $\of\in\mathcal O_M$), the weight is $[(v_{e_1(\nf_j)}-v_{e_2(\nf_j)})\cdot\nu_j]_-$ in all cases. We write it as $\rho_j\cdot(1+|v_{e_1(\nf_j)}-v_{e_2(\nf_j)}|)$, with
\[ \rho_j:=\frac{[(v_{e_1(\nf_j)}-v_{e_2(\nf_j)})\cdot\nu_j]_-}{1+|v_{e_1(\nf_j)}-v_{e_2(\nf_j)}|}, \qquad\int_{\mathbb S^{d-1}}\rho_j\d \nu_j\le c_d . \]
The factors $1+|v_{e_1}-v_{e_2}|$ are bounded by \eqref{eq.weight_est_1}.
 
In \cite{DHM24} the variable of a $\{3\}$ molecule is the apsidal line
$\omega$. For fixed $\omega$, the map $(v_{e_1},v_{e_2})\mapsto(v_{e'_1},v_{e'_2})$
is linear, which gives (9.110). Here the variable in \eqref{eq.intmini_3} is
the impact direction $\nu$, and for fixed $\nu$ this map is not linear.
Changing $\nu$ to $\omega$ would replace the weight by the cross section $b$,
which is not bounded under our assumptions. We therefore replace (9.110) by
the following estimate.

\begin{proposition}[Modified Proposition 9.9 in \cite{DHM24}.]\label{prop.modify_9.9} Fix the times of the atoms of the molecules in $\Sf$ ($t_{12}$ for binary atoms, $t_\of$ for special ones), and fix all variables of the molecules not in $\Sf$. Fix $X_e\in[1,|\log\eps|^{C^*}]\cap2^{\mathbb Z}$ for $e\in\Ec_{\rm end}$. Then
    \begin{multline}\label{eq.vol_est}
    \int\prod_{\Mb'_j\in\Sf}\rho_j\d \nu_j\d v_{e_2(\nf_j)} \prod_{e\in\Ec_{\rm end}}\mathbf 1_{|v_e|\le X_e} \prod_{e}\mathbf 1_{|v_e|\le|\log\eps|^{C^*}}\\
    \le\prod_{e\in\Ec_{\rm end}}(3X_e)^d\cdot C^{|\Mb|+|\mathcal O_M|}\cdot |\log\eps|^{C^*(K+\Xf+\Bf)},
    \end{multline}
    where all other velocities are the functions of the integration variables given by the collision rules, and $K$ is the number of molecules of $\Mb'$ not
    in $\Sf$, not counting special molecules formed from a C-constituent. 
\end{proposition}
  \begin{proof} We follow the exploration in the proof of Proposition 9.9 of \cite{DHM24}, with three changes.
\begin{enumerate}
\item Let $\Gc $ consist of the ends of $\Mb$ and the C-edges of the molecules not in $\Sf$. A molecule has at most two atoms, and its C-edges belong to binary atoms, to C-constituents or multi-particle C-atoms (counted in $\Xf$), or to M-atoms (counted in $\Bf$).   Hence $|\mathcal G\setminus\Ec_{\rm end}|\le C(K+\Xf+\Bf)$.
\item Starting from any $e_2(\nf)$, explore to the fixed end $e^\sharp$ created in the same cutting (Definition \ref{def.cutting}; $v_e=v_{e^\sharp}$). Stop if $e^\sharp\in\Gc$, or if $e^\sharp=e_1(\nf_j)$ for some $\Mb'_j\in\Sf$. If $e^\sharp$ is an O-fixed end of a molecule not in $\Sf$, continue from its serial free end, which carries the same velocity (see Definition \ref{def.reg}).
\item The starting edges are $e'_1(\nf_j),e'_2(\nf_j)$ for C-atoms, and $e'_2(\nf_j)$ for O-atoms and special atoms; for an O-atom we also start from $e'_1(\nf_j)$ if $e_1(\nf_j)$ is reached.
\end{enumerate}
As in \cite{DHM24}, every exploration terminates, because each step $e\mapsto e^\sharp$ moves to a strictly later molecule in $\prec_{\rm cut}$. 
The fixed end of an auxiliary simple pair is never created at an extremal atom of an ov-segment, so it is never reached.

We say $\nf_i$ is the parent of $\nf_j$ if an exploration from $\nf_i$ reaches $e_1(\nf_j)$. This defines a rooted forest, whose nodes are molecules of $\Sf$, and different trees are called components. It is clear that the number of components is no more than $|\Sf| \leq |\Mb| + |\Oc|$. 
The integral in \eqref{eq.vol_est} factorizes over components, since their explorations end at different edges of $\Gc$. Fix a component and its root value $u:=v_{e_1}$ at a root ($V_\pb(t_\of)$ for a special root).
\begin{enumerate}
\item For O-atoms and special atoms, integrate $\nu_j$ first, for all other variables fixed; this gives at most $c_d$. These atoms then pass $v_{e_2}$ (and, for O-atoms, $v_{e_1}$) on unchanged.
\item For a C-atom, to apply \eqref{eq.transport}, we still need to integrate respect to the fixed variable $u$ of the root. For $|\delta|\leq1$, shift $v_{e_1(\nf_j)}$ and all $v_{e_2(\nf_j)}$ of the component by $\delta$. Every velocity of the component shifts by $\delta$ (Galilean invariance), every $\rho_j$ is unchanged, and the velocities reached in $\Gc$ shift by $\delta$. Hence, the integral with $v_{e_1(\nf_j)}$ is at most the integral with $u+\delta$ in which the balls $\{|v_e|\le X_e\}$, $\{|v_e|\le|\log\eps|^{C^*}\}$ are enlarged by $1$. Averaging over $\delta$ makes $v_{e_1(\nf_j)}$ an integration variable, at the cost of a factor $|B_1|^{-1}$ (volume of unit ball in $\Rb^d$).
\item Process the C-atoms of the component from the root down. At each of them, both incoming velocities are integration variables, and the rest of
  the integrand depends on them only through the outgoing velocities. By \eqref{eq.transport}, $\int\rho_j\d \nu_j\d v_{e_1}\mathrm dv_{e_2}$ becomes $c_d\int\mathrm dv_{e'_1}\mathrm dv_{e'_2}$.
\end{enumerate}
At the end, only the velocities reached in $\mathcal G$ remain, each integrated once over its enlarged ball. This gives $(3X_e)^d$ for ends of $\Mb$, and $(3|\log\eps|^{C^*})^d$ for each of the at most $C(K+X+B)$ other edges of $\mathcal G$. Moreover, the cost of constant is at most $c_d^{|\Sf|}|B_1|^{-\#\text{components}} \leq C^{|\Mb|+ |\Oc|}$.
\end{proof}
\textbf{Proof Part 4.} \emph{Use in (9.65)--(9.75).} Based on part 2 and part 3, we can argue in the same way as \cite{DHM24} to obtain the desired estimates. Note that in Proposition \ref{prop.comb_est}, we have an extra factor $C^{\Nf}$, we claim that $\Nf \leq (d+1)\Xf + \mathsf{N} \Bf + |\Oc_M| + |\Oc_C|$, where $\Oc_C$ is the set of O-connections between an edge constituent and a C-constituent. Indeed, the particle lines of a multi-particle atom, M-atom or MO-atom $\nf$ are the particle lines of its constituents. Each C-constituent (or multi-particle C-atom) is counted in $\Xf$ and has at most $d+1$ particle lines; each M-constituent has at most $\mathsf{N}$ particle lines and is counted in $\Bf$; and by Proposition \ref{prop.mol_axiom} \eqref{it.axiom1_CO}, each edge constituent is O-connected to exactly one atom constituent, so the number of edge constituents is at most $|\Oc_M|+|\Oc_C|$. To control the factor $C^{|\Oc_C|}$, after the operation sequence of Proposition \ref{prop.comb_est} we further partially cut each remaining CO- or MO-atom into its constituents, by partially cutting its C-constituents (Definition \ref{def.partial_cutting} \eqref{it.partial_cutting1}). This does not change $\#_{\mathrm{bad}}$ (a C-constituent of a CO-atom without fixed end is cut as a bad atom, and the other pieces have the fixed end of an auxiliary simple pair), nor the initial numbers $\#_{\mathrm{cycle}}$ and $\#_{\mathrm{fixed}}$. Now each O-connection in $\Oc_C$ gives a special atom, formed by an edge constituent and an auxiliary simple pair, which is O-connected to a C-constituent marked as short interaction; its time variable lies in a interval of length $\eps^{\upsilon_0}$ around a fixed time (see \eqref{eq.associated_dist_CO_multi}). If this special atom has one fixed end, it is good by Definition \ref{def.good_normal}, and we gain a factor $\eps^{\upsilon_0}$, which absorbs the corresponding factor $C$. Otherwise it has no fixed end or at least two fixed end, and the number of such molecule is no more than $\#_{\mathrm{bad}}+\#_{\mathrm{cycle}}+\#_{\mathrm{fixed}}\leq|H|+C\#_{\mathrm{good}}+\rho+C\Bf$. Hence the factor $C^{|\Oc_C|}$ is acceptable. The factor $C^{|\Oc_M|}$ is absorbed by the factor $(C\gamma^{-1/2}\tau)^{|\Oc_M|}$ in \eqref{eq.weight_est_2}, as $\tau$ is small. 

And \eqref{eq.vol_est} and \eqref{eq.weight_est_2} give (9.75) with $\prod_{e\in \Ec_{\mathrm{end}}} X_e^d$ replaced by $\prod_{e\in \Ec_{\mathrm{end}}}(3X_e)^d$, and with the extra factors $C^{\Nf}|\log\eps|^{C^*(\Xf+\Bf)}$ and $(C\gamma^{-1/2}\tau)^{|\mathcal O_M|}$. The rest is (9.77)--(9.82) of \cite{DHM24}.

Now we consider those special $\{3_2\}$ atoms come from a C-constituent, or a crossing atom.
\item Special $\{3_2\}$ molecules formed from a C-constituent, and crossing $\{3_2\}$ molecules. These are not in $\Sf$, and we estimate them in their
natural position.
\begin{enumerate}
\item Let $\Mb'_j$ be a special $\{3_2\}$ molecule formed by partially cutting a C-constituent. By \eqref{eq.associated_dist_CO_multi}, $t_\of$ must be supported on an interval of length $\eps^{\upsilon_0}$. Hence it must be a good molecule. 
\item A crossing $\{3_2\}$ molecule makred as short interaction must satisfies $\|I_{\Mb'_j}\|_{L^\infty\to L^\infty}\leq \eps^{\upsilon_0}|\log\eps|^{C^*}$ due to the restriction on time variable. 
\end{enumerate}
\end{proof}
  
\begin{lemma}\label{lem.liouville}
Let $w=v_1-v_2\neq0$. For $\nu\in\mathbb S^{d-1}$ with $w\cdot\nu<0$ let $\omega=\omega(w,\nu)$ be the apsidal direction, so that $v_1'=v_1-(w\cdot\omega)\omega$, $v_2'=v_2+(w\cdot\omega)\omega$. Then
$\nu\mapsto\omega$ is a bijection from $\{w\cdot\nu<0\}$ onto $\{w\cdot\omega<0\}$, and
\begin{equation}\label{eq.b_def}
(w\cdot\nu)_-\,\mathrm d\nu=b(|w|,\omega)\,\mathrm d\omega,\qquad
\int_{\{w\cdot\omega<0\}}b(|w|,\omega)\,\mathrm d\omega=c_d|w|,
\end{equation}
where $b$ depends only on $|w|$ and $|\hat w\cdot\omega|$, and $c_d=\int(e\cdot\nu)_-\,\mathrm d\nu$.
Consequently, for measurable $h\ge0$ on $\mathbb R^{2d}$,
\begin{equation}\label{eq.transport}
\int\frac{(w\cdot\nu)_-}{1+|w|}\,h(v_1',v_2')\,\mathrm d\nu\,\mathrm dv_1\mathrm dv_2
\le c_d\int h(v_1',v_2')\,\mathrm dv_1'\mathrm dv_2'.
\end{equation}
\end{lemma}
\begin{proof}
The bijection is Proposition \ref{prop.coll.kernel} (it is monotone since $\partial_{\Jc_0}\Theta>0$); the first
identity in \eqref{eq.b_def} defines $b$, and the second is the first one integrated over $\nu$.
For fixed $\omega$, $(v_1,v_2)\mapsto(v_1',v_2')$ is the hard-sphere map: a linear involution with unit
Jacobian, preserving $|w|$ and mapping $\{w\cdot\omega<0\}$ onto $\{w'\cdot\omega>0\}$. By Tonelli,
the left side of \eqref{eq.transport} equals
\[
\int\mathrm d\omega\int_{\{w\cdot\omega<0\}}\frac{b(|w|,\omega)}{1+|w|}h(v')\,\mathrm dv
=\int\frac{h(v')}{1+|w'|}\int_{\{w'\cdot\omega>0\}}b(|w'|,\omega)\,\mathrm d\omega\,\mathrm dv'
=c_d\int\frac{|w'|}{1+|w'|}h(v')\,\mathrm dv',
\]
using $b(|w'|,-\omega)=b(|w'|,\omega)$.
\end{proof}

\section{Cutting algorithm I: reduction to UD molecules}\label{sec.cuttingI}
In this section and Section \ref{sec.maincr}, we will construct the cutting algorithm in Proposition \ref{prop.comb_est}. The goal of this section is to prove several cases of Proposition \ref{prop.comb_est}.

\subsection{The UP algorithm}
In this subsection, we introduce an algorithm {\bf UP}, which is a modified version of UP algorithm in \cite{DHM24} to deal with the boxes and multi-particle atoms. We first introduce the following algorithm to deal with MO-atoms and CO-atoms.
\begin{definition}[The algorithm {\bf PM} and {\bf PM1}]\label{def.alg_pm}
    Let $\mf$ be an MO-atom, a CO-atom or a BO-atom (by viewing each box as a whole), with atom constituents $\nf_1,...,\nf_q$, and edge constituents $e_1,...,e_r$. These constituents are connected by joining O-connections, and each edge constituent is O-connected with exactly one atom constituent. Moreover, recall that the cluster graph is a tree. We define the algorithm \textbf{PM} as follows:
    \begin{enumerate}
        \item Choose a constituent of $\mf$ as follows: if some atom constituent has a fixed end, choose one of them; otherwise, if some edge constituent has a fixed end, choose one of them; otherwise, arbitrarily choose one atom constituent.
        \item\label{it.alg_pm_2} If this constituent is an atom $\nf_j$, then we partially cut it; otherwise, if this constituent is an edge $e_i$, let $\nf_j$ be the atom constituent to which $e_i$ is O-connected (cf. Proposition \ref{prop.mol_axiom} \eqref{it.axiom1_CO}), and partially cut the pair $\{\nf_j,e_i\}$ as free.
        \item Then proceed to the next atom constituent (if any) that is O-connected to an atom constituent already cut, partially cut it together with the newly added simple pair created when partially cut hte previous one. Repeat until all atom constituents of $\mf$ have been cut.
\end{enumerate} 
We further define the algorithm {\bf PM1}, which is similar to the algorithm {\bf PM}, but with the following modification: If \(\mf\) has exactly one C-constituent and all its remaining constituents are edges, then we simply cut \(\mf\) as free, instead of partially cut the atom constituent as free in \eqref{it.alg_pm_2}.
\end{definition}

With the {\bf PM} algorithm  defined, we introduce the {\bf UP} algorithm as follows.
\begin{definition}[The algorithm \textbf{UP}]\label{def.alg_up} Define the following cutting algorithm \textbf{UP}. It takes as input any regular molecule $\Mb$, which has \textbf{no C-top fixed end} (which means that, in any C-atom, CO-atom, or MO-atom, no C-top edge is a fixed end. \textbf{C-bottom fixed ends} can be understood in the same way). For any such $\Mb$, we define the cutting sequence as follows:
\begin{enumerate}
\item\label{it.alg_up_1} If $\Mb$ contains any fully fixed atom, then cut it as free, and repeat until there is no such atom left.
\item\label{it.alg_up_2} Consider the set of all remaining atoms with at least one fixed end in $\Mb$ (i.e. those atoms that have not been cut). We arbitrarily choose the lowest atom in this set \emph{in the sense that viewing each box as a whole}\footnote{For example, if there are two lowest atoms $\nf_1$ and $\nf_2$, where $\nf_1$ is a ancestor of some atom in box $B$, and $\nf_2$ is a descendant of some atom in $B$, then we choose $\nf_2$ instead of $\nf_1$, although $\nf_2$ is not a descendant of $\nf_1$.}. If $\Mb$ has no such atom, then choose a lowest atom $\nf$ in $\Mb$. Let $S_\nf$ be the set of descendants of $\nf$ (including $\nf$); this $\nf$ and $S_\nf$ will be fixed until all atoms in $S_\nf$ have been cut.
\item \label{it.alg_up_3}
Denote $A := S_\nf$. Within $A$, repeatedly choose a highest remaining atom $\mf$ (the highest atom is chosen in the sense by viewing each box as a whole as above)

\begin{enumerate}
\item If $\mf$ is not a M-atom or a MO-atom, proceed as follows.

If $\mf$ is a $\{3_2\}$ atom and if we cut $\mf$ as free, then either a non-M (i.e. not an M-atom or MO-atom) ov-parent $\mf^{+}$ or  ov-child $\mf^{-}$ of $\mf$ becomes fully fixed, then cut the pair $\{\mf,\mf^{+}\}$ or $\{\mf,\mf^{-}\}$ as free; Otherwise, cut $\mf$ as free. Note that it is possible that $\mf$ has a double bond with $\mf^\pm$.

Repeat this procedure until either a M-atom or a MO-atom is reached or no atom remains in $A$. For the later case, we return to \eqref{it.alg_up_1}.
\item If $\mf$ is an M-atom (or an MO-atom), this means $\mf$ (or at least a constituent of $\mf$) belongs to a box) proceed to \eqref{it.alg_up_4}.
\end{enumerate}
\item \label{it.alg_up_4}
Suppose that $\mf$ is an M-atom or an MO-atom. 
\begin{enumerate}
\item\label{it.alg_up_BM2} If $\mf$ is an MO-atom, and none of the M-constituents of $\mf$ containing any fixed end, we partially cut $\mf$ by {\bf PM}. Then we go to \eqref{it.alg_up_3}.
\end{enumerate}
 Now we assume that $\mf$ has an M-constituent $\pf$ with fixed ends.
\begin{enumerate}\setcounter{enumii}{1} 
\item\label{it.alg_up_BM1} If $\pf$ has C-bottom fixed ends, then $\mf$ must be the atom $\nf$ chosen in \eqref{it.alg_up_2}. Then the constituent $\pf$ must be fully fixed (otherwise $\nf$ will not be the lowest atom with fixed end), and we partially cut it as free and go back to \eqref{it.alg_up_2}. 
 
Now we assume that $\pf$ contains no bottom C-fixed end.
\item\label{it.alg_up_BM3} If every C-top edges of $\pf$ are fixed ends, then we partially cut $\pf$ as free, and go to \eqref{it.alg_up_3}.

\item\label{it.alg_up_BM4} Suppose $\pf$ originally (before joining O-connections) belongs to a box $B$. If at least one C-top edge of $\pf$ is not a fixed end, let $\nf'$ be the lowest atom in $B$ with at least one fixed end. Then we return to \eqref{it.alg_up_3} and replace $A$ by $A \cup S_{\nf'}$.
\end{enumerate}
\end{enumerate}

We also define the dual algorithm \textbf{DOWN}, by reversing the notions of parent/child, lowest/highest, top/bottom etc.
\end{definition} 
\begin{definition}\label{def.alg_up1}
We define an analogous algorithm, \({\bf UP1}\), as a modification of \({\bf UP}\): \begin{enumerate}
    \item\label{it.alg_up1_1} in the whole process, we view each box \emph{as a whole};
    \item\label{it.alg_up1_2} in step \eqref{it.alg_up_3} of Definition \ref{def.alg_up}, we cut only one atom each time instead of a possible two-atom molecule. If $\mf$ is a CO- or BO-atom chosen in Definition \ref{def.alg_up} \eqref{it.alg_up_3}, we apply {\bf PM1} to $\mf$.
    \item\label{it.alg_up1_3} After each atom has been cut, we further cut each box by {\bf UP} or {\bf DOWN} (depending on the position of fixed ends) into elementary molecules.
\end{enumerate}
We also introduce algorithm {\bf UP2}, which is identical to {\bf UP1}, except that we apply {\bf PM} to $\mf$ instead of applying {\bf PM1} in \eqref{it.alg_up1_2}.
\end{definition}
These different algorithms will be used by the context. The next proposition provides a key property of the \textbf{UP} algorithm. 
\begin{proposition}
\label{prop.alg_up} Let $\Mb$ be any molecule as in Definition \ref{def.alg_up}. We may assume $\Mb$ is connected (otherwise consider each component of $\Mb$). Then after applying algorithm \textbf{UP} to $\Mb$ (and same for \textbf{DOWN}):
\begin{enumerate}
\item\label{it.up_proof_1} Fix $\nf$ chosen in Definition \ref{def.alg_up} \eqref{it.alg_up_2}. Let the set of atoms $\nf'$ in Definition \ref{def.alg_up} \eqref{it.alg_up_BM4} be $G$. 
We have the following monotonicity property (MONO): throughout the process of \textbf{UP}, the set $A \backslash\{\nf \cup G\}$ will not have any \emph{C-bottom} fixed end, and $(\Mb\backslash A)\cup\{\nf \cup G\}$ will not have any \emph{C-top} fixed end. In particular, the molecule $\Mb$ is cut into elementary molecules, and $\#_{\{3_23_2B\}}=\#_{\{4_2 (2l)_l\}}=0$, and all $\{3_2 (l+2)_l\}$ molecules have a double bond. Moreover, for each M-atom (or MO-atom) generated in the cutting process, it is either cut as fully fixed, or cut with no fixed end, together with a simple pair which has O-connection with it.
\item\label{it.up_proof_2} We have the number of bad molecules is no more than $1$, and equal to $1$ if and only if $\Mb$ is full.  
\end{enumerate}
In (\ref{it.up_proof_3})--(\ref{it.up_proof_MO}) below we assume $\Mb$ has no fully fixed atoms.
\begin{enumerate}
\setcounter{enumi}{2} 
\item\label{it.up_proof_3} If $\Mb$ contains a cycle, a box or contains at least two fixed ends, then either there exists a good molecule which contains a binary atom, or
$$\sigma_0 := \sum_{l\geq 2}\#_{\{3_2 (l+1)_l \}} + \sum_{l\geq 2}\#_{\{3_2 (l+2)_l \}}+ \sum_{l\geq 3,l<k\leq 2l}\#_{\{k_l\}}  > 0.$$
 If not (i.e. $\Mb$ is a tree, has no box, and contains at most one fixed end), then all elementary molecules are all $\{(2l -1)_l\}$ with at most one exception which might be a $\{(2l)_{l}\}$ molecule.
\item\label{it.up_proof_4} If $\Mb$ contains at most one fixed end \emph{and} contains a cycle, or $\Mb$ has a multi-particle atom, then we can split into at most $2$ sub-cases such that for each sub-case we have $\#_{\mathrm{good}}\geq 1$.
\item\label{it.up_proof_MO} Suppose $\Mb$ has a single MO-atom $\nf$, which is formed by joining an O-connection between an M-atom $\nf_1$ and a simple pair. Moreover, all other ends at $\nf$ are free ends. Then we have 
    \begin{equation}\label{it.est_M}
        \Ic_\Mb (Q) \leq C \eps^\upsilon | \log \eps |^{C^*} \cdot \| Q \|_{L^\infty}.
    \end{equation} 
\end{enumerate}
\end{proposition}
\begin{proof}The proof basically follows the proof of  Proposition 12.2 in \cite{DHM24}.

\textbf{Proof of \eqref{it.up_proof_1} and \eqref{it.up_proof_2}} We first prove the (MONO) property.  Note that (MONO) is true initially, as $\Mb$ has no C-top fixed end. Moreover, there exists no fixed end in $(S_G \cup S_\nf) \backslash \{\nf\}$ by our choice of $\nf$ in Definition \ref{def.alg_up}. Therefore, it suffices to show that (MONO) is preserved. First note that cutting a fully fixed atom in Definition \ref{def.alg_up} \eqref{it.alg_up_1} will only generate bottom fixed ends at its ov-parents, so the absence of the C-top fixed end in the whole molecule is preserved.

    Next, suppose $\nf$ is chosen in Definition \ref{def.alg_up} \eqref{it.alg_up_2}. If $G = \emptyset$, then (MONO) preserved by the same argument as Proposition 12.2 in \cite{DHM24}. Now we assume $G \neq \emptyset$. It is clear that $G \cap S_\nf = \emptyset$ by the choice of $G$. 

    First, we choose the highest remaining atom $\mf \in A$. If we cut $\mf$ as free, this does not create any C-bottom fixed end in $A \backslash  (\{\nf\} \cup G)$, and does not create any C-top fixed end in $(\Mb \backslash A)\cup \{\nf\}  $. Moreover, each time we add an atom $\nf' \in G$, it has only C-bottom fixed end as it does not belong to $A$. Therefore, in the whole process, $(\Mb \backslash A)\cup (\{\nf\}\cup G)$ has no C-top fixed end.
    Moreover,  when we choose $\nf'$, $S_{\nf'} \backslash \{\nf'\}$ (here $S_{\nf'}$ denotes the descendants of $\nf'$ that have not been cut) has no bottom fixed end. Indeed, $S_{\nf'} \backslash \{\nf'\}$ has no bottom fixed end in the box $B$, as we choose $\nf'$ to be the lowest atom with fixed end in $B$. It has no bottom fixed end outside $B$, this is because when we choose $\mf \in A$ in Definition \ref{def.alg_up} \eqref{it.alg_up_3}, we choose the highest atom in $A$ by viewing each box as a whole. But if there exists an bottom fixed end outside $B$, then we must have chosen some $\mf$ which is a descendant of the box, but at the moment not all atoms in the box have been cut, contradiction. 
    

   The same can be proved for operation of cutting $\{\mf, \mf^+\}$ or $\{\mf, \mf^-\}$. Here, cutting $\{ \mf,\mf^+\}$ (the case $\{ \mf,\mf^+\}$ is similar) may create bottom fixed ends at ov-parents of $\mf$ and $\mf^+$, but none of them belong to $A  \backslash (\{\nf\} \cup G)$ because $\mf$ is the highest remaining atom. As for $\{\mf^+\}$, it does not belong to $A$, so it has only C-bottom fixed ends. Since it would become fully fixed after cutting $\mf$, we know that any other ov-child of $\mf^+$ must belong to the maximal ov-segment (there may be one or two such ov-segments, depending on whether there exists a double bond between $\{\mf,\mf^+ \}$) containing $\{\mf,\mf^+ \}$. Viewed in the molecule, cutting $\{\mf,\mf^+\}$ is equivalent to first cut $\mf$ then cut $\mf^+$. After cutting $\mf$, cutting $\mf^+$ does not create any C-top fixed end (hence only create C-bottom fixed ends in $\Mb \backslash A$). This shows (MONO) is preserved in the whole process of Definition \ref{def.alg_up} \eqref{it.alg_up_3}.

Moreover, the statements for $\{3_2\,3_2\,B\}$ and $\{4_2\,(2l)_l\}$ follow directly from the (MONO) property. We now prove the statements for the M-atom and the $\{3_2\,k_l\}$ molecules.

Consider a $\{3_2k_l\}$ molecule with $l\geq 3$, and write it as $\{\nf_1,\nf_2\}$, where $\nf_1$ is the $\{3_2\}$ atom. We show that it satisfies the condition in Definition \ref{def.elementary}: when $\nf_1$ is cut as free, $\nf_2$ remains regular. By Definition \ref{def.alg_up}(\ref{it.alg_up_2}), the $3_2$ atom must be chosen as $\mf$, so the other atom is $\mf^\pm$. If it is $\mf^+$, it has no fixed top end and remains regular when $\mf$ is cut as free. If it is $\mf^-$, it has no fixed bottom end and likewise remains regular. Moreover, in either case, $\nf_2$ becomes fully fixed after $\nf_1$ is cut as free. If there exists double bond between $\nf_1$ and $\nf_2$, then it is a \{$3_2 (k+2)_k$\} molecule, otherwise, it is a \{$3_2 (k+1)_k$\} molecule
 
Now we prove the statement on M-atoms and MO-atoms.  When we (partially) cut an MO-atom in Definition \ref{def.alg_up} (\ref{it.alg_up_BM2}), then it is exactly formed by joining an O-connection between a M-atom and a simple pair, with the only fixed end on the simple pair. If it is cut in Definition \ref{def.alg_up} (\ref{it.alg_up_BM1}), then it must be chosen in Definition \ref{def.alg_up} (\ref{it.alg_up_2}). Assuming it is originally belongs to a box $B$, then at this time each atom in the bottom part of $B$ must have been cut, as $B$ is connected by our construction. Then, this M-atom must be cut as fully fixed. If it is cut in Definition \ref{def.alg_up} (\ref{it.alg_up_BM3}) or (\ref{it.alg_up_BM4}), then by our construction of the algorithm, it must be cut as fully fixed when it is cut. 

\textbf{Proof of (\ref{it.up_proof_3}).} Here and in the proof of (\ref{it.up_proof_4}) below we assume that $\Mb$ has no fully fixed atoms. We claim that if $\sigma_0 = 0$ and $\#_{\mathrm{good}}=0$, then necessarily also $\#_{\{\mathrm{fufi}\}}=0$, where $\#_{\{\mathrm{fufi}\}}$ denotes the number of fully fixed atoms. Assume the contrary, consider the first fully fixed atom $\mf$ when it is cut.

We first assume that $\mf$ is not an M-atom or an atom constituent of any MO-atom, and consider the cutting operation that turns $\mf$ into fully fixed, which must be cutting a single atom $\pf$. By assumption $\pf$ is not fully fixed and it must be a binary atom (since no multi-particle atom is generated, and no M-atom or MO-atom has been cut, as each M-atom is either cut as a good molecule, or cut as a fully fixed atom). Therefore, it should be cut together with either $\mf$ or some other atom, contradiction.

Now assume $\mf$ is a non-M constituent of an MO-atom $\mf'$. By Definition \ref{def.alg_up} \eqref{it.alg_up_4}, if $\mf$ is cut as in Definition \ref{def.alg_up} \eqref{it.alg_up_BM2}, then the M-constituent in $\mf'$ must be cut as a good molecule. If $\mf$ is cut in Definition \ref{def.alg_up} \eqref{it.alg_up_BM1}- \eqref{it.alg_up_BM4}, then it must be cut after some M-constituent, which is either a fully fixed atom or a good molecule. 

Now assume $\mf$ is an M-atom or an M-constituent of an MO-atom. We claim that when cutting $\mf$, there must exist at least one $\{3_2 k_l\}$ atom or multi-particle atom or a fully fixed atom before cutting $\mf$. Then the desired conclusion follows once this claim is established. Denote the box containing $\mf$ by $B$. At the time of cutting $\mf$, all C-atoms and C-constituents (which are C-atoms belong to the box before joining O-connections) in the (either the bottom part or the top part of) $B$ must have been cut (as $\mf$ is cut as fully fixed). We assume that all the top C-atoms have been cut (the other case is similar). Denote all C-atoms, and those C-constituent belongs to the top part of $B$ to be set $\Ac$ (then $\Ac$ is exactly the top part of $B$ before adding any O-connections in $B$). We may assume that all atoms in $\Ac$ are binary atoms (otherwise a multi-particle atom of $\Ac$ has been cut before $\mf$). Consider the quantity
\begin{equation}
    \sigma_\Ac : = \#_{\mathrm{bo / fr}} - 3 |\Ac| 
\end{equation}
where $\#_{\mathrm{bo / fr}}$ denotes the number of maximal ov-segments plus the number of free ends at atoms in $\Ac$. We consider the increment of the quantity $\sigma_\Ac$. Note that by construction, every atoms in $\Ac$ belong to the same layer, and each cluster graph in each layer is a tree.

Cutting an atom outside $\Ac$ does not increase the value of $\sigma_\Ac$ (this step may turn some free ends of $\Ac$ into fixed ends, and may turn a bond into fixed ends in $\Ac$.) Cutting a $\{2_2\},\{3_2\}, \{3_2 3_2\}$ and $\{4_2\}$ molecule in $\Ac$ increases the value of $\sigma_\Ac$ by $1,,0,1$ and $-1$. In view of the first and final values of $\sigma_A$, we have 
\begin{equation}\label{counting.box}
    \#_{\{4_2\}}-(\#_{\{2_2\}}+\#_{\{3_23_2\}})=  1 -\#_{\mathrm{cycle}}(\Ac) - \#_{\mathrm{fixed}} ,
\end{equation}
where $\#_{\mathrm{fixed}}$ denotes the number of fixed ends in $\Ac $ that produced by cutting atoms not in $\Ac$, and $\#_{\mathrm{cycle}}(\Ac)$ denotes the number of cycles in $\Ac$. By construction of the box, we know that  $\#_{\mathrm{cycle}}(\Ac) \geq 1$ and $ \#_{\mathrm{bad}} =1$ if and only if $\#_{\mathrm{fixed}} = 0$. Hence $(\#_{\{2\}}+\#_{\{3_23_2\}}) \geq 1$. Hence before cutting $\mf$ there must exists either a \{$3_2 3_2$\} molecule or a \{$2_2$\} molecule; contradiction to the assumption of $\mf$. Hence, we have proved the claim.
 
Next, let $\Mb$ be a connected molecule. Define $\sigma:=\#_{\mathrm{bo/fr}}- \sum_{l \geq 2} (2 l -1) |\Mb_l|$ and consider its increment under cutting operations, where $\#_{\mathrm{bo/fr}}$ is the number of bonds plus the number of free ends (not counting fixed ends), and $|\Mb_l|$ is the number of (remaining) $l$-atoms in $\Mb$. It is straightforward to check that cutting a $\{k_l\}$ atom increases the value of $\sigma$ by $ (2l -1) - k$, cutting a $\{3_2 (l+1)_l\}$ increases the value by $l - 1$, cutting a $\{3_2 (l+2)_l\}$ (with a double bond) increases the value by $l-1$. We also claim that partially cut a $k_l$ atom as free from a ${(k_0)}_{(l_0)}$ atom increases the value of $\sigma$ by $ (2l -1) - k$. Indeed, considering a partially cutting generates $i$ atoms and each is a $(k_j)_{(l_j)} \ (1\leq j \leq i)$ atom, we know that
\begin{equation} 
        \sum_{j \geq 1} (k_j - 1) + k = k_0, \qquad  \sum_{j} (l_j - 1) + l = l_0 .
\end{equation}
Hence,
\begin{equation}
    \Delta \sigma = \sum_{j \geq 1} (k_j - (2 l_j -1)) - \big(k_0 -(2 l_0 -1)\big) = (2l -1) - k.
\end{equation}
In view of the initial and final values of $\sigma$, we then get  
\begin{equation}\label{eq.incr}
\begin{aligned}
    -\sum_{2\leq l \leq k\leq 2l} \#_{\{k_l\}}(2 l - k -1) -& \sum_{l \geq 2}(l - 1) ( \#_{\{3_2 (l+1)_l \}} + \#_{\{3_2 (l+2)_l\}}) = \sigma_{\mathrm{init}}\\
   & = \sum_{l \geq 2} |\Mb_l| - \#_{\mathrm{bond}} - \#_{\mathrm{fixed}} = 1 - \#_{\mathrm{cycle}} - \#_{\mathrm{fixed}},
\end{aligned}
\end{equation}
 where $\sigma_{\mathrm{init}}$ is the initial value of $\sigma$ (with final value $\sigma_{\mathrm{fin}}=0$), $\#_{\mathrm{cycle}}$ is the number of cycles (including cycles inside boxes) initially, in $\Mb$ and $\#_{\mathrm{fixed}}$ is the number of fixed ends (initially) in $\Mb$. Moreover, we further obtain that
 \begin{equation} 
    \sum_{2 \leq l\leq k\leq 2l-1}\#_{\{k_l\}}(2 l - k -1) + (l -1)\sum_{l\geq 2} ( \#_{\{3_2 (l+1)_l\}} + \#_{\{3_2 (l+2)_l\}}) = \#_{\{(2l)_l\}} + \#_{\mathrm{cycle}} + \#_{\mathrm{fixed}} -1 .
\end{equation}

By \eqref{it.up_proof_2},  if $\#_{\mathrm{fixed}}\geq 1$, the number of bad molecule (i.e. $\#_{(2l)_l}$) is $0$, if $\#_{\mathrm{fixed}}=0$ then $\#_{\{(2l)_l\}} =1$, and if $\#_{\mathrm{fixed}}> 0$ then  $\#_{\{(2l)_l\}} =0$.  We also note that if there is a fully fixed atom, then either there exists a good binary atom, or $\sigma_0 \geq  1$. Therefore, if $\#_{\mathrm{fixed}} \geq 2$, or has a cycle or a box, we always have $\sigma_0 + \sum_{l\geq 2} \#_{\{l_l\}} > 0$. Moreover, since $\sum_{l\geq 2}\#_{\{l_l\}}\geq 1$ implies there exists a good binary atom, or $\sigma_0\geq 1$, we conclude the result.

If $\Mb$ is a tree and has at most one fixed end, then we have
\begin{equation} 
        \sum_{2 \leq l\leq k\leq 2l-1}\#_{\{k_l\}}(2 l - k -1) + (l -1)\sum_{l\geq 2} ( \#_{\{3_2 (l+1)_l\}} + \#_{\{3_2 (l+2)_l\}})  = 0,
\end{equation}
Therefore, every atom must be cut by itself and each contains exactly one fixed end. This proves (\ref{it.up_proof_3}).

\textbf{Proof of (\ref{it.up_proof_4}).} Consider the first $\{k_l\}$ $(2 <l < k\leq 2l)$ atom, or $3_2 (l+1)_l$ atom or $3_2 (l+2)_l$ atom. Consider the first such atom, we assume it is a $\{3_2 3_2\}$ atom (so it is a $\{3_2 3_2A\}$ atom by \eqref{it.up_proof_1}, and other cases are the same). There are fixed ends $e_1$ at the $\{3_2\}$ atom and $e_2$ to be the last fixed end at the $k_l$ atom, which are created when cutting previous $\pf_1$ and $\pf_2$ as free ($\pf_1$ and $\pf_2$ may be the same), and both of them are $\{3_2\}$ atoms, denote by $\Xc_j$ and $\Yc_j$ (Note that no fully fixed atom exists at this time). 

We perform the splitting by considering the dichotomy:
\begin{equation*}
    \text{either}\quad  |v_{e_1'} - v_{e_2'}| \leq \eps^{1 - C_\Phi\upsilon} \quad \text{or} \quad \!|v_{e_1'} - v_{e_2'}|\geq \eps^{1 - C_\Phi \upsilon}.
\end{equation*}
Then, in the former case, this implies either $\Xc_j$ or $\Yc_j$ is good $\{3_2\}$ molecule, in the later case, we have $\{3_2 (l+1)_l\}$ is a good molecule.

 We assume it is a $\{3_2 (l+1)_l\}$ or $\{3_2 (l+2)_l\}$ atom. If it is a $\{3_2 4_2\}$ atom, then it must be a good molecule. If not, there are fixed ends $e_1$ at the $\{3_2\}$ atom and $e_1',...,e_j'$ at the other atom.

Now assume that $\max_{1 \leq i \leq l-1} (|x_1-x'_i|,|v_1-v'_i|) \sim \mu_1$ for some $\eps \lesssim \mu_1 \lesssim 1$. Then, we perform the splitting by considering 
\begin{equation*}
     \text{either}\quad  \mu_1 \cdot \mu_2 \leq \eps^{1- C_\Phi \upsilon} \quad \text{or}\quad \!  \mu_1 \cdot \mu_2 \geq \eps^{1 - C_\Phi \upsilon}.
\end{equation*}
Then in the former case, consider all atoms that cut before the two-atom molecule and is connected to it. There must exist at least one good $\{3_2\}$ molecule among these atoms. In the latter case, we have that this $\{3_2 k_l \}$ ($k =l+1$ or $l+2$) molecule is good.
 
\textbf{Proof of \eqref{it.up_proof_MO}.} Let the particle lines in M-atom $\nf_1$ be $\pb_1,...,\pb_n$ with $n \leq \mathsf{N}$, and the corresponding particle line of the simple pair is $\pb_{n+1}$. Hence $\pb_{n+1}$ overlaps with some $\pb_i$. Here $n \leq \mathsf{N}$. Then by \eqref{eq.associated_dist_M} and \eqref{eq.associated_dist_MO}, we have 
 \begin{equation}
    \begin{aligned}
        \Ic_\Mb (Q) 
        =& \eps^{-n(d - 1)} \int_{\Rb^{2 d (n + 1)}} \mathbbm{1}_{\Mf_{\nf_1}}(\vz_{\nf_1}) \cdot\mathbbm{1}_{\mathrm{overlap}}(z_{i},z_{n+1}) (Q) \d \vz_{p(\Mb)}\\
        =& \eps^{-n(d - 1)} \int_{\Rb^{2 d (n + 1)}} \mathbbm{1}_{\Mf_{\nf}}(\vz_\nf) (Q) \d \vz_{p(\Mb)}  ,
    \end{aligned}
 \end{equation}
 where $\mathbbm{1}_{\mathrm{overlap}}(z_{i},z_{n+1})$ denotes that particles $\pb_i$ and $\pb_{n+1}$ satisfies the pre-collision condition at some time $t_{i(n+1)}$. And $\mathbbm{1}_{\Mf_\nf}$ is defined by $$ \mathbbm{1}_{\Mf_\nf} := \mathbbm{1}_{\Mf_{\nf_1}}(\vz_{\nf_1}) \cdot\mathbbm{1}_{\mathrm{overlap}}(z_{i},z_{n+1}).$$ 
 
The idea is to further expand the $\Mf$, which also has lots of collisions (by our choice of $G(q)$, bottom part, top part of the box in Definitions \ref{def.pathological} and \ref{def.box}), to gain from the cycle inside the M-atom. We denote $\Cc(\Mf_\nf)$ to be the set of all smallest boxed sub-molecules of $\Mf_\nf$, such that every particle lines in $\Mf_\nf$ are connected, it contains at least one cycle. It is clear that, the number of such molecule is bounded by a constant $C=C(\mathsf{N})$. Such sub-molecules exist: by Part 1 of the proof of Proposition \ref{prop.layerrec3}, the bottom and top parts of the box have at most $G(n-1)$ atoms each, so the middle part has at least $G(n)-2G(n-1)-2$ atoms; a connected sub-molecule satisfying (iii) of Definition \ref{def.pathological} with fewer than $n$ particle lines has fewer than $G(n-1)$ atoms (otherwise it would have been made into a box earlier), and a connected sub-molecule with $n$ particle lines and no cycle has at most $n-1$ atoms. Since $G(n)-2G(n-1)-2>2n\,G(n-1)+2n$, the middle part contains a connected sub-molecule with all $n$ particle lines and a cycle. We take such an $\Mb'$ formed by the collisions in a time interval which contains the overlap time of $\pb_i$ and $\pb_{n+1}$. We can rewrite the indicator function $\mathbbm{1}_{\Mf_\nf}$ in terms of
\begin{equation*}
    \mathbbm{1}_{\Mf_\nf} \leq \sum_{\Mb' \in \Cc(\Mb)} \mathbbm{1}_{\Mf^+} \cdot  \mathbbm{1}_{\Mb'} \cdot  \mathbbm{1}_{\Mf^-}  
\end{equation*}
where the event $\Mf^+$ and $\Mf^-$ are the set of all possible molecules in $\Mf_\nf \setminus \Mb'$. Moreover, there exists exactly one fixed end at some O-atom or O-edge of some CO-atom in $\Mb'$. Since $n \leq \mathsf{N}$, the number of molecules in $\Cc(\Mb)$ has an upper bound $C(\mathsf{N})$. Hence, by the $L^1$ preserve property of $\Hc$, we have
\begin{equation*}
    \begin{aligned}
        \Ic_\Mb (Q) \leq
        C(\mathsf{N}) \|Q\|_{L^\infty} \sup_{\Mb' \in \Cc(\Mb)} \eps^{-n (d-1)}  \int_{\Rb^{2d (n+1)}} \mathbbm{1}_{\Mb'} (\vz_{p(\Mb)})  \d \vz_{p(\Mb)}. 
    \end{aligned}
\end{equation*}
Note that $\Mb'$ has exactly one fixed end, and a cycle. Then we apply {\bf UP} to $\Mb'$, and then by \eqref{it.up_proof_4}, we can gain at least one good molecule, or an MO-molecule with only one fixed end on the simple pair. Hence, we have either $\eps^{-n (d-1)}\int_{\Rb^{2d (n+1) }} \mathbbm{1}_{\Mb'}  \d \vz_{p(\Mb)}  \leq \eps^\upsilon | \log \eps|^{C^*}$, or 
$$ \Ic_\Mb(Q) \leq \|Q\|_{L^\infty} \eps^{-n' (d-1)}\int_{\Rb^{2d (n'+1) }} \mathbbm{1}_{\Mb'}  \d \vz_{p(\Mb')},$$ where $\Mb'$ is the set of all possible molecule correspond to the MO atom $\nf'$ with $n' < n$. (formed by a smaller box with $n'$ particles inside $\Mb'$ and the simple pair). In latter case, we can argue above until obtain the desired factor $\eps^\upsilon$ (by induction in $n$; since a box has at least $3$ particles, this stops after at most $\mathsf{N} $ steps).
\end{proof}

\begin{proposition}\label{prop.box_multiple_atom}
    Proposition \ref{prop.comb_est} is true in the case when $\Xf +\Bf \geq 10 \cdot\rho$.
\end{proposition}
\begin{proof}
In this case, we simply apply the {\bf UP2} algorithm to $\Mb$. We first claim that $\#_{\mathrm{bad}} \leq |H|\leq \rho$. Indeed, here we present an argument which will be used many times later in different variants, which we refer to as the \textbf{(PL) argument} below. Start from any component $W$ of $\Mb$; by Proposition \ref{prop.mol_axiom} (\ref{it.axiom4}), this $W$ must intersect a particle line $\pb\in p(\Mb_{\ell'})$ for some $\ell'$. By Proposition \ref{prop.mol_axiom} (\ref{it.axiom5}), this $\pb$ is equal to or connected to some $\pb'\in r(\Mb_{\ell'+1})$ via $\Mb_{\ell'+1}$; in either case $\pb'$ also intersects $W$ (note that each particle line can intersect at most one component of $\Mb$, see Remark \ref{rem.layer_interval}). We then replace $\pb$ by $\pb'$ and repeat the same argument to get $\pb''\in r(\Mb_{\ell'+2})$ and so on, in the end we see that $W$ must intersect a particle line $\boldsymbol{q}\in r(\Mb_{\ell})=H$. Since each $\boldsymbol{q}$ intersects at most one component of $\Mb$, this shows that $\#_{\mathrm{comp}}$ is bounded by the number of particle lines in $H$.

Then, by the same arguments as in Proof of Proposition \ref{prop.alg_up} \eqref{it.up_proof_4} and counting the increasing of $\sigma$, we know that with at most $\rho + |H|$ exceptions, each atom has exactly one fixed end. By algorithm {\bf UP2}, we know that each box is either cut by itself, or together with a simple pair. For each such box with exactly one fixed end, we can obtain at least one good molecule by Proposition \ref{prop.alg_up} after applying {\bf UP} algorithm to it. Moreover, for each atom or constituent that contributes to $\Xf$, if it is a multi-particle C-atom before cutting, it still is after {\bf UP2}. If it is a C-constituent of a CO- or BO-atom, and if it is the first constituent partially cut in this atom, then it may not be cut as a multi-particle atom (if it is a C-constituent with two particles). But such atom must belong to a CO- or BO-atom with at least two atom constituents (otherwise, we will not apply partial cutting to it, see Definition \ref{def.alg_pm}). Moreover, all other constituents contributes in $\Xf$ will be cut as multi-particle atom. Therefore, we can obtain at least $(\Xf + \Bf - \rho - |H|)/2$ good molecules and at most $|H|$ bad molecules, which proves Proposition \ref{prop.comb_est} by our assumption.
\end{proof}
Below, we assume that $\Xf+ \Bf \leq 10 \rho$, it suffices to proof that 
\begin{equation}\label{eq.overall_alg1}
    (\upsilon/2) \cdot \#_{\mathrm{good}} + (d-1) (|H| - \#_{\mathrm{bad}}) \geq 20 (C^*_{13})^{-1} \cdot \rho.
\end{equation}
\begin{proposition}\label{prop.case1} Proposition \ref{prop.comb_est} is true in the case when $|H|\geq (C_{12}^*)^{-1}\cdot\rho$.
\end{proposition}
\begin{proof}
    The proof is basically the same as Proposition 12.3 in \cite{DHM24}. 
\end{proof}
\subsection{Strong degeneracies}\label{sec.strdeg} From now on we only need to consider the case $|H|<(C_{12}^*)^{-1}\cdot\rho$ in Proposition \ref{prop.comb_est}, and the goal is to prove (\ref{eq.overall_alg}). Now it suffices to prove
\begin{equation}\label{eq.case2_2}
    (\upsilon/2) \cdot \#_{\mathrm{good}} - d \cdot \#_{\mathrm{bad}} \geq 20 (C^*_{13})^{-1} \cdot \rho.
\end{equation}
In this subsection we treat the case of sufficiently many strong degeneracies (Definition \ref{def.strdeg} below), see Proposition \ref{prop.case2}. For technical reasons we also need the extra notion of \textbf{primitivity} in Definition  \ref{def.strdeg} below. 
\begin{definition}\label{def.strdeg}
Let $\Mb$ be a molecule as in Proposition \ref{prop.mol_axiom}. Let $\{\nf,\nf'\}$ be a pair of atoms in $\Mb$ satisfying the following conditions follow:
\begin{enumerate}[{(i)}]
\item They are ov-adjacent and in the same layer ($\ell[\nf]=\ell[\nf']$).
\item \label{def.strdeg_3}If both $\{\nf,\nf'\}$ are binary atoms: there exist edges $(e,e')$ at $\nf$ and $\nf'$ respectively, such that cutting $\{\nf,\nf'\}$ as free in $\Mb$ and turning $(e,e')$ into fixed ends results in a \{$3_2 3_2A$\} molecule, as in Definition \ref{def.elementary}; 
\item\label{def.strdeg_5} If one of $\{\nf,\nf'\}$ is not binary atom (say $\nf$) and the other is a binary atom (say $\nf'$), and $\nf$ is not a M-atom. There exist edges $(e_1,e_2,...e_j)$ at $\nf$ and $e'$ at $\nf'$ respectively, such that cutting $\{\nf,\nf'\}$ as free in $\Mb$ and turning $(e_1,e_2,...e_j)$ and $e'$ into fixed ends results in a $\{3_2 (l+1)_l\}$ molecule or a $\{3_2 (l+2)_l\}$ with a double bond, as in Definition \ref{def.elementary};
 
\item\label{def.strdeg_7} For \eqref{def.strdeg_3}, we have the restriction that $|x_{e}-x_{e'}|\leq\varepsilon^{1- C_\Phi \upsilon}$ and $|v_{e}-v_{e'}|\leq\varepsilon^{1- C_\Phi \upsilon}$. And for \eqref{def.strdeg_5}, we have the restriction $\max(|x_{e_i} -x_{e'}| , |v_{e_i}-v_{e'}|) \lesssim \eps^{1 - C_\Phi \upsilon}$ for every $i$. 
 \end{enumerate}
Then we define a pair of atoms $\{\nf,\nf'\}$ in $\Mb$ to be \textbf{strongly degenerate}.

We also define an atom pair $\{\nf,\nf'\}$ in $\Mb$ to be \textbf{primitive}, if they are in the same layer ($\ell[\nf]=\ell[\nf']$), $\nf$ is ov-parent of $\nf'$, and there \emph{do not exist} C-atoms or CO-atoms $\qf_j\,(1\leq j\leq r)$ with $r\geq 1$, such that $\qf_{j+1}$ is ov-child of $\qf_j$ along the C-bond of $\qf_j$ for each $j$ with $\qf_0=\nf$ and $\qf_{r+1}=\nf'$.
\end{definition} 
\begin{proposition}\label{prop.case2} Suppose we make the further restriction that $\Mb$ contains at least $(C_{12}^*)^{-1}\rho$ disjoint strongly degenerate primitive pairs (by inserting the corresponding indicator function $\mathbbm{1}_{\mathrm{str.deg.prm}}$), then Proposition \ref{prop.comb_est} is true.
\end{proposition}
\begin{proof}First assume there exists at least $\Af:=(C_{12}^*)^{-1}\rho$ disjoint strongly degenerate  primitive pairs. By splitting into at most $|\log\varepsilon|^{C^*\Af}$ sub-cases (i.e. decomposing $1$ into the corresponding indicator functions), we may fix one choice of these pairs. By using pigeonhole principle and replacing $\Af$ with $\Lf^{-1}\cdot\Af$, we may restrict all these atom pairs $\{\nf_j,\nf_j'\}$ to a single layer, say $\Mb_{\ell'}$. Define $\Qb$ to be the union of all connected components of $\Mb_{\ell'}$ that contain at least one strongly degenerate primitive pair $\{\nf_j,\nf_j'\}$. We then
\begin{enumerate}[{(i)}]
\item Cut $\Qb$ as free, then cut $\Mb_{>\ell'}$ as free from $\Mb\backslash\Qb$ and cut it into elementary molecules using \textbf{UP};
\item Then cut $\Mb_{<\ell'}\cup (\Mb_{\ell'}\backslash \Qb)$ into elementary molecules using \textbf{DOWN}.
\end{enumerate} 
It is easy to see that this can always be done, and the contribution of these steps (i.e. coming from atoms outside $\Qb$) to $\#_{\mathrm{bad}}$ is at most $2|H|$. This because the number of components of $\Mb_{>\ell'}$ is $\leq|H|$ by repeating the (PL) argument in the proof of Proposition \ref{prop.case1} (i.e. starting with a particle line intersecting one component and iteratively applying Proposition \ref{prop.mol_axiom} (\ref{it.axiom5})), and moreover each \emph{full} component of $\Mb_{<\ell'}\cup (\Mb_{\ell'}\backslash \Qb)$ \emph{after the rest has been cut} must also be a full component in the original molecule $\Mb$ (see Remark \ref{rem.full_cut}) and there are at most $|H|$ of them, so the upper bound on $\#_{\mathrm{bad}}$ follows from Proposition \ref{prop.alg_up} (\ref{it.up_proof_2}).

Now we only need to cut $\Qb$ into elementary molecules. For each connected component $\Qb'$ of $\Qb$, consider the clusters in $\Qb'$ as defined in Proposition \ref{prop.mol_axiom} (\ref{it.axiom3}) (i.e. ov-components of the set of C-atoms in $\Qb'$); denote these clusters by $\Qb_j$, then $\Qb'$ is formed by joining these $\Qb_j$ with O-atoms. Denote by $X$ the number of components of $\Qb$, and by $Y$ the total number of $\Qb_j$ that are cyclic (i.e. $\rho(\Qb_j)>0$) or $\Qb_j$ contain at least one box. Moreover, denote by $Z$ the number of $\Qb_j$ that contains a box (note also that each cluster contains at most $2 \Gamma$ boxes). Depending on the sizes of $X$, $Y$ and $Z$, we may then choose to perform one of the following three cutting sequences.

\textbf{Choice 1}: in this choice we will gain from $X$. For each component $\Qb'$ of $\Qb$, choose one pair $\{\nf_j,\nf_j'\}$ in it, then cut $\{\nf_j,\nf_{j}'\}$ as free, which is a good \{$4_2 (2l)_l$\} molecule due to the restrictions in Definition \ref{def.strdeg} \eqref{def.strdeg_7}. Recall that $S_{\nf_j}$ is the set of descendants of $\nf_j$. We next prove that $\nf_j'$ has no ov-parent $\qf\neq\nf_j$ that is a C-atom in $S_{\nf_j}$ or a CO-atom connecting $\nf_j'$ along  the C-bond. In fact, if such a $\qf$ exists, then $\qf$ is a descendant of $\nf_j$, so there exist a path from $\nf_j$ to $\qf$ by iteratively taking children. This path forms a cycle together with the ov-segments from $\qf$ to $\nf_j'$ and from $\nf_j'$ to $\nf_j$. As the cycle belongs to a single layer $\Mb_{\ell'}$, by Proposition \ref{prop.mol_axiom} (\ref{it.axiom2}) it must be formed by ov-segments connecting C-atoms or CO-atoms along the C-bonds; in particular, by considering the path between $\nf_j$ and $\qf$, we see that the $\qf_j$ must exist as in Definition \ref{def.strdeg}, so $\{\nf_j,\nf_j'\}$ is not primitive, contradicting our assumption that $\{\nf_j,\nf_j'\}$ is a strongly degenerate primitive pair.

Now, we have proved that the above $\qf$ does not exist; this implies that, after cutting $\{\nf_j,\nf_j'\}$ as free, there is no C-top fixed end in $\Qb'\backslash S_{\nf_j}$ (recall that $\nf_j'\in S_{\nf_j}$), and no C-bottom fixed end in $S_{\nf_j}$, in particular the monotonicity property (MONO) in Proposition \ref{prop.alg_up} is true. Since this property is preserved under \textbf{UP}, we can then repeat the steps in Definition \ref{def.alg_up} (\ref{it.alg_up_3}) and subsequently apply \textbf{UP} to the rest of $\Qb'$ to cut it into elementary molecules. By applying this to each component $\Qb'$ of $\Qb$, we obtain $\#_{\mathrm{good}}\geq X$ and $\#_{\mathrm{bad}}=0$.

\textbf{Choice 2}: in this choice we will gain from $Y$. For each component $\Qb'$ of $\Qb$, let $\Qb_j$ be the clusters in $\Qb'$ defined as above, which are joined by O-connections to form the full component $\Qb'$. For convenience, we allow some $\Qb_j$ to be empty ends (where joining an edge with an empty end simply means creating an O-atom in the interior of this edge, and drawing two serial free ends at this O-atom). We may reorder them as follows: first choose $\Qb_1$ arbitrarily, and each time choose one remaining cluster $\Qb_j$ that is joined with some $\Qb_{j'}\,(j'<j)$ by some O-connection (such $\Qb_j$ must exist because $\Qb'$ is connected). 

In this case, as the cluster graph in the same layer is a tree by Proposition \ref{prop.mol_axiom} (\ref{it.axiom2}), we have reordered the $\Qb_j$ such that each $\Qb_j$ is joined with \emph{exactly one} $\Qb_{j'}\,(j'<j)$ by \emph{exactly one} O-connection. Suppose that the atom \(\pf_2\) contains this O-connection. We then cut $\Qb_1$ as free (note that \(\Qb_1\) contains \(\pf_2'\), if it exists, where \(\pf_2'\) denotes the C-constituent of \(\pf_2\) belonging to \(\Qb_1\)), which creates the \emph{unique} atom with fixed end in an atom in $\Qb_2$ (as $\pf_2$ is the only atom joining $\Qb_2$ with $\Qb_1$); then we cut $\Qb_2\cup\{\pf_2\} $ as free, after which we have the unique atom with fixed end $\pf_3$ in an ov-segment connecting two C-atoms in $\Qb_3$ (as $\pf_3$ is the only atom joining $\Qb_3$ with $\Qb_1\cup \Qb_2$), and so on. Note that when we cut each $\Qb_j\cup\{\pf_j\}$, it always contains \emph{at most one} atom with fixed end, so by Proposition \ref{prop.alg_up} (\ref{it.up_proof_4}), we can subsequently cut it into elementary molecules using \textbf{UP} or \textbf{DOWN} depending on whether the unique fixed end is bottom or top, and get a good molecule for each $\Qb_j$ that is \emph{cyclic}. 

By applying the above process to each component $\Qb'$ of $\Qb$, we obtain $\#_{\mathrm{good}}\geq Y$ and $\#_{\mathrm{bad}}\leq X$ (total number of components).  

\textbf{Choice 3}: in this choice we will gain from $\Af$. For each component $\Qb'$ of $\Qb$, we cut it as free using \textbf{UP}, but with the modification that we cut one atom each time (so there will be no elementary molecules with two atoms). For each $\Qb'$, we define the quantity $\sigma=\#_{\mathrm{bo/fr}}-\sum_{l \geq 2}(2l-1)|\Qb'_l|$ as in the proof of Proposition \ref{prop.alg_up} (\ref{it.up_proof_3}). By the same arguments there, we get
\begin{equation} 
   \sum_{2 \leq l \leq k \leq 2l} \#_{\{k_l\}}(2l-k-1)  = \#_{\mathrm{cycle} } - \#_{\mathrm{comp}}
\end{equation} (note that initially $\Qb'$ has no fixed ends and there exists no two-atom molecule). Since $\#_{\{(2l)_l\}} = \#_{\mathrm{comp}} = X$, we have
\begin{equation}
    \sum_{l \leq k \leq 2l -1 }\#_{\{k_l\}}(2l-k-1) =  \#_{\mathrm{cycle}} \leq   C (\mathsf{N}) \Gamma\cdot Y.
\end{equation}
We then have
\[(\textrm{number of molecules that are not } \{(2l-1)_l\}) \leq X+  C (\mathsf{N}) \Gamma\cdot Y.\] In particular, by excluding this many exceptions, we see that the number of pairs $\{\nf_j,\nf_j'\}$ such that \emph{both $\nf_j$ and $\nf_j'$ are cut as $\{(2l-1)_l\}$ molecules} is at least $\Lf^{-1}\cdot\Af-2(X+  C \cdot   \Gamma Y)$.

For each such pair $\{\nf_j,\nf_j'\}$, if one of the atoms is a multi-particle atom, then it is a good molecule. If both are $\{3_2\}$ atoms, suppose $\nf_j$ is cut before $\nf_j'$, we know that $\{\nf_j'\}$ must be good due to the restriction $|v_{e_j}-v_{e'_j}|\leq\varepsilon^{1-C_\Phi \upsilon}$ and Definition \ref{def.good_normal}. Here when $\nf_j'$ is O-atom, the edge $e_j'$ cannot be serial with the bonds between $\nf_j$ and $\nf_j'$, otherwise $\{\nf_j,\nf_j'\}$ will not be regular (and thus not a \{$3_2 3_2A$\} molecule) after turning $e_j'$ into a fixed end. In this way we get $\#_{\mathrm{good}}\geq \Lf^{-1}\cdot\Af-C(\mathsf{N})(X+\Gamma Y +\Gamma Z)$ and $\#_{\mathrm{bad}}\leq X$.

Summing up, by choosing one of the \textbf{Choices 1--3}, we get
\begin{equation}\label{eq.case2_5} (\upsilon/2)\cdot (\#_{\mathrm{good}})-d\cdot\#_{\mathrm{bad}}\geq\max\big((\upsilon/2)X,(\upsilon/4)(Y+Z)-d\cdot X,(\Lf^{-1}\upsilon/2)\Af-C(\mathsf{N})(X+\Gamma (Y + Z))\big)\geq (C_{2}^*)^{-1}\Af,
\end{equation} which proves (\ref{eq.case2_2}) as $\Af=(C_{12}^*)^{-1}\rho$.
\end{proof}
Note that Proposition \ref{prop.case2} deals with the case where there exist many strongly degenerate primitive pairs. As for the non-primitive pairs, it turns out that they always lead to good molecules if they are cut as \{$3_2 3_2A$\}, as shown in the following lemma.
\begin{lemma}\label{lem.alg_up_ex}
For each two-atom molecule $\{\nf,\nf'\}$, which is either a $\{3_2 3_2 A\}$ molecule, a $\{3_2 (k+1)_{k}\}$ molecule, or a $\{3_2 (l+2)_{l}\}$ molecule, is obtained \emph{in any cutting sequence after any deletion of O-connection}, but $\{\nf,\nf'\}$ \emph{originally belong to the molecule} (not newly formed in partially cutting or deleting). If the atom pair $\{\nf,\nf'\}$ is \emph{not primitive in $\Mb$ before any operation}, then by splitting into 2 sub-cases, we can either make this molecule good, or cut it into a good molecule and a normal molecule.
\end{lemma}
\begin{proof} The proof is the same as Lemma 12.6 in \cite{DHM24}.
\end{proof}

\begin{reduct}\label{red.1} With Propositions \ref{prop.case1} and \ref{prop.case2} and Lemma \ref{lem.alg_up_ex} proved, from now on, we may always assume in Proposition \ref{prop.comb_est} that $|H|\leq (C_{12}^*)^{-1}\rho$ and the number of disjoint strongly degenerate primitive pairs is at most $(C_{12}^*)^{-1}\rho$.

Using also Lemma \ref{lem.alg_up_ex}, we see that for each $\{3_2 (l+1)_l\}$ ($\{3_2 3_2 A\}$ if $l =2$) molecule, if both atoms belong to the original molecule, then it must be good with at most $(C_{12}^*)^{-1}\rho$ exceptions (these exceptions are when the corresponding atom pair is both strongly degenerate and primitive; if it is not strongly degenerate then we get a good molecule by Definition \ref{def.good_normal}, if it is not primitive we get a good molecule by Lemma \ref{lem.alg_up_ex}). But those $\{3_2 (l+1)_l\}$ with at least one atom not in the original molecule, we can only regard it as a normal molecule, unless there are extra restrictions.

As such, we see that Proposition \ref{prop.comb_est} will follow if we can prove that
\begin{equation}\label{eq.overall_alg_3} (\upsilon/2)\cdot(\#_{\mathrm{good}})-d\cdot(\#_{\mathrm{bad}})\geq (C_{11}^*)^{-1}\cdot\rho
\end{equation}
under the assumption that \textbf{all \{$3_2 3_2A$\}, $\{3_2 (l+1)_l\}$ and $\{3_2 (l+2)_l\}$ molecules are good, if both atoms originally belong to the  molecule $\Mb$, rather than created by deleting or partial cutting.} This will be assumed throughout the rest of the proof.
\end{reduct}

\subsection{Layer refinement and the large $\Rf$ case} \label{sec.layer_refine} In this subsection, we introduce the layer refining process and prove the next case of Proposition \ref{prop.comb_est} when $\Rf$ is large, namely Proposition \ref{prop.comb_est_case3}. In layer refinement, we {\bf view each box as a whole}, i.e. regard as one atom and do not count the cycle inside the box. Then the result follows directly from Section 12.3 of \cite{DHM24}.
\begin{lemma}\label{lem.layer_refine}  We regard each box as a whole in Layer refinement.  Let $\Mb$ be a connected molecule with $\rho(\Mb)\leq\Gamma$, see Definition \ref{it.recollision_number_box}. Then we can divide $\Mb$ into at most $\Gamma+1$ disjoint subsets $\Mb_{(k)}$ for $1\leq k\leq\gamma\leq\Gamma+1$, such that (i) each $\Mb_{(k)}$ is a tree, and (ii) if $k<k'$ then no atom in $\Mb_{(k)}$ can be parent of any atom in $\Mb_{(k')}$. Note that, unlike the original layer defined by time intervals, there is \emph{no time inequality} between atoms in different thin layers within the same original layer, other than those coming from parent-child relations.
\end{lemma}

\begin{definition}[Layer refining]
\label{def.layer_refine} Let $\Mb$ be a molecule formed by layers $\Mb_{\ell'}$, as in Proposition \ref{prop.mol_axiom}. We define the following operation of \textbf{layer refining} to turn $\Mb_{\ell'}$ into a new layered molecule $\widetilde{\Mb}_{\ell'}$. Here the atoms and edges of $\widetilde{\Mb}_{\ell'}$ are the same as $\Mb_{\ell'}$ except some O-connections being deleted as in Definition \ref{def.delete} (for convenience, we will not distinguish between $\Mb_{\ell'}$ and $\widetilde{\Mb}_{\ell'}$ below), but the layering of the atoms are different: the atoms in $\widetilde{\Mb}_{\ell'}$ are divided into at most $\Gamma+2$ disjoint subsets, each of them indicated by a refined layer number. We will call these \textbf{thin layers}. The layer refining process is defined as follows:
\begin{enumerate}
\item\label{it.thin_1} Recall all the clusters (i.e. C-molecules) in $\Mb_{\ell'}$. We repeatedly delete O-connections in $\Mb_{\ell'}$, as long as each deletion disconnects two \emph{cyclic} clusters $\Qb_j$ (i.e. these two clusters belong to different components after deletion), until no more O-atoms can be deleted. For each connected component of $\Mb_{\ell'}$ that does not contain a cyclic cluster, we assign thin layer $(\ell',0)$ to all its atoms; note that with this definition, the thin layer $(\ell',0)$ plays a special role different from those thin layers $(\ell',k)$ for $k\geq 1$.
\item\label{it.thin_2} Apply Lemma \ref{lem.layer_refine} to each component $\Qb$ of $\Mb_{\ell'}$ that contains a cyclic cluster (we are guaranteed to have $\rho(\Qb)\leq\Gamma$, which will be proved in Proposition \ref{prop.layer_refine_2} below), to divide it into $\Qb_k$, where $1\leq k\leq \Gamma+1$. Then, for each atom in $\Mb_{\ell'}$, assign it thin layer $(\ell',k)$ if it belongs to $\Qb_k$ for some component $\Qb$.
\end{enumerate}

The thin layers defined above are linearly ordered by lexicographic ordering (where we first compare $\ell'$ and then compare $k$); if layer $\ell'$ is not refined, we identify it with $(\ell',0)$. Note that the original layer $\ell'$ precisely contains all the thin layers $(\ell',k)$. For any $\zeta$ and $\xi$, we define $\Mb_\zeta^T$ to be the set of atoms with thin layer $\zeta$; the associated notions $\Mb_{(\zeta:\xi)}^T$, $\Mb_{[\zeta:\eta]}^T$ etc. are defined similarly.
\end{definition}
In Proposition \ref{prop.layer_refine_2} below we state and prove some properties of the molecule $\Mb$ after layer refining, which corresponds to Proposition 12.9 in \cite{DHM24}.
\begin{proposition}
\label{prop.layer_refine_2} If we perform the layer refining for $\ell'$ as in Definition \ref{def.layer_refine}, then
\begin{enumerate}
\item\label{it.layer_refine_21} The number of O-connections deleted in the process is at most $\Rf_{\ell'}:=\rho(\Mb_{\ell'})$.
\item\label{it.layer_refine_21+} After deleting the O-connections, each component $\Qb$ of $\Mb_{\ell'}$ that contains a cyclic cluster satisfies $\rho(\Qb)\leq\Gamma$, in particular Lemma \ref{lem.layer_refine} is applicable in Definition \ref{def.layer_refine}.
\item\label{it.layer_refine_22} Each $\Mb_{(\ell',k)}^T$ (where $1\leq k\leq \Gamma+1$) is a forest, and has most $\Rf_{\ell'}$ connected components.
\item\label{it.layer_refine_23} $\Mb_{(\ell',0)}^T$ is a forest, and the number of components of $\Mb_{(\ell',0)}^T$ is not more than the number of components of $\Mb_{\ell'}$ before refining.
\end{enumerate}

Now consider the result of refining any set of layers in $\Mb$, which we still denote by $\Mb$ for convenience. Then we have the followings
\begin{enumerate}
    \setcounter{enumi}{4}
\item\label{it.layer_refine_24} For $\zeta<\xi$, no atom in $\Mb_\zeta^T$ can be parent of any atom in $\Mb_\xi^T$.
\item\label{it.layer_refine_25} Each connected component of $\Mb_{(\ell',0)}^T$ (where $\ell'\leq\ell$) must intersect some particle line $\pb\in r(\Mb_{\ell'})$.
\item\label{it.layer_refine_26} With at most $2\Rf_{\ell'}$ exceptions, for each particle line $\pb\in p(\Mb_{\ell'})$, either $\pb$ intersects a component of $\Mb_{(\ell',0)}^T$, or it intersects a component of $\Mb_{(\ell',k)}^T$ for some $1\leq k\leq\Gamma+1$, or $\pb\in r(\Mb_{\ell'})$.
\item\label{it.layer_refine_27} Recall $\Rf=\sum_{\ell'}\Rf_{\ell'}$. For each $\ell'<\ell$, with at most $2\Rf$ exceptions, each particle line $\pb\in r(\Mb_{\ell'})$ is connected to a cycle or another particle line $\pb'\in r(\Mb_{\ell'})$ via the subset $\Mb_{<(\ell'+1,0)}^T$.
\end{enumerate}
\end{proposition}
 \begin{proof}
    The proof is as Proposition 12.9 in \cite{DHM24}, by viewing each box as an atom. The only change is here we delete O-connections instead of O-atoms, which does not affect the proof.
 \end{proof}
Before stating and proving Proposition \ref{prop.comb_est_case3} below (i.e. the case $\Rf\gtrsim\rho$), we need one more auxiliary result concerning a special case of \textbf{UP}. The intuition will be explained in the proof of Proposition \ref{prop.comb_est_case3}.
\begin{proposition}
\label{prop.alg_up_recl} Suppose $\Mb$ is a full molecule divided into two subsets $\Mb_U$ and $\Mb_D$, such that (i) $\Mb_U$ is a forest (by viewing each box as a whole), and is in the same layer, and (ii) no atom in $\Mb_D$ is parent of atom in $\Mb_U$. Define $\#_{\mathrm{comp}(\Mb)}$ to be the number of components of $\Mb$, and define $\#_{\mathrm{cyc.multi}(\Mb_D)}$ to be the number of components of $\Mb_D$ that contains a cycle, or a box, or a multi-particle atom within this component. Then, after applying the algorithm \textbf{UP}, but with the following modifications:
\begin{enumerate} 
    \item\label{it.cbm_1} When we choose the lowest atom $\nf$ with fixed ends in Definition \ref{def.alg_up} \eqref{it.alg_up_2}, we first choose the atom $\nf$ in $\Mb_D$.
    \item\label{it.cbm_2} When we cut some $\{3_2\}$ atom $\nf$ in $\Mb_U$, together with those possible $\nf^\pm$ in Definition \ref{def.alg_up} \eqref{it.alg_up_3}, we first choose those $\nf^-$ in a component of $\Mb_D$ with cycle, multi-particle atom or box, and cut $\{\nf,\nf^-\}$ as free.
    \item\label{it.cbm_3} For each ${k_l}$ atom (with $l \geq 3$) with only C-constituents, we further cut its constituents by \textbf{PM}
\end{enumerate}
We have that $\#_{\{3_2 3_2B\}}=\#_{\{4_2(2l)_{l}\}}=0$, and
\begin{equation}\label{eq.alg_up_recl}
    \#_{\mathrm{bad}} \leq \#_{\mathrm{comp}(\Mb)},\qquad \#_{\mathrm{good}}\geq 1/(10d) \cdot \big(\#_{\mathrm{cyc.multi(\Mb_D)}} - C(\mathsf{N})\Bf(\Mb_U) - C \Xf(\Mb_U)\big),
\end{equation}
with the assumption that all $\{3_2 k_l\}$ molecules are good (not including special \{$3_2 k_l$\} molecules).
\end{proposition}
\begin{proof} We only prove the statement about $\#_{\mathrm{good}}$, as the statements about $\#_{\{3_23_2B\}}$, $\#_{\{4_2(2l)_{l}\}}$ and $\#_{\mathrm{bad}}$ follow from Proposition \ref{prop.alg_up}.
    
\textbf{Proof part 1.} Note that at this stage, each special \{$3_2 k_l$\} molecule \emph{is not regarded as a good molecule}. To deal with this, we prove the following claim: if $A$ is a component of $\Mb_D$ that contains a cycle, or has a multi-particle atom, or has a box, and no atom in $A$ belongs to a good molecule, then at least one atom in $A$ must be fully fixed when it is cut or cut as a normal $\{k_l\}$ atom $(k\geq l+1 \geq4$).
    
To prove the claim, we may think each cycle as formed by ov-segments, and consider the endpoints of these ov-segments. Consider the first cutting step at which either a multi-particle atom in $A$ is involved or every endpoint of some cycle of \(A\) has been cut. We know that, at this stage, no partial cutting happens in $A$. Indeed, if partial cutting happens in Definition \ref{def.alg_up} \eqref{it.alg_up_BM2}, then we get a good MO-molecule by Proposition \ref{prop.alg_up} \eqref{it.up_proof_MO}. If partial cutting happens in Definition \ref{def.alg_up} \eqref{it.alg_up_BM1} or \eqref{it.alg_up_BM3}, then the top/bottom part $B$ has been cut, then there must exist a cycle whose endpoints have already been cut, as the top part and bottom part of a box have at least one cycle.

Assume such first step involves a multi-particle atom $\pf$. If there is no fully fixed atom and good atom in $A$, $\pf$ must be cut as a normal multi-particle atom, or cut together with some special $\{3_2\}$ atom $\qf'$ that appears by partial cutting some atom $\qf \in \Mb_U$. Consider the later case, and denote $\nf$ to be the lowest atom with fixed end chosen in Definition \ref{def.alg_up} \eqref{it.alg_up_2}. Then there must be two different paths from $\pf$ to $\nf$. Indeed, when we choose $\nf$, no atom in $\Mb_D$ has fixed ends (by \eqref{it.cbm_1}), but $\pf$ has fixed ends when cutting $\pf'$. So there exists a path other than the ov-segment connecting $\pf$ and $\qf$. Moreover, $\qf$ also have a path only in $\Mb_U$ to $\nf$ as $\qf$ has no fixed end when we choose $\nf$, (because we choose the lowest atom with fixed in the sense by viewing each box as a whole). So there exists a path from $\nf \to \pf \to \qf \to \nf$. This forms a cycle, and there exists a O-bond at $\qf$ which is not serial to another O-edge inside the cycle. Hence, they are not in the same layer by Proposition \ref{prop.mol_axiom} \eqref{it.axiom2}. By assumption all atoms in $\Mb_U$ belongs to the same layer, hence $\pf$ and $\qf'$ are in different layers, and $\{\qf',\pf\}$ is a good atom; contradiction. 


Now we assume such first step involves a cycle $\theta$. By the above analysis, we know that there exists no partial cutting in $A$ until all atoms in $\theta$ has been cut. Let $\pf$ be the one among the endpoints of $\theta$ that is cut \emph{after all the other endpoints}. Then if $\pf$ is cut together with some atom $\qf$, but forms a normal $\{3_2k_l\}$ molecule. Then this $\qf$ must be in $\Mb_U$ as no partial cutting happens in $A$ so far (otherwise $\{\pf,\qf\}$ forms a good molecule). By the same arguments as above, we can prove that $\pf$ and $\qf$ are in different layers; hence it is a good molecule. Therefore, $\pf$ must be cut by itself. Moreover, $\pf$ is ov-adjacent to other endpoints, which are cut before $\pf$. Therefore, $\pf$ must be cut as fully fixed as there is no good molecule in $A$ (if $\pf$ is a binary atom, then it has at least two fixed ends; if $\pf$ is a multi-particle atom, then it must be fully fixed otherwise it is a good molecule), hence prove the claim.

\textbf{Proof part 2.} Now, assume that $A$ has no good molecule, and suppose $\pf\in A$ is the first fully fixed atom or normal multi-particle atom in $A$ when it is cut. 

\emph{Case 1.} We first assume $\pf$ is a fully fixed atom. By the above analysis we know that $\pf$ cannot be an M-atom (or originally belong to a MO-atom, otherwise either the top/bottom part of a box has been cut, or there exists a good M-atom as in Proposition \ref{prop.alg_up} \eqref{it.up_proof_MO}). Now consider the cutting operation that turns $\pf$ into fully fixed. Then this cutting operation must involve an atom $\qf$ that is ov-adjacent to $\pf$; note that this $\qf$ belongs to either $A$ or $\Mb_U$. If $\qf\in A$, then $\qf$ cannot be fully fixed when it is cut (as $\pf$ is the first); this $\qf$ also cannot be bad when it is cut (because we cut a bad atom only when every atom has no fixed end, and $\qf$ must be a $\{4_2\}$ atom due to the absence of good molecule in $A$. Cutting a bad $\{4_2\}$ atom $\qf$ cannot turn any atom $\pf$ without fixed end into fully fixed, cf. Proposition \ref{prop.mol_axiom} (\ref{it.axiom1})). This $\qf$ must be a binary atom when it is cut. Note that by Definition \ref{def.alg_up} \eqref{it.alg_up_3}, $\qf$ cannot be cut by itself as it is ov-adjacent to $\pf$. If $\qf$ is cut with an atom $\nf' \in A$, then it must be a good molecule due to no partial cutting happens in $A$. If it is cut with a special atom $\nf' \in \Mb_U$, but they are in different layers by the same arguments as above, which is also a good molecule.

Now we have $\qf \in \Mb_U$. This $\qf$ may be formed in some partial cutting, or originally belongs to the $\Mb_U$. Define $\nf$ to be the atom chosen in Definition \ref{def.alg_up} \eqref{it.alg_up_2}. It is clear that $\qf$ is not a MO-atom, otherwise the bond connecting $\pf$ and $\qf$ must be a O-bond at $\pf$. Note that in Definition \ref{def.alg_up} \eqref{it.alg_up_4}, we do not break any such bond when we cut an MO-atom.

Now assume $\qf$ is chosen as $\mf$ in Definition \ref{def.alg_up} \eqref{it.alg_up_3}. If $\qf$ is a $\{k_l\}$ \ $(2 < l <  k \leq 2l )$ atom, then it must be cut by itself. If $\qf \neq \nf $ has several constituents, then at most one constituent of $\qf$ has fixed ends (as $\Mb_U$ is a tree in a single layer, by viewing each box as a whole), hence we always have $k \geq 2l-(d+1)$ (as each constituent of a multi-particle atom has at most $d+1 $ particles). Hence, we can gain excess $\eps^{\sigma_1 \upsilon_0/2}$, where $\sigma_1 = \min \big(k-l, l-2 \big) \geq \max(1, l-(d+1)) $. Moreover, this step may turn at most $l$ atoms in $\Mb_D$ into fully fixed, and those components may do not contain any good molecule. If $\qf= \nf$, then its C-fixed end must be bottom, and cutting $\qf$ turns at most $k-l$ atoms in $\Mb_D$ into fully fixed, and those components may do not contain any good molecule. After partially cutting $\qf$ as in \eqref{it.cbm_3} (each special $\{3_2\}$ atom is good, thanks to the support on time variable), we see that gaining each good molecule costs at most $2d$ components in $\Mb_D$.

If $\qf$ is chosen as $\mf$ in Definition \ref{def.alg_up} \eqref{it.alg_up_3} and is a $\{3_2\}$ atom. It must be cut together with another atom in $\Mb_D$ by our modification of the algorithm. This molecule must also be good, because if one of the atoms is a special atom, then the these two atoms must be in different layers, by the same argument as above. Such a good molecule can be obtained by at most $2$ components in $\Mb_D$. If $\qf$ is already a fully fixed atom, then it has no bottom fixed end (as it is connected to an atom in $A \subset \Mb_D$) and it must be a lowest C-atom of some box (it cannot be a CO-atom, as each CO-atom in $\Mb_U$ cannot be cut as fully fixed as $\Mb_U$ is in the same layer) as $\Mb_U$ is a tree by viewing box as a whole. This costs at most $\mathsf{N} \cdot\#_{\mathrm{box}(\Mb_U)}$ components in $\Mb_D$. 
 
Finally, if $\qf$ is cut with some atom $\qf' \in \Mb_U$, then it must be chosen as $\mf^-$ in Definition \ref{def.alg_up} \eqref{it.alg_up_3} (it can not be chosen as $\mf^+$; otherwise $\qf'$ would be the ov-child of $\qf$ along the ov-segments containing $\pf \in A$, and we would cut $\mf$ with $\pf$ instead with $\qf$ by \eqref{it.cbm_2}). We claim that it must be a good molecule. Indeed, if it is a normal molecule, then either $\qf$ or $\qf'$ (say $\qf'$) is formed by partial cutting some MO-atom $\qf''$. Note that $\qf$ has no fixed end when we choose $\nf$, and becomes a $\{3_2\}$ atom when cutting $\qf'$. Then there must exists a cycle containing $\pf$, such the ov-segment connecting $\pf$ and $\qf''$ is inside the cycle, while the bond at $\qf''$ serial to this bond does not belong to this cycle; contradicts the fact that $\Mb_U$ is in the same layer. Such good molecule can be obtained by at most $d+1$ components in $\Mb_D$, as $\qf$ cannot be a CO-atom as $\Mb_U$ is in a single layer and each C-atom has at most $d+1$ particles.

\emph{Case 2.} In this case, we assume that $\pf$ is a normal $\{(l+1)_l\}$ ($k > l \geq 3$) atom. Then there exist fixed ends $e_i$, $e_j$ at $\pf$ such that  $|v_i- v_j| \leq \eps^{\upsilon_0/2}$ at $\pf$. Consider the atoms (say $\pf_i$ and $\pf_j$, and it is possible that $\pf_i = \pf_j $ due to double bond) that is connected to $\pf$ by $e_i,e_j$ before cutting. They cannot be both in $A$, otherwise at least one of them would be a good $\{3_2\}$ molecule (as $\pf$ is the first normal multi-particle atom, all atom cut before in $A$ must be binary atom). If at least one of them belong to $\Mb_U$ (say $\pf_i$). Assume that $\pf_i$ is cut later than $\pf_j$. With at most  $C(\mathsf{N})\Bf(\Mb_U) + C \Xf(\Mb_U)$ exceptions, $\pf_i$ is a binary atom. We claim it cannot be fully fixed. Indeed, as it has a ov-child in $A$, its bottom end cannot be all fixed ends. As it has at most one path to the $\nf$ chosen in Definition \ref{def.alg_up} \eqref{it.alg_up1_2}, its top end cannot be all fixed ends. Moreover, the restriction $|v_i - v_j|\leq \eps^{\upsilon_0/2}$ makes $\pf_i$ be a good $\{3_2\}$ molecule. Hence the result. Therefore $\pf_i$ must belong to a good molecule (either a good $\{3_2\}$ molecule, or a good two-atom molecule), and each such good molecule cost at most $10d$ components in $\Mb_D$.

As a result, combining Case 1 and Case 2, we have the total number of good molecule  is at least
$$1/(10d) \cdot \big(\#_{\mathrm{cyc.multi(\Mb_D)}}  - C(\mathsf{N})\Bf(\Mb_U) - C \Xf(\Mb_U)\big).$$ 
 
\end{proof}

We are now ready to prove Proposition \ref{prop.comb_est}, in the case when $\Rf + \Bf + \Xf \geq (C_{10}^*)^{-1}\rho$.
\begin{proposition}
\label{prop.comb_est_case3}
Proposition \ref{prop.comb_est} is true, if $\Rf +\Bf + \Xf \geq (C_{10}^*)^{-1}\rho$.
\end{proposition}
\begin{proof} With Proposition \ref{prop.alg_up_recl}, we can prove the proposition by similar arguments as in Proposition 12.11 in \cite{DHM24}. We only outline the idea of the proof. We first choose a maximal $\ell_1\in[1:\ell]$ such that
\begin{equation}\label{eq.comb_est_case31}
    \begin{aligned}
       & {\Bf}_{\ell_1} +\Xf_{\ell_1} + \Rf_{\ell_1}>(C_{6}^*)^{-\ell_1}\cdot (C_{10}^*)^{-1}\rho;\\
       &  {\Bf}_{\ell'}+\Xf_{\ell'}+\Rf_{\ell'}\leq (C_{6}^*)^{-\ell'}\cdot (C_{10}^*)^{-1}\rho< (C_{6}^*)^{-1}\cdot ({\Bf}_{\ell_1} +\Xf_{\ell_1}+ \Rf_{\ell_1}),\quad\forall \ell_1<\ell'\leq\ell.
    \end{aligned}
\end{equation} 

Consider all the components of $\Mb_{\ell_1}$ (i.e. layer $\ell_1$) that contain at least one cyclic cluster, a box, or a multi-particle atom, let them be $\Xb_j\,(1\leq j\leq N)$ where $N$ is the total number of such components, and let their union be $\Xb$. If $N< (C_{5}^*)^{-1}\cdot\big(\Rf_{\ell_1} +{\Bf}_{\ell_1}+\Xf_{\ell_1}\big)$, then we 
\begin{enumerate}[{(i)}]
\item Cut $\Xb$ as free, then cut $\Mb_{>\ell_1}$ as free from $\Mb\backslash\Xb$ and cut it into elementary molecules using \textbf{UP};
\item Then cut $\Mb_{<\ell_1}\cup (\Mb_{\ell_1}\backslash \Xb)$ into elementary molecules using \textbf{DOWN}.
\end{enumerate} 
Then by repeating (PL) argument, the number of bad molecules in this step is at most $2 |H|$. We consider the following two cases.

If $\Xf_{\ell_1} \geq C^*_2(\Rf_{\ell_1} + \Bf_{\ell_1})$, we apply {\bf UP2} algorithm.
After this, we can gain at least $\frac{\Xf_{\ell_1}}{3}$ multi-particle atoms. Indeed, for each CO-atom with at least two atom constituents, the first C-constituent may be cut by itself without any O connections, and such C-constituent may be a binary C-atom, which is not a good molecule, but such CO-atom must have at least another atom constituent, and all other C-constituents are cut as multi-particle atoms. If there is no other C-constituents, then there exists a box constituent, and the number of such BO-atom is no more than $\Bf_{\ell_1} \leq (C_2^*)^{-1} \Xf_{\ell_1}$. Moreover, apart from $\Rf_{\ell_1} +C\Bf_{\ell_1}$ exceptions, those multi-particle atoms are cut with at most one fixed end. Hence, we can gain at least $\frac{\Xf_{\ell_1}}{4}$ good molecules and at most $N$ bad molecules, which easily implies \eqref{eq.case2_2}, in view of $N \leq (C^*_5)^{-1}\cdot (\Rf_{\ell_1} + \Bf_{\ell_1}+ \Xf_{\ell_1}) $ and $\Xf_{\ell_1} \geq C^*_2(\Rf_{\ell_1} + \Bf_{\ell_1})$.

Now if $\Xf_{\ell_1} \leq C^*_2(\Rf_{\ell_1} + \Bf_{\ell_1})$, we argue as \textbf{Choice 2} in the proof of Proposition \ref{prop.case2}. The same proof of Proposition \ref{prop.case2} implies that each cluster in $\Xb_j$ that has a cycle or a box provides one good molecule. So for the contribution of the cuttings involving $\Xb$, we have $\#_{\mathrm{good}}\geq (2\Gamma)^{-1}\cdot \big(\Rf_{\ell_1} + \Bf_{\ell_1} \big) $ and $\#_{\mathrm{bad}}\leq N$, which easily implies (\ref{eq.case2_2}), in view of $N< (C_{5}^*)^{-1}\cdot\Rf_{\ell_1}$, and (\ref{eq.comb_est_case31}).

Now assume $N\geq (C_{5}^*)^{-1}(\Rf_{\ell_1} + \Bf_{\ell_1}+ \Xf_{\ell_1})$ (in particular we must have $\ell_1<\ell$ as $|H|\leq (C^*_{12})^{-1} \rho$). We then perform layer refining (Definition \ref{def.layer_refine}) for each layer $\ell_1<\ell'\leq\ell$. After that, for each thin layer $\zeta\geq(\ell_1,0)$ (where we identify the unrefined layer $\ell_1$ with $(\ell_1,0)$), define
\begin{equation}\label{eq.comb_est_case32}
     N_\zeta:=\textrm{number of components of $\Mb_{[(\ell_1,0):\zeta]}^T$ that contains some $\Xb_j\,(1\leq j\leq N)$}.
\end{equation}
 Note that for $\zeta=(\ell_1,0)$ we have $N_\zeta=N$, and for $\zeta$ being the highest thin layer of $\Mb$, we have 
\begin{equation}\label{eq.comb_est_case32.5}
    N_\zeta\leq |H|+(2\Gamma+4)\sum_{\ell'>\ell_1}\Rf_{\ell'}<C_{2}^*(C_{6}^*)^{-1}\cdot (\Rf_{\ell_1} + \Bf_{\ell_1}+ \Xf_{\ell_1});
\end{equation}
this is proved by repeating the (PL) argument in the proof of Proposition \ref{prop.box_multiple_atom} but using Proposition \ref{prop.layer_refine_2} (\ref{it.layer_refine_22}), (\ref{it.layer_refine_25}) and (\ref{it.layer_refine_26}).

Now, using (\ref{eq.comb_est_case32.5}), we see that there exists a (maximal) thin layer $\zeta\geq(\ell_1,0)$, which is not the highest thin layer of $\Mb$, such that
\begin{equation}\label{eq.comb_est_case33}
     \begin{aligned}
        & N_\zeta\geq (C_{5}^*)^{n((\ell_1,0))-n(\zeta)-1}\cdot (\Rf_{\ell_1} + \Bf_{\ell_1}+ \Xf_{\ell_1});\\
        & N_{\zeta'}<(C_{5}^*)^{n((\ell_1,0))-n(\zeta')-1}\cdot (\Rf_{\ell_1} + \Bf_{\ell_1}+ \Xf_{\ell_1})\leq (C_{5}^*)^{-1}\cdot N_\zeta,\quad \forall \zeta'>\zeta,
     \end{aligned}
\end{equation} where $n(\zeta)$ is the number of thin layers below or equal to $\zeta$. Let $\zeta^+$ be the lowest thin layer above $\zeta$, then define $\Yb$ to be the union of all components of $\Mb_{[(\ell_1,0):\zeta^+]}^T$ that contains some $\Xb_j\,(1\leq j\leq N)$, and write it as $\Yb=\Yb_U\cup\Yb_D$, where $\Yb_U=\Yb\cap \Mb_{\zeta^+}^T$ and $\Yb_D=\Yb\cap\Mb_{[(\ell_1,0):\zeta]}^T$. It is easy to verify that the assumptions of Proposition \ref{prop.alg_up_recl} are satisfied, and that
\begin{equation}\label{eq.comb_est_case34} \#_{\mathrm{comp}(\Yb)}\leq N_{\zeta^+},\qquad \#_{\mathrm{cyc.multi}(\Yb_D)}\geq N_\zeta
\end{equation} (note that each $\Xb_j$ already contains a cycle, box or multi-particle atom, and so does any component that contains some $\Xb_j$; therefore each of the $N_\zeta$ components of $\Mb_{[(\ell_1,0):\zeta]}^T$ containing some $\Xb_j$ will contribute to $\#_{\mathrm{cyc.comp}(\Yb_D)}$).

Then, as above, we 
\begin{enumerate}[{(i)}]
\item Cut $\Yb$ as free, then cut $\Mb_{>\zeta^+}^T$ as free from $\Mb\backslash\Yb$  and subsequently cut it using \textbf{UP};
\item Then cut $\Mb_{<(\ell_1,0)}^T\cup (\Mb_{[(\ell_1,0):\zeta^+]}^T\backslash \Yb)$ into elementary molecules using \textbf{DOWN}.
\end{enumerate} It is easy to see that this can always be done (as there is no bond between $\Yb$ and $\Mb_{[(\ell_1,0):\zeta^+]}^T\backslash \Yb$), and the contribution of these steps to $\#_{\mathrm{bad}}$ is at most 
\begin{equation}\label{eq.comb_est_case35}
2|H|+8(\Gamma+1)\sum_{\ell'>\ell_1}\Rf_{\ell'}<C_{2}^*(C_{6}^*)^{-1}\cdot( \Rf_{\ell_1} +\Bf_{\ell_1} +\Xf_{\ell_1});
\end{equation}
this follows from the upper bounds on the number of components of $\Mb_{>\zeta^+}^T$ and $\Mb$, which is proved by repeating the above (PL) arguments which uses Proposition \ref{prop.layer_refine_2} (\ref{it.layer_refine_22}), (\ref{it.layer_refine_25}) and (\ref{it.layer_refine_26}).

Finally, let $C_0 =C(\mathsf{N})$,
\begin{enumerate}
    \item if $\Bf(\Yb_U) +\Xf(\Yb_U)\geq  N_\zeta / (2C_0)$, we first cut $\Yb_U$ as free, and then apply \textbf{UP}, but with the modification that \emph{we view box as a single atom} and we cut one atom each time. After that, each box has at most one fixed end since $\Yb_U$ is a tree. We then apply \textbf{UP} or \textbf{DOWN} to each box; moreover, each constituent contribute to $\Xf$ must be a good multi-particle atom (has at most one fixed end), as $\Yb_U$ is a forest.
    \item otherwise, we cut $\Yb$ into elementary molecules using \textbf{UP},
\end{enumerate}
In the first case, for the contribution of the cuttings involving $\Yb_U$, we have 
\begin{equation}\label{eq.case_MU_box}
     \#_{\mathrm{good}} \geq C^*_2 (C_7^*)^{-1}\cdot (\Rf_{\ell_1}+\Bf_{\ell_1} + \Xf_{\ell_1}),\qquad\#_{\mathrm{bad}}\leq N_{\zeta^+}. 
\end{equation}
In the second case, with at most $2 \sum_{\ell'> \ell_1} \Rf_{\ell'}$ exceptions (this is because in the layer refinement, there may produce some strongly degenerate primitive pairs, but the number of such pair does not exceed $2 \sum_{\ell'\geq \ell_1} \Rf_{\ell'}$), each non-special $\{3_2 k_l\}$ molecule is good. By Proposition \ref{prop.alg_up_recl}, \eqref{eq.comb_est_case34} and $\Bf(\Yb_U) +\Xf(\Yb_U) < N_\zeta / (2C_0)$, we have 
\begin{equation}\label{eq.case35.5}\#_{\{\mathrm{good}\}}\geq N_\zeta/(20d) ,\qquad\#_{\mathrm{bad}}\leq N_{\zeta^+}.
\end{equation} By putting together (\ref{eq.comb_est_case31}), (\ref{eq.comb_est_case33}), (\ref{eq.comb_est_case35}) and (\ref{eq.case35.5}), we get
\begin{equation}\label{eq.comb_est_case37}
(\upsilon/2)\cdot(\#_{\mathrm{good}})-d\cdot(\#_{\mathrm{bad}})\geq (\upsilon/(80d))\cdot N_\zeta-d\cdot (N_{\zeta^+}+C_2^*(C_6^*)^{-1}\cdot\Rf_{\ell_1})
\geq (\upsilon/160)\cdot N_\zeta\geq (C_{10}^*)^{-2}\cdot\rho,
\end{equation} which proves (\ref{eq.case2_2}). The proof is now complete.
\end{proof}

To deal with the special atoms, we introduce the definition of short degeneracy. 
\begin{definition}\label{def.short_degeneracy}
    Let $\Mb$ be a full molecule without any operation. An atom pair $\{\nf,\nf'\}$ is called \textbf{short degeneracy}, if $\nf$ and $\nf'$ are in the same layer, and are ov-adjacent, such that either is an MO-atom or CO-atom or binary atom, with at least one is an MO-atom or CO-atom. Moreover, the ov-segment that connects $\{\nf,\nf'\}$ is an O-bond at the MO-atom or CO-atom. Moreover, $|t_1 - t_2|\leq \eps^\upsilon$, where $t_1 ,t_2$ denote the initial time of the particle associate with the ov-segment connecting $\nf$ and $\nf'$, gets involved in $\nf$ and $\nf'$. 
\end{definition}
\begin{proposition}\label{prop.short_degeneracy}
Proposition \ref{prop.comb_est} is true if there are at least $(C^*_{9})^{-1} \rho$ disjoint short degenerate pairs, (note that one CO-atom may belong to many degenerate pairs). 
\end{proposition}
\begin{proof}
    By using pigeonhole principle, there exists a layer $\ell'$ with at least $\Cf:=\Lf^{-1} \cdot (C^*_{9})^{-1} \rho$ degeneracies. We define $\Zb$ to be the union of components of $\Mb_{\ell'}$ that contain at least one short degeneracy. We then  
\begin{enumerate}[{(i)}]
\item Cut $\Zb$ as free, then cut $\Mb_{>\ell'}$ as free from $\Mb\backslash\Zb$ and cut it into elementary molecules using \textbf{UP};
\item Then cut $\Mb_{<\ell'}\cup (\Mb_{\ell'}\backslash \Zb)$ into elementary molecules using \textbf{DOWN}.
\end{enumerate} 
It is easy to see that this can always be done, and the contribution of these steps (i.e. coming from atoms outside $\Zb$) to $\#_{\mathrm{bad}}$ is at most $2|H|$ by repeating (PL) argument as above. Moreover, the number of components is no more than $(C^*_{10})^{-1} \rho$, as the number of boxes and multi-particle atoms is no more than $(C^*_{10})^{-1} \rho$.

Now we apply \textbf{UP1} to cut each component of $\Zb$.
By counting the quantity $\sigma = \#_{\mathrm{bo/fr}} - \sum (2l-1) \#_{k_l}$, and similar arguments as in Proposition \ref{prop.case2}, we know that the number of molecules that are not $\{3_2\}$ is at most $\#_{\mathrm{comp}} + \Xf + C(\mathsf{N}) \Bf + \Rf $. With at most these exceptions, for each short degenerate pair, there generates two $\{3_2\}$ atoms (at least one of them is a special $\{3_2\}$ atom). Moreover, at least one of these two atoms is a good molecule due to the restriction on time. We get 
\begin{equation}
    (\upsilon/2)\#_{\mathrm{good}} - d \cdot \#_{\mathrm{bad}} \geq \upsilon/2 \Cf -  (\Xf + C \Bf + \Gamma \cdot \Rf + C|H|)  \geq (C^*_2)^{-1} \cdot \Cf
\end{equation}
which proves \eqref{eq.case2_2}.   
\end{proof}
Using Propositions \ref{prop.case2}--\ref{prop.short_degeneracy}, we see that each $\{3_2 3_2 A\}$ and $\{3_2 k_l \}$ molecules are good, with at most $2(C^*_{9})^{-1}\rho$ exceptions. Indeed, for any such two-atom molecule, if both atoms already exist in the original molecule, then by Proposition \ref{prop.case2}, we know these atoms are good with at most $(C^*_{12})^{-1}\rho$ exceptions. If one of the atoms is formed by deleting, then it must originally contribute to $\Xf$, and the number of such molecules is no more than $(C^*_{10})^{-1}\rho$ by Proposition \ref{prop.comb_est_case3}. If one of atoms are formed by partial cutting, then by Proposition \ref{prop.short_degeneracy}, the number of such molecules is no more than $(C^*_{9})^{-1}\rho$.
\begin{reduct}\label{red.0} As such, we see that Proposition \ref{prop.comb_est} is true if we can prove that 
\begin{equation}\label{eq.overall_alg_4}
    (\upsilon/2) \cdot (\#_{\mathrm{good}}) - d \cdot (\#_{\mathrm{bad}}) \geq (C^*_{8})^{-1} \cdot \rho
\end{equation}
by treating \textbf{all \{$3_23_2A$\} and $\{3_2 k_l\}$ ($k=l+1$ or $k=l+2$, depends on whether there exists a double bond), including those special atoms, are good}.
\end{reduct}
\subsection{Layer selection and weak degeneracies}\label{sec.layer_select} From now on we will focus on proving (\ref{eq.overall_alg_4}) under the assumptions in Reduction \ref{red.0}. The goal of this subsection is to reduce this result to Proposition \ref{prop.case5} by applying the procedure of \textbf{layer selection}. Throughout this section, we regard each box as a whole.
\begin{definition}[Layer selection]
\label{def.layer_select} Let $\Mb$ be a molecule as in Proposition \ref{prop.mol_axiom} and view each box as a single atom. Assume that $|H|\leq (C_{12}^*)^{-1}\rho$ as in Reduction \ref{red.1} and $\Rf+ \Bf + \Xf\leq (C_{10}^*)^{-1}\rho$, and the assumptions in Reduction \ref{red.0}. Then we will perform the following steps to determine the values of two \textbf{thin layers} $\zeta_U$ and $\zeta_D$ (in one case it is possible that $\zeta_D$ is absent and we only determine $\zeta_U$). Recall that $r(\Mb_{\ell})=H$, and that we denote $s_{\ell'}=|r(\Mb_{\ell'})|$.

\begin{enumerate}
\item\label{it.layer_select_1} Choose $ 1 \leq \ell_1\leq \ell$  such that
\begin{equation}\label{eq.layer_select_1}s_{\ell_1-1}\geq (C_{7}^*)^{-(\ell_1-1)}\cdot\rho;\qquad s_{\ell'}<(C_{7}^*)^{-\ell'}\cdot\rho\leq (C_{7}^*)^{-1}
\cdot s_{\ell_1-1},\quad\forall \ell_1\leq\ell'\leq\ell,\end{equation}
 and refine layer $\Mb_{\ell_1}$ as in Definition \ref{def.layer_refine}. Define also $s_{\ell_1-1}:=\rho'$.
\item\label{it.layer_select_2} After refining, consider all particle lines $\pb\in r(\Mb_{\ell_1-1})$ (there are $\rho'$ of them), for any $0\leq\gamma\leq \Gamma+1$, define $u_\gamma^*$ to be the number of these $\rho'$ particle lines that \emph{do not} contain any atom in $\Mb_{(\ell_1,k)}^T$ for any $k\leq\gamma$ (denote $u_{-1}^*:=s_{\ell_1-1}$). Then choose $0\leq\gamma\leq\Gamma+1$ such that
\begin{equation}\label{eq.layer_select_2}
    u_{\gamma-1}^*\geq (C_{6}^*)^{-\gamma}\cdot \rho',\quad u_\gamma^*<(C_{6}^*)^{-\gamma-1}\cdot \rho'\leq (C_{6}^*)^{-1}\cdot u_{\gamma-1}^*.
\end{equation} Fix \textbf{the thin layer $\zeta_U:=(\ell_1,\gamma)$}. Define also $u_{\gamma-1}^*-u_\gamma^*:=\rho''$. If $\ell_1=1$, then we stop the process with $\zeta_D$ absent. 

\item\label{it.layer_select_3} If $\ell_1 >1$, refine the layers $\Mb_{\ell'}$ for all $\ell'\leq\ell_1-1$, and consider all the particle lines $\pb\in r(\Mb_{\ell_1-1})$ that does not contain any atom in $\Mb_{(\ell_1,k)}^T$ for any $k\leq\gamma-1$ but does contain an atom in $\Mb_{(\ell_1,\gamma)}^T$; there are $\rho''$ of them, and define the set of them to be $S$. Let $\zeta_0$ be the highest thin layer below $\zeta_U$. For each thin layer $\zeta\leq\zeta_0$, define $v_\zeta^*$ to be the number of particle lines in $S$ that are either \emph{connected to a cycle} via $\Mb_{[\zeta:\zeta_0]}^T$, or \emph{are connected to another particle line} $\pb'\in r(\Mb_{\ell_1-1})$ via $\Mb_{[\zeta:\zeta_0]}^T$ (this $\pb'$ need not be in $S$). Now, if for the lowest thin layer $(1,0)$ we have $v_{(1,0)}^*\leq (C_{5}^*)^{-1}\cdot\rho''$, then stop the process with $\zeta_D$ absent; otherwise \textbf{choose a thin layer} $\zeta=\zeta_D$ such that
\begin{equation}\label{eq.layer_select_3}v_{\zeta_D}^*\geq (C_{5}^*)^{-n(\zeta_D)}\cdot\rho'',\qquad v_{\zeta}^*<(C_{5}^*)^{-n(\zeta)}\cdot \rho''\leq(C_{5}^*)^{-1}\cdot v_{\zeta_D}^*,\quad \forall \zeta_D<\zeta\leq\zeta_0,\end{equation} where $n(\zeta)$ is the number of thin layers below or equal to $\zeta$. Define also $v_{\zeta_D}^*:=\rho'''$. In the case when $\zeta_D$ is absent, we also define $\rho''':=\rho''$.
\end{enumerate}
\end{definition} 
\begin{proposition}
\label{prop.layer_select} The construction process in Definition \ref{def.layer_select} is well-defined. Moreover, in Definition \ref{def.layer_select}, consider all the connected components of $\Mb_{(\zeta_D:\zeta_U)}^T$ (or $\Mb_{<\zeta_U}^T$ if $\zeta_D$ is absent) that (i) intersect exactly one particle line in $S$ (as defined in Definition \ref{def.layer_select} (\ref{it.layer_select_3})), and (ii) does not intersect any other particle line in $r(\Mb_{\ell_1-1})$, and (iii) does not contain a cycle. Denote the union of these components by $Y$. Then, each of the following quantities is bounded by $10(C_{5}^*)^{-1}\cdot\rho'''$:
\begin{enumerate}
\item\label{it.layer_select_new1} The number of components of $\Mb$, $\Mb_{>\zeta_U}^T$, $\Mb_{>\zeta_D}^T$, $\Mb_{\zeta_U}^T$, and $\Mb_{[\zeta_D:\zeta_U]}^T$ and $\Mb_{(\zeta_D:\zeta_U]}^T$;
\item\label{it.layer_select_new2} The number of components of $\Mb_{\zeta_D}^T$ that are not connected to $\Mb_{\zeta_U}^T\cup\Mb_{(\zeta_D:\zeta_U)}^T$ by a bond.
\item\label{it.layer_select_new3} The number of components of $\Mb_{(\zeta_D:\zeta_U)}^T\backslash Y$;
\end{enumerate}
\end{proposition} 
\begin{proof}
    This is similar to Proposition 12.14 in \cite{DHM24}.
\end{proof}
Now, once the thin layers $\zeta_U$ and $\zeta_D$ are chosen by Definition \ref{def.layer_select}, we need to cut as free all the exceptional components in $\Mb_{(\zeta_D:\zeta_U)}^T\backslash Y$. Together with a few auxiliary operations treating other layers (and cutting components of $\Mb_{\zeta_D}^T$ that have no connections to $\Mb_{\zeta_U}^T\cup Y$ which is purely for convenience), this forms the \textbf{Option 1} in Definition \ref{def.layer_cutting} below, which is the main option to be executed in the rest of the proof. On the other hand, the \textbf{Option 2} in Definition \ref{def.layer_cutting} involves cutting $\Mb_{\zeta_D}^T$ first, and will be used in the \emph{weakly degenerate} case (Definition \ref{def.weadeg}, Proposition \ref{prop.comb_est_case4}). 
\begin{definition}\label{def.layer_cutting} Recall the layer selection process in Definition \ref{def.layer_select} and the definition of $Y$ in Proposition \ref{prop.layer_select}. Define also the set $\Vb$ to be the union of all components of $\Mb_{\zeta_D}^T$ that are connected to $\Mb_{\zeta_U}^T\cup Y$ by at least one bond, and let $G$ be the number of those components of $\Vb$ that are connected to $\Mb_{\zeta_U}^T\cup \Mb_{(\zeta_D:\zeta_U)}^T$ by at least two bonds. Below we define two options of cutting sequences to be applied to $\Mb$.
\begin{enumerate}
\item \textbf{Option 1}. We perform the following cutting operations, by treating each box as a whole (note that if $\zeta_D$ is absent, then the steps here involving $\Mb_{\zeta_D}^T$ and $\Mb_{<\zeta_D}^T$ will be skipped). 
 
\begin{enumerate}[{(i)}]
\item We cut $\Mb_{[\zeta_D:\zeta_U]}^T$ as free, then cut $\Mb_{>\zeta_U}^T$ as free from $\Mb_{>\zeta_U}^T\cup\Mb_{<\zeta_D}^T$ and cut it into elementary molecules using \textbf{UP}, then cut $\Mb_{<\zeta_D}^T$ into elementary molecules using \textbf{DOWN}.
\item We cut as free all the components of $\Mb_{(\zeta_D:\zeta_U)}^T$  (or $\Mb_{<\zeta_U}^T$ if $\zeta_D$ is absent) not in $Y$, and cut each of them into elementary molecules using \textbf{UP}. If any atom in $\Mb_{\zeta_U}^T$ or $\Mb_{\zeta_D}^T$ becomes fully fixed (including those atoms in boxes), then we cut it as free, and repeat until there is no such atom in $\Mb_{\zeta_U}^T$ and $\Mb_{\zeta_D}^T$. 
\item We cut as free all components of $\Mb_{\zeta_D}^T$ that are not connected to $\Mb_{\zeta_U}^T\cup Y$ by a bond, and cut each of them into elementary molecules using \textbf{DOWN}. 
\end{enumerate}
After the cutting sequence in \textbf{Option 1}, define the remaining molecule to be $\Mb_{UD}$  (or $\Mb_{ U}$ if $\zeta_D$ is absent); it contains parts of $\Mb_{\zeta_U}^T$, $\Mb_{\zeta_D}^T$ and $Y$ and will be studied in Section \ref{sec.maincr}.

\item \textbf{Option 2}. We perform the following cutting operations by viewing each box as a whole.
 
\begin{enumerate}[{(i)}]
\item With $\Vb$ defined as above, we cut $\Mb_{\zeta_U}^T\cup\Mb_{(\zeta_D:\zeta_U)}^T\cup \Vb$ as free, then cut $\Mb_{>\zeta_U}^T$ as free from $\Mb_{<\zeta_D}^T\cup \Mb_{>\zeta_U}^T\cup (\Mb_{\zeta_D}^T\backslash\Vb)$ and cut it into elementary molecules using \textbf{UP}, then cut $\Mb_{<\zeta_D}^T\cup (\Mb_{\zeta_D}^T\backslash\Vb)$ into elementary molecules using \textbf{DOWN}.
\item If a component of $\Vb$ has only one bond connecting to $\Mb_{\zeta_U}^T\cup \Mb_{(\zeta_D:\zeta_U)}^T$, then cut it as fixed, otherwise cut it as free. Then cut each component into elementary molecules using \textbf{DOWN}, by viewing each box as a whole, and then cut $\Mb_{\zeta_U}^T\cup \Mb_{(\zeta_D:\zeta_U)}^T$ into elementary molecules using \textbf{UP}.
\end{enumerate}
\end{enumerate}
\end{definition}

\begin{proposition}\label{prop.layer_cutting} Recall $G$ defined in Definition \ref{def.layer_cutting}. For \textbf{Options 1 and 2} in Definition \ref{def.layer_cutting}, we have the followings.

\begin{enumerate}
\item\label{it.layer_cutting_1} The contribution of cuttings in \textbf{Option 1} to $\#_{\mathrm{bad}}$ is at most $100 (C_{5}^*)^{-1}\cdot \rho'''$. If $\zeta_D$ is absent, then after \textbf{Option 1}, we can choose at least $(\rho'''/2)$ particle lines $\pb \in S$ such that 
\begin{enumerate}
\item\label{it.layer_cutting_1i} there exists a subset $T(\pb)$ of $\Mb_U\backslash \Mb_{\zeta_U}^T=Y$ which is either empty or a component of $Y$; 
\item\label{it.layer_cutting_1ii} if $T(\pb)\neq\varnothing$ then it is a tree and has only one bond connecting to $\Mb_{\zeta_U}^T$, which is along $\pb$; 
\item\label{it.layer_cutting_1iii} consider the set $\Ec_0$ of all bottom free ends $e$ at atoms in $T(\pb)$ (or along $\pb$ if $T(\pb)=\varnothing$) for all the chosen $\pb$, then for at least  $\rho'''/4$ chosen particle lines $\pb$, there exist $e\in \Ec_0$ at an atom in $T(\pb)$ or along $\pb$, such that an initial link exists between $e$ and $e'$, where $e'\neq e$ either belongs to $\Ec_0$ or has an atom in some component already cut in \textbf{Option 1}.
\end{enumerate}  

\item\label{it.layer_cutting_2} If $\zeta_D$ exists, after \textbf{Option 2}, $\Mb$ is cut into elementary molecules with $\#_{\{3_23_2B\}}=\#_{\{4_2 (2l)_l\}}=0$. We also have $\#_{\mathrm{bad}}\leq 100(C_{5}^*)^{-1}\cdot\rho'''+G$ (note here $G$ equals the number of components of $\Vb$ that are cut as free in \textbf{Option 2} (ii)). Moreover, each atom in $\Vb$ that does not belong to any box has at most one fixed end.

\item\label{it.layer_cutting_3} If $\zeta_D$ exists, after \textbf{Option 1}, the molecule $\Mb_{UD}$ consists of parts of $\Mb_{\zeta_U}^T$, parts of $\Mb_{\zeta_D}^T$, and $Y$. Define $\Mb_U$ to be the union of $\Mb_{\zeta_U}^T$ and $Y$, and all the components in $\Mb_{\zeta_D}^T$ that are full (after \textbf{Option 1}) and have exactly one bond connecting to $\Mb_{\zeta_U}^T\cup Y$. Let also $\Mb_D:=\Mb_{UD}\backslash \Mb_U$ (which is a subset of $\Mb_{\zeta_D}^T$). Then we have the following properties for $\Mb_U$ and $\Mb_D$:
\begin{enumerate}
\item\label{it.layer_cutting_3i} Each of $\Mb_U$ and $\Mb_D$ is a forest by viewing each box as a whole, and has no fully fixed atoms. No atom in $\Mb_D$ is a parent of any atom in $\Mb_U$, and $\Mb_U$ (resp. $\Mb_D$) has no C-top (resp. C-bottom) fixed end.
\item\label{it.layer_cutting_3ii} The number of full components of $\Mb_U$ is at most $10(C_{5}^*)^{-1}\cdot\rho'''$. The number of components of $\Mb_D$ is $\geq G$. Each component in $\Mb_D$ either has at $\geq 2$ bonds connecting to $\Mb_U$, or has one bond connecting to $\Mb_U$ and one atom with fixed ends.
\item\label{it.layer_cutting_3iii} Define the set of \textbf{UD connections} to be the set of bonds connecting an atom in $\Mb_U$ and an atom in $\Mb_D$. We say two UD connections are \textbf{connected via} $\Mb_D$, if the $\Mb_D$ endpoints of them belong to the same connected component of $\Mb_D$. Then, there exist at least $(\rho'''/4)$ many UD connections, such that each UD connection is either connected to another UD connection via $\Mb_D$, or has its $\Mb_D$ endpoint belonging to component of $\Mb_D$ that contains an atom with fixed ends.
\end{enumerate}
\end{enumerate}
\end{proposition}
\begin{proof}
    The proof is similar to Proposition 12.16 in \cite{DHM24}. We list the changes. All cuttings are done by viewing each box as a whole; when an atom of a box becomes fully fixed in Option 1 (ii), the whole box is considered, and it is cut into elementary molecules as in Definition \ref{def.alg_up1} \eqref{it.alg_up1_3}, which creates no bad molecule. The O-connections deleted in the layer refining (Proposition \ref{prop.layer_refine_2}) do not affect the forest property in (\ref{it.layer_cutting_3iii}), as deleting O-connections does not create cycles. The rest of the proof is unchanged.
\end{proof}
Now we can prove the next case of Proposition \ref{prop.comb_est}, namely Proposition \ref{prop.comb_est_case4} below, which contains the initial link case, and the weak degenerate case (where there are many weak degeneracies in $\Mb_D$ as defined in Definition \ref{def.weadeg}) below. In the weakly degenerate case we choose \textbf{Option 2} in Definition \ref{def.layer_cutting}, which allows us to gain from the weak degeneracy conditions at the $\{3_2\}$ and \{$4_2$\} molecules in $\Mb_D$.
\begin{definition}\label{def.weadeg}
Define the atom pair $\{\nf,\nf'\}$ to be \textbf{weakly degenerate} if 
\begin{enumerate}
    \item they are both binary C-atoms;
    \item they are ov-adjacent; and
    \item we have the restriction (by indicator functions) $|t^-_{\nf}-t^-_{\nf'}|\leq \varepsilon^{\upsilon}$.
\end{enumerate} 
A single atom $\nf$ is called \textbf{weakly degenerate} if it is a binary C-atom, and at least one of the following conditions holds:
\begin{enumerate}
    \item  we have the restriction $|v_{e_i}-v_{e_j}|\leq\varepsilon^{C_\Phi \upsilon}$ for two distinct edges $(e_i,e_j)$ at $\nf$;
    \item  for two bottom edges $e_i,e_j$, we have 
    $$\frac{|(v_{e_i}-v_{e_j})\wedge \nu_{ij}|}{|v_{e_i}-v_{e_j}|}\leq \eps^\upsilon,$$
     where $\nu_{ij} = \eps^{-1} (x_{e_i} -x_{e_j} + t^-_\nf(v_{e_i} -v_{e_j}))$.
     \item we have the restriction $|t^-_\nf - \ell'\tau|\leq \eps^{\upsilon_0}$ for $\ell'= \ell[\nf]$ or $\ell[\nf]+1$.  
\end{enumerate}
\end{definition}
\begin{proposition} 
\label{prop.comb_est_case4}
Suppose $\Mb$ satisfies the assumptions in Reduction \ref{red.1}, and we apply the layer selection process in Definition \ref{def.layer_select}. Then, Proposition \ref{prop.comb_est} is true if $\zeta_D$ is absent.

 If $G\leq  (C_{3}^*)^{-1}\cdot\rho'''$ (with $G$ defined in Proposition \ref{prop.layer_cutting}) and we make the restriction that there exist at least $C_{2}^*)^{-1}\cdot\rho'''$ weakly degenerate atoms and atom pairs in $\Vb$ (by inserting an indicator function $\mathbbm{1}_{\mathrm{wea.deg}}$),  then Proposition \ref{prop.comb_est} is also true.
\end{proposition}
\begin{proof} The proof is the same as Proposition 12.18 as \cite{DHM24}.
\end{proof}
Now, summarizing Propositions \ref{prop.case1}, \ref{prop.case2}, \ref{prop.comb_est_case3} and \ref{prop.comb_est_case4}, we see that, for the rest of the proof of Proposition \ref{prop.comb_est}, we may assume that the assumptions in Reduction \ref{red.1} are satisfied, that $\zeta_D$ in Definition \ref{def.layer_select} exists, and that either $G\geq (C_{3}^*)^{-1}\cdot\rho'''$ with $G$ defined in Proposition \ref{prop.layer_cutting} or we can restrict the number of weakly degenerate atoms and pairs in $\Vb$ to be less than $(C_{2}^*)^{-1}\cdot\rho'''$ (by inserting the indicator function $1-\mathbbm{1}_{\mathrm{wea.deg}}$ opposite to the one defined in Proposition \ref{prop.comb_est_case4}).

We then apply \textbf{Option 1} in Definition \ref{def.layer_cutting} and reduce to the molecule $\Mb_{UD}$ that satisfies (3a)--(3c) of Proposition \ref{prop.layer_cutting}. In particular the number of components of $\Mb_D$ satisfies $\#_{\mathrm{comp}(\Mb_D)}\geq G$ by (3b) of Proposition \ref{prop.layer_cutting}. Using also the upper bound on the contribution to $\#_{\mathrm{bad}}$ of \textbf{Option 1} in Definition \ref{def.layer_cutting} (proved in Proposition \ref{prop.layer_cutting}), we see that Proposition \ref{prop.comb_est} (i.e. (\ref{eq.overall_alg_3})) now follows from the following
\begin{proposition}
\label{prop.case5} Let $\Mb_{UD}$ be a molecule satisfying the properties (3a)--(3c) in Proposition \ref{prop.layer_cutting} with some $\rho'''\geq (C_{8}^*)^{-1}\cdot\rho$. Moreover, assume that either the number of components of $\Mb_D$ satisfies $\#_{\mathrm{comp}(\Mb_D)}\geq (C_{3}^*)^{-1}\cdot\rho'''$, or (we can restrict by some indicator function such that) the number of weakly degenerate atoms and pairs in $\Mb_D$ is less than $(C_{2}^*)^{-1}\cdot\rho'''$. Then in either case, we can cut $\Mb_{UD}$ into elementary molecules such that $\#_{\{4_2 (2l)_l\}}=0$ and each \{$3_23_2$B\} molecule is good, and moreover
\begin{equation}\label{eq.comb_est_case51} (\upsilon/2)\cdot(\#_{\mathrm{good}})-d\cdot(\#_{\mathrm{bad}})\geq (C_{4}^*)^{-1}\cdot \rho''',
\end{equation} where we again understand that all \{$3_2 3_2A$\} molecules are good as in Reduction \ref{red.1}.
\end{proposition}

\section{Cutting algorithm II: analyzing UD molecules}\label{sec.maincr} The goal of this section is to prove Proposition \ref{prop.case5} for the two-layer molecule $\Mb_{UD}$, by introducing the various algorithms. To simplify notation in this section, if a molecule $\Mb_{UD}=\Mb_U\cup\Mb_D$ satisfies that (a) each of $\Mb_U$ and $\Mb_D$ is a forest (by viewing each box as a whole) and has no fully fixed atoms, and (b) no atom in $\Mb_D$ is parent of any atom in $\Mb_U$, then we say it is a \textbf{UD molecule}. If in addition $\Mb_U$ has no C-top fixed end and $\Mb_D$ has no C-bottom fixed end, then we say $\Mb_{UD}$ is \textbf{canonical}.
\begin{reduct} Our proof below relies on the properties (3a)--(3c) in Proposition \ref{prop.layer_cutting}. For convenience, we record them here:
\label{red.2}
\begin{enumerate}[{(i)}]
\item Each of $\Mb_U$ and $\Mb_D$ is a forest by viewing each box as a whole, and has no fully fixed atoms. No atom in $\Mb_D$ is a parent of any atom in $\Mb_U$, and $\Mb_U$ (resp. $\Mb_D$) has no C-top (resp. C-bottom) fixed end.
\item The number of full components of $\Mb_U$ is at most $10(C_{5}^*)^{-1}\cdot\rho'''$. Each component in $\Mb_D$ either has at least two bonds connecting to $\Mb_U$, or has one bond connecting to $\Mb_U$ and one atom with fixed end.
\item There exist at least $(\rho'''/4)$ many UD connections, such that each UD connection is either connected to another UD connection via $\Mb_D$, or has its $\Mb_D$ endpoint belonging to component of $\Mb_D$ that contains an atom with fixed end.
\end{enumerate}
\end{reduct}
\begin{remark}\label{rem.BM}
    Since $\Xf + \Bf \leq (C^*_{10})^{-1} \rho \leq (C^*_5)^{-1}\rho'''$. For simplicity, we call every box constituent or C-constituent of a BO- or CO-atom a BM-constituent, and every B-atom or multi-particle C-atom a \textbf{BM-atom}. By a slight abuse of terminology, we refer to all these objects as BM-atoms, even when they are constituents rather than atoms. Denote $\#_{\mathrm{BM}}=\#_{\mathrm{BM}}(\Mb)$ by the number of these constituents and atoms, then we have $\#_{\mathrm{BM}} \leq (C^*_{10})^{-1} \rho$.
\end{remark}
\subsection{The \textbf{2CONNUP} algorithm}\label{sec.2connup}
We start by introducing the \textbf{2CONNUP} algorithm, in which we view each box as a whole. 
\begin{definition}[The algorithm \textbf{2CONNUP}]
\label{def.alg_2connup}
Define the following cutting algorithm \textbf{2CONNUP}. It takes as input any canonical UD molecule $\Mb_{UD}$, such that $\Mb_D$ has no full component that has no bond connecting to $\Mb_U$. We view each box as a whole. Define $\Mb_{\mathrm{2conn}}^1$ to be the set of atoms $\nf\in\Mb_D$ that are $\{3_2\}$ atoms and ov-adjacent to one atom $\mf\in\Mb_U$, such that no other $\Mb_D$ atom on the ov-segment between $\nf$ and $\mf$ has fixed end. Define $\Mb_{\mathrm{2conn}}^2$ to be the set of binary atoms $\nf\in\Mb_D$ that are ov-adjacent to two atoms in $\Mb_U$ along two different particle lines, and $\Mb_{\mathrm{2conn}}=\Mb_{\mathrm{2conn}}^1\cup\Mb_{\mathrm{2conn}}^2$. For any such $\Mb=\Mb_{UD}$, we define the cutting sequence as follows:
\begin{enumerate}
\item\label{it.alg_2connup_1} If any (maximal) ov-segment contains at least $A:=32\upsilon^{-2}$ atoms in $\Mb_{\mathrm{2conn}}^2$, then we choose $\nf_0$ to be any one of these atoms, and cut it as free. Repeat until no such ov-segment exists.
\item\label{it.alg_2connup_2} If $\Mb$ contains any fully fixed atom $\nf$, then cut it as free, and repeat until there is no such atom left.
\item\label{it.alg_2connup_3} Choose a lowest atom with fixed end $\nf\in\Mb_U$ among all remaining such atoms in $\Mb_U$ (or a lowest atom in $\Mb_U$ if $\Mb_U$ has no atom with fixed end). Let $S_\nf$ be the set of descendants of $\nf$ in $\Mb_U$; this $\nf$ and $S_\nf$ will be fixed until the end of (\ref{it.alg_2connup_4}).
\item\label{it.alg_2connup_4} Now, starting from $\nf$, each time choose a highest atom $\mf$ in $S_\nf$. If $\mf$ is a $\{3_2\}$ atom, and $\mf$ is ov-adjacent to an atom $\pf\in \Mb_{\mathrm{2conn}}$ that also is a $\{3_2\}$ atom, then cut $\{\mf,\pf\}$ as free; if is a CO-atom or BO-atom, we partially cut $\mf$ by {\bf PM}. Otherwise, we cut $\mf$ as free. Repeat until all atoms in $S_\nf$ has been cut. Then go to (\ref{it.alg_2connup_2}).
\item\label{it.alg_2connup_5} Repeat (\ref{it.alg_2connup_2})--(\ref{it.alg_2connup_4}) until all atoms in $\Mb_U$ has been cut. Then apply algorithm \textbf{DOWN} to cut the rest of $\Mb_D$ into elementary molecules. After every atom have been cut, we apply {\bf UP} algorithm to each B- and BO-atom and cut it into elementary molecules.
\end{enumerate}
We also define the dual algorithm \textbf{2CONNDN}, by reversing the roles of $\Mb_U$ and $\Mb_D$ and the notions of parent/child etc., just as how we define \textbf{DOWN} from \textbf{UP} in Definition \ref{def.alg_up}.
\end{definition}
 
\begin{proposition}
\label{prop.alg_2connup} For the algorithm \textbf{2CONNUP}, we have $\#_{\{3_23_2B\}}=\#_{\{4_2 (2l)_l\}}=0$, and
\begin{equation}\label{eq.alg_2connup_1}
(\upsilon/2)\cdot\#_{\{3_23_2A\}}-d\cdot \#_{\mathrm{bad}}\geq C^{-1}(|\Mb_{\mathrm{2conn}}^1|+|\Mb_{\mathrm{2conn}}^2|)-C\cdot \big(\#_{\mathrm{fucp}(\Mb_U)} + \mathsf{N} \#_{\mathrm{BM}}\big),
\end{equation} where $\#_{\mathrm{fucp}(\Mb_U)}$ is the number of \emph{full components} of $\Mb_U$ (note the difference with $\#_{\mathrm{comp}(\Mb_U)}$, i.e. the number of \emph{all components}). Recall also the notion of constant $C$ in Definition \ref{def.notation} (\ref{it.defC*}). Similar results hold for the dual algorithm \textbf{2CONNDN}, if we reverse the roles of $\Mb_U$ and $\Mb_D$.
\end{proposition}
\begin{proof} The proof is essentially the same as Proposition 13.2 in \cite{DHM24}. For an atom $\mf$ in Definition \ref{def.alg_2connup} \eqref{it.alg_2connup_4}, suppose it is a B-, BO- or CO-atom or a multi-particle C-atom, and it is also ov-adjacent to a $\{3_2\}$ atom in $\Mb_{\mathrm{2conn}}$. Then $\mf$ is (partially) cut as free by itself, and hence cannot gain a $\{3_2 3_2 A\}$ atom. On the other hand, with at most one exception (i.e., the case that, when we consider $\mf$, each C- or B-constituent has no fixed end in Definition \ref{def.alg_pm}), every edge constituent of a BO- or CO-atom becomes a special $\{3_2\}$ atom after the partial cutting. Therefore, the number of such ov-segments containing an atom in $\Mb_{\mathrm{2conn}}$, together with either the exceptional edge constituent described above (if any) or a C- or B-constituent is at most $(\mathsf{N}+1)  \#_{\mathrm{BM}}$. Hence (13.3) in \cite{DHM24} should be 
    \begin{equation}\label{eq.alg_2connup_3}
        \#_{\mathrm{good}}   \geq \frac{1}{2}\max\bigg(|\Pc_1|-|\Lc_1| - (\mathsf{N}+1)  \#_{\mathrm{BM}},\,|\Mb_{\mathrm{2conn}}^1| - (\mathsf{N}+1)  \#_{\mathrm{BM}},\,A^{-1}|\Pc_2|-  (\mathsf{N}+1)  \#_{\mathrm{BM}}\bigg)-\#_{\mathrm{exc}}. 
    \end{equation}
The remainder of the proof proceeds similarly, yielding the desired result.
\end{proof}
Now we can prove the first case of Proposition \ref{prop.case5}, namely Proposition \ref{prop.comb_est_case6}. This corresponds to the case that there exist (many components of $\Mb_D$ or) many 2-connections in $\Mb_U$.
\begin{proposition}
\label{prop.comb_est_case6} Proposition \ref{prop.case5} is true if either $\#_{\mathrm{comp}(\Mb_D)}\geq (C_{3}^*)^{-1}\cdot\rho'''$, or if $\#_{\mathrm{comp}(\Mb_D)}\leq (C_{3}^*)^{-1}\cdot\rho'''$ and $|\overline{\Mb}_{\mathrm{2conn}}^2|\geq (C_{2}^*)^{-1}\cdot\rho'''$, where $\overline{\Mb}_{\mathrm{2conn}}^2$ is defined as in Definition \ref{def.alg_2connup}, but with the roles of $\Mb_U$ and $\Mb_D$ reversed.
\end{proposition}
\begin{proof} Our goal is to prove that
\begin{equation}\label{eq.case51_copy} (\upsilon/2)\cdot(\#_{\mathrm{good}})-d\cdot(\#_{\mathrm{bad}})\geq (C_{4}^*)^{-1}\cdot \rho'''.
\end{equation}

\textbf{Proof part 1.} Suppose $\#_{\mathrm{comp}(\Mb_D)}\geq (C_{3}^*)^{-1}\cdot\rho'''$. Recall the definition of $\Mb_{\mathrm{2conn}}^j$ for $j\in\{1,2\}$ in Definition \ref{def.alg_2connup}. Consider all components of $\Mb_D$, if at least half of them contain an atom in $\Mb_{\mathrm{2conn}}:=\Mb_{\mathrm{2conn}}^1\cup \Mb_{\mathrm{2conn}}^2$, then we have $|\Mb_{\mathrm{2conn}}|\geq \#_{\mathrm{comp}(\Mb_D)}/2$. In this case we apply \textbf{2CONNUP} to $\Mb_{UD}$, by Proposition \ref{prop.alg_2connup} we have
\begin{equation}\label{eq.comb_est_case61}
(\upsilon/2)\cdot(\#_{\mathrm{good}})-d\cdot(\#_{\mathrm{bad}})\geq C^{-1}|\Mb_{\mathrm{2conn}}|-C\cdot \big(\#_{\mathrm{fucp}(\Mb_U)} + \#_{\mathrm{BM}} \big)\geq\big(C^{-1}(C_{3}^*)^{-1}-C(C_{5}^*)^{-1}\big)\cdot\rho'''
\end{equation} (where the upper bound of $\#_{\mathrm{fucp}(\Mb_U)}$ follows from Reduction \ref{red.2} (ii)), which implies (\ref{eq.case51_copy}).

If at least half of the components of $\Mb_D$ do not contain an atom in $\Mb_{\mathrm{2conn}}$, then we shall cut $\Mb_U$ as free and cut it into elementary molecules using \textbf{UP}. Then at least $(4 C_{3}^*)^{-1}\cdot\rho'''$ components in $\Mb_D$ does not contain any BM-atom or BM-constituents. Then, for each of these components $X\subseteq \Mb_D$, it will not contain any BM-atom, or $\{2_2\}$ atom (because this component does not contain any atom in either $\Mb_{\mathrm{2conn}}^1$ or $\Mb_{\mathrm{2conn}}^2$); moreover, this $X$ may be divided into new components due to the broken ov-segments, but each new component will contain at least one $\{3_2\}$ atom (cf. Remark \ref{rem.full_cut}), and at least one new component will contain at least two $\{3_2\}$ atoms. This last fact is because at least two atoms in $X$ will become $\{3_2\}$ atoms after cutting $\Mb_U$ by Reduction \ref{red.2} (ii); by considering a path in $X$ between these two atoms before cutting and whether or where this path breaks after cutting, we can obtain a component of $X$ after cutting that contains two $\{3_2\}$ atoms (as the atom where the path breaks must have a fixed end). 

We then cut each resulting component of $\Mb_D$ into elementary molecules using \textbf{DOWN}. By Proposition \ref{prop.alg_up}, and note that all $\{3_2 3_2 A\}$ molecules are treated as good, we get that $\#_{\mathrm{good}}\geq \#_{\mathrm{comp}(\Mb_D)}/2$ and $\#_{\mathrm{bad}}\leq\#_{\mathrm{fucp}(\Mb_U)}$ (this is because $\Mb_D$ has no full component after cutting $\Mb_U$, due to Remark \ref{rem.full_cut} and Reduction \ref{red.2} (ii)), so (\ref{eq.case51_copy}) follows in the same way as (\ref{eq.comb_est_case61}).

\textbf{Proof part 2.} Assume $\#_{\mathrm{comp}(\Mb_D)}\leq (C_{3}^*)^{-1}\cdot\rho'''$ and $|\overline{\Mb}_{\mathrm{2conn}}^2|\geq (C_{2}^*)^{-1}\cdot\rho'''$. We first cut as free each full component in $\Mb_U$ that has no bond connecting to $\Mb_D$ (and then cut it into elementary molecules using \textbf{UP}), to make $\Mb_{UD}$ satisfy the assumptions of Definition \ref{def.alg_2connup} after reversing the roles of $\Mb_U$ and $\Mb_D$. Then we apply \textbf{2CONNDN} to $\Mb_{UD}$. Proposition \ref{prop.alg_2connup} then gives that (note $\#_{\mathrm{fucp}(\Mb_D)}\leq\#_{\mathrm{comp}(\Mb_D)}$)
\begin{multline}\label{eq.comb_est_case62}
(\upsilon/2)\cdot(\#_{\mathrm{good}})-d\cdot(\#_{\mathrm{bad}})\geq C^{-1}\cdot|\overline{\Mb}_{\mathrm{2conn}}^2|-C\cdot(\#_{\mathrm{comp}(\Mb_D)}+\#_{\mathrm{fucp}(\Mb_U)} + \#_{\mathrm{BM}})\\
\geq\big(C^{-1}(C_{2}^*)^{-1}-C(C_{3}^*)^{-1}\big)\cdot\rho''',
\end{multline} which implies (\ref{eq.case51_copy}). Here, we used the bound on $\#_{\mathrm{fucp}(\Mb_U)}$ provided by Reduction \ref{red.2} (ii).
\end{proof}
\begin{reduct}
\label{red.3} With Proposition \ref{prop.comb_est_case6}, in the remaining proof of Proposition \ref{prop.case5}, we may now assume for $\Mb_{UD}$, in addition to Reduction \ref{red.2} (i)--(iii):
\begin{enumerate}
\item That $\#_{\mathrm{comp}(\Mb_D)}\leq (C_{3}^*)^{-1}\cdot\rho'''$, and the number of weakly degenerate atoms and pairs in $\Mb_D$ is less than $(C_{2}^*)^{-1}\cdot\rho'''$. 
\item We have $|\overline{\Mb}_{\mathrm{2conn}}^2|\leq (C_{2}^*)^{-1}\cdot\rho'''$, with the set $\overline{\Mb}_{\mathrm{2conn}}^2$ defined in Proposition \ref{prop.comb_est_case6}.
\end{enumerate}
\end{reduct}
\subsection{The \textbf{3COMPUP} algorithm}\label{sec.3compup}
In this subsection, we introduce the \textbf{3COMPUP} algorithm, which is analog to 3COMPUP algorithm in subsection 13.2 in \cite{DHM24}, but including boxes and multi-particle atoms. Throughout this subsection, we always regard each box as a whole. We first introduce the following lemma from Lemma 13.4 in \cite{DHM24}.
\begin{lemma}\label{lem.cutconnected} Let $\Mb$ be a molecule that is a tree, by viewing each box as a whole, and $S\subseteq\Mb$ is an \emph{ov-connected} set. Suppose we cut $S$ as free from $\Mb$. Then, for each atom $\qf\in\Mb\backslash S$ that is \emph{ov-adjacent} to an atom in $S$ before cutting, there is a unique component $X_\qf\subseteq \Mb\backslash S$ after cutting that contains $\qf$. This $X_\qf$ is in \emph{bijection} with $\qf$, and when viewed in $X_\qf$, the operation of cutting $S$ creates \emph{only one} fixed end in $X_\qf$ at $\qf$. 
\end{lemma}
\begin{proof}
    This is Lemma 13.4 in \cite{DHM24}; the proof only uses that $\Mb$ is a tree and that cutting an ov-connected set breaks each maximal ov-segment at the atoms of $S$, which also holds here by viewing each box as one atom. If $S$ is (partially) cut at a CO- or MO-atom, the O-connections deleted in the partial cutting only create simple pairs (an O-fixed end serial with an O-free end) at constituents outside $S$, and the ov-segment through such a simple pair is not broken; hence partial cutting does not create a second fixed end in $X_\qf$ which would break the conclusion.
\end{proof}
\begin{lemma}\label{lem.3comp_aux} Suppose $\Mb$ is a regular molecule that is a tree, by viewing each box as a whole, contains no C-top fixed end, and fully fixed atoms. Let the set of $\{3_2\}$ atoms in $\Mb$ be $\Mb_{(3)}$ and assume the number of $\{3_2\}$ atoms in $\Mb$ is at least $2$. Define a set $A\subseteq\Mb$ to be \textbf{no-bottom} if any child of any atom in $A$ must also be in $A$. Then, there always exists a connected, no-bottom subset $A\subseteq \Mb$ with either of two cases:
    \begin{enumerate}
    \item\label{it.3comp_aux1} there exists exactly two $\{3_2\}$ atoms in $A$;
    \item\label{it.3comp_aux2} there exists more than two $\{3_2\}$ atoms in $A$, but $A=S_\pf$ for some BM-atom $\pf$.
\end{enumerate}
\end{lemma}
\begin{proof} For any $\pf\in\Mb$, let $S_\pf$ be the set of descendants of $\pf$, then it is always connected and no-bottom. If there exists $\pf$ such that $|S_\pf \cap \Mb_{(3)}| \geq 2$, then choose a lowest such $\pf$. Since $\pf$ has no C-top fixed end, we know that one of the following holds:
     \begin{enumerate}
    \item $\pf$ is a $\{3_2\}$ atom and has one descendant $\pf_1$ which is a $\{3_2\}$ atom; 
    \item $\pf$ is a $\{4_2\}$ atom and has two descendants $\pf_1$, $\pf_2$ which are $\{3_2\}$ atoms; or
    \item $\pf$ is a BM-atom, and has descendants $\pf_1,...\pf_j$ which are $\{3_2\}$ atoms, and $|S_{\pf_j} \cap \Mb_{(3)}| \leq 1$ for each $j$.
\end{enumerate}
In the first two cases, we have $|S_{\pf} \cap \Mb_{(3)}| = 2$; and in the last case, we have $|S_{\pf} \cap \Mb_{(3)}| \geq 2$.

Now assume $|S_{\pf} \cap \Mb_{(3)}| \leq 1$ for each $\pf \in \Mb$. If $|S_{\pf} \cap \Mb_{(3)}| = 1$, we denote the unique $\{3_2\}$ by $\qf = \qf (\pf)$. Now consider all atom pairs $(\pf,\pf')$ such that $\qf(\pf)$ and $\qf(\pf')$ both exist and \emph{are different}, and one such pair $(\pf,\pf')$ that are connected by \emph{a shortest path} in $\Mb$. Let this path be $(\pf_0,\pf,\cdots,\pf_j=\pf')$, then for each $0<i<j$ we must have either $S_{\pf_i}\cap \Mb_{(3_2)}=\varnothing$ or $\qf(\pf_i)=\qf(\pf)$ (otherwise the pair $(\pf,\pf_i)$ corresponds to a shorter path, contradiction). Now let $A=\cup_{i=0}^j S_{\pf_i}$, it is clear that $A$ is connected, no-bottom and $A\cap \Mb_{(3_2)}=\{\qf(\pf),\qf(\pf')\}$, as desired.
\end{proof}
\begin{definition}[The \textbf{3COMPUP} algorithm]\label{def.3comp_alg} Let $\Mb$ be a molecule which is a forest, by viewing each box as a whole, contains no fully fixed atoms. Assume that each component of $\Mb$ contains at least two $\{3_2\}$ atoms. Let the set of $\{3_2\}$ atoms in $\Mb$ be $\Mb_{(3)}$. Define the following cutting sequence:
\begin{enumerate}[{(1)}] 
\item \label{it.3comp_1} Apply Lemma \ref{lem.3comp_aux} to one component of $\Mb$ and select a connected no-bottom subset $A$ with at least two $\{3_2\}$ atoms. This set is fixed until the end of \eqref{it.3comp_5}.
\item \label{it.3comp_3} If $A$ is as in Lemma \ref{lem.3comp_aux} \eqref{it.3comp_aux1}, then choose a lowest atom $\nf$ with at least one fixed end in $A$ that has not been cut. Let $S_\nf$ be the set of descendants of $\nf$ (note that $S_\nf\subseteq A$). 

If $A$ is as in Lemma \ref{lem.3comp_aux} \eqref{it.3comp_aux2}, then choose a lowest atom $\nf$ with at least one fixed end \emph{that is a descendant of $\pf$ along the C-bond of $\pf$ in Lemma \ref{lem.3comp_aux} \eqref{it.3comp_aux2}} in $A$ that has not been cut. If no such atom exists, then we choose $\pf$ as $\nf$. Let $S_\nf$ be the set of descendants of $\nf$, but $S_\pf$ be the set of descendants of $\pf$ \emph{along the C-bonds at $\pf$} for the case $\nf=\pf$. This $\nf$ and $S_\nf$ are fixed until the end of \eqref{it.3comp_4}.
\item \label{it.3comp_4} Starting from $\nf$, each time choose a highest atom $\mf$ in $S_\nf$ that has not been cut. If $\mf$ is a $\{3_2\}$ and has an ov-parent $\mf^+$ or ov-child $\mf^-$ that is a $\{3_2\}$ atom, then cut $\{\mf,\mf^\pm\}$ as free. If $\mf$ is a CO- or BO-atom, we apply \textbf{PM} to it. Otherwise, we just cut it as free. Repeat until all atoms in $S_\nf$ have been cut. 
\item \label{it.3comp_4+} If $\Mb$ has any fully fixed atom, then cut it as free; repeat until $\Mb$ has no such atoms. Then go to \eqref{it.3comp_3}.
\item \label{it.3comp_5} When all atoms in $A$ have been cut, we subsequently cut all remaining fully fixed atoms in $\Mb$ as in \eqref{it.3comp_4+}, and then cut all components of $\Mb$ with only one $\{3_2\}$ atom using \textbf{UP1}.
 \item\label{it.t3comp_6}Go to \eqref{it.3comp_1} and select the next $A$ from one of the remaining components of $\Mb$, and so on.
\end{enumerate}
\end{definition}
 
\begin{proposition}\label{prop.3comp} Let $\Mb$ be as in Definition \ref{def.3comp_alg}. Define $\#_{\mathrm{3ovcp}(\Mb)}$ to be the number of ov-components of $\Mb_{(3)}$, and $\#_{\mathrm{comp}(\Mb)}$ to be the number of components of $\Mb$. Then $\Mb$ is cut into elementary molecules after \textbf{3COMPUP} (Definition \ref{def.3comp_alg}), and we have
\begin{equation}\label{eq.alg_3ovcpup_1}\#_{\{3_23_2B\}}=\#_{\{4_2 (2l)_l\}}=\#_{\mathrm{bad}}=0;\qquad
\#_{\{3_23_2A\}} + 2 \mathsf{N} \cdot \#_{\mathrm{BM}}(\Mb) \geq \frac{1}{10}\cdot\#_{\mathrm{3ovcp}(\Mb)}+\frac{1}{5}\cdot\#_{\mathrm{comp}(\Mb)}.
\end{equation}
\end{proposition}
\begin{proof} We first explain the choice of $\nf$ and $S_\nf$ in Definition \ref{def.3comp_alg} \eqref{it.3comp_3} when $A$ is as in Lemma \ref{lem.3comp_aux} \eqref{it.3comp_aux2}. The highest atom $\pf \in A$ may have lots of O-bonds, and along each O-bond there may be a \{$3_2$\} atom which is a descendant of $\pf$, and each \{$3_2$\} atom may form a $\Mb_{(3)}$ component. If we naively choose the lowest atom in $A$, we may cut all those $\{3_2\}$ atoms by itself, and we cannot gain from them. This is undesirable, as we expect each $\Mb_{(3)}$ component to contribute a good molecule. To resolve this, when each descendant along the C-bond of $\pf$ has no fixed end, we choose $\pf$ as $\nf$ and (partially) cut it.
    
The (MONO) property is preserved since $A$ is no-bottom. If $\nf$ is chosen as $\pf$ in Lemma \ref{lem.3comp_aux} \eqref{it.3comp_aux2}, (partially) cut it does not create any C-top fixed ends in $A \backslash S_\nf$. In fact, for those newly created components in $A \backslash S_\nf$, it originally has an O-connection with $\nf$, and partially cut $\nf$ only create O-fixed ends at those components. Moreover, $\#_{\{3_23_2B\}}=\#_{\{4_2 (2l)_l\}}=\#_{\mathrm{bad}}=0$ follows by the same proof as in Proposition \ref{prop.alg_up}. We now focus on the proof of (\ref{eq.alg_3ovcpup_1}). Note that the algorithm \textbf{3COMPUP} consists of finitely many loops where one loop is formed by steps \eqref{it.3comp_1}--\eqref{it.3comp_5} in Definition \ref{def.3comp_alg}, and (say) turns $\Mb$ into $\Mb'$. Define $\lambda(\Mb)$ to be 
\begin{equation}
\Lambda(\Mb) =  \frac{1}{10}\cdot\#_{\mathrm{3ovcp}(\Mb)}+\frac{1}{5 }\cdot\#_{\mathrm{comp}(\Mb)} - 2\mathsf{N}\cdot \#_{\mathrm{BM}}(\Mb),
\end{equation}
then we only need to consider each single loop \eqref{it.3comp_1}--\eqref{it.3comp_5} in Definition \ref{def.3comp_alg}, and prove that:
\begin{itemize}
\item[($\clubsuit$)] If $\Mb$ is a forest, has no top fixed end and contains no fully fixed atom, with each component containing at least two $\{3_2\}$ atoms, then the result $\Mb'$ after executing the loop \eqref{it.3comp_1}--\eqref{it.3comp_5} in Definition \ref{def.3comp_alg} satisfies the same property. Moreover, for this loop we have
\begin{equation}\label{eq.3comp_2}\#_{\{3_23_2\}} \geq\Lambda(\Mb)-\Lambda(\Mb').
\end{equation}
\end{itemize} Note that adding up (\ref{eq.3comp_2}) leads to the desired lower bound of $\#_{\{3_23_2\}}$, thus it suffices to prove ($\clubsuit$) for a given loop. We divide the proof into 4 parts.

\textbf{Proof part 1.} We prove that each component of $\Mb'$ satisfies the desired properties. This is obvious, as the lack of top fixed ends follows from the same proof as Proposition \ref{prop.alg_up}, and the lack of fully fixed atoms and components with only one $\{3_2\}$ follows from the cuttings done in Definition \ref{def.3comp_alg} \eqref{it.3comp_5}.

\textbf{Proof part 2.} Now we need to prove (\ref{eq.3comp_2}). We start by making some preparations. Note that $A$ is an ov-connected subset of $\Mb$; by Lemma \ref{lem.cutconnected}, each atom $\mf\in \Mb\backslash A$ that is ov-adjacent to an atom $\nf\in A$ by a pre-loop ov-segment corresponds to a unique component $X_\mf$ after cutting the atoms in $A$ as free (the loop will in fact cut as free more atoms). Note also that $\mf$ must be an ov-parent of $\nf$ as $A$ has no bottom.

Let $\Mb_{(3)}$ be the set of $\{3_2\}$ atoms in $\Mb$ (before executing the loop), and consider all ov-components of $\Mb_{(3)}$ that intersect $A$ or are ov-adjacent to $A$ by an ov-segment (we call them \textbf{good}). We know $A$ has at least two atoms in $\Mb_{(3)}$; apart from the ov-components that contain these atoms, for any atom in any other good ov-component, it must belong to $X_\mf$ for a unique $\mf$, and this $\mf$ must be a $\{3_2\}$ atom, with the corresponding $\nf$ in $A$ being either a BM-atom or a $\{4_2\}$ atom. Let $B$ be the union of all good ov-components and $A$. Note that each atom in $B$ has been cut after steps \eqref{it.3comp_1}--\eqref{it.3comp_5} as they become fully fixed after cutting $A$. 



\textbf{Proof part 3.} Next we obtain a lower bound on $\#_{\mathrm{good}}$. We first consider the case that $A$ is as the first case in Lemma \ref{lem.3comp_aux}. We claim that, with at most $ \#_{\mathrm{BM}}(A) +1$ exceptions, each binary atom in $A$ has exactly one fixed end when it is cut. In this case, $A$ contains \emph{exactly} two $\{3_2\}$ atoms. Suppose that a binary atom $\qf \in A$ is cut as a $\{2_2\}$ atom. Then, initially, either \(\qf\) is a \(\{3_2\}\) atom and has at least one descendant with a fixed end, or it has two distinct descendants with fixed ends along two distinct particle lines. In fact, if $\qf$ is a \{$3_2$\} atom initially, then it must have a descendant with fixed ends, otherwise we will choose $\qf$ as $\nf$ in Definition \ref{def.3comp_alg} \eqref{it.3comp_3}, and it will be cut as a $\{3_2\}$ atom; contradiction. If it is a \{$4_2$\} atom initially, it has no top fixed end when it is cut; otherwise it belongs to $S_\nf \backslash \nf$ for some $\nf$ chosen in Definition \ref{def.3comp_alg} \eqref{it.3comp_3}, but $\qf$ has at most one parent in $S_\nf$ as $A$ is a tree, it has at most one fixed end when it is cut. Then $\qf$ must be cut in Definition \ref{def.3comp_alg} \eqref{it.3comp_4+}. Consider the moment that $\qf$ becomes a $\{3_2\}$ atom (since $A$ is a tree, cutting each atom produces at most one fixed end at each atom). It must have a descendant with fixed end, otherwise we will choose $\qf$ as the lowest atom with fixed end; hence the result. Moreover, the number of such atoms is no more than $ \#_{\mathrm{BM}}(A)+1$ because $A$ is a tree and there are at most $ \#_{\mathrm{BM}}(A)+2$ atoms with fixed ends. 

Now let $\lambda$ be the number of good ov-components of $\Mb_{(3)}$ (i.e. those forming $B$). With at most $2$ exceptions (i.e., those in $A$), each such ov-component must contain some $\mf$ which is a $\{3_2\}$ atom and is an ov-parent of an atom $\nf\in A$; consider the maximal ov-segment $\sigma$ containing $\mf$ and $\nf$. Note that each $\sigma$ corresponds to at most one good ov-component (if $\sigma$ is given, then any $\{3_2\}$ atoms on $\sigma$ must be ov-adjacent to each other and thus belong to the same ov-component), so the number of such ov-segments $\sigma$ is at least $\lambda - 2$. Moreover, we (partially) cut each BM-atom in $A$ breaks at most $\mathsf{N}$ such ov-segments (as we partially cut any O-bonds at the BM-atom only when it is a simple pair, and cutting a simple pair does not break any such ov-segments). We claim, with at most $(\mathsf{N} + 2) \cdot  \#_{\mathrm{BM}}(A) + 2 $ exceptions, the first atom $\qf \in \sigma$ that is cut, must belong to a $\{3_2 3_2\}$ atom. Indeed, if $\qf$ is cut as a BM-atom, then the corresponding ov-segment must have a C-bond at $\qf$, and there are at most $  \mathsf{N}\cdot  \#_{\mathrm{BM}}(A) $ such ov-segments. If $\qf$ is cut as $\{2_2\}$ atom, and it is originally a binary atom in $A$, then there are at most $2  \#_{\mathrm{BM}}(A)+ 2$ such ov-segments. Note that $\qf$ cannot be cut as a special $\{2_2\}$ atom; otherwise, $\sigma$ either initially is a simple pair at some BM-atom, or $\qf$ is not the first atom cut in $\sigma$. With these exceptions, each $\qf$ must be a $\{3_2\}$ atom when it is cut. At this time $\mf$ (which is the ov-parent of $\qf$ along $\sigma$) also is a $\{3_2\}$ atom, so by Definition \ref{def.3comp_alg} \eqref{it.3comp_3}, we know that $\qf$ must belong to a $\{3_2 3_2\}$ molecule (which may or may not be $\{\qf,\mf\}$). Since each $\{3_2 3_2 \}$ molecule can be so obtained from at most three $\sigma$ (if we cut a \{$3_2 3_2$\} molecule with both atoms in $A$, and they are initially $4_2$ atoms in $A$, then such \{$3_2 3_2$\} may have three different edges connecting three different $\sigma$), we thus obtain that 
\begin{equation*}
    \#_{\{3_2 3_2\}} \geq \frac{\lambda - 2 - (\mathsf{N} +2 )\cdot  \#_{\mathrm{BM}}(A) - 2 }{3}.
\end{equation*}
Moreover, if $ \#_{\mathrm{BM}}(A) = 0$, we can show that $\#_{\mathrm{good}} \geq 1$ in the same way as Proposition \ref{prop.alg_up} \eqref{it.up_proof_3} as $A$ has two $\{3_2\}$ atoms. Hence,
\begin{equation}\label{eq.3compc}
    \#_{\{3_2 3_2\}} + \mathsf{N}\cdot  \#_{\mathrm{BM}}(A)  \geq \max \Big(\frac{\lambda  - 4 }{3}, 1\Big).
\end{equation}


Now assume $A$ is as the second case in Lemma \ref{lem.3comp_aux}. It is clear that $\Mb$ is formed by joining O-connections among several components at $\pf$. We denote by $\Mb_0$ to be the components that contains atom constituent at $\pf$ (so it has C-bonds at $\pf$). Denote $\Mb_j$ to be the  components with only edge constituent at $\pf$ ($\Mb_j$ does not contain $\pf$). Denote by $A^j_0$ to be the components in $(A \backslash \pf) \cap \Mb_0$, $A_0  = \cup_{j} A^j_0$, and $\Mb_0^j$ to be the components containing $A^j_0$, and denote by $A_j = \Mb_j \cap A$. Then for each $A^j_0$ and $A_j$, it has at most one $\{3_2\}$ atom. First, for each $A^j_0$ (may cut before or later than $\pf$, if it is cut later than $\pf$, then it has exactly one fixed end), by the same argument as above, we have 
\begin{equation*}
     \#_{3_2 3_2}(\Mb^j_0) \geq  \frac{\lambda(\Mb^j_0)  - (\mathsf{N}+2)\cdot   \#_{\mathrm{BM}}(A^j_0) - 4 }{3},
\end{equation*} 
which implies 
\begin{equation} 
     \#_{3_2 3_2}(\Mb^j_0)  +  \mathsf{N}  \#_{\mathrm{BM}}(A^j_0) \geq  \frac{\lambda(\Mb^j_0)  - 4 }{3}  .
\end{equation} 
We take a summation with respect to $j$, and notice that $j\leq \mathsf{N}\cdot  \#_{\mathrm{BM}}(\pf)$ (here $ \#_{\mathrm{BM}}(\pf)$ denotes the number of atom constituents in $\pf$), obtain
\begin{equation}\label{eq.3compa}
     \#_{3_2 3_2}(\Mb_0)  +  \sum_j \mathsf{N} \cdot \#_{\mathrm{BM}}(A^j_0) \geq  \frac{ \sum \lambda(\Mb^j_0) - 4 \mathsf{N} \cdot  \#_{\mathrm{BM}}(\pf)}{3}.
\end{equation} 
Moreover, for each $A_j\ (j\geq 1)$ with $\lambda (\Mb_j) \geq 1$, it has a special $\{3_2\}$ atom after cutting $\pf$. By the same arguments as above, we have 
\begin{equation*}
     \#_{3_2 3_2}(\Mb_j) + \mathsf{N} \cdot  \#_{\mathrm{BM}}(A_j) \geq  \frac{\lambda(\Mb_j)  -  4 }{3},
\end{equation*} 
and $\#_{3_2 3_2}(\Mb_j) \geq 1$ if $ \#_{\mathrm{BM}}(\Mb_j) = 0$ (as $\lambda(\Mb_j) \geq 1,$ alongside the special $\{3_2\}$ atom, there are at least two $\{3_2\}$ atoms when we cut $B \cap \Mb_j$). Hence, 
\begin{equation}\label{eq.3compb}
     \#_{3_2 3_2}(\Mb_j)  +  \mathsf{N}  \#_{\mathrm{BM}}(A_j) \geq \max\Big(\frac{\lambda(\Mb_j) - 4}{3}, 1\Big) \geq  \frac{\lambda(\Mb_j)  + 2 }{9} . 
\end{equation} 
Moreover, for those good components in $\Mb_0 \backslash A$, they are cut as $\{2_2\}$ after (partially) cutting $\pf$, but there are at most $\mathsf{N} \cdot  \#_{\mathrm{BM}}(\pf)$ such good components. Note that $\lambda \leq \sum \lambda(A^j_0) + \sum \lambda(A_j) + \mathsf{N} \cdot  \#_{\mathrm{BM}}(\pf)$ and $ \#_{\mathrm{BM}}(A) =   \sum \#_{\mathrm{BM}}(A^j_0) + \sum \#_{\mathrm{BM}}(A_j) +  \#_{\mathrm{BM}}(\pf)$. As a consequence, by \eqref{eq.3compa}, \eqref{eq.3compb} and the fact that $ \#_{\mathrm{BM}}(\pf) \geq 1$, we have 
\begin{equation}\label{eq.3comp_3}
    \begin{aligned}
         \#_{\{3_2 3_2\}}+ 2 \mathsf{N} \cdot  \#_{\mathrm{BM}}(A)& = \sum_{j\geq 0} \#_{3_2 3_2}(\Mb_j) + 2\mathsf{N} \cdot \big(\sum_{j}  \#_{\mathrm{BM}}(A_j) + \sum  \#_{\mathrm{BM}}(A^j_0)  +  \#_{\mathrm{BM}}(\pf)\big) \\
          &\geq \sum_j \frac{\lambda(\Mb_j) + 2}{9 } + \frac{\sum\lambda(\Mb^j_0) - 4 \mathsf{N}  \#_{\mathrm{BM}}(\pf)}{3} + 2 \mathsf{N}  \#_{\mathrm{BM}}(\pf)\\
         &\geq \frac{\lambda + 2 }{9}.
    \end{aligned}
\end{equation}






\textbf{Proof part 4.} Now consider the components of $\Mb'$ and ov-components of $\Mb_{(3)}'$ (set of $\{3_2\}$ atoms in $\Mb'$). Define $B'$ as the union of all atoms in $B$ and all fully fixed ends (including possible BM-atoms) cut in Definition \eqref{def.3comp_alg} \eqref{it.3comp_4+}. Then $\Mb'$ is formed from $\Mb$ by cutting $B'$ as free and then removing all components with only one $\{3_2\}$ atom. By applying Lemma \ref{lem.cutconnected} to $B'$, we see that each component $X_\qf$ of $\Mb'$ corresponds to a unique atom $\qf\in \Mb\backslash B'$ that is ov-adjacent to an atom in $B'$. Let $\mu$ be the number of such atoms (i.e. corresponding to those new components with at least two $\{3_2\}$ atoms after cutting $B'$), then \begin{equation}\label{eq.3comp_4}\#_{\mathrm{comp}(\Mb)}-\#_{\mathrm{comp}(\Mb')}=1-\mu.\end{equation}

Finally, consider the set of $\{3_2\}$ atoms, i.e. $\Mb_{(3)}$ and $\Mb_{(3)}'$, and their ov-components. Note that for each new component $Y$ formed after cutting $B'$ that contains only one $\{3_2\}$ atom, it will be removed and will not occur in $\Mb'$; however, if $Y$ originally has a $\{3_2\}$ atom before any cutting, and still have only one $\{3_2\}$ atom after cutting $B'$, then $Y$ must have a BM-atom which is connected to $B'$ (otherwise $Y$ will have at least two $\{3_2\}$ atoms). Denote $Y'$ to be the union of all such component $Y$, and ignore those components in $Y \backslash Y'$ as they initially do not contribute to $\#_{\mathrm{3covp(\Mb)}}$. The number of components in $Y'$ is at most $\#_{\mathrm{BM}}(Y')$. 
The ov-components of $\Mb_{(3)}$ includes those $\lambda$ ov-components forming $B$, the $\Mb_{(3)}$ in $B' \backslash B$ and $Y'$, and the other ov-components which are contained in some $X_\qf$ defined above. Moreover, the number of $\Mb_{(3)}$ in $B' \backslash B$ is no more than $\mathsf{N} \cdot  \#_{\mathrm{BM}}(B' \backslash B)$. Indeed, after cutting atoms in $B$, those fully fixed atoms must be BM-atoms (otherwise they will be $\{3_2\}$ that ov-adjacent to $B$ before cutting $B$, and will be contained in $B$ by the definition of $B$). Cutting each such BM-atom turns at most $\mathsf{N}$ ov-components in $\Mb_{(3)}$ into fully fixed. 

By Lemma \ref{lem.cutconnected}, when viewed in each $X_\qf$, the effect of cutting $B'$ exactly turns one free end into fixed end. Therefore the part of $\Mb_{(3)}'$ within $X_\qf$ is formed by \emph{adding at most one $\{3_2\}$ atom} into the existing part of $\Mb_{(3)}$ within $X_\qf$. This new $\{3_2\}$ atom may \emph{decrease} $\#_{3\mathrm{ovcp}(X_\qf)}$ (i.e. the number of ov-components of the set of $\{3_2\}$ atoms in $X_\qf$) by joining existing ov-components together, but this decrease is \emph{at most 2} because this new $\{3_2\}$ atom can join at most 3 ov-components into 1. From this, and accounting for the $\lambda$ ov-components removed when cutting $B' \cup Y'$, we deduce that
\begin{equation}\label{eq.3comp_5}
    \#_{3\mathrm{ovcp}(\Mb)}-\#_{3\mathrm{ovcp}(\Mb')}\leq \lambda + 2 \mu +  \#_{\mathrm{BM}}(B'\backslash B) \cdot \mathsf{N} +  \#_{\mathrm{BM}}(Y')  
\end{equation}
Putting together \eqref{eq.3compc}, (\ref{eq.3comp_3})--(\ref{eq.3comp_5}), and using the definition of $\Lambda(\Mb)$, we deduce that
\begin{equation}\label{eq.3comp_6}
\begin{aligned}\Lambda(\Mb)-\Lambda(\Mb')&=\frac{1}{5 } \big(\#_{\mathrm{comp}(\Mb)}-\#_{\mathrm{comp}(\Mb')}\big)+ \frac{1}{10  }\big(\#_{3\mathrm{ovcp}(\Mb)}-\#_{3\mathrm{ovcp}(\Mb')}\big) - 2\mathsf{N} \cdot\big( \#_{\mathrm{BM}}(\Mb) -  \#_{\mathrm{BM}}(\Mb')\big)
\\&\leq \frac{2 + \lambda }{10 } - 2\mathsf{N}\cdot  \#_{\mathrm{BM}}(A) \leq \#_{\{3_2 3_2\}} .
\end{aligned}\end{equation} This proves (\ref{eq.3comp_2}) and completes the proof of Proposition \ref{prop.3comp}.
\end{proof}
\subsection{Pre-processing}\label{sec.preprocess}
In this subsection we introduce the \textbf{pre-processing} step. We start with the notion of \emph{ov-proper}.
\begin{definition}\label{def.ov_proper} Let $\Mb$ be a molecule, we say it is \textbf{ov-proper}, if it is a forest, contains only $\{3_2\}$ and $\{4_2\}$ atoms (hence no boxes and multi-particle atoms), and the ov-distance (Definition \ref{def.ov_connect}) between any two $\{3_2\}$ atoms is at least $3$. We say $\Mb$ is \textbf{proper} if $\Mb$ is ov-proper and is a C-molecule.
\end{definition}
Note that $\Mb_D$ may have fixed ends as the result of the cuttings in Definition \ref{def.layer_cutting} and might not be ov-proper. The idea of pre-processing is to turn $\Mb_D$ into ov-proper again by cutting as free a suitable set $S$. We start by stating Lemma \ref{lem.3_comp_property} below.
\begin{lemma}
\label{lem.3_comp_property} Let $\Mb_{UD}$ be a UD molecule such that $\Mb_D=\varnothing$, and $\Mb_U$ does not have C-top fixed end. Suppose $A\subseteq \Mb_U$ is ov-connected within $\Mb_U$ (cf. Definition \ref{def.ov_connect}) such that each atom in $A$ is a $\{3_2\}$ atom. Then $A$ contains a \emph{unique} lowest atom $\nf$, such that (i) each atom $\mf\in A$ is obtained from $\nf$ by iteratively taking ov-parents in $A$, and (ii) each ov-child of $\mf\in A$ is either one of the terms preceding $\mf$ in the ov-parent sequence defined in (i), or an O-atom between two consecutive terms preceding $\mf$, or an BM-atom between two consecutive terms preceding $\mf$ with the bond connected to $A$ is an O-bond at the BM-atom, or an ov-child of $\nf$.

Similarly, if we assume $\Mb_D\neq \varnothing$, $A$ is ov-connected within $\Mb_U$, each atom of $A$ is binary C-atom, and each $\{4_2\}$ atom of $A$ has at most one ov-child in $A$, then conclusion (i) is still true.
\end{lemma}
\begin{proof} This statement is intuitive. We refer the reader to Lemma 13.9 in \cite{DHM24} for the detailed argument.
\end{proof}

\begin{definition}\label{def.W}
    Let $\Mb$ be any UD molecule. For any BM-atom $\qf$ in $\Mb_U$, define \(W_\qf\) to be a set of atoms \(\nf\in\Mb_D\) satisfying the following property. The atom \(\nf\) is an ov-child of \(\pf\in\Mb_U\) along an ov-segment \(\sigma\), where \(\sigma\) contains a C-bottom edge at \(\pf\). There exists another ov-child \(\nf'\in\Mb_D\) of \(\pf\), lying on a different ov-segment \(\sigma'\neq\sigma\), such that \(\nf\) and \(\nf'\) are connected by a path in \(\Mb_D\). No internal atom of this path belongs to either \(\sigma\) or \(\sigma'\). Define $W = \cup_{\qf \in \Mb_U} W_\qf$.
\end{definition}
\begin{proposition}\label{prop.w}
    Let $W$ be defined in Definition \ref{def.W}. Then we have $ |W| \leq 2 \mathsf{N} \cdot \#_{\mathrm{BM}}$, where $\#_{\mathrm{BM}}$ denotes the number of BM-atoms in $\Mb_U$.
\end{proposition}
\begin{proof}
     It suffices to prove that each atom constituent of the BM-atom $\qf\in\Mb_U$ contributes at most \(2\mathsf{N}\) atoms to \(W\). Without loss of generality, we may assume $\qf$ has exactly one atom constituent. Denote the set of ov-segments in Definition \ref{def.W} by $\sigma_\pf=\{\sigma_j, \ ( 1 \leq j \leq k \leq \mathsf{N})\}$. Now we claim $|W_\qf| \leq 2k - 2$. We prove it by induction. If $k = 1$ then by definition we have $|W_\qf| = 0$. Now assume it is true for some $ k \leq \mathsf{N} -1$, and now we consider $k+1$. If $\sigma_{k+1}$ is not connected with any ov-segment in $\sigma_\pf \backslash \sigma_{k+1}$, then the result is trivial by induction. Now assume that $\sigma_{k+1}$ is connected with some $\sigma_j$. Before adding the ov-segment $\sigma_{k+1}$ in $\Mb_D$, there exists $i\leq k$ different components in $\Mb_D$ that contain at least one ov-segment in $\sigma_\pf$. After adding $\sigma_{k+1}$, each of these ov-components intersect exactly one atom with $\sigma_{k+1}$ (because $\Mb_D$ is a tree). This contributes at most $2 i$ atoms in $W_\qf$. Let the number of atoms in $W_\qf$ in each component that contain an ov-segment in $\sigma_\pf$ (before adding $\sigma_{k+1}$) is $a_j$ $(1 \leq j\leq i)$, and the number of ov-segments in $\sigma_\pf$ be $b_j$. Therefore, we have by induction
     \begin{equation}
        |W_\qf| \leq \sum_{1\leq j\leq i}  a_j + 2i  \leq \sum_{1 \leq j\leq i} (2 b_j -2) + 2i \leq 2 k =2(k+1) -2.
     \end{equation}
     Hence, we conclude the result by induction.
\end{proof}
\begin{remark}
    Let us explain the choice of $W$. In the algorithm {\bf MAINUD} below, each time we (partially) cut an atom in $\Mb_U$, we expect that there exists at most one fixed end in each component in $\Mb_D$. But due to BM-atoms, if we cut a BM-atom in $\Mb_U$, those atoms in $W$ will have fixed end and they are connected. Hence, in the function \textbf{SELECT} below, we choose those atoms in $W$ and will cut them as free. After that, each time we cut an atom in $\Mb_U$, there generates exactly one fixed end in each component in $\Mb_D$.
\end{remark}

\begin{definition}[The function \textbf{SELECT}]
\label{def.func_select} Let $\Mb$ be any UD molecule, such that $\Mb_D$ has no fully fixed atoms, and has no C-bottom fixed end. Let $X$ be the set of multi-particle C-atoms, C-constituents in CO-atoms, B-atoms or B-constituents in $\Mb_D$,  $Z$ be the set of $\{3_2\}$ atoms in $\Mb_D$, and $Y$ be a subset of the set of $\{4_2\}$ atoms in $\Mb_D$. Moreover, let \(W\) be defined as in Definition \ref{def.W}.

Given $(\Mb_D,W,X,Z,Y)$, we define a set $S:=\textbf{SELECT}(\Mb_D,W,X,Z,Y)$ as follows:
\begin{enumerate}
\item\label{it.func_select_1} Consider all the ov-components of $W\cup X\cup Z\cup Y$ in $\Mb_D$ (this means we only consider the O-atoms in $\Mb_D$ when defining ov-connectedness in Definition \eqref{def.ov_connect}), which is a finite collection of disjoint subsets of $\Mb_D$.
\item\label{it.func_select_2} If any two of these subsets, say $A$ and $B$, have the shortest ov-distance (within $\Mb_D$) which is at most $4$, then we choose the shortest path between an atom in $A$ and an atom in $B$, which can be divided into at most $4$ ov-segments. Write this path as $(\nf_0,\cdots,\nf_q)$, where $q\leq 4$ and $\nf_0\in A$, $\nf_q\in B$, and each $(\nf_j,\nf_{j+1})$ is an ov-segment. Then replace the two sets $A$ and $B$ by a single set $A\cup B\cup\{\nf_1,\cdots,\nf_{q-1}\}$.
\item\label{it.func_select_4} Repeat (\ref{it.func_select_2}) until this can no longer be done (i.e. when only one set is left, or any two sets have of distance at least 5). Then output $S:=\textbf{SELECT}(\Mb_D,W,X,Z,Y)$ as the union of all the current sets.
\end{enumerate}
\end{definition}
\begin{proposition}
\label{prop.func_select} Let the function \textbf{SELECT} be as in Definition \ref{def.func_select}. Then we have the followings:
\begin{enumerate}
\item\label{it.func_select_new1} After (partially) cutting $S:=\textbf{SELECT}(\Mb_D,W,X,Z,Y)$ from $\Mb$ as free, then the set $\Mb_D\backslash S$ in the resulting molecule $\Mb\backslash S$ is ov-proper (Definition \ref{def.ov_proper}). Note that $\Mb_D\backslash S$ may have both C-top and C-bottom fixed ends, however this is compensated by the good property of being ov-proper.

\item\label{it.func_select_new2} Consider all the maximal ov-segments in $\Mb$ that contains at least one atom in $S$ and one atom in $\Mb_U$. Then the number of these maximal ov-segments does not exceed $10(\#_{\mathrm{ovcp}(Z)}+|Y|+|W|) + \mathsf{N} |X|$, where $\#_{\mathrm{ovcp}(Z)}$ is the number of ov-components of $Z$ (i.e. the set of all $\{3_2\}$ atoms in $\Mb_D$), where the ov-components are determined by considering $\Mb_D$ only.
\end{enumerate}
\end{proposition}
\begin{proof} \textbf{Proof of (\ref{it.func_select_new1}).} At the end of the construction in Definition \ref{def.func_select} (\ref{it.func_select_4}), we get a number of disjoint subsets, such that each one is ov-connected and the ov-distance between each two of them is at least $5$, and $S$ is the union of all these subsets. Note also that the sub-molecule $\Mb_D\backslash S$ of the molecule formed by cutting $S$ as free from $\Mb$, is equal to the molecule formed by cutting $S$ as free from $\Mb_D$.

Now suppose after cutting $S$ as free, there exist two $\{3_2\}$ atoms $\nf,\nf'\in\Mb_D\backslash S$ of ov-distance at most $2$; the case of a $\{2_2\}$ atom in $\Mb_D\backslash S$ can be treated similarly. Since $\nf,\nf'\not\in S$, they must be $\{4_2\}$ atoms initially in $\Mb_D$. As they become $\{3_2\}$ atoms after cutting $S$ as free, we know that for each of $\nf$ and $\nf'$, it must be ov-adjacent to an atom in $S$. Assume $\nf$ is ov-adjacent to $\mf\in S$, and $\nf'$ is ov-adjacent to $\mf'\in S$. If $\mf$ and $\mf'$ belong to two different subsets $A$ and $B$ defined above, then the ov-distance between $\mf$ and $\mf'$ (and hence the ov-distance between $A$ and $B$) is at most $1+2+1=4$, contradicting the assumption about these subsets.

Therefore, we know that $\mf$ and $\mf'$ belong to the same set $A$ defined above. But then there exists a path $\mu_1:\mf\to\nf\to\nf'\to\mf'$ (composed of ov-segments), and another path $\mu_2$ between $\mf$ and $\mf'$, which is composed of ov-segments connecting atoms in $S$ (because $\mf$ and $\mf'$ belong to the same ov-component $A$ of $S$). The path $\mu_1$ cannot be expressed in this way (because the atoms $\nf$, $\nf'$ separating the different ov-segments, as well as the possible atom ov-adjacent to both $\nf$ and $\nf'$, all belong to $\Mb\backslash S$), unless this whole path belongs to a single ov-segment connecting $\mf$ and $\mf'$ (meaning the atoms $\nf$, $\nf'$ etc., all belong to this single ov-segment connecting $\mf$ and $\mf'$), but in this latter case the path $\nf\to\nf'$ would no longer exist after cutting $S$ as free (see Definition \ref{def.cutting}). Consequently, we know that the paths $\mu_1$ and $\mu_2$ \emph{cannot coincide}, but then they will form a cycle, which contradicts the fact that $\Mb_D$ is a forest.

\textbf{Proof of (\ref{it.func_select_new2}).}
In each step in Definition \ref{def.func_select} (\ref{it.func_select_2}) we add at most $3$ atoms to the union of the current subsets. Moreover, the number of such steps that can be carried out cannot exceed the initial number of subsets (as each step reduces the number of subsets by at least 1), which is the number of ov-components of $Z\cup Y \cup X$, and is at most $\#_{\mathrm{ovcp}(Z)}+|W|+|X|+ |Y|$. This means that the number of newly added atoms in the construction of $S$ (starting from $W\cup X \cup Y\cup Z$) is at most $3(\#_{\mathrm{ovcp}(Z)}+ |W|+|X|+|Y|)$.

Therefore, we only need to consider the required maximal ov-segments that contains an atom in $Z$. In each such maximal ov-segment, we have an atom $\mf$ belonging to an ov-component $V$ of $Z$, which has an ov-parent in $\Mb_U$. By (the dual version of) Lemma \ref{lem.3_comp_property} (ii), this ov-parent of $\mf$, which does not belong to $\Mb_D$, must also be an ov-parent of the unique highest atom $\nf$ of $V$; therefore this maximal ov-segment must contain $\nf$. As each ov-component $V$ contains a unique highest atom $\nf$, and each atom belongs to at most two maximal ov-segments, we conclude that the number of such maximal ov-segments does not exceed $2\cdot\#_{\mathrm{ovcp}(Z)}$ (as $\#_{\mathrm{ovcp}(Z)}$ is the number of ov-components of $Z$ in $\Mb_D$), as desired. Moreover, cutting atoms in $X$ breaks at most $\mathsf{N} |X|$ the required maximal ov-segments, and cutting an atom not in $Z \cup X$ breaks at most 2 the required maximal ov-segments.

\end{proof}
\subsection{The \textbf{MAINUD} algorithm}\label{sec.mainud}
In this section we introduce the last ingredient, which is also the central algorithm in our proof, namely the \textbf{MAINUD} algorithm, which is a modification of the MAINUD algorithm of \cite{DHM24}, due to the appearance of BM-atoms. In Definition \ref{def.alg_maincr} below, we present the full version of this algorithm. 
 
\begin{definition}[The algorithm \textbf{MAINUD}] \label{def.alg_maincr}
Define the following cutting algorithm \textbf{MAINUD}. It takes as input any UD molecule $\Mb_{UD}$ that satisfies the following assumptions:
\begin{enumerate}[{(i)}]
\item $\Mb_U$ contains no C-top fixed end, and $\overline{\Mb}_{\mathrm{2conn}}^2=\varnothing$, where $\overline{\Mb}_{\mathrm{2conn}}^2$ is defined as in Proposition \ref{prop.comb_est_case6}.
\item $\Mb_D$ is ov-proper (Definition \ref{def.ov_proper}), and does not have any full component that is not connected to $\Mb_U$ by a bond.
\item\label{it.BM-no-twobond} For each BM-atom $\pf\in\Mb_U$, each component of $\Mb_D$ has at most one ov-segment connecting to $\pf$ which contains a C-bottom edge at $\pf$.
\end{enumerate}

 For any such $\Mb$, we define the sets $X_0$ and $X_1$ in $\Mb_U$ as follows:
\begin{itemize}
 \item Define $\nf\in X_0^0$ if $\nf$ is a $\{3_2\}$ C-atom and has one ov-child in $\Mb_D$, or $\nf$ is a BM-atom, and each bottom edge is either fixed end, or belongs to an ov-segment connected to $\Mb_D$.
 \item Define $\nf\in X_0^q$, if $\nf\in X_0^{q-1}$, or $\nf$ is a $\{3_2\}$ C-atom and has one ov-child in $X_0^{q-1}\cup\Mb_D$, or if $\nf$ is a $\{4_2\}$ C-atom and has two ov-children in $X_0^{q-1}\cup\Mb_D$ along two ov-segments in two different particle lines, or if $\nf$ is a BM-atom, and each bottom edge is either a bottom fixed end, or belongs to an ov-segment connected to an atom in $X_0^{q-1}\cup\Mb_D$. Then define $X_0=\cup_{q\geq 0}X_0^q$.
 \item With $X_0$ fixed, define $\nf\in X_1$ if and only if $\nf\not\in X_0$, and $\nf$ either is a $\{3_2\}$ atom or a $\{4_2\}$ atom and has one ov-child in $X_0\cup\Mb_D$, or $\nf$ is a BM-atom (including every constituent), has either a fixed end, or (at least) one ov-child in $X_0\cup\Mb_D$.
 \end{itemize}
 
Note these sets are fixed throughout the cutting sequence below. Define the following cutting sequence:
\begin{enumerate}
\item\label{it.alg_maincr_2} If $\Mb_U$ contains any fully fixed atom $\nf$, then cut it as free, and repeat until there is no such atom left. 
\item\label{it.alg_maincr_3} Consider the set of all (remaining) atoms that either belong to $X_0\cup X_1$, or has at least one fixed end; choose $\nf$ to be a lowest atom in this set (if this set is empty then choose $\nf$ to be any lowest atom). Let $S_\nf$ be the set of descendants of $\nf$ in $\Mb_U$, and fix $\nf$ and $S_\nf$ until the end of (\ref{it.alg_maincr_7}).
\item\label{it.alg_maincr_4} Starting from $\nf$, choose a highest atom $\pf$ in $S_\nf$; if $\pf$ is a binary atom and is ov-adjacent to a $\{3_2\}$ atom $\qf\in\Mb_D$, then cut $\{\pf,\qf\}$ as free, and cut $\qf$ as free from $\{\pf,\qf\}$ if $\pf$ is a $\{4_2\}$ atom. If $\pf$ is a BM-atom, then each time we choose a constituent $\pf'$ of $\pf$ as in \textbf{PM}.  If $\pf'$ is ov-adjacent to some $\{3_2\}$ $\qf_1,...,\qf_j$ atoms in $\Mb_D$, then we (partially) cut $\{\pf',\qf_1,...,\qf_j\}$ as free, then cut each $\qf_1,...,\qf_j$ as free. Otherwise (partially) cut $\pf$ as free.
\item\label{it.alg_maincr_5} If any two $\{3_2\}$ atoms $(\rf,\rf')$ in $\Mb_D$ becomes ov-adjacent, then choose $\{\rf,\rf'\}$ such that this ov-segment contains the fewest bond, and cut $\{\rf,\rf'\}$ as free. Repeat until no such instances exist.
\item\label{it.alg_maincr_6} If the ov-distance of any two $\{3_2\}$ atoms $(\rf,\rf')$ in $\Mb_D$ becomes 2 (within $\Mb_D$), say $\rf$ and $\rf'$ are both ov-adjacent to some $\rf''\in\Mb_D$, then cut $\{\rf,\rf',\rf''\}$ as free and then cut $\rf$ as free from $\{\rf,\rf',\rf''\}$. Go to (\ref{it.alg_maincr_5}). Repeat until $\Mb_D$ becomes ov-proper.
\item\label{it.alg_maincr_7} When $\Mb_D$ becomes ov-proper, go to (\ref{it.alg_maincr_4}). Repeat until all atoms in $S_\nf$ have been cut.
\item\label{it.alg_maincr_8} When all atoms in $S_\nf$ have been cut, go to (\ref{it.alg_maincr_2}). Repeat until all atoms in $\Mb_U$ have been cut.
\item\label{it.alg_maincr_9} When all atoms in $\Mb_U$ have been cut, choose any $\{3_2\}$ atom in $\Mb_D$ and cut it as free. If the ov-distance between any two $\{3_2\}$ atoms becomes at most 2, go to (\ref{it.alg_maincr_5}) and repeat (\ref{it.alg_maincr_5})--(\ref{it.alg_maincr_6}) until $\Mb_D$ becomes proper again. Then choose and cut the next $\{3_2\}$ atom in $\Mb_D$, and so on.
\end{enumerate}
\end{definition}
\begin{proposition}
\label{prop.alg_maincr} Let $\Mb$ be any UD molecule that satisfies (i)--(ii) in Definition \ref{def.alg_maincr}. Consider the number of UD connections that are either (i) connected to some other UD connection via $\Mb_D$, or (ii) connected to a component of $\Mb_D$ that contains $\{3_2\}$ atoms. We call such UD connections \textbf{good} and let the number of them be $\#_{\mathrm{UD}}$.

Let $\#_{\mathrm{fucp}(\Mb_U)}$ be the number of full components of $\Mb_U$. Let $\overline{X}_1$ be all binary atoms in $X_1$ defined in Definition \ref{def.alg_maincr}, and let $\#_{\mathrm{ovcp}(\overline{X}_1)}$ be the number of ov-components of the set $\overline{X}_1$. Then, for the algorithm \textbf{MAINUD}, we have $\#_{\{4_2 (2l)_l\}}=0$, and each \{$3_23_2$B\} molecule must have both its atoms in $\Mb_D$, and that
\begin{equation}\label{eq.alg_maincr}
\begin{aligned}
    \#_{\{3_23_2A\}}+\#_{\{3_23_2B\}} &\geq\frac{1}{5}\cdot(\#_{\mathrm{UD}}-\#_{\mathrm{fucp}(\Mb_U)}-\#_{\mathrm{ovcp}(X_1)} - \mathsf{N}\cdot \#_{\mathrm{BM}}),\\
    \#_{\mathrm{bad}} &\leq \#_{\mathrm{fucp}(\Mb_U)}+\#_{\mathrm{ovcp}(\overline{X}_1)} + \#_{\mathrm{BM}}.
\end{aligned}
\end{equation}
\end{proposition}
\begin{proof} We divide the proof into 4 parts.

\textbf{Proof part 1.} We prove the (MONO) property as in Proposition \ref{prop.alg_up}: at any time in the cutting sequence, there is no C-top fixed end in $(\Mb_U\backslash S_\nf)\cup\{\nf\}$ and no C-bottom fixed end in $S_\nf\backslash\{\nf\}$. In fact, when we choose $\nf$ as in Definition \ref{def.alg_maincr} (\ref{it.alg_maincr_3}), every atom $\pf\in S_\nf\backslash\{\nf\}$ must have no fixed end. In the process of Definition \ref{def.alg_maincr} (\ref{it.alg_maincr_4})--(\ref{it.alg_maincr_7}), if we only consider the cutting within $\Mb_U$, then the desired property follows from the same proof as in Proposition \ref{prop.alg_up}. Now suppose an arbitrary set of atoms in $\Mb_D$ are also cut in the process, then this will only create C-bottom fixed ends, and only at C-atoms $\tf$ that are ov-adjacent to some atom in $\Mb_D$. By definition we must have $\tf\in X_0\cup X_1$, and thus $\tf\not\in S_\nf\backslash\{\nf\}$ by our choice, so it does not affect the desired (MONO) property. Similarly, for any molecule $\{\pf,\qf\}$ cut in (\ref{it.alg_maincr_4}), the $\pf$ must be either be O-atom or $\pf=\nf$. If this molecule is a $\{3_2 3_2\}$ molecule, then it must be a $\{3_2 3_2 A\}$ (as $\pf$ has not C-top fixed end). 

\textbf{Proof part 2.} We prove that each atom in $\Mb_D$ belongs to a $\{3_2\}$ or $\{3_2 3_2\}$ molecule (which can be either \{$3_2 3_2$A\} or \{$3_23_2$B\}). In fact, note first that initially $\Mb_D$ is ov-proper (thus has no $\{2_2\}$ atom), and does not have full components not connected to $\Mb_U$ by a bond (thus will not have any full component after cutting $\Mb_U$ as free), so the cuttings in Definition \ref{def.alg_maincr} (\ref{it.alg_maincr_9}) will not create any $\{4_2\}$ molecule. Next, note that $\Mb_D$ initially is a forest, and will remain a forest in the whole cutting process as cutting an atom will not create any cycle. Now consider each cutting step in Definition \ref{def.alg_maincr} (\ref{it.alg_maincr_4}) involving $\pf\in\Mb_U$, at which time $\Mb_D$ must be proper by the construction in Definition \ref{def.alg_maincr} (\ref{it.alg_maincr_5})--(\ref{it.alg_maincr_7}).

Each $\pf$ (or $\pf'$ in Definition \ref{def.alg_maincr} \eqref{it.alg_maincr_4}) may have ov-segment connecting to $\Mb_D$. If no such ov-segment exists, then cutting $\pf$ does not affect $\Mb_D$. In particular, if $\pf$ is a binary atom, then it has \emph{at most one} such ov-segment, because $\overline{\Mb}_{\mathrm{2conn}}^2=\varnothing$ by Definition \ref{def.alg_maincr} (i). The effect of the cutting step in Definition \ref{def.alg_maincr} (\ref{it.alg_maincr_4}) \emph{on $\Mb_D$}, has the following three possibilities:
\begin{enumerate}[{(a)}]
\item\label{it.stepa1} If $\pf$ has one ov-segment connecting to $\Mb_D$. If we cut $\pf$ as free (which means all atoms in $\Mb_D$ ov-adjacent to $\pf$ are $\{4_2\}$ atoms), then this step breaks $\Mb_D$ into connected components $X_\qf$ which are in bijection with atoms $\qf\in\Mb_D$ ov-adjacent to $\pf$, and this cutting creates only one fixed end at $\qf$ when viewed within $X_\qf$, turning it from a \{$4_2$\} atom to a $\{3_2\}$ atom.
\item\label{it.stepb1} If $\pf$ has one ov-segment connecting to $\Mb_D$. If we cut $\{\pf,\qf\}$ as free (which means $\qf\in\Mb_D$ is ov-adjacent to $\pf$ which is a $\{3_2\}$ atom), then when viewed in $\Mb_D$, this step equivalent to cutting $\qf$ as free in $\Mb_D$. 
\item\label{it.stepc1} If $\pf$ (or $\pf'$) is a BM-atom, and it has two or more ov-segments connecting to $\Mb_D$. By assumption \ref{it.BM-no-twobond} in Definition \ref{def.alg_maincr}, each component in $\Mb_D$ has \emph{at most one ov-segment} connecting to $\pf$.  Therefore, viewed in $\Mb_D$, the cutting involving $\pf$ or $\pf'$ is either equal to step \eqref{it.stepa1} or step \eqref{it.stepb1}. 
\end{enumerate}
 Here \eqref{it.stepa1} follows from the same proof in Lemma \ref{lem.cutconnected} (i.e. $\qf$ in bijection with $X_\qf$, in which only one fixed end is created at $\qf$) as $\Mb_D$ is a forest. In \eqref{it.stepb1}, the effect of cutting $\{\pf,\qf\}$ as free to $\Mb_D$ is equivalent to cutting $\{\qf\}$ as free, because $\pf$ has no other ov-segment connecting to other atoms in $\Mb_D$ (apart from those on the ov-segment containing $\qf$). In \eqref{it.stepc1}, the effect of cutting $\pf$ as free in $\Mb_D$ is turning a free end into fixed end. By the definition of $W$ in Definition \ref{def.W}, after cutting $W$, each component in $\Mb_D$ has at most one atom ov-connected to $\pf$, then view within each component $X_\qf$, the effect of cutting $\pf$ is the same as \eqref{it.stepa1}. If $\pf$ is cut together with some $\{3_2\}$ atom in $\Mb_D$, then in the component containing this $\{3_2\}$, the cutting is equivalent to \eqref{it.stepb1}.

Now, regarding the effects described in \eqref{it.stepa1}--\eqref{it.stepc1}, we have the following lemma.
\begin{lemma}\label{lem.mainudlem} Suppose $\Mb_D$ is proper at some point of the algorithm, and becomes non-proper after \eqref{it.stepa1}, \eqref{it.stepb1} or \eqref{it.stepc1}. Then, $\Mb_D$ will once more become proper after applying Definition \ref{def.alg_maincr} (\ref{it.alg_maincr_5})--(\ref{it.alg_maincr_6}), via cuttings of \{$3_2$\} and \{$3_23_2$\} molecules only, and \emph{without seeing any $\{2_2\}$ atom}.
\end{lemma}
\begin{proof} This is the same as Lemma 13.14 in \cite{DHM24} as $\Mb_D$ now has no BM-atom, and the effect of \eqref{it.stepc1} is equivalent to \eqref{it.stepa1} or \eqref{it.stepb1} viewed within each newly formed component in $\Mb_D$.
\end{proof}

\textbf{Proof part 3.} We prove that the number of bad atoms $\nf\in\Mb_U$ cut in Definition \ref{def.alg_maincr} (\ref{it.alg_maincr_4}) is at most  $\#_{\mathrm{fucp}(\Mb_U)}+\#_{\mathrm{ovcp}(X_1)} + \#_{\rm BM}$. In fact, such $\nf$ must be the atom chosen in Definition \ref{def.alg_maincr} (\ref{it.alg_maincr_3}) (the $\pf\in S_\nf\backslash\{\nf\}$ in (\ref{it.alg_maincr_4}) will not be cut as a bad molecule because it has an ov-parent in $S_\nf$ which must be cut before $\pf$). Moreover, there are two possibilities when this $\nf$ can be bad: either the component of $\Mb_U$ containing $\nf$ contains no fixed end, or $\nf$ is the lowest remaining atom in $X_0\cup X_1$ and has no fixed end. The former case leads to at most $\#_{\mathrm{fucp}(\Mb_U)}$ choices for $\nf$.

In the latter case, $\nf$ must have no fixed end and must not have any ov-child in $X_0\cup X_1$ (otherwise this ov-child must have been cut before $\nf$ as $\nf$ is the lowest remaining atom in $X_0\cup X_1$, so $\nf$ would have at least one fixed end), which means that $\nf \in X_0$ if and only if it is a BM-atom, by definition of $X_0$ (and that $\overline{\Mb}_{\mathrm{2conn}}^2=\varnothing$). If $\nf\in X_1$ and is a binary atom, then it must belong to one of the initial ov-components of $\overline{X}_1$, say $A$. Since $\nf$ has no ov-child in $A$, by Lemma \ref{lem.3_comp_property} (i), it must be the unique lowest atom in $A$ (here $A$ satisfies the assumption of Lemma \ref{lem.3_comp_property}, because each atom in $X_1$ has two particles, and each $\{4_2\}$ atom in $X_1$ must have an ov-child in $X_0\cup\Mb_D$, thus has at most one ov-child in $X_1$), which leads to at most $\#_{\mathrm{ovcp}(\overline{X}_1)}$ choices for $\nf$ due to the same number of choices of $A$. Finally, if $\nf\in X_0\cup X_1$ is a BM-atom without fixed end, this gives at most $\#_{\mathrm{BM}}$ further choices. Putting together, this finishes Part 3 of the proof.

\textbf{Proof part 4.} We now complete the proof of (\ref{eq.alg_maincr}). The bound for $\#_{\mathrm{bad}}$ follows directly from Part 3, as no bad atom can be cut in Definition \ref{def.alg_maincr} (\ref{it.alg_maincr_5})--(\ref{it.alg_maincr_6}), or Definition \ref{def.alg_maincr} (\ref{it.alg_maincr_9}) (as shown in Part 2). To prove the lower bound for $\#_{\{3_2 3_2\}}$ (which is $\#_{\{3_2 3_2 A\}}+\#_{\{3_2 3_2 B\}}$), note that with at most $\mu:=\#_{\mathrm{fucp}(\Mb_U)}+\#_{\mathrm{ovcp}(\overline{X}_1)} + \mathsf{N} \cdot \#_{\mathrm{BM}}$ exceptions, each $\{3_2\}$ molecule in $\Mb_D$ created in Definition \ref{def.alg_maincr} (\ref{it.alg_maincr_5})--(\ref{it.alg_maincr_6}) is always paired with a \{$3_23_2$\} molecule (the exceptions correspond to those atoms $\pf$ in Definition \ref{def.alg_maincr} (\ref{it.alg_maincr_4}) are $\{4_2\}$ atoms, or BM-atoms, which may cause $\qf$ to be cut as a $\{3_2\}$ molecule without a pairing $\{3_2 3_2\}$ molecule), so the number of atoms in $\Mb_D$ cut in these steps is at most $3\cdot\#_{\{3_23_2\}}+\mu$.

Now consider the $\#_{\mathrm{UD}}$ good UD connections and the corresponding atoms $\rf\in \Mb_D$. If any such $\rf$ is not cut in Definition \ref{def.alg_maincr} (\ref{it.alg_maincr_5})--(\ref{it.alg_maincr_6}), then in the molecule $\Mb_D$ after Definition \ref{def.alg_maincr} (\ref{it.alg_maincr_8}), this $\rf$ must be a $\{3_2\}$ atom (it cannot be $\{4_2\}$ as it belong to a UD connection, and cannot be $\{2_2\}$ because $\Mb_D$ is ov-proper by Part 2), and must belong to a component of $\Mb_D$ that contains at least two $\{3_2\}$ atoms (by considering the pre-algorithm path connecting $\rf$ to a $\{3_2\}$ atom or another UD connection, and the possible place where a bond of this path is broken; note that a broken bond creates a fixed end and thus $\{3_2\}$ atom). We collect all these atoms $\rf$ (which are not cut in Definition \ref{def.alg_maincr} (\ref{it.alg_maincr_5})--(\ref{it.alg_maincr_6})), with the number of them being $\geq \#_{\mathrm{UD}}-(3\cdot\#_{\{3_23_2\}}+\mu)$ by the above proof, and consider all the components of $\Mb_D$ (after Definition \ref{def.alg_maincr} (\ref{it.alg_maincr_8})) containing them; let the number of $\{3_2\}$ atoms in each of these components be $m_j\geq 2$, then after Definition \ref{def.alg_maincr} (\ref{it.alg_maincr_8}) we have that (note that each component in $\Mb_D$ will have at least one $\{3_2\}$ atom after cutting $\Mb_U$, see Remark \ref{rem.full_cut}):
\begin{equation}\label{eq.alg_maincr_2}
\#_{\{3_2\}} - \#_{\mathrm{comp}}\geq\sum_j(m_j-1) \geq\frac{1}{2}\sum_j m_j\geq\,(\#_{\mathrm{UD}}-3\cdot\#_{\{3_23_2\}}-\mu)/2,
\end{equation}
 where $\#_{\{3_2\}}$ and $\#_{\mathrm{comp}}$ are the number of $\{3_2\}$ atoms and number of components respectively. Then, consider the increment of the left-hand side of (\ref{eq.alg_maincr_2}) in each cutting in Definition \ref{def.alg_maincr} (\ref{it.alg_maincr_9}). Using Definition \ref{def.cutting}, and noting that each such cut can only create one fixed end at one $\{4_2\}$ atom (adding only one $\{3_2\}$ atom) in each resulting (extra) component, we can verify that this is invariant after cutting one $\{3_2\}$ atom, and decreases by 1 after each cutting in Definition \ref{def.alg_maincr} (\ref{it.alg_maincr_5})--(\ref{it.alg_maincr_6}). Since this value should be come $0$ after all atoms of $\Mb_D$ have been cut, we conclude that \[\#_{\{3_23_2\}}\geq \frac{1}{2}(\#_{\mathrm{UD}}-3\cdot\#_{\{3_23_2\}}-\mu),\] and therefore
\[\#_{\{3_23_2\}}\geq\frac{1}{5}\cdot(\#_{\mathrm{UD}}-\mu).\]This completes the proof.
\end{proof}
\subsection{Finishing the proof}\label{sec.finish} We now finish the proof of Proposition \ref{prop.case5}. This part corresponds to subsection 13.5 in \cite{DHM24}. The proof is divided into $6$ stages. Following from Reductions \ref{red.2} and \ref{red.3}, we assume that 
\begin{itemize}
\item Each of $\Mb_U$ and $\Mb_D$ is a forest and has no fully fixed atoms. No atom in $\Mb_D$ is a parent of any atom in $\Mb_U$, and $\Mb_U$ (resp. $\Mb_D$) has no C-top (resp. C-bottom) fixed end.
\item The number of BM-atoms in $\Mb_U \cup \Mb_D$ is at most $(C^*_5)^{-1} \rho'''$.
\item The number of full  components of $\Mb_U$ is $\#_{\mathrm{fucp}(\Mb_U)}\leq 10(C_{5}^*)^{-1}\cdot\rho'''$. Each component in $\Mb_D$ has either at least two bonds connecting to $\Mb_U$, or one bond connecting to $\Mb_U$ and one atom with fixed ends.
\item The number of good UD connections (i.e. UD connections that are connected to either another UD connection or an atom with fixed end in $\Mb_D$, see Proposition \ref{prop.alg_maincr}) satisfy $\#_{\mathrm{UD}}\geq \rho'''/4$.
\item The number of components of $\Mb_D$ is $\#_{\mathrm{comp}(\Mb_D)}\leq (C_{3}^*)^{-1}\cdot\rho'''$, and the number of weakly degenerate atoms and pairs in $\Mb_D$ is less than $(C_{2}^*)^{-1}\cdot\rho'''$.
\item We have $|\overline{\Mb}_{\mathrm{2conn}}^2|\leq (C_{2}^*)^{-1}\cdot\rho'''$, with the set $\overline{\Mb}_{\mathrm{2conn}}^2$ defined in Proposition \ref{prop.comb_est_case6}.
\end{itemize}
\begin{proof}[Proof of Proposition \ref{prop.case5}]
\textbf{Stage 1.} We cut each atom in $\overline{\Mb}_{\mathrm{2conn}}^2$ as free. If any atom in $\Mb_{UD}$ becomes fully fixed, we also cut it as free, until $\Mb_{UD}$ contains no such atoms. Note that the contribution of these operations to $\#_{\mathrm{bad}}$ is at most $|\overline{\Mb}_{\mathrm{2conn}}^2|  \leq 2 (C_{2}^*)^{-1}\cdot\rho'''$. Moreover, after \textbf{Stage 1}, we still have 
\begin{equation}\label{eq.comb_est_case71}\#_{\mathrm{UD}}\geq \rho'''/4-2|\overline{\Mb}_{\mathrm{2conn}}^2| \geq \rho'''/5.
\end{equation}
After \textbf{Stage 1}, $\Mb_{UD}$ is still a canonical UD molecule, and now $\overline{\Mb}_{\mathrm{2conn}}^2=\varnothing$. The number of full components of $\Mb_U$ and $\Mb_D$ still satisfy that $\#_{\mathrm{fucp}(\Mb_U)}+\#_{\mathrm{fucp}(\Mb_D)}\leq 4(C_{3}^*)^{-1}\cdot \rho'''$.

\textbf{Stage 2.} Now consider all the components of $\Mb_D$ that contain at least two $\{3_2\}$ atoms, and all the $\{3_2\}$ atoms contained in one such component. Let the set of all such $\{3_2\}$ atoms be $\overline{Z}_0$, and let $\overline{Y}$ be the union of all weakly degenerate atoms and atom pairs in $\Mb_D$ (see Proposition \ref{prop.comb_est_case4}). Let the number of ov-components of $\overline{Z}_0$ in $\Mb_D$ be $\#_{\mathrm{ovcp}(\overline{Z}_0)}$ as in Propositions \ref{prop.3comp} and \ref{prop.func_select}. 

Now, if $\#_{\mathrm{ovcp}(\overline{Z}_0)}\geq (10/\upsilon)^4(C_{2}^*)^{-1}\cdot \rho'''$, then we cut $\Mb_D$ as free. For each component of $\Mb_D$ containing at least two $\{3_2\}$ atoms, we cut it into elementary molecules using \textbf{3COMPDN} (the dual of \textbf{3COMPUP} in Definition \ref{def.3comp_alg}); for other components we cut them into elementary molecules using \textbf{DOWN}. Then cut $\Mb_U$ into elementary molecules using \textbf{UP}.

In this case, by Proposition \ref{prop.3comp} we get $\#_{\{3_23_2B\}}=\#_{\{4_2 (2l)_l\}}=0$, and
\begin{equation}\label{eq.comb_est_case72}
\begin{aligned}
        \#_{\{3_2 3_2 A\}}\geq (1/10 )\cdot\#_{\mathrm{ovcp}(\overline{Z}_0)} - \mathsf{N} \cdot \#_{\mathrm{BM}} \geq 100  \upsilon^{-4} \cdot (C_{2}^*)^{-1}\cdot\rho''',\\
    \#_{\mathrm{bad}}\leq \#_{\mathrm{fucp}(\Mb_U)}+\#_{\mathrm{fucp}(\Mb_D)}+(C_{2}^*)^{-1}\cdot\rho'''\leq 2\cdot(C_{2}^*)^{-1}\cdot \rho''',
\end{aligned}
\end{equation}
where we used Proposition \ref{prop.3comp} and the fact that $\#_{\mathrm{BM}} \leq (C^*_{5})^{-1} \rho'''$ by Remark \ref{rem.BM}. This already proves (\ref{eq.comb_est_case51}) (recall that all \{$3_23_2A$\} molecules are treated as good).

\textbf{Stage 3.} Now suppose $\#_{\mathrm{ovcp}(\overline{Z}_0)}\leq (10/\upsilon)^4(C_{2}^*)^{-1}\cdot \rho'''$. Define $\overline{Z}$ to be the union of $\overline{Z_0}$, together with the $\{3_2\}$ atoms in those components of $\Mb_D$ that contain one $\{3_2\}$ atom and at least one atom in $W$, $\overline{Y}$ or at least one BM-atom. Then for the set $\overline{Z}$ we have $\#_{\mathrm{3ovcp}(\overline{Z})}\leq \#_{\mathrm{ovcp}(\overline{Z}_0)} + |W| +|\overline{Y}| + \#_{\mathrm{BM}}$. Note also $|\overline{Y}|\leq 10(C_{2}^*)^{-1}\cdot\rho'''$ by Reduction \ref{red.3}.

We then define $S:=\textbf{SELECT}(\widetilde{\Mb}_D,W, X,\overline{Y}, \overline{Z})$ as in Definition \ref{def.func_select}, where $\widetilde{\Mb}_D$ contains only those components of $\Mb_D$ that contain an atom in $W \cup X \cup \overline{Y}\cup \overline{Z}$. We then cut $S$ as free and subsequently cut it into elementary molecules using \textbf{DOWN}. If any atom in $\Mb_{U}$ becomes fully fixed, we also cut it as free, until $\Mb_{U}$ contains no such atoms.

By Proposition \ref{prop.func_select}, we know that $\Mb_{UD}$ is still a UD molecule with $\overline{\Mb}_{\mathrm{2conn}}^2=\varnothing$, and $\Mb_D$ becomes ov-proper after cutting $S$. Moreover, since $W\subseteq S$, after cutting $S$ the molecule satisfies (iii) in Definition \ref{def.alg_maincr}: if a BM-atom $\pf\in\Mb_U$ had two ov-segments containing C-bottom edges at $\pf$ and reaching the same component of $\Mb_D$, then the first atoms of $\Mb_D$ on them belong to $W_\pf$ (Definition \ref{def.W}), so they are cut, and these ov-segments are broken. Now $\Mb_D$ may contain C-bottom fixed ends, but $\Mb_U$ still does not contain C-top fixed end. We know that the number of ov-segments corresponding to UD connections $e$, that are broken in the process of cutting $S$ as free is at most
\begin{equation}\label{eq.comb_est_case74}10\cdot(3\cdot\#_{\mathrm{ovcp}(\overline{Z})} + |W| + |\overline{Y}| ) + \mathsf{N} |X| \leq C\cdot (C_{2}^*)^{-1}\cdot\rho''',
\end{equation} by Proposition \ref{prop.w}, therefore we still have
\begin{equation}\label{eq.comb_est_case75}\#_{\mathrm{UD}}\geq \rho'''/5-C\cdot (C_{2}^*)^{-1}\cdot\rho'''\geq \rho'''/6,
\end{equation} after \textbf{Stage 3}. After cutting $S$, there may be some components in $\Mb_D$ not connected to $\Mb_U$, note that those components of $\Mb_D$ cannot be full.

\textbf{Stage 4.} Now the molecule $\Mb_{UD}$ satisfies the assumptions (i)--(ii) in Definition \ref{def.alg_maincr}, and has $\#_{\mathrm{UD}}\geq \rho'''/6$. Moreover, $\Mb_D$ is ov-proper and does not contain any weakly degenerate atom or atom pair, and does not contain any BM-atom. Next we perform the following operations, and repeat until this can no longer be done:
\begin{itemize}
\item If a subset $A\subseteq\Mb_U$ is not connected to $\Mb_{UD}\backslash A$ by any bond, then we cut it into elementary molecules using \textbf{UP}.
\item If a subset $A\subseteq \Mb_U$ has no fixed ends, and is connected to $\Mb_{UD}\backslash A$ by only one bond (this bond can be connected to either $\Mb_U\backslash A$ or $\Mb_D$), then we cut $A$ \emph{as fixed} and subsequently cut it into elementary molecules using \textbf{UP} or \textbf{DOWN} (depending on the direction of this only bond).
\end{itemize}
By the same argument as in the corresponding part of \cite{DHM24}, after \textbf{Stage 4}, we have that no subset $A\subseteq\Mb_U$ as described above, and we still have
\begin{equation}\label{eq.comb_est_case755}
\#_{\mathrm{UD}}\geq \rho'''/6-20 (C_5^*)^{-1}\cdot\rho'''\geq \rho'''/7.
\end{equation}
The assumptions (i)--(ii) in Definition \ref{def.alg_maincr} are still preserved (if needed, we can further remove any full component of $\Mb_D$ that is not connected to $\Mb_U$ by a bond, and cut it into elementary molecules \emph{following Definition \ref{def.alg_maincr} (\ref{it.alg_maincr_9})}; the number of bad molecules produced is at most $\#_{\mathrm{fucp}(\Mb_D)}\leq (C_3^*)^{-1}\cdot\rho'''$ and is negligible). Now define the sets $X_0$, $X_1$ and $\overline{X}_1$ as in Definition \ref{def.alg_maincr}, and consider the value of $\#_{\mathrm{ovcp}(\overline{X}_1)}$ as in Proposition \ref{prop.alg_maincr}.

\textbf{Stage 5.} If $\#_{\mathrm{ovcp}(\overline{X}_1)}\geq 10^{-6}\upsilon^2\cdot\rho'''$, then we cut $\Mb_D$ as free and cut it into elementary molecules \emph{following Definition \ref{def.alg_maincr} (\ref{it.alg_maincr_9})} (except that if a component of $\Mb_D$ is full, then we first cut any $\{4_2\}$ atom as free). If any atom in $\Mb_U$ becomes fully fixed then we cut it as free, and repeat until no such atom exists in $\Mb_U$. This produces elementary molecules (with the \{$3_23_2$B\} molecule only in $\Mb_D$ and \{$2_2$\} molecules only in $\Mb_U$), with the number of \{$4_2$\} molecules bounded by $\#_{\mathrm{fucp}(\Mb_D)}$; moreover all \{$3_23_2$B\} molecules are good due to absence of weakly degenerate atoms and atom pairs.

By definitions of $X_j$ in Definition \ref{def.alg_maincr}, it is clear that all the atoms in $X_0$ have been cut, and all the atoms in $X_1$ now have fixed ends and are not fully fixed, and all the other C-atoms or BM-atoms now have no fixed end. In addition, we claim that after the above cuttings, with \emph{at most $\#_{\mathrm{BM}}$ exceptions}, each atom $\nf\in X_1$ must belong to a component of $\Mb_U$ with at least two $\{3_2\}$ atoms. In fact, let the component of $\Mb_U$ containing $\nf$ be $A$. Then before the cuttings, $A$ must have one bond connecting to $\Mb_{UD}\backslash A$, and either has another such bond or an atom with fixed end (due to the operations in \textbf{Stage 4}). In either case, after the cutting operations, these two bonds (or one bond and one atom with fixed ends) will lead to two fixed ends in $A$ (or a fixed end at an O-atom). If $A$ has at most one $\{3_2\}$ atom, then there must exist at least one BM-atom in $A$. The number of such components is no more than $\#_{\mathrm{BM}}$, and the number of $\Mb_{(3)}$ components in these components is no more than $\#_{\mathrm{BM}}$. 
 
Now consider all the components of $\Mb_U$ that contains at least two $\{3_2\}$ C-atoms in $X_1$ after the above cutting operations, (i.e. ignore the other components, which can always be cut into elementary molecules using \textbf{UP}, in view of absence of C-top fixed ends in $\Mb_U$); denote the union of these components by $\widetilde{\Mb}_U$. Recall $\overline{X}_1$ in Proposition \ref{prop.alg_maincr}, and define $X_2 :=\overline{X}_1 \cap \widetilde{\Mb}_U $. In $\Mb_{UD}$ before cutting, the number of ov-components of $X_2$ is at least $\mu:=\#_{\mathrm{ovcp}(\overline{X}_1)} - \#_{\mathrm{BM}}\geq 10^{-6}\upsilon^2\cdot\rho'''$. After the cutting operations, the set of $\{3_2\}$ atoms in $\widetilde{\Mb}_U$ belongs to $X_2 \cup X_O$, where $X_O$ is the set of $\{3_2\}$ atoms in $\widetilde{\Mb}_U$ that are O-atoms. Now note that, the number of ov-components of $X_2\cup X_O$ is no less than the number of ov-components of $X_2$ (which is obviously no less than $\mu$, because these cuttings can only break ov-segments but not create them). This is because for each $\{3_2\}$ O-atom, it must have two edges (free ends or bonds) that are serial, and one free end that is serial with one fixed end; in other words one cannot change to a different ov-segment by going through this O-atom. Therefore, if two atoms in $X_2$ are ov-connected within $X_2\cup X_O$, they must also be ov-connected within $X_2$.

Then, we can apply \textbf{3COMPUP} to $\widetilde{\Mb}_U$, and use Proposition \ref{prop.3comp} to get
\begin{equation}\label{eq.comb_est_case75-}
\#_{\{3_23_2A\}}\geq \frac{1}{10}\cdot\mu - 2 \mathsf{N}\cdot \#_{\mathrm{BM}} \geq C^{-1}\cdot\rho''',
\end{equation}
\begin{equation}\label{eq.comb_est_case75+}\#_{\mathrm{bad}}\leq \#_{\mathrm{fucp}(\Mb_U)}+\#_{\mathrm{fucp}(\Mb_D)}+C(C_{2}^*)^{-1}\cdot\rho'''\leq C(C_{2}^*)^{-1}\cdot \rho''',
\end{equation} which already proves (\ref{eq.comb_est_case51}) (with all \{$3_2 3_2A$\} molecules treated as good).

\textbf{Stage 6.} Finally, if $\#_{\mathrm{ovcp}(\overline{X}_1)}\leq 10^{-6}\upsilon^2\cdot\rho'''$, then we apply \textbf{MAINUD} to cut $\Mb_{UD}$ into elementary molecules. By Proposition \ref{prop.alg_maincr}, and note that all \{$3_23_2$B\} molecules are good (as they must have both atoms in $\Mb_D$, and $\Mb_D$ does not contain weakly degenerate atoms or atom pairs; note also that all the \{$3_23_2$B\} molecules must have both atoms being C-atoms, because the fixed end at any O-atom cannot be serial with the bond), we get that
\begin{equation}\label{eq.comb_est_case76}
\#_{\mathrm{good}}\geq \frac{1}{5}\cdot\big(\#_{\mathrm{UD}}-\#_{\mathrm{fucp}(\Mb_U)}-\#_{\mathrm{ovcp}(X_1)} - \mathsf{N} \#_{\mathrm{BM}})\geq\frac{1}{40}\cdot\rho'''
\end{equation}
\begin{equation}\label{eq.comb_est_case76+}\qquad \#_{\mathrm{bad}}\leq \#_{\mathrm{fucp}(\Mb_U)}+\#_{\mathrm{ovcp}(X_1)}+ C(C_{2}^*)^{-1}\cdot \rho'''\leq 10^{-5}\upsilon^2\cdot \rho''',
\end{equation} which proves (\ref{eq.comb_est_case51}). This completes the proof of Proposition \ref{prop.case5}, and thus proves Proposition \ref{prop.comb_est}, in the case when $\Mb\in\Fc_{\Lambda_\ell}$.

\textbf{The finale: the case when $\Mb\in\Fc_{\Lambda_\ell}^{\mathrm{err}}$.} Finally, we prove Proposition \ref{prop.comb_est} when $\Mb\in\Fc_{\Lambda_\ell}^{\mathrm{err}}$. First, note that we can cut $\Mb$ using \textbf{UP} which leads to $\#_{\mathrm{good}}\geq 0$ and $\#_{\mathrm{bad}}\leq|H|$ (see Part 1 of the proof of Proposition \ref{prop.case1}). In view of the extra term $C^*|H|$ on the right hand side of (\ref{eq.overall_alg}), we see that (\ref{eq.overall_alg}) is true if $\rho\lesssim C^*|H|$.

Now, suppose $\rho\gg C^*|H|(\gg C_{12}^*|H|)$, then we can repeat exactly the same proof in Sections \ref{sec.cuttingI}--\ref{sec.maincr}. The only place where Proposition \ref{prop.mol_axiom} (\ref{it.axiom6}) is used (beyond Proposition \ref{prop.case1} which is already taken cared of), is in layer selection (Definition \ref{def.layer_select}). However, in this case we always have $\ell_1-1\leq\ell-1$, so Proposition \ref{prop.mol_axiom} (\ref{it.axiom6}) remains true for the set $r(\Mb_{\ell_1-1})$, therefore the proof goes exactly the same way as in the $\Mb\in\Fc_{\Lambda_\ell}$ case.
\end{proof}
 \section{Proof of Proposition \ref{prop.est_fa}}\label{sec.f_Ac} In this section we prove Proposition \ref{prop.est_fa}. 
\begin{proof}[Proof of Proposition \ref{prop.est_fa}]
The proof is similar to Section 14 of \cite{DHM24} in the hard sphere case.  We list the changes.
\begin{enumerate}[(i)]
\item As in (14.3)--(14.4) of \cite{DHM24}, we write $f^\Ac(\ell\tau)=f^{\Ac,\mathrm{led}}+f^{\Ac,\mathrm{sub}}$. Here $f^{\Ac,\mathrm{led}}$ is the sum over $\Mb\in\Tc_{\Lambda_\ell}$ with $\rho(\Mb)=0$ (recall that $\rho$ does not count cycles inside boxes, Definition \ref{def.boxed_molecule} \eqref{it.recollision_number_box}) and $|\Mb|_p<\Lambda_\ell$, and with no multi-particle atom, no box and no top crossing atom. The term $f^{\Ac,\mathrm{sub}}$ contains the terms with $\rho(\Mb)>0$ or $|\Mb|_p\geq\Lambda_\ell$, or with a multi-particle atom, a box or a top crossing atom, the terms with $\Mb\in\Tc^{\mathrm{err}}_{\Lambda_\ell}$, and $f^{\Ac,\mathrm{err}}$ in \eqref{eq.fAterm}. The term $f^{\Ac,\mathrm{sub}}$ is put into $(f^\Ac)_{\mathrm{rem}}$.

We now explain the gain for the terms in $f^{\Ac,\mathrm{sub}}$. As in \cite{DHM24}, the top end of the root particle line is the unique fixed end, and all other ends of $\Mb$ are free. The molecule may not be regular (it may have an O-fixed end which is serial with an O-bond), and we make it regular again as in Point 4 of Section 14.3 of \cite{DHM24}. Let $\Xf$ and $\Bf$ be the numbers of multi-particle C-atoms or C-constituents and of boxes in $\Mb$.

If $\Xf+\Bf\geq10\rho(\Mb)$, we cut $\Mb$ using \textbf{UP2}, starting from the root, and view each box as a single atom. The number of fixed ends created in this process equals the number of bonds, which is $|\Mb|-1+\rho(\Mb)$, and at most one atom is cut with no fixed end (Proposition \ref{prop.alg_up} \eqref{it.up_proof_2}), so at most $\rho(\Mb)$ atoms are cut with at least two fixed ends. Hence at least $\Xf+\Bf-\rho(\Mb)\geq\frac{9}{10}(\Xf+\Bf)$ multi-particle atoms and boxes are cut with at most one fixed end. Each of them is either a good molecule, or can be further cut into elementary molecules with at least one good molecule (Proposition \ref{prop.est_single} \eqref{it.est_kl} and Proposition \ref{prop.alg_up} \eqref{it.up_proof_4}). Hence we gain $\eps^{c(\Xf+\Bf)}$. Otherwise, if $\Xf+\Bf<10\rho(\Mb)$, we gain $\eps^{c\rho(\Mb)}$ from the cycles, as in Section 14.3 of \cite{DHM24}. In both cases we gain $\eps^{c(\rho(\Mb)+\Xf+\Bf)}$, which absorbs the combinatorial factor $|\log\eps|^{C^*(\rho(\Mb)+\Xf+\Bf)}$ from Proposition \ref{prop.layerrec3}.

If $\rho(\Mb)=0$, $\Xf=\Bf=0$, $|\Mb|_p<\Lambda_\ell$ and $\Mb\notin\Tc^{\mathrm{err}}_{\Lambda_\ell}$, then $\Mb$ is in $f^{\Ac,\mathrm{sub}}$ because it has a top crossing atom $\nf$. This atom is marked as short interaction, so $|t^-_\nf-\ell\tau|\leq\eps^{\upsilon_0}$, and the integral in $t^-_\nf$ gives an extra factor $\eps^{\upsilon_0}|\log\eps|^{C^*}$.

In all cases, the weights in $\mathrm{Bol}^{\beta_\ell}$ are controlled as in \cite{DHM24}: since $\Mb$ has no bottom crossing atom (Definition \ref{def.set_T_F} \eqref{it.set_T}; this also holds for $\Tc^{\mathrm{err}}_{\Lambda_\ell}$), the particles of $\Mb$ have no potential energy at time $(\ell-1)\tau$, so by conservation of energy and $\Phi\geq0$ (Assumption \ref{ass.phi1}), the kinetic energy of the root at time $\ell\tau$ is at most the total kinetic energy of the incoming particles, even though the root may still be interacting at time $\ell\tau$. Hence these terms are $O(\eps^{\theta_0})$ in $\mathrm{Bol}^{\beta_\ell}$.
\item The terms in $f^{\Ac,\mathrm{led}}$, i.e. the trees with only binary C-atoms (and O-atoms) and without crossing atoms, give the Duhamel iterates of the Boltzmann equation \eqref{eq.boltzmann_omega} with the cross-section $b$: by Lemma \ref{lem.liouville}, $[(v-v_1)\cdot\nu]_-\,\d\nu=b(|v-v_1|,\omega)\,\d\omega$, and for fixed $\omega$ the map $(v,v_1)\mapsto(v',v_1')$ is the hard sphere map, with unit Jacobian. Compared with hard spheres, the outgoing positions in \eqref{eq.relation_collision} are shifted by at most $C\eps+\Delta t\,|v_{e_1}-v_{e_2}|\leq C\eps$ by Proposition \ref{prop.upperbond.interaction}; using the bound on $\nabla_x f$ in \eqref{eq.boltzmann_decay_4}, this shift is also absorbed into $(f^\Ac)_{\mathrm{rem}}$.
\item The remaining estimates in Section 14 of \cite{DHM24} use the collision kernel only through $b\geq0$, $b(|w|,-\omega)=b(|w|,\omega)$ and $\int b(|w|,\omega)\,\d\omega=c_d|w|$, which hold here by Lemma \ref{lem.liouville}, so they apply without change.
\end{enumerate} 
\end{proof}

\section{Proof of the error term estimates}\label{sec.error} In this section we prove Proposition \ref{prop.cumulant_error}. The key ingredient here is a cutting algorithm specifically designed for the layer $\Mb_\ell\in\Fc_{\Lambda_\ell}^{\mathrm{trc.err}}$ for the molecule $\Mb\in\Fc_{\vLambda_\ell}^{\mathrm{trc.err}}$. This is stated in Propositions \ref{prop.comb_est_extra} \ref{prop.comb_est_extra2} below. These two Propositions play the role of Proposition 15.1 in \cite{DHM24}.

\begin{proposition}\label{prop.comb_est_extra} Let $\Gamma$ be a large absolute constant depending only on the dimension $d$ and $C_\Phi$. Let $\Mb$ be a connected boxed molecule without O-connections, with $|\Mb|\leq |\log\varepsilon|^{C^*}$ and $\Gamma<\rho(\Mb)\leq 2\Gamma$, in which every C-atom is marked as short interaction, and has at most $d+1$ particles. Then there exists an operation sequence (formed by cutting and splitting) which produces at most $C^{|\Mb|}|\log\varepsilon|^{C^*}$ sub-cases, such that in each sub-case, $\Mb'$ consists of good, normal, and bad elementary molecules, and the numbers $\#_{\mathrm{good}}$ and $\#_{\mathrm{bad}}$ of good and bad molecules satisfy that
\begin{equation}
\label{eq.extra_case_1} (\upsilon/2)\cdot(\#_{\mathrm{good}})-d\cdot (\#_{\mathrm{bad}})\geq \upsilon/2,
\end{equation}
or there exists at least $2 d^2$ good molecules which are multi-particle atoms, each of then have excess $\eps^{\upsilon_0/2}$, and one bad molecule (note that the total excess $\eps^{d^2 \cdot \upsilon_0} \leq \eps^{d - 1 + \upsilon}$).
\end{proposition}
\begin{proof} The proof will occupy Section \ref{sec.last_proof}. 
\end{proof}

\begin{proposition}\label{prop.comb_est_extra2}  Let $\Mb$ be a connected molecule with only C-atoms, with $|\Mb|\leq |\log\varepsilon|^{C^*}$ and either of the following
\begin{enumerate}
    \item  $\Mb_\ell$ has $2 \Gamma$--$4 \Gamma$ boxes and $\rho(\Mb_{\ell}) \leq \Gamma$. 
    \item\label{it.prop9.2_2}  $\Mb_\ell$ has an atom that is marked as long interaction or an atom with more than $d+1$ particles. 
\end{enumerate}
    Then there exists a cutting sequence, which produces at most $C^{|\Mb|}|\log\varepsilon|^{C^*}$ sub-cases, such that in each sub-case, $\Mb'$ consists of good, normal, and bad elementary molecules, and the numbers $\#_{\mathrm{good}}$ and $\#_{\mathrm{bad}}$ of good and bad molecules satisfy 
\eqref{eq.extra_case_1}, or (in case (\ref{it.prop9.2_2})) the following excess form holds: there is exactly one bad molecule, which has excess $\eps^{d-1+\upsilon}$ (Definition \ref{def.excess}; i.e. $\|\Ic(Q)\|_{L^\infty}\leq\eps^\upsilon|\log\eps|^{C^*}\|Q\|_{L^\infty}$ for this molecule), and all other molecules are good or normal. 
\end{proposition}
\begin{proof}
    If $\Mb_\ell$ has $2 \Gamma$--$4 \Gamma$ boxes. Note that $\rho(\Mb_\ell) \leq \Gamma$. By applying {\bf UP1} algorithm, the number of atoms with more than $1$ fixed end is less than $\Gamma$. Hence, apart from $\Gamma$ boxes, other boxes has at most one fixed end. We then apply {\bf UP} to those boxes. Hence, $\#_{\mathrm{good}} \geq \Gamma $ and the number of bad molecule is $1$.

    If $\Mb_\ell$ has an atom that is marked as long interaction or an atom with more than $d+1$ particles, then we cut $\Mb_{\ell}$ by {\bf UP}, but first choose the atom marked as long interaction. By applying Proposition \ref{prop.est.d+2}, there exists exactly one bad atom with excess $\eps^{d-1+ \upsilon}$ (see Definition \ref{def.excess}).
\end{proof}

\begin{proposition}\label{prop.est.d+2}
    Let $\Mb$ be a C-molecule with only one atom $\nf$ with no fixed end. Assume that it is marked as long interaction (see Definition \ref{def.top_reduction}), or has more than $d+1$ particles, then $\Ic_\Mb(Q) \leq \eps^{\upsilon}\cdot |\log \eps|^{C^*} \cdot \|Q\|_{L^\infty}$.
\end{proposition}
    \begin{proof}
         We denote $t^-$ and $t^+$ to be the moment that the collision happens and completes (in the sense of Definition \ref{def.admissible}. First we assume this atom is a binary atom (so it is marked as long interaction). By the same calculation as in Proposition \ref{prop.est_single} \eqref{it.est_3}, we have
        \begin{equation*}
            \Ic_\Mb(Q) \leq \eps^{-(d-1)} \int [(v_1- v_2)\cdot \nu_{12}]_- \d \nu_{12} \d v_2 \d t_{12} \d v_1 \d x_1.
        \end{equation*}
        Recall that in Definition \ref{def.top_reduction}, there are two cases that an atom is marked as long interaction. The first case is $|v_1 -v_2| \leq \eps^{1-\upsilon_0}$ at time $(\ell'[\nf] -1)\tau$, then $\d v_1 \d v_2$ in $\Ic(Q)$ can provide a factor $\eps^{d- d \upsilon_0} \leq \eps^{\upsilon}$, which meets our requirement. For the second case, the interaction duration of $\nf$ is at least $\eps^{\upsilon_0}$. By Proposition \ref{prop.upperbond.interaction}, we have $|v_1 -v_2| \cdot (t^+ - t^-) \leq C \eps$. Then, we also have $|v_1 - v_2| \leq \eps^{1 - \upsilon_0}$, and we can obtain the same estimate.

        Now assume the collision involves $3$ particles $\pb_1,\pb_2$ and $\pb_3$, assume the collision duration $\Delta t \sim \mu'$ for some $\eps \lesssim \mu' \lesssim 1$. We assume particle $\pb_1$ first collides with $\pb_2$, then $\pb_3$ collides with $\pb_1$. Therefore,
        \begin{equation*}
            \Ic_{\Mb}(Q) \leq \eps^{-(d-1)} \int [v_1 - v_2]_- \cdot [v_3 -V_1(t_{13})]_-\d \nu_{12} \d \nu_{13} \d v_1 \d v_2 \d v_3 \d t_{12} \d t_{13}.
        \end{equation*} We divide into the following cases:
        If at $t^+$ or $t^-$, we have $\sum_{i=1}^3|x_i - x_c + t(v_i-v_c)|^2 \sim \mu^2$ for some $\eps \lesssim \mu \lesssim 1$ (assume $t^-$ by symmetry, and $\sum_{i=1}^3|x_i - x_c + t^+ (v_i-v_c)|^2\lesssim \mu$). We assume $|x_3 -x_1 + t^-(v_3 - v_1)| \sim \mu$ (at $t^-$ we have $|x_1 -x_2 + t^- (v_1 -v_2)| = \eps$). Since $|v| \leq |\log \eps|^{C^*}$,  the collision between $\pb_1$ and $\pb_2$ lasts at least $\mu \cdot |\log \eps|^{-C^*} $. By Proposition \ref{prop.upperbond.interaction}, we have $|v_1 -v_2| \lesssim \eps \cdot \mu^{-1} $ ($\pb_1$ and $\pb_2$ must keep interacting until $\pb_3$ gets into the collision, which takes at least $\mu \cdot |\log \eps|^{-C^*}$),  and we can gain the volume $(\eps\cdot \mu^{-1})^d  $ in the integration in $\d v_1 \d v_2$. Moreover, by applying Proposition \ref{prop.upperbond.interaction} to particles $\pb_1,\pb_2$ and $\pb_3$, we have $\sum_{i=1}^3 |v_i - (v_1 + v_2  + v_3)/3| \cdot \mu' \lesssim \mu$, we thereby gain from $\d v_3$ with volume $(\mu/\mu')^{d}$. Finally, we note that $|t_{13} - t_{12}| \lesssim \mu'$, and we can gain $\mu'$ from $\d t_{13}$. As a consequence, we totally gain $(\eps\cdot \mu^{-1})^d \cdot(\mu/\mu')^{d} \cdot  \mu \lesssim \eps^d \cdot (\mu')^{-(d-1)}$, as $\mu \leq \mu' \cdot |\log \eps|^{C^*}$. In particular, if $|v_i-v_j| \leq \eps^{1-\upsilon_0}$ or $\mu' \geq \eps^{\upsilon_0}$, we can gain at least $\eps^{d-1+ \upsilon}$, which meets our requirement.
        
        Now assume the collision involves $4$ particles, and the collision duration is $\mu''$. By the similar arguments as in the proof of Proposition \ref{prop.est_single} \eqref{it.est_kl}, we have an extra factor $\eps^{-(d-1)}$ in $\Ic_{\Mb}(Q)$. We assume that at time $t^-$, $\pb_1$ collides with $\pb_2$, then $\pb_3$ collides with $\pb_1$, and $\pb_4$ collides with $\pb_1$ (other collision orders can be proved in the similar way). We assume that $\sum_{i=1}^4|x_i - x_c + t^+ (v_i -v_c) | \lesssim \sum_{i=1}^4|x_i - x_c + t^- (v_i -v_c)| \sim \mu'$, where $x_c,v_c$ denotes the average position and velocity. Then we may assume that $|x_4 - x_1 + t^- (v_4 -v_1)| \sim \mu'$.  Before particle $\pb_4$ gets into the collision, the interaction duration of $\pb_1,\pb_2\pb_3$ is at least $\mu' |\log \eps |^{-C^*}$. By the same arguments as above, applied to  $\pb_1,\pb_2,\pb_3$ on $[t^-,s_4]$, where $s_4$ is the time when $\pb_4$ enters the collision, we gain $\eps^d(\mu')^{-(d-1)}$ from the variables of $\pb_1,\pb_2,\pb_3$. Here, if the spread of $\{\pb_1,\pb_2,\pb_3\}$ at $s_4$ is larger than at $t^-$, we use the time-reversed argument, taking as variables the state of $\pb_1,\pb_2,\pb_3$ at $s_4$, which is legitimate by the $L^1$-preservation of the flow.  
        Moreover, we can gain $\mu''$ from $\d t_4$ as  $|t_{14} - t^-| \leq \mu''$, and by applying Proposition \ref{prop.upperbond.interaction} to particles $\pb_1,\pb_2,\pb_3,\pb_4$, we gain volume $(\mu'/\mu'')^{d}$ in $\d v_4$. In total, we gain volume $ \eps^d\cdot {\mu''}^{-(d-2)}$. In particular, if $\nf$ is marked as long interaction, we can also gain the desired factor.

        Next, we consider the collision with $4$ particles, but with a different decoration: $\pb_1$ collides with $\pb_2$, and $\pb_3$ collides with $\pb_4$, then $\pb_2$ collides with $\pb_3$. We assume at $t^-$, the first collision is between $\pb_1$ and $\pb_2$. We assume that $\sum_{i=1}^4|x_i - x_c + t_1 (v_i -v_c)| \sim \mu'$. Then we may assume that $|x_4 - x_1 + t^- (v_4 -v_1)| \sim \mu'$, and $|x_1 - x_2 +t_{12} (v_1 - v_2)| \leq \eps$. Then the collision duration of $\pb_1$ and $\pb_2$ is at least $\mu'\cdot |\log \eps|^{-C^*}$. By Proposition \ref{prop.upperbond.interaction}, we have $|v_1- v_2 | \lesssim \eps \cdot \mu'^{-1}$, which provides a factor $(\eps \cdot \mu'^{-1})^d$ in $\d v_1 \d v_2$. Moreover, we have $|t_{14} - t^-|\leq \mu''$ and $|t_{34} - t^-| \leq \mu''$. Integrating these two time variables provides a factor $\mu'^2$. On the other hand, by applying Proposition \ref{prop.upperbond.interaction}, we have $(\sum_{i=1}^4|v_1 -v_c| )\cdot \mu'' \lesssim \mu'$, which provides a factor $(\mu' \cdot (\mu'')^{-1})^d$ in $\d v_3 \d v_4$. As a result, we can also gain a factor $\eps^d \cdot \mu''^{-(d-2)}$.

        By performing the above analysis repeatedly, we see that for a collision with $l$ particles and collision duration $\mu$, we can gain a factor $\eps^d \cdot (\mu)^{-(d+2-l)}$. The estimate for the case that $\nf$ is marked as long interaction can be obtained accordingly. Moreover, for a collision involves at least $d+2$ particles we can gain volume $\lesssim \eps^d$. This completes the proof.
    \end{proof}

We now prove Proposition \ref{prop.cumulant_error} assuming Proposition \ref{prop.comb_est_extra}.
\begin{proof}[Proof of Proposition \ref{prop.cumulant_error} assuming Proposition \ref{prop.comb_est_extra}] With Propositions \ref{prop.comb_est_extra} and \ref{prop.comb_est_extra2}, the proof of Proposition \ref{prop.cumulant_error} follows similar arguments as Section 15 in \cite{DHM24}.
\end{proof}
\subsection{The proof of Proposition \ref{prop.comb_est_extra}} \label{sec.last_proof}
Now we prove Proposition \ref{prop.comb_est_extra}. Although we do not have Burago--Ferleger--Kononenko \cite{BFK98} theorem for short-range potentials, this can be compensated by introducing boxes. 
 We first define the \textbf{transversal} subsets $A$ of $\Mb$, and prove some properties.
\begin{definition}
\label{def.trans} Let $\Mb$ be a full C-molecule, and $A\subseteq \Mb$ is an atom set. Throughout this definition, we view each box as a whole. We say $A$ is transversal, if we can decompose $\Mb\backslash A$ into two disjoint subsets $A^+$ and $A^-$, such that no atom in $A^-$ is parent of any atom in $A\cup A^+$, and no atom in $A$ is parent of any atom in $A^+$. For a box in $\Mb$, all atoms of it must either be in $A$, $A^+$ or $A^-$.

For any transversal set $A$, define the set $X_0(A)$ such that, an atom $\nf\in X_0(A)$ if and only if $\nf\not\in A$ and $\nf $ has $l$ (here and below $l$ denotes the number of particle lines of $\nf$) bonds connected to $A$. Inductively define $X_q(A)$ such that, $\nf\in X_q(A)$ if and only if $\nf\not\in A$, and either $\nf\in X_{q-1}(A)$, or $\nf$ has $l$ bonds connected to atoms in $X_{q-1}(A)\cup A$. Define $X(A)=\cup_{q\geq 0}X_q(A)$.
\end{definition} 
\begin{proposition}
\label{prop.trans}
Let $A$ be a connected transversal subset and $\Mb$ is a connected molecule.
\begin{enumerate}
\item\label{it.trans1} We can always choose $(A^+,A^-)$ in Definition \ref{def.trans} such that each component of $A^+$ has at least one bond connected to $A$.
\item\label{it.trans2} Define the sets $X_q^\pm (A)$ in the same way as $X_q(A)$, but with the sentence `` $\nf$ has $l$ bonds connected to $A$'' replaced by ``$\nf$ has $l$ ($l$ denotes the number of particle lines of $\nf$) children (or $l$ parents) in $X^\pm(A)=X(A) \cap A^\pm$''.
Then we have $X_q^\pm(A)=X_q(A)\cap A^\pm$. Define $X^\pm(A)=\cup_{q\geq 0}X_q^\pm(A)$, then $X^\pm(A)=X(A)\cap A^\pm$.
\item\label{it.trans3} The set $X(A)\cup A$ is also connected and transversal. Moreover any atom in $\Mb\backslash (X(A)\c A)$ with $l$ particles has at most $l-1$ bond with atoms in $X(A)\cup A$, so $X_0(X(A)\cup A)=\varnothing$.
\item\label{it.trans4} Recall $\rho(A)$ defined in Definition \ref{def.boxed_molecule} \eqref{it.recollision_number_box}. Then there exists a connected transversal set $B\supseteq A$, obtained by adding atoms $\nf_1,\nf_2,\cdots$ in this order, such that each $\nf_i$ has at least one and at most $l_{\nf_i} -1$ bonds to $A \cup \{\nf_1,\dots,\nf_{i-1}\}$, all of them bottom edges of $\nf_i$ (or all top edges), and such that either $B=\Mb$ or $X_0(B) \neq \emptyset$. In particular, we have $\rho(B) \leq \rho(A) + \sum_i (l_{\nf_i} - 2)$.
\end{enumerate}
\end{proposition}
\begin{proof} The proof of \eqref{it.trans1}--\eqref{it.trans3} follows from the same arguments as Proposition 15.4 in \cite{DHM24}. 



\textbf{Proof of (\ref{it.trans4}).} Suppose $A$ is connected and transversal; we may assume $A\neq\Mb$ and $X_0(A)=\varnothing$ (otherwise choose $B=A$). Since $\Mb$ is connected, there must exist an atom $\nf\in \Mb\backslash A$ that is either a parent or child of an atom in $A$. We may assume it is a parent, and then choose a lowest atom $\nf\in \Mb\backslash A$ among these parents. Since $X_0(A)=\varnothing$, we know that an atom $\nf$ with $l_\nf$ particle lines has at most $l_\nf-1$ bonds with atoms in $A$; let $A_1=A\cup\{\nf\}$, then $A_1$ is connected, and $\rho(A_1) \leq \rho(A) + l_{\nf}-2$

Now we claim that $A_1$ is transversal. To see this, define $C$ to be the set of atoms $\mf\in A^+\backslash\{\nf\}$ that are descendants of $\nf$, and decompose $\Mb\backslash A_1=A_1^+\cup A_1^-$, where $A_1^+=A^+\backslash (C\cup\{\nf\})$ and $A_1^-=A^-\cup C$. It is clear by checking all cases that no atom in $A_1$ can be parent of any atom in $A_1^+$, and no atom in $A^-$ can be a parent of any atom in $A_1\cup A_1^+$. If an atom $\mf\in C$ is parent of some atom $\pf\in A_1^+$, then $\pf\in A^+\backslash C$, but $\mf\in C$ and $\pf$ is child of $\mf$, so we also have $\pf\in C$, contradiction. Finally if $\mf\in C$ is parent of $\pf\in A_1=A\cup \{\nf\}$, clearly $\pf\in A$, but then $\mf$ is also parent of an atom in $A$, contradicting the lowest assumption of $\nf$. 

Now we know that $A_1$ is also connected transversal; by replacing $A$ by $A_1$ and repeating the above discussion, we eventually will reach some $B$ such that either $B=\Mb$ or $X_0(B)\neq\varnothing$, as desired.
\end{proof}

\begin{definition}\label{def.overline_X0}
    Let $\Mb$ be a connected C-molecule and $A$ be a subset (not necessarily transversal). We view each atom individually (i.e. not view box as a whole).  We define $\overline{X}^+_0(A)$ to be the set of atoms with $l $ (the number of particle line at the atom) children in $A$ (including those contained in box); $\overline{X}^-_0(A)$ is defined similarly. We also define $X$ Note that $\overline{X}^+_1(A)$ to be the set of atoms with at least one bond in $A$. Note that $\overline{X}^\pm_1(A)$ and $\overline{X}^\pm_1(A)$ may be different from ${X}^\pm_0(A)$ and $X^\pm_1(A)$defined in \eqref{def.trans} if $A$ is a transversal set.
\end{definition}
Now we define the algorithm \textbf{TRANSDN}, which is analogous to the \textbf{2CONNDN} algorithm and is applicable in the case when $|\overline{X}_0^+(A)|$ (Definition \ref{def.overline_X0}) is large.
\begin{definition}[Algorithm \textbf{TRANSDN}]
\label{def.alg_transup}
Suppose $\Mb$ is connected full C-molecule and $A$ is a connected transversal subset. We may choose $(A^+,A^-)$ as in Proposition \ref{prop.trans} (\ref{it.trans1}). Define the following algorithm:
\begin{enumerate}
\item\label{it.alg_transup_1} Throughout all steps in this algorithm, \emph{we treat each box in $A$ as a whole but atoms in box that belongs to $X^+(A)$ individually.} If $A$ contains any fully fixed atom $\nf$, then cut it as free, and repeat until there is no such atom left.
\item\label{it.alg_transup_2} Choose a \emph{highest} atom $\nf$ in the set of all atoms with fixed end in $A$ (or a highest atom in $A$, if $A$ does not contain such atom). Let $Z_\nf$ be the set of ancestors of $\nf$ in $A$.
\item \label{it.alg_transup_3} Starting from $\nf$, choose a \emph{lowest} atom $\mf$ in $Z_\nf$ that has not been cut. If $\mf$ is a box, then cut it as free, and then cut it into elementary molecules by {\bf UP} or {\bf DOWN}. If $\mf$ is a $\{3_2\}$ atom, and has a parent $\pf\in \overline{X}^+_0(A)$ that becomes fully fixed after cutting $\mf$, then cut $\{\mf,\pf\}$ as free. In all remaining cases, cut only $\mf$ as free. Repeat until all atoms in $Z_\nf$ have been cut. Then go to (\ref{it.alg_up_1}).
\item \label{it.alg_transup_4} Repeat (\ref{it.alg_transup_1})--(\ref{it.alg_transup_3}) until all atoms in $A$ have been cut. Then cut the remaining part of $A^+$ as free and cut it into elementary components using \textbf{UP}, then cut $A^-$ into elementary components using \textbf{DOWN}.
\end{enumerate}
We also define the dual algorithm \textbf{TRANSUP} by reversing the notions of parent/child etc. Note that we replace $A^+$ and $\overline{X}^+_0(A)$ by $A^-$ and $\overline{X}^-_0(A)$, and also replace $A^+$ by $A^-$ in Proposition \ref{prop.trans} (\ref{it.trans1}). 
\end{definition}

\begin{proposition}
\label{prop.alg_transup} For the algorithm \textbf{TRANSDN} (and same for \textbf{TRANSUP}), we have $\#_{\{3_23_2B\}}=\#_{\{4_2 (2l)_l\}}=0$ and $\#_{\mathrm{bad}}=1$, and $\#_{\mathrm{good}}\geq \frac{|\overline{X}^+_0(A)|}{2}-\rho(A) - \mathsf{N}\cdot \#_{\mathrm{box}}(A)- \#_{\rm mp} - 1 $, under the assumption that all two-atom molecules are good.
\end{proposition}
\begin{proof} It is clear, in the same way as in the proof of Proposition \ref{prop.alg_up}, that during the process, there is no bottom fixed end in $(A\backslash Z_\nf)\cup\{\nf\}$ or in $A^-$, and there is no top fixed end in $Z_\nf\backslash\{\nf\}$. Moreover $\#_{\mathrm{bad}}=1$ because no component of the remaining part of $A^+$ after Definition \ref{def.alg_transup} (\ref{it.alg_transup_3}) can be full when it is first cut in Definition \ref{def.alg_transup} (\ref{it.alg_transup_4}) (thanks to Proposition \ref{prop.trans} (\ref{it.trans1})), and no component of $A^-$ can be full when it is first cut in Definition \ref{def.alg_transup} (\ref{it.alg_transup_4}) (thanks to $\Mb$ being connected).

To prove the lower bound for $\#_{\mathrm{good}}$, note that only one atom in $A$ is a bad molecule (as $A$ is connected). Let the number of atoms in $A$ that belongs to fully fixed molecule be $q$, then by considering the increment of the same quantity $\sigma:=\#_{\mathrm{bo/fr}}-\sum_{l\geq 2}(2l -1)|\Mb_l|$ and applying the same arguments in the proof of Proposition \ref{prop.alg_up} (but restricted to $A$), and noticing that each atom viewed in $A$ is cut by itself, and initially all atoms in $A$ have no fixed end, we deduce that
\[\sum_{2\leq l \leq k\leq 2l}(2l-k-1)\#_{\{k_l\}} =\rho(A)-1,\]
and $\#_{\{(2l)_l\}}=1$, where $\rho(A)$ does not count cycles inside boxes. Therefore, we have $q\leq\rho(A)$. Then, for each $\nf\in \overline{X}^+_0(A)$, there exists atoms or boxes in $A$ which cut later than any other atoms connected to $\nf$. There are at least $\frac{|\overline{X}^+_0(A)| -\mathsf{N}  \cdot \#_{\mathrm{box}(A)} - (d+1)\#_{\rm mp}|}{2}$ such binary atoms. For each such atom, it is either cut by itself (if it is a \{$2_2$\} atom, $\{4_2\}$ atom), or cut together with $\nf$. Due to all two-atom molecules are good, we have
$$\#_{\mathrm{good}}\geq \frac{|\overline{X}^+_0(A)|}{2}-\rho(A) - \mathsf{N}\cdot \#_{\mathrm{box}}(A)- \#_{\rm mp} - 1,$$
which meets our requirement.

\end{proof}

\begin{definition}[The function \textbf{SELECT2}]
\label{def.func_select2} Let $A$ be a connected molecule with only C-atoms, no bottom fixed end, and no fully fixed atoms. Let $W$ be the set of atoms in $A$ that are connected to the same atom in $X^+(A) \backslash \overline{X}^+_0(A)$.\footnote{This corresponds to those multi-particle atoms in $X^+(A)$, with more than one bonds connected to $A$.}
Let $X$ be the set of B-atoms. Let $Z$ be the set of atoms in $A$ with fixed ends, or are multi-particle atoms, or are involved in weak degeneracies, and let $Y$ be a subset of atoms in $A \setminus (X\cup Z)$ such that $A$ becomes a forest after removing the atoms in $Y\cup X$. We define the function $\textbf{SELECT2}=\textbf{SELECT2}(A,W,X,Z,Y)$ as follows. 
\begin{enumerate}
\item\label{it.func_select2_1} Consider all the components of $W\cup X\cup Z\cup Y$ in $A$, which is a finite collection of disjoint subsets of $A$.
\item\label{it.func_select2_2} Since $A$ only has C-atoms, the shortest ov-distance between sets (Definition \ref{def.ov_connect}) coincides with the standard definition of distance on a graph. If any two of the subsets in (\ref{it.func_select2_1}), say $U$ and $V$, have the shortest distance (within $A$) which is at most $4$, then we choose one shortest path between an atom in $U$ and an atom in $V$, and let the atoms on this path be $\nf_j$. Then replace the two sets $U$ and $V$ by a single set which is $U\cup V\cup \{\nf_j\}$.
\item\label{it.func_select2_4} Repeat (\ref{it.func_select2_2}) until this can no longer be done. Next, if a single subset $U$ contains two atoms (which may be the same) that are connected by a path of length at most $4$ with none of the intermediate atoms belonging to $U$ or any other subset, then add the atoms on this path to $U$. If this causes the shortest distance between two subsets to be $\leq 4$, then proceed to (\ref{it.func_select2_2}) and repeat it as above.
\item\label{it.func_select2_5} When no scenario in (\ref{it.func_select2_2}) or (\ref{it.func_select2_4}) occurs, we output $S:=\textbf{SELECT2}(A,W,X,Z,Y)$ as the union of all the current sets.
\end{enumerate}
\end{definition}

\begin{proposition}\label{prop.func_select2}  After applying {\bf SELECT2}, the molecule $A$ becomes a proper forest (Definition \eqref{def.ov_proper}) after cutting $S$ as free. Consider all the bonds connecting an atom in $S$ and $X^+(A)$. Then the number of these bonds does not exceed $|S|\leq 2 \mathsf{N}\cdot |X|+ 10d \cdot(|W| + |Y|+|Z|+\rho(A))$;
\end{proposition} 
\begin{proof}
    The proof is essentially the same as Proposition \ref{prop.func_select}. The fact that $\Mb_D$ becomes proper after cutting $S$ as free in both function due to the same reason in the proof of Proposition \ref{prop.func_select}. 
    
    Now we proof the upper bound for the bonds connecting $\Mb_U$ and $S$ in function {\bf SELECT2}.
    Simply note that each single step in Definition \ref{def.func_select2} \eqref{it.func_select2_2} and \eqref{it.func_select2_4} adds at most $3$ atoms to $S$. Each single step in Definition \ref{def.func_select2} \eqref{it.func_select2_2} decrease the number of subsets by $1$, and does not decrease $\rho(S)$. Each single step in Definition \ref{def.func_select2} \eqref{it.func_select2_4} does not change the number of subsets, but \emph{adds one more cycle to $S$}; in particular adding this cycle increasess the value of $\rho(S)$ by at least one. Putting together, we know that the value of
    \[\gamma := (\text{number of sets in } S)-\rho(S)\] decreases by at least $1$ in each step of Definition \ref{def.func_select2} \eqref{it.func_select2_2} and \eqref{it.func_select2_4}. Since we always have $\rho(S)\leq \rho(A)$ because $S\subseteq A$ (again, it is easy to see that adding one atom or arbitrarily many atoms does not decrease the $\rho$ value), we know that the initial value of $\gamma$ is at most $|Y|+|Z|$ and the final value of $\gamma$ is at least $-\rho(A)$. This means that the total number of steps in Definition \ref{def.func_select2} \eqref{it.func_select2_2} and \eqref{it.func_select2_4} is bounded by \[|W|+|X|+|Y|+|Z|+\rho(A).\]
    Moreover, cutting an atom in $Y\cup Z$ breaks at most $2$ the desired bonds, and cutting an atom in $X$ breaks at most $\mathsf{N}$ the desired bonds, cutting an atom in $\Mb_D\setminus (X\cup Y\cup Z)$ breaks at most $2$ the desired bonds. Hence, the total bonds connecting $\Mb_U$ and $S$ is no more than $2 \mathsf{N} \cdot |X| + 10d(|W|+ |Y|+|Z|+\rho(A))$.
\end{proof}
Now we can define the algorithm \textbf{MAINTRUP1} and \textbf{MAINTRUP2}, which are the main algorithms in this section. For each of them, it has three options, \textbf{Option 1}, which is applied only in the (simple) weakly degenerate case, and essentially trivially cuts most atoms in $A$ as \{$3_2$\} molecules, so one can gain from extra restrictions from weak degeneracies. In \textbf{Option 2}, we deal with the case when $|W|$ is large (Definition \ref{def.func_select2} and Definition \ref{def.func_select3}. In \textbf{Option 3}, which is the main part, we basically combine \textbf{Stage 1} (cutting atoms in $\overline{X}^+_0(A)$) and \textbf{Stage 3} (pre-processing) in the proof of Proposition \ref{prop.case5} at the end of Section \ref{sec.maincr}, and the \textbf{MAINUD} algorithm, but wih the roles of $\Mb_U$ and $\Mb_D$ played by $X^+(A)$ and $X$ respectively.
\begin{definition}[Algorithm \textbf{MAINTRUP1}] 
\label{def.alg_maintrup}
Suppose $\Mb$ is a connected full molecule and $A$ is a connected transversal subset. We may choose $(A^+,A^-)$ as in Proposition \ref{prop.trans} \eqref{it.trans1}. Let the number of bonds connecting an atom in $X^+(A)$ and an atom in $A$ be $\#^+_{\mathrm{conn}}$. We always view each box in $A$ as a whole, and treat each atom in box of $A^\pm$ individually. Moreover, by splitting into at most $C^{|A|+\rho(A)}$ sub-cases, we can identify a set of weakly degenerate atoms and atom pairs in $A$ in the sense of Definition \ref{def.weadeg}; let $X$ be the B-atoms in $A$ and let $Y_1\subseteq A \setminus X$ be the atoms involved in these weak degeneracies. We define the following algorithm, which contains three Options that we can choose at the beginning.

In \textbf{Option 1} we do the followings:
\begin{enumerate}
\item Cut $A$ as free, then cut it into elementary components using \textbf{UP1}, by viewing each box as a whole. Then, cut each box by {\bf UP} or {\bf DOWN}, depending on whether the box has a top or bottom fixed end.
\item Then cut $A^+$ as free and cut it into elementary components using \textbf{UP}, then cut $A^-$ into elementary molecules using \textbf{DOWN}.
\end{enumerate}
In \textbf{Option 2} we do the followings:
\begin{enumerate}
    \item Cut $\overline{X}^+_0(A)$ as free, then we cut $X^+(A)$ into elementary components using \textbf{UP}, but only cut one atom each time. Then cut each component in $A+ \backslash X^+(A)$ by \textbf{UP} and each component in $A \cup A^-$ by \textbf{DOWN}.
\end{enumerate}

In \textbf{Option 3} we do the followings:
\begin{enumerate}
\item\label{it.alg_maintrup_1}  Cut all atoms in $\overline{X}^+_0(A)$ as free. If any atom in $A$ becomes fully fixed, also cut it as free until $A$ has no such atom.
\item\label{it.alg_maintrup_2}  Define $S=\textbf{SELECT2}(A,W,X,Z,Y)$ as in Definition \ref{def.func_select2}, then cut $S$ as free and cut it into elementary components using \textbf{DOWN}. 
\item\label{it.alg_maintrup_3} If not all atoms in $X^+(A)$ have been cut, then choose a lowest atom $\nf$ in $X^+(A)$ that has not been cut. If $\nf$ is  a binary atom and is adjacent to an atom $\pf\in A$ that is a $\{3_2\}$ atom, then cut $\{\nf,\pf\}$ as free; if $\nf$ is adjacent to some $\pf_{1},\dots,\pf_j \in A$ that is a $\{3_2\}$ atom, then we cut $\{\nf,\pf_1,\dots,\nf_j\}$ as free, and further cut $\pf_1,\dots,\pf_j$ as free in this order, if $\nf$ will not be fully fixed after the cutting; otherwise cut $\nf$ as free. If $A$ becomes non-proper then repeat the steps in Definition \ref{def.alg_maincr} (\ref{it.alg_maincr_5})--(\ref{it.alg_maincr_6}) until $A$ becomes proper again.
\item Repeat (\ref{it.alg_maintrup_3}) until all atoms in $X^+(A)$ have been cut. Then repeat the step in Definition \ref{def.alg_maincr} (\ref{it.alg_maincr_9}) until all atoms in $A$ have been cut.
\item\label{it.alg_maintrup_5} Finally, cut (the remaining parts of) $A^+$ as free and cut it into elementary components using \textbf{UP}, and then cut $A^-$ into elementary components using \textbf{DOWN}.
\end{enumerate}

We define the dual algorithm \textbf{MAINTRDN1} in the same way (so $\#_{\mathrm{conn}}^+$ is replaced by $\#_{\mathrm{conn}}^-$ etc.).
\end{definition}
\begin{proposition}
\label{prop.alg_maintrup} Assume all two-atom molecules and multi-particle atoms are good. In \textbf{Option 1} of \textbf{MAINTRUP1} in Definition \ref{def.alg_maintrup} (same for \textbf{MAINTRDN1}), we have
\begin{equation}\label{eq.alg_maintrup_1}
\#_{\mathrm{bad}}=1;\qquad \#_{\mathrm{good}}\geq\frac{1}{10}\cdot |Y_1|-\rho(A)-1.
\end{equation}
In \textbf{Option 2} of of \textbf{MAINTRUP1} in Definition \ref{def.alg_maintrup}, we have  \begin{equation}\label{eq.alg_maintrup_W}
\#_{\mathrm{bad}}\leq |\overline{X}^+_0(A)|;\qquad \#_{\mathrm{good}}\geq \frac{|W|}{2d}.
\end{equation}
 In \textbf{Option 3}   of \textbf{MAINTRUP1} in Definition \ref{def.alg_maintrup}, we have  all \{$3_23_2$B\} molecules are good, and moreover
\begin{equation}\label{eq.alg_maintrup_2}
\begin{aligned}
\#_{\mathrm{good}}&\geq\frac{1}{10}\cdot\big(\#_{\mathrm{conn}}^+ - 2 \mathsf{N} \cdot\#_{\mathrm{box}}(A) - 10^4 d^3 (|W| + |Y_1|+\rho(A)+|\overline{X}^+_0(A)| + \#_{\mathrm{mp}}(A))\big),\\
\#_{\mathrm{bad}}&\leq |W| + |Y_1|+\rho(A)+|\overline{X}^+_0(A)| +\#_{\mathrm{box}}(A) + \#_{\mathrm{mp}}(A).
\end{aligned}
\end{equation}
\end{proposition} 
\begin{proof}
The case of \textbf{Option 1} is easy; $\#_{\mathrm{bad}}=1$ because $A$ is connected and $A^+$ satisfies the assumption in Proposition \ref{prop.trans} \eqref{it.trans1}. Moreover, each atom in $A$ belongs to a $\{k_l\}$ molecule, and the total number of fully fixed atoms is no more than $\rho(A)$ (see the proof of Proposition \ref{prop.alg_transup}). For any atom $\nf$ that belongs to a $\{3_2\}$ molecule, if it is a weakly degenerate atom, or if it is part of a weakly degenerate pair and is cut after the other atom $\nf'$ of the pair, then $\{\nf\}$ must be a good component by Definition \ref{def.good_normal} \eqref{it.good_1}. This gives at least
\[
|Y_1|/10-\rho(A)-1
\]
good components, as desired. 

In the case of \textbf{Option 2}, denote the set of atoms in $\Mb_U$ that connects at least one atom in $W$ by $W'$. Note that $|W'|\leq |W|\leq (d+1)|W'|$ Note that each atom in $W'$ must be a multi-particle atom, as it does not belong to $\overline{X}^+_0(A)$ but has at least two bonds connecting $A$. Moreover, each atom in $X^+(A) \backslash A$ must be an ancestor of some atom in $\overline{X}^+_0(A)$, but the definition of $X^+(A)$. Therefore, each atom in $W'$ has only bottom fixed end when it is cut (as it has descendants with bottom fixed ends). Note that only those atoms in $\overline{X}^+_0(A)$ are cut as bad molecules as $X^+(A)$ is connected. Moreover, each atom in $W'$ cannot be cut as fully fixed as it has descendant in $A$, which is cut after all atoms in $X^+(A)$ have been cut, so every atom in $W'$ is cut as a good molecule.

In the case of \textbf{Option 3}, note that in the whole process, there is no top fixed end in $A^+$ and no bottom fixed end in $A^-$ by Definition \ref{def.trans} (each atom in $X^+(A)$ has descendants in $X^+_0(A)$, which have been cut). Moreover, any binary atom $\nf$ chosen in this step either has no bond connecting to $A$, or has exactly one bond connecting to $A$ (because all atoms in $\overline{X}^+_0(A)$ have been cut in Definition \ref{def.alg_maintrup} \eqref{it.alg_maintrup_1}). Also any \{$3_2 3_2$\} molecule $\{\nf,\mathfrak p\}$ cut in this way must be a \{$3_2 3_2$A\} molecule, and no full component can be cut (hence no bad molecule created) in Definition \ref{def.alg_maintrup} \eqref{it.alg_maintrup_5}. Therefore, all \{$3_2 3_2$B\} molecule must have both atoms in $A$, and the only bad molecules created are those created in Definition \ref{def.alg_maintrup} \eqref{it.alg_maintrup_1}(contributing at most $|\overline{X}^+_0(A)| $ many bad molecules) and Definition \ref{def.alg_maintrup} \eqref{it.alg_maintrup_2} (contributing at most $\#_{\mathrm{box}}(A)+\rho(A) + \#_{\mathrm{mp}}(A) + |Y_1|$ many bad molecules), hence the upper bound
\[
\#_{\mathrm{bad}}\le \#_{\mathrm{box}}(A) + \#_{\mathrm{mp}}(A) + |Y_1|+\rho(A)+|\overline{X}^+_0(A)|.
\]

Now we prove the lower bound on $\#_{\mathrm{good}}$. First, we note that each multi-particle atom $\nf$ cut in \eqref{it.alg_maintrup_3} together with an atom in $A$ must belong to a good molecule, and $\nf$ is cut together with at most $d$ $\{3_2\}$ atoms in $A$ (as $\nf \notin X^+_0(A)$). Note that after Definition \ref{def.alg_maintrup} \eqref{it.alg_maintrup_1}, there is no bottom fixed end in $A$ (and hence none in $S$), so we can cut $S$ into elementary components using \textbf{DOWN} as in Definition  \ref{def.alg_maintrup} \eqref{it.alg_maintrup_2}. By Proposition \ref{prop.func_select2}, we know that after Definition \ref{def.alg_maintrup} \eqref{it.alg_maintrup_2} is finished, $A$ will become a forest which is proper, and contains no weakly degenerate atoms or atom pairs.

Next we prove an upper bound on the bonds connecting $S$ and $X^+(A)$; cutting an atom in $\overline{X}^+_0(A)$ breaks at most $d+1$ such bonds. Let $Z_0$ be the set of atoms in $A$ that is cut in Definition \ref{def.alg_maintrup} \eqref{it.alg_maintrup_1}. After this step, denote atoms with fixed ends in $A$ by $Z_1$. Then
\[
|Z_1|\le (d+1) \bigl(|\overline{X}^+_0(A)| + Z_0\bigr)
\]
(because the atoms with fixed ends in $Z$ must be child to one of the atoms in $\overline{X}^+_0(A)$ or $Z_0$ that has been cut). By Definition \ref{def.alg_maintrup} \eqref{it.alg_maintrup_1}, we know that each atom with $2 \leq l\leq d+1$ particles in $Z_0$ must have $l$ parents with each of them in either $Z_0$ or $\overline{X}^+_0(A)$. This then leads to
\[
\rho(Z_0)\ge |Z_0|-(d+1)|\overline{X}^+_0(A)|
\]
(by counting the total number of top bonds at all atoms of $Z_0$, note that each bond is counted only once, and the number of such bonds that do not connect two atoms in $Z_0$ is at most $d|\overline{X}^+_0(A)|$). Since also $\rho(Z_0)\le \rho(A)$ (because $Z_0\subseteq A$), we know
\[
|Z_0|\le \rho(A)+ (d+1) |\overline{X}^+_0(A)|.
\]
Thus, $|Z_1| \leq 2d^2 \overline{X}^+_0(A)+ 2 d \rho(A)$ and hence $|Z| \leq 2d^2 \overline{X}^+_0(A) + |Y_1 + \#_{\mathrm{mp}}(A)|$. By Proposition \ref{prop.func_select2}, the number of bonds connecting $S$ and $X^+(A)$ is no more than 
\[
  2 \mathsf{N} \cdot \#_{\mathrm{box}}(A) + (10d)^3  \cdot \bigl(\rho(A)+|Y_1|+|\overline{X}^+_0(A)| + \#_{\mathrm{mp}}(A)\bigr).
\]
Now, after cutting $\overline{X}^+_0(A)\cup Z_0\cup S$ as free, the number of bonds connecting an atom in $A$ to an atom in $X^+(A)$ is still at least
\[
\begin{aligned}
    (\#_{\mathrm{conn}}^+)' \geq & \#_{\mathrm{conn}}^+
    - 2 \mathsf{N} \cdot \#_{\mathrm{box}}(A)  -  10^4 d^3  \cdot \bigl(\rho(A)+|Y_1|+|\overline{X}^+_0(A)| + \#_{\mathrm{mp}}(A) \bigr)  
    \end{aligned}
\]

At this point we can apply the same arguments as in the proof of Proposition \ref{prop.alg_maincr} (especially using Part 2 of the proof and \eqref{eq.alg_maincr_2}) to show that
\[ \#_{\{3_23_2A\}}+\#_{\{3_23_2B\}} \ge \frac{1}{10} \big( (\#_{\mathrm{conn}}^+)' - 2 \mathsf{N}\cdot \#_{\mathrm{box}}(A) - 10^4 d^3\bigl(\rho(A)+|Y_1|+|\overline{X}^+_0(A)|\bigr) \big). \]
Here the role of UD connections is played by the $(\#_{\mathrm{conn}}^+)'$ bonds connecting an atom in $A$ to an atom in $X^+(A)$. Note that there exists no M-atom in $X^+(A)$ has a bond to atom in $A$, thanks to our choice of $X(A)$ by viewing each box as a whole. This completes the proof.
\end{proof}

Now, we can finally prove Proposition \ref{prop.comb_est_extra}.
\begin{proof}[Proof of Proposition \ref{prop.comb_est_extra}] {\bf Proof part 1.} First assume $\Mb$ contains at least one strongly degenerate primitive pair $(\mf,\nf)$, then we can cut $\{\mf,\nf\}$ as free and cut the rest of $\Mb$ into elementary components using either \textbf{UP} or \textbf{DOWN} (cf. proof of Proposition \ref{prop.case2}, \textbf{Choice 1}). This creates one good \{$4_2 (2l)_l$\} molecule and no bad molecule (as $\Mb$ is connected), so $\#_{\mathrm{good}}\geq 1$ and $\#_{\mathrm{bad}}=0$, which implies (\ref{eq.extra_case_1}). Moreover, if $\Mb$ contains any multi-particle atom with $|v_{e_i} - v_{e_j} |\leq \eps^{1 - \upsilon}$, where $e_i,e_j$ are both top/bottom edges at the multi-particle atom, then we cut it as free and cut the rest of $\Mb$ by using either \textbf{UP} or \textbf{DOWN} and (\ref{eq.extra_case_1}) also holds. Below we assume each multi-particle atom is good, by Proposition \ref{prop.est_single} \eqref{it.est_kl2}.
    
Now assume $\Mb$ contains no strongly degenerate primitive pair, then we can show (by splitting into two sub-cases) that each \{$3_2 3_2A$\} molecule is either good or can be cut into one \{$2_2$\} and one good \{$3_2$\}-atom. The proof is similar to Lemma \ref{lem.alg_up_ex}; suppose this \{$3_2 3_2$A\} molecule is $\{\nf,\nf'\}$ where $\nf$ is parent of $\nf'$, then we may assume this is strongly degenerate and $|t_\nf-t_{\nf'}|\leq \varepsilon^{\upsilon}$ otherwise this is already good by Definition \ref{def.good_normal}. Now our assumption implies that $\{\nf,\nf'\}$ must be non-primitive; in the absence of O-atoms, this implies that there exists another chain of atoms going from $\nf$ to $\nf'$ by iteratively taking children. If any atom $\qf$ in this chain is cut before $\{\nf,\nf'\}$, then the inequality $t_\nf\geq t_{\qf}\geq t_\nf'$ allows us to argue as in the proof of Lemma \ref{lem.alg_up_ex} to get a good molecule; if not, this means that at the time $\{\nf,\nf'\}$ is cut, $\nf$ still has two children (i.e. $\nf'$ and a child on the chain), and $\nf'$ still has two parents (i.e. $\nf$ and a parent on the chain), which means that $\nf$ must have a top fixed end and $\nf'$ must have a bottom fixed end, but this contradicts that $\{\nf,\nf'\}$ is \{$3_2 3_2$A\} molecule. Moreover, each \{$3_2 k_l$\} molecule can be viewed as good, by the same argument.

Therefore, from now on, we will treat \{$3_2 3_2A$\} molecules and \{$3_2 k_l$\} molecules as good. We start by choosing a connected transversal subset $A_1$ of $\Mb$ and either $A_1=\Mb$ or $X_0(A_1)\neq\varnothing$; the existence of such $A_1$ follows in the same way as in Proposition \ref{prop.trans} (\ref{it.trans4}) starting from a single atom with no children.  If $A_1\neq\Mb$, let $X_1=X(A_1)$, then $A_1\cup X_1$ is connected transversal by Proposition \ref{prop.trans} (\ref{it.trans3}), so we can find $A_2\supseteq A_1\cup X_1$ starting from $A_1\cup X_1$, using Proposition \ref{prop.trans} (\ref{it.trans4}). If $A_2\neq \Mb$, then let $X_2=X(A_2)$ and find $A_3$ by Proposition \ref{prop.trans} (\ref{it.trans4}), then define $X_3=X(A_3)$ and so on.

Recall the function $G(q)=(2q+1)\cdot \binom{q}{2}\cdot(32q^{3/2})^{q^2}$. Define $K := (2d+1)/\upsilon$, $D:= 8K^{2K+3}$, and define $N_1 = D^D$, $N_{j+1} = (60d N_j)^{60d}$ for $j\geq 1$, and $\mathsf{N} = N_D$. We also define $K_1:= D\cdot G(\mathsf{N})$, and $K_{j+1}=\big(G(K_j)\cdot G(\mathsf{N})\big)^{60d}$ for $j\geq 1$. We will assume $\Gamma>(60d)^{60d} K_D$, and consider several cases.

{\bf Proof part 2.} Suppose $A_j$ exists and $X_0(A_j)\neq\varnothing$ for all $j\leq D$. By construction in Proposition \ref{prop.trans} (\ref{it.trans4}), we can list the elements of $B_j:=A_j\backslash (A_{j-1}\cup X_{j-1})$ as $\nf_1^j,\cdots ,\nf_{|B_j|}^j$ such that $\nf_1^j$ has at most $l_{\nf_1^j}-1$ (here and below $l_\nf$ denotes the number of particles of an atom $\nf$, note that $\nf$ may be a box, see Definition \ref{def.trans}) bonds and at least one bond with $A_{j-1}\cup X_{j-1}$ and $\nf_i^j$ has at most $l_{\nf_i^j}-1$ bonds and at least one bond with with $A_{j-1}\cup X_{j-1}\cup\{\nf_1^j,\cdots,\nf_{i-1}^j\}$. Similarly we can list the atoms in $B_1:=A_1$. Choose also an element $\mf_j\in X_0(A_j)$ for each $j$, then this $\mf_j$ is adjacent to $l_{\mf_j}$ atoms in $A_j$, with at least one being in $B_j$. We then perform the following cutting sequence:

Starting from $j=1$, for each $j$, cut the atoms $(\nf_i^j:i\leq |B_j|)$ as free in the increasing order of $i$. However, if $\nf_i^j$ is a \{$3_2$\} atom, and is adjacent to $\mf_{j}$ which becomes a fully fixed atom after cutting $\nf^j_i$, then we cut $\{\nf_i^j,\mf_{j}\}$ instead of $\nf_i^j$. After all atoms in $B_j$ have been cut, we then cut the remaining atoms in $X_j$ as free (starting from the lowest ones; each will be fully fixed when we cut it), and proceed with $B_{j+1}=A_{j+1}\backslash (A_j\cup X_j)$ and so on.

In this way, each atom in $B_j$ cannot be fully fixed when it is cut. It is clear that for each $\mf_j$, if it is cut as a fully fixed atom, then $\nf^j_i$ cut before $\mf_j$ must be a multi-particle atom or a box, and it cannot be a fully fixed atom. So $\nf^j_i$ is a good molecule (if $\nf_j$ is a box, we further cut it by {\bf UP1} or {\bf DOWN1}, which also gives a good molecule by the same argument as in the proof of Proposition \ref{prop.alg_up} \eqref{it.up_proof_3}). Otherwise, $\mf_j$ must be cut with some $\nf_i^j$ which is a good molecule, because $\mf_j$ will have no top (or bottom) fixed ends when it is cut, assuming $\mf_j\in X_0^+(A_j)$ (or $\mf_j\in X_0^-(A_j)$). This produces at least $D$ many good molecules, while $\#_{\mathrm{bad}}=1$, which implies (\ref{eq.extra_case_1}).

{\bf Proof part 3.} Suppose $A_j=\Mb$ for some $j<D$, then $\rho(A_j)\geq\Gamma$ for this $j$. If $\rho(A_1) \geq K_1$, by construction $A_1$ is constructed from a single node by adding atoms one by one. We can list the elements in $A_1$ in the order that they are added in $A_1$ as $\nf_1,\cdots, \nf_j$. Note that $\Gamma <\rho(A_1) \leq \sum_{i}(l_{\nf_i} -2 )\leq j \cdot \mathsf{N}$ (as we view each box as a whole when we construct in Proposition \ref{prop.trans}), the number of multi-particle atoms and boxes is at least $\Gamma /\mathsf{N}$. By the same argument in (1), starting from $\nf_1$ and cut each atom in order, each multi-particle atom or box can provide a good molecule. This meets the requirement in Proposition \ref{prop.comb_est_extra}.

Now we assume $\rho(A_1) \leq K_1$, and we choose the smallest $j$ such that $\rho(A_{j+1})\geq K_{j+1}$. Let $A=A_j$, then we have
\[\rho(A)<K_j,\quad K_{j+1}\leq \rho(A_{j+1}),\]

\emph{Case 1.} We first assume that $\#_{\mathrm{conn}}^+ \geq (60d)^{60d}\cdot K_j^2$, we then
    \begin{enumerate}
        \item Run algorithm {\bf UP1} to $A$ if $\#_{\mathrm{box}}(A) \geq (10d)^{10d} \cdot K_j$, or the number of multi-particle atom $\#_{\mathrm{mp}}(A)\geq (10d)^{10d} \cdot K_j$, then cut $A^+$ by {\bf UP}, and then cut $A^-$ by {\bf DOWN};
        \item Run algorithm {\bf TRANSUP} and apply Proposition \ref{prop.alg_transup} if $$\#_{\mathrm{box}}(A)\  \mathrm{and}\  \#_{\mathrm{mp}}(A)\leq (10d)^{10d} \cdot K_j \quad \text{and}\quad\! \overline{X}^+_0(A) \geq (20d)^{20d} \cdot K_j^{4/3};$$
        \item Run algorithm {\bf MAINTRUP1, Option 1} and apply Proposition \ref{prop.alg_maintrup} if $$\#_{\mathrm{box}}(A)\  \mathrm{and}\  \#_{\mathrm{mp}}(A)\leq (10d)^{10d} \cdot K_j \quad \text{and}\quad\! \overline{X}^+_0(A) \leq (20d)^{20d} \cdot K_j^{4/3};$$
        and $|Y_1| \geq (20d)^{20d} \cdot K_j$;
        \item Run algorithm {\bf MAINTRUP1, Option 2} and apply Proposition \ref{prop.alg_maintrup} if $$\#_{\mathrm{box}}(A)\  \mathrm{and}\  \#_{\mathrm{mp}}(A)\leq (10d)^{10d} \cdot K_j \quad \text{and}\quad\! \overline{X}^+_0(A) \leq (20d)^{20d} \cdot K_j^{4/3};$$
        and $|W| \geq (30d)^{30d} \cdot K_j^{5/3}$;
        \item Run algorithm {\bf MAINTRUP1, Option 3} and apply Proposition \ref{prop.alg_maintrup} if $$\#_{\mathrm{box}}(A)\  \mathrm{and}\  \#_{\mathrm{mp}}(A)\leq (10d)^{10d} \cdot K_j \quad \text{and}\quad\! \overline{X}^+_0(A) \leq (20d)^{20d} \cdot K_j^{5/3};$$
         and $|Y_1| \leq (20d)^{20d} \cdot K_j$, $|W| \leq (30d)^{30d} \cdot K_j^{5/3}$. 
    \end{enumerate}
In any case it is easy to verify that \eqref{eq.extra_case_1} is true. If $\#_{\mathrm{conn}}^- \geq (60d)^{60d}\cdot K_j^2$, we can also argue as above and apply the corresponding dual algorithm.

\emph{Case 2.} We then assume that $\#_{\mathrm{conn}}^-$ and $\#_{\mathrm{conn}}^+ \leq (60d)^{60d}\cdot K_j^2$. We first assume there exists at least $\upsilon^{-1} d^2+1$ multi-particle atoms or boxes (note that in Proposition \ref{prop.trans} we view each box as a whole) in $A_{j+1} \setminus (A_j \cup X(A_j))$. We cut $A_{j+1}$ as free and then apply $\textbf{UP1}$ to $A_{j+1} \setminus (A_j \cup X(A_j))$. By Proposition \ref{prop.trans}, there will be no fully fixed atom in $A_{j+1} \setminus (A_j \cup X(A_j))$. For each multi-particle, either it is cut as good molecule, or the atom cut before it is a good molecule with excess $\eps^{\upsilon_0/2}$. As a result, we gain enough good molecule, with $\#_{\rm bad} =1$. After, we cut the rest components by {\bf UP} or {\bf DOWN}.

Then we have $\rho(A_{j+1}) \leq  \rho(X(A_j)\cup A_j) + d^2 \upsilon^{-1}\mathsf{N}$ (note that the quantity $\rho$ does not count the cycles inside the boxes)
\[ \rho(A_j) + \rho(X^+(A)) +\rho(X^-(A)) + \#_{\mathrm{conn}}^-+\#_{\mathrm{conn}}^+ \geq \rho(A_{j+1}) - d^2 \upsilon^{-1}\mathsf{N} \geq K_{j+1}/2   \]
which implies that 
\[\rho(X^+(A)) +\rho(X^-(A)) \geq  (G(K_j)\cdot G(\mathsf{N}))^{30d}.\]
 We now claim there must exists a component in $X^+(A)$ or $X^-(A)$ with at least $\mathsf{N}$ particles. Indeed, assume that each component in $X(A)$ has less than $\mathsf{N}$ particles. By our construction of box, for a component with $q\leq \mathsf{N}$ particles, there exists no more than $G(q)$ atoms (by viewing box as a whole). Therefore, the number of atoms in $X(A)$ is no more than $2 \cdot \frac{G(\mathsf{N})}{\mathsf{N}} \cdot \big((60d)^{60d}\cdot  K_j\big)$ atoms, which is less than $\rho(X(A))$, contradicts the assumption of $X(A)$. We may assume this component $\Xc \in X^+(A)$. The rest of the proof of Case 2 is given in Section \ref{subsec_MUMD} below, which completes the proof of Proposition \ref{prop.comb_est_extra}.
\end{proof}

\subsection{Proof of Case 2 in Proof part 3 of Proposition \ref{prop.comb_est_extra}}\label{subsec_MUMD}
In this subsection we finish the proof of Case 2 in Proof part 3 of Proposition \ref{prop.comb_est_extra}. We use the notations $A$, $\Xc$ there, where $A$ is a transversal set and $\Xc$ is a component in $ X^+(A)$ with at least $\mathsf{N}$ particles. Recall the constants $K = (2d+1)/\upsilon$, $D= 8K^{2K+3}$, $N_1 = D^D$, $N_{j+1} = (60d N_j)^{60d}$ for $j\geq 1$, and $\mathsf{N} = N_D$. We start with the definition of canonical cycle.
 
\begin{definition}[Canonical cycle]\label{def.canonical_cycle}
     Let $\Mb$ be a C-molecule, we construct the canonical cycle as follows. We first view each box as a whole.
     
     Starting from a smallest cycle of $\Mb$, we view it as a whole, namely as a single node. If this new node is contained in another cycle, then we choose the smallest cycle containing this node, and again view the whole cycle as a whole. We repeat this procedure until the resulting node is not contained in any cycle. The set of all atoms in $\Mb$ contained in the resulting node is called a \textbf{canonical cycle}.
     
     We then start from another smallest cycle in $\Mb$ which is not contained in any previously constructed canonical cycle, and repeat the same procedure. Finally, if a box is not contained in any canonical cycle constructed above, we also \emph{regard this box as a canonical cycle}. 
\end{definition}
\begin{proposition}\label{prop.canonical_cycle}
   The canonical cycles constructed in Definition \ref{def.canonical_cycle} are independent of the choices of the smallest cycles. After viewing each canonical cycle as a whole, the resultant graph is a tree.
\end{proposition}
\begin{proof}
We first prove that the construction is independent of the choices made during the procedure. Starting from a smallest cycle, the construction repeatedly adds cycles containing the node already constructed. Consequently, the resulting canonical cycle contains every cycle that can be connected to the initial cycle through a finite sequence of cycles, where two consecutive cycles have at least one common atom.

Indeed, suppose that the procedure terminates before one of these cycles has been included. Consider the first cycle in such a sequence which has not been included. The preceding cycle has already been viewed as part of the constructed node, so the next cycle contains this node. This contradicts the stopping condition. Conversely, every cycle added during the construction contains the node constructed in the preceding step. Hence it is connected to the initial cycle through such a sequence of cycles. Therefore, the final canonical cycle depends only on the initial collection of connected cycles, and not on the choices made during the construction.

Moreover, starting from any other smallest cycle contained in the same collection produces the same canonical cycle. Thus, two canonical cycles constructed by Definition \ref{def.canonical_cycle} are either identical or disjoint.

We now view every canonical cycle as a whole. Suppose that the resultant graph contains a cycle. This cycle either contains a node corresponding to a canonical cycle or consists only of atoms not contained in any canonical cycle. In the first case, the construction of canonical cycle would not have terminated, since the corresponding node is still contained in a cycle. In the second case, this cycle would have been used to construct another canonical cycle. Both cases give a contradiction. Therefore, the resultant graph contains no cycle. Since $\Mb$ is connected, the resultant graph is a tree.
\end{proof}
  
We construct the following $\Mb^j_D$ inductively.
\begin{definition}[Construction of $\Mb^j_D$]\label{def.construct_MD}
{\bf Construction of $\Mb^{0}_D$.} With $A$ and $\Xc$ defined above, choose the highest connected subset of $A$ and denote it by $\Ac_0$, such that it contains every particle in $A$ and there exists no atom in $A\backslash \Ac_0$ is a parent of some atom in $\Ac_0$ (we call $\Ac_0$ the \emph{highest part} of $A$). Define $\Mb_D^0$ by $\Ac_0$ and here and below we always define $\Mb^j_U$ by $\Xc$ for $j\geq 0$.

{\bf Construction of $\Mb^{j+1}_D$ from $\Mb^{j}_D$.} Suppose $\Mb^j_D$ contains at most $K$ (recall that $K = (2d+1)/\upsilon$) canonical cycles, and there exists a canonical cycle in $\Mb^j_D$ with at least $N_{D-j-1}$ bonds connecting $\Mb_U$. We view each canonical cycle in $\Mb^j_D$ as a single node, and also view each atom which does not belong to any canonical cycle as a node. By Proposition \ref{prop.canonical_cycle}, the resultant graph is a tree. We say a node is a parent of another node, if the node has an atom with a child in another node.

\emph{Step 1.} Choose the smallest subtree which contains all canonical cycles in $\Mb^j_D$ having a bond with $\Mb_U$, and add all ancestors in $\Mb^j_D$ of its nodes. Denote the set of atoms in the resulting nodes by $\Ac_{j+1}$, and let $S:=\Ac_{j+1}$.

\emph{Step 2.} Consider the lowest nodes of $\Ac_{j+1}$. We say that $\Ac_{j+1}$ is of type \eqref{it.type1} if at least one of its lowest nodes is a canonical cycle; otherwise, we say it is of type \eqref{it.type2}. We apply one of the following depending on the type of $\Ac_{j+1}$.
\begin{enumerate}[(1)]
    \item\label{it.type1} If at least one lowest node of $\Ac_{j+1}$ is a canonical cycle, then for each such lowest canonical cycle $\Cc$, we repeat the following steps until none of them applies:
    \begin{enumerate}[(a)]
        \item\label{it.defMj_a} if an atom (or a box) $\nf\in\Cc$ is a lowest atom of $S$, and $S\setminus\{\nf\}$ is still connected and contains every particle line in $S$, then we replace $S$ by $S\setminus \{\nf\}$;
        \item\label{it.defMj_b} otherwise, if a box $B\subseteq\Cc$ is a lowest atom of $S$, then we replace $B$ by its top part in $S$ (i.e. the M-atom and the bottom part of $B$ are removed). From now on we view the atoms in this top part individually;
        \item\label{it.defMj_c} otherwise, if a C-atom $\qf\in\Cc$ is a lowest atom of $S$ and is contained in a cycle in $S$, and \eqref{it.defMj_c} has not been used for $\Cc$, then we replace $S$ by the component of $S\setminus\{\qf\}$ with the most bonds connecting $\Mb_U$.
    \end{enumerate}
    After each step we go back to \eqref{it.defMj_a}. In particular, after \eqref{it.defMj_c} we continue with \eqref{it.defMj_a} and \eqref{it.defMj_b}, and they may again be applied to any canonical cycle in this set. Note that \eqref{it.defMj_c} is used at most once for each $\Cc$, also when the remaining part of $\Cc$ contains smaller canonical cycles; these smaller canonical cycles are only treated in later rounds.
    \item\label{it.type2} Otherwise, each lowest node of $\Ac_{j+1}$ is a single C-atom which is not contained in any cycle (a box not contained in any cycle is a canonical cycle). We choose one such atom $\pf$, remove it, and replace $S$ by the component of $S\setminus\{\pf\}$ with the most bonds connecting $\Mb_U$.
\end{enumerate}
 The components of $S\setminus\{\qf\}$ in \eqref{it.defMj_c}, or of $S\setminus\{\pf\}$ in type \eqref{it.type2}, which are not chosen are called \emph{dropped pieces}.

\emph{Step 3.} Define $\Mb^{j+1}_D:=S$ and $R_j = \Mb^{j-1}_D\setminus \Mb^{j}_D$.
\end{definition}
 
In Step 1 we keep all canonical cycles having bonds with $\Mb_U$, the paths between them and their ancestors, so that $\Mb^{j+1}_D$ has no parent in $\Mb^j_D \setminus\Mb^{j+1}_D$. In Step 2 \eqref{it.defMj_a} and \eqref{it.defMj_b} is exactly choosing the highest part of the lowest canonical cycle. Step \eqref{it.defMj_c} is used only when the lowest atom is a multi-particle atom.

\begin{proposition}\label{prop_MUMD}
     Let $\Mb_U = \Xc$ and $\Mb^j_D$ ($j\geq0$) be constructed as above. Then they satisfy the following properties:
    \begin{enumerate}[(i)]
        \item\label{it.prop_MUMD2} each child of an atom in $\Xc$ which is not in $\Xc$ belongs to $\Ac_0$, and it is the highest atom of its particle line in $A$;
        \item\label{it.prop_MUMD1} $\Mb_U$ and $\Mb^j_D$ are connected, and no atom in $A\setminus\Mb^j_D$ is a parent of some atom in $\Mb^j_D$. In particular, no atom outside $\Mb^j_D$ has both ancestor and descendant in $\Mb^j_D$; 
        \item\label{it.prop_MUMD3} no atom in $\Mb_U$ is a descendant of some atom in $\Mb^j_D$, and there exists no M-atom in $\Mb_U$ adjacent to an atom in $\Mb^j_D$;
        \item\label{it.prop_MUMD4} the number of bonds between $\Mb_U$ and $\Mb^0_D$ is at least $\mathsf{N}$, and for $j\geq1$ the number of bonds between $\Mb_U$ and $\Mb^j_D$ is at least $N_{D-j}/(2(d+1))^{K+1}$;
        \item\label{it.prop_MUMD5} for $j\geq1$, each component of $R_j$ is adjacent to $\Mb^j_D$, and each bond between them goes from an atom in $\Mb^j_D$ to its child in $R_j$;
        \item\label{it.prop_MUMD6} for $j\geq1$, $R_j$ consists of the following atoms:
        \begin{enumerate}[(a)]
            \item the components of $\Mb^{j-1}_D\setminus\Ac_j$. They contain no canonical cycle of $\Mb^{j-1}_D$ having a bond with $\Mb_U$;
            \item the atoms removed in \eqref{it.defMj_a}, and the M-atoms and bottom parts removed in \eqref{it.defMj_b}. None of them is a child of an atom in $\Mb_U$;
            \item the atoms $\qf$ removed in \eqref{it.defMj_c} and $\pf$ removed in \eqref{it.type2}. Each of them is a C-atom with at most $d+1$ parents, and there are at most $K$ of them;
            \item the dropped pieces. A dropped piece has no bond with $\Mb^j_D$ or with a dropped piece created after it, and each parent of an atom in a dropped piece is either in this piece, or not in $\Ac_0$.
        \end{enumerate}
    \end{enumerate}
    Hence, $\Mb^j_D\cup\Mb_U$ is a UD-molecule.
\end{proposition}
\begin{proof}

    Note that there exists no fixed end in $A\cup X(A)$. We first prove \eqref{it.prop_MUMD2}. Since $\Xc$ is a component in $X^+(A)$, so its descendant outside $\Xc$ must belong to $A$. As $\Ac_0$ is the highest part of $A$ with all particle lines in $A$, each child in $\Xc$ outside $\Xc$ must belong to $\Ac_0$. 

    Now we prove \eqref{it.prop_MUMD1} by induction. $\Mb_U=\Xc$ is connected, as it is a component of $X^+(A)$. For $j=0$, \eqref{it.prop_MUMD1} follows from the choice of $\Ac_0$. Assume it is true for $j$. In Step 1, $\Ac_{j+1}$ is connected, and it contains all parent nodes of its nodes, so no atom in $\Mb^j_D\setminus\Ac_{j+1}$ is a parent of an atom in $\Ac_{j+1}$. In \eqref{it.defMj_a}, we remove a lowest atom of $S$, which keeps this property, and $S$ is still connected. In \eqref{it.defMj_b}, as $B$ is a lowest atom of $S$, each bond between $B$ and $S\setminus B$ goes to the top part of $B$; moreover the top part is connected and contains all particle lines of $B$ (Definition \ref{def.box}). So this property is also kept. In \eqref{it.defMj_c} and \eqref{it.type2}, we remove a lowest atom $\qf$ (or $\pf$) of $S$ and the chosen component is also connected; for an atom in the chosen component, its parents are in $S$ (by the choice of $\Ac_{j+1}$), and there exists no atom in $A \setminus S$ has child in $S$. In particular, for an atom in $(A\setminus S)\cup A^-$, it has only ancestors in $S$ and for an atom in $A^+$, it has only descendants in $S$. This property follows for each $\Mb^j_D$. 

    For \eqref{it.prop_MUMD3}, as $\Mb^j_D\subseteq A$ and $\Xc\subseteq A^+$, no atom in $\Xc$ is a descendant of an atom in $\Mb^j_D$. The statement for M-atoms follows from Definition \ref{def.trans}, as we view each box in $A$ and $X(A)$ as a whole, and an M-atom is only adjacent to atoms in the same box. Moreover, \eqref{it.prop_MUMD5} follows from \eqref{it.prop_MUMD1} and the connectedness of $\Mb^{j-1}_D$.

    Now we prove \eqref{it.prop_MUMD6}. (a) is clear since all canonical cycles having a bond with $\Mb_U$ are in $\Ac_j$. For (b), let $\nf\in S$ be a child of an atom in $\Xc$. By \eqref{it.prop_MUMD2}, $\nf$ is the highest atom of its particle line in $A$. If $\nf$ is a lowest atom of $S$, then it is the only atom of its particle line in $S$; otherwise, the atoms of this particle line between $\nf$ and another atom of it in $S$ are ancestors of an atom in $S$, so they are in $S$ by the proof of \eqref{it.prop_MUMD1}, and $\nf$ has a child in $S$. So removing $\nf$ loses a particle line, and $\nf$ is not removed in \eqref{it.defMj_a}. Moreover, the bonds from atoms outside a box to this box go to its top part, which is kept in \eqref{it.defMj_b}. (c) is clear, since each C-atom has at most $d+1$ particle lines, \eqref{it.defMj_c} is used at most once for each lowest canonical cycle, and in \eqref{it.type2} we only remove one atom. For (d), let $P$ be a dropped piece, which is a component of $S\setminus\{\qf\}$ (or $S\setminus\{\pf\}$) that is not chosen. The chosen component contains $\Mb^j_D$ and all dropped pieces created after $P$, and it has no bond with $P$. A parent of an atom in $P$ which is in $\Mb^{j-1}_D$ is in $S$, adjacent to this atom and different from $\qf$ (as $\qf$ is a lowest atom of $S$), so it is in $P$; and by \eqref{it.prop_MUMD1} there is no such parent in $\Ac_0\setminus\Mb^{j-1}_D$.

    Finally we prove \eqref{it.prop_MUMD4}. For $j=0$, each particle line of $\Xc$ is a particle line of $A$, so the lowest atom of $\Xc$ on it has a child in $A$, which is in $\Ac_0$ by \eqref{it.prop_MUMD2}. Different particle lines give different bonds, and $\Xc$ has at least $\mathsf{N}$ particles. Now consider the construction of $\Mb^{j+1}_D$. At the beginning, $S$ contains a canonical cycle with at least $N_{D-j-1}$ bonds connecting $\Mb_U$. By the proof of (b) in \eqref{it.prop_MUMD6}, the atoms at the ends of these bonds are not removed in \eqref{it.defMj_a} and \eqref{it.defMj_b}. In each Step \eqref{it.defMj_c}, or in \eqref{it.type2}, the removed atom has at most $d+1$ bonds connecting $\Mb_U$, and $S\setminus\{\qf\}$ has at most $d+1$ components, as each of them contains a parent of $\qf$. So the number of bonds connecting $\Mb_U$ and $S$ decreases at most by a factor $2(d+1)$. As there are at most $K$ such steps (the number of canonical cycles is no more than $K$) for each $j$, we get \eqref{it.prop_MUMD4}.
\end{proof}

When part of $\Mb$ is cut, we will cut the rest of $\Mb$ by the following lemma.
\begin{lemma}\label{lem.case2_rest}
    Let $A$ be a transversal set in a connected C-molecule $\Mb$. Suppose the set of atoms that have been cut is $E\cup F$, where $E\subseteq A$ such that no atom in $\Mb \setminus E$ has both ancestor and descendant in $E$, and $F\subseteq A^+$ such that each child of an atom in $F$ which has not been cut is in $A$ and is not an ancestor of an atom in $E$. We first cut the ancestors of atoms in $E$ that belong to $A$ by \textbf{UP}, then cut the rest of $A^+$ by \textbf{UP}, and finally cut all the rest by \textbf{DOWN}. Then the regularity is kept, and there is no bad molecule in this process (if $E\cup F\neq\varnothing$).
\end{lemma}
\begin{proof}
    An ancestor of an atom in $E$ which belongs to $A$ and is adjacent to an atom in $E$ must be a parent of it. By assumption it is not a child of an atom in $F$, and no atom in $A$ is a child of an atom in $A^+\setminus F$ that has been cut. So in the first step there are only bottom fixed ends. After that, an atom in $A^+$ adjacent to an atom that has been cut must be a parent of it, since no atom in $A$ is a parent of an atom in $A^+$, and the children of atoms in $F$ in $A^+$ have been cut. In the last step, an atom in $A$ which has not been cut is not a parent of an atom that has been cut (otherwise it is an ancestor of an atom in $E$, or a parent of an atom in $A^+$), and no atom in $A^-$ is a parent of an atom in $A\cup A^+$. So there are only top fixed ends. As $\Mb$ is connected, each component that we cut is adjacent to an atom that has been cut, so there is no bad molecule.
\end{proof}

\begin{lemma}\label{lem.case2_manycc}
    If $\Mb^j_D$ contains at least $K+1$ canonical cycles, then there exists a cutting sequence with $\#_{\rm good}\geq K$ and $\#_{\rm bad}=1$.
\end{lemma}
\begin{proof}
    We view each canonical cycle in $\Mb^j_D$ as a single node. By Proposition \ref{prop.canonical_cycle}, $\Mb^j_D$ becomes a tree, and two adjacent nodes are connected by exactly one bond. We then cut $\Mb^j_D$ by {\bf UP} but viewing each canonical cycle as a whole. After that, except the first node that is cut, each node has exactly one fixed end. We apply \textbf{UP} or \textbf{DOWN} to cut each node into elementary molecules. For each canonical cycle, by Proposition \ref{prop.alg_up} \eqref{it.up_proof_4}, we get a good molecule in each canonical cycle. There are at least $K$ canonical cycles, and we can gain at least $K$ good molecules and one bad molecule. Finally, by Proposition \ref{prop_MUMD} \eqref{it.prop_MUMD1}, we can cut the rest of $\Mb$ by Lemma \ref{lem.case2_rest} with $E=\Mb^j_D$ and $F=\varnothing$.
\end{proof}

For $i\geq1$, let $\Sigma_i$ be the set of atoms $\mf\in\overline{X}^+_0(\Mb^{i-1}_D)$ (recall that $\overline{X}^+_0$ is defined in Definition \ref{def.overline_X0}), such that $\mf$ has at least one child in $R_i$, and no child of $\mf$ in $R_i$ belongs to a canonical cycle of $\Mb^{i-1}_D$.
\begin{lemma}\label{lem.case2_sigma}
    Suppose $|\Sigma_i|\geq(d+1)K$ for some $i\geq1$. We cut $\Mb^i_D$ as free and cut it into elementary molecules by \textbf{UP}. Then we cut $R_i$ by \textbf{DOWN} but cut one atom each time; however, if we cut an atom $\pf\in R_i$ which is a $\{3_2\}$ atom, and all other children of some $\mf\in\Sigma_i$ have been cut, then we cut $\{\pf,\mf\}$ together. Finally we cut the rest of $\Mb$ by Lemma \ref{lem.case2_rest}, with $E=\Mb^{i-1}_D$ and $F$ the set of atoms in $\Sigma_i$. Then we have $\#_{\rm good}\geq K$ and $\#_{\rm bad}=1$.
\end{lemma}
\begin{proof}
    By Proposition \ref{prop_MUMD} \eqref{it.prop_MUMD1} and \eqref{it.prop_MUMD5}, cutting $\Mb^i_D$ gives the only bad molecule, and when we start to cut $R_i$, there are only top fixed ends in $R_i$. All children in $F$ will be cut as $F \subset \overline{X}^+_0(\Mb^{i-1}_D)$. So we can apply Lemma \ref{lem.case2_rest} in the last step.

    For $\mf\in\Sigma_i$, all children of $\mf$ are in $\Mb^{i-1}_D$, and at least one of them is in $R_i$. So the last child $\pf$ of $\mf$ that is cut belongs to $R_i$, and $\pf$ does not belong to any canonical cycle of $\Mb^{i-1}_D$.
    We claim that $\pf$ is not fully fixed when it is cut. Indeed, in \textbf{DOWN}, after we choose an atom $\nf$, we cut all the ancestors of $\nf$ which have not been cut, starting from the lowest one. If a child of $\pf$ is cut before $\pf$ in the step of some $\nf$, then $\pf$ is also an ancestor of $\nf$, and it is cut in the same step. So if two bottom edges of $\pf$ are fixed ends when it is cut, there are two different paths of children from $\pf$ to $\nf$, which form a cycle containing $\pf$ in $\Mb^{i-1}_D$. This contradicts that $\pf$ does not belong to any canonical cycle. Moreover, the top edge of $\pf$ connecting $\mf$ is not a fixed end, as $\mf$ has not been cut. Hence $\pf$ is not fully fixed when it is cut. If $\pf$ is a multi-particle atom, it is a good molecule. If $\pf$ is a binary atom, then by the regularity and the above argument, it has exactly one fixed end, so it is a $\{3_2\}$ atom and is cut together with $\mf$, which gives a $\{3_2 3_2A\}$ or $\{3_2 k_l\}$ molecule, and it is good. Each $\pf$ is the last child of at most $d+1$ atoms in $\Sigma_i$, so we get at least $K$ good molecules.
\end{proof}

{\bf Induction Round $j$} ($0\leq j\leq D-2$). Suppose $\Mb^j_D$ has been constructed, which satisfies Proposition \ref{prop_MUMD}. We consider the following cases.
\begin{enumerate}
    \item If the number of canonical cycles in $\Mb^j_D$ is $\geq K+1$, then we apply Lemma \ref{lem.case2_manycc}.
    \item If the number of multi-particle atoms outside those canonical cycles in $\Mb_D^j$ is $\geq d^2$, then we cut $\Mb^j_D$ as free and apply {\bf UP1} to $\Mb_D^j$. By counting the quantity $\sigma$, each such multi-particle atom must have at most one fixed end, so it is a good molecule with excess $\eps^{\upsilon_0}$ (by Proposition \ref{prop.est_single} \eqref{it.est_kl}), and the total gain is $\eps^{\upsilon_0 \cdot d^2} \leq \eps^{(d-1) + \upsilon}$. After that, we cut the rest of $\Mb$ by Lemma \ref{lem.case2_rest} with $E=\Mb^j_D$ and $F=\varnothing$. Note that $\#_{\mathrm{bad}} = 1$ in the whole process.
    \item If $j\geq1$ and $|\Sigma_i|\geq(d+1)K$ for some $1\leq i\leq j$, then we apply Lemma \ref{lem.case2_sigma}.
    \item Otherwise, if each canonical cycle in $\Mb^j_D$ has less than $N_{D-j-1}$ bonds connecting $\Mb_U$, then we
    \begin{enumerate}[(i)]
        \item run algorithm {\bf TRANSUP1} and apply Proposition  \ref{coro.alg_transup1} if $|\overline{X}^+_0(\Mb^j_D)| \geq (15d)^{15d} \cdot (N_{D-j-1})^2$;
        \item run algorithm {\bf MAINTRUP2, Option 1} and apply Proposition \ref{prop.alg_maintrup2} if $|Y_1(\Mb^j_D)| \geq (30d)^{30d} \cdot( N_{D-j-1})^2$;
        \item run algorithm {\bf MAINTRUP2, Option 2} and apply Proposition \ref{prop.alg_maintrup2} if $|\overline{X}^+_0(\Mb^j_D)| \leq (15d)^{15d} \cdot (N_{D-j-1})^2$ and $|W(\Mb^j_D)| \geq (30d)^{30d} \cdot (N_{D-j-1})^3$;
        \item run algorithm {\bf MAINTRUP2, Option 3} and apply Proposition \ref{prop.alg_maintrup2} if $|\overline{X}^+_0(\Mb^j_D)| \leq (15d)^{15d} \cdot (N_{D-j-1})^2$, $|W(\Mb^j_D)|\leq  (30d)^{30d} \cdot (N_{D-j-1})^3$, and $|Y_1(\Mb^j_D)| \leq (30d)^{30d} \cdot (N_{D-j-1})^2$.
    \end{enumerate}
    After each of the above algorithms, we cut $\Ac_0\setminus\Mb^j_D$ by \textbf{DOWN}, and then cut the rest of $\Mb$ by Lemma \ref{lem.case2_rest} with $E=\Ac_0$ and $F=\Mb_U$. By Proposition \ref{prop_MUMD} \eqref{it.prop_MUMD4}, $\#^+_{\rm conn}\geq N_{D-j}/(2(d+1))^{K+1}$, which is much larger than $(30d)^{30d}(N_{D-j-1})^3+ DK^2$, as $N_{D-j}=(60dN_{D-j-1})^{60d}$ and $N_1=D^D$. So in each case \eqref{eq.extra_case_1} is true.
    \item Otherwise, there exists a canonical cycle in $\Mb^j_D$ with at least $N_{D-j-1}$ bonds connecting $\Mb_U$. We construct $\Mb^{j+1}_D$ as above, and go to Round $j+1$.
\end{enumerate}

{\bf Completing the argument.} Now the remaining case is that in each Round $j$ ($0\leq j\leq D-2$) we construct $\Mb^{j+1}_D$. In particular, each $\Mb^j_D$ contains at most $K$ canonical cycles. We first note the following facts, which follow from Definition \ref{def.canonical_cycle}. Let $\Cc$ be a canonical cycle of a connected molecule $\Mb'$.
\begin{enumerate}[(F1)]
    \item\label{it.F1} Each cycle in $\Mb'$ which contains an atom in $\Cc$ is contained in $\Cc$. There exists no atom in $\Mb' \setminus \Cc$ that has both ancestor and descendant in $\Cc $; if $\Cc$ is not a single box, each atom (or box) in $\Cc$ is contained in a cycle in $\Cc$, and if $\Cc$ is a single box, then each atom inside the box is also contained in a cycle in $\Cc$.
    \item\label{it.F2} If $\Mb''\subseteq\Mb'$, then each canonical cycle of $\Mb''$ is contained in a canonical cycle of $\Mb'$. Each canonical cycle of $\Mb'$ which is contained in $\Mb''$ and contains no opened box is also a canonical cycle of $\Mb''$.
\end{enumerate}
If Round $j$ is of type \eqref{it.type1} in Definition \ref{def.construct_MD}, and $\Cc$ is a lowest canonical cycle of $\Ac_{j+1}$, we denote $H_j(\Cc):=\Cc\cap\Mb^{j+1}_D$ (hence $H_j(\Cc)$ has no parent in $\Cc \setminus H_j(\Cc)$, by the construction of $\Mb^{j+1}_D$).

\begin{lemma}\label{lem.case2_trim}
    \begin{enumerate}[(i)]
        \item\label{it.trim1} Let Round $j$ be of type \eqref{it.type1}, and $\Cc$ be a lowest canonical cycle of $\Ac_{j+1}$. Then $H_j(\Cc)\neq\Cc$, and the canonical cycles of $\Mb^{j+1}_D$ which contain an atom in $H_j(\Cc)$ are exactly the canonical cycles of $H_j(\Cc)$.
        \item\label{it.trim2} Let $\Cc$ be a canonical cycle of $\Mb^j_D$, which is not a lowest canonical cycle of $\Ac_{j+1}$ in a Round $j$ of type \eqref{it.type1}. If $\Cc$ has an atom in $\Mb^{j+1}_D$, then $\Cc\subseteq\Mb^{j+1}_D$, and $\Cc$ is a canonical cycle of $\Mb^{j+1}_D$.
        \item\label{it.trim3} There are at least $(D-1-K)/(K+1)$ rounds of type \eqref{it.type1}.
    \end{enumerate}
\end{lemma}
\begin{proof}
    \eqref{it.trim1} A cycle in $\Mb^{j+1}_D$ which contains an atom in $H_j(\Cc)$ gives a cycle in $\Mb^j_D$ which has an atom in $\Cc$ (by viewing the boxes opened in Round $j$ as a whole), so it is contained in $\Cc$ by (F\ref{it.F1}), hence in $H_j(\Cc)$. This gives the second statement. If $H_j(\Cc)=\Cc$, then $\Cc$ is not a single box (otherwise it is removed in \eqref{it.defMj_a} or opened in \eqref{it.defMj_b}), and \eqref{it.defMj_c} is not used for $\Cc$. But a lowest atom of $\Cc$ is a lowest atom of $S$ (as $\Cc$ is a lowest node of $\Ac_{j+1}$), it is contained in a cycle in $\Cc\subseteq S$ by (F\ref{it.F1}), and it is not a box; so one of \eqref{it.defMj_a}--\eqref{it.defMj_c} applies, which is a contradiction.  

    \eqref{it.trim2} The atoms of $\Cc$ are not removed in \eqref{it.defMj_a}--\eqref{it.defMj_c} of Definition \ref{def.construct_MD}, and $\Cc$ does not contain the atom $\pf$ in \eqref{it.type2}, which is not contained in any cycle. All other atoms are removed as a whole component of $\Mb^j_D\setminus\Ac_{j+1}$, or as a whole component of $S\setminus\{\qf\}$ (or $S\setminus\{\pf\}$). Since $\Cc$ is connected, it is contained in one of them. Then we use (F\ref{it.F2}).

    \eqref{it.trim3}  
    Let $a_j$ be the number of canonical cycles of $\Mb^j_D$ having a bond with $\Mb_U$, then $a_j\leq K$. Consider a Round $j$ of type \eqref{it.type2}. By Definition \ref{def.construct_MD}, the canonical cycles of $\Mb^{j+1}_D$ are canonical cycles of $\Mb^j_D$. Moreover, each leaf node of $\Ac_{j+1}$ is a canonical cycle having a bond with $\Mb_U$. So each component in $S \setminus \{\pf\}$ has at least one canonical cycle, so $a_{j+1}\leq a_j-1$. Hence among any $K+1$ consecutive rounds, at least one is of type \eqref{it.type1}, and the number of rounds of type \eqref{it.type1} among Rounds $0,\dots,D-2$ is at least $(D-1-K)/(K+1)$.
\end{proof}

Now we choose a sequence of cycles inside cycles, as follows. Each time, we choose the canonical cycle which contains the most lowest canonical cycles cut in later rounds.
\begin{proposition}\label{prop.case2_chain}
    There exist rounds $r_1<r_2<\cdots<r_{2K+1}$ of type \eqref{it.type1}, and canonical cycles $\Cc_1\supsetneq\Cc_2\supsetneq\cdots\supsetneq\Cc_{2K+2}$, such that for $1\leq i\leq 2K+1$, $\Cc_i$ is a lowest canonical cycle of $\Ac_{r_i+1}$, and $\Cc_{i+1}$ is a canonical cycle of $H_{r_i}(\Cc_i)$.
\end{proposition}
\begin{proof}
    For a canonical cycle $\Cc$ of $\Mb^k_D$, let $n_k(\Cc)$ be the number of pairs $(l,\Cc')$, where $l\geq k$, Round $l$ is of type \eqref{it.type1}, $\Cc'$ is a lowest canonical cycle of $\Ac_{l+1}$, and $\Cc'\subseteq\Cc$. By Lemma \ref{lem.case2_trim} \eqref{it.trim3}, there are at least $E:=(D-1-K)/(K+1)$ such pairs with $l\geq0$. By (F\ref{it.F2}), each $\Cc'$ is contained in a canonical cycle of $\Mb^0_D$, and there are at most $K$ of them. So we can choose one of them, denoted by $\Cc_1$, with $n_0(\Cc_1)\geq E/K$.

    Suppose $\Cc_i$ is a canonical cycle of $\Mb^k_D$ with $n_k(\Cc_i)\geq1$ (where $k=0$ for $i=1$, and $k=r_{i-1}+1$ for $i\geq2$). Let $r_i$ be the first round $l\geq k$ with a pair $(l,\Cc')$ such that $\Cc'\subseteq\Cc_i$. In Rounds $ k \leq l <r_i$, $\Cc_i$ is not a lowest canonical cycle of $\Ac_{l+1}$ of type \eqref{it.type1}, and it intersects $\Mb^{l+1}_D$, so by Lemma \ref{lem.case2_trim} \eqref{it.trim2}, $\Cc_i$ is a canonical cycle of $\Mb^{r_i}_D$ ($\Cc_i$ is unchanged in Rounds $k \leq l <r_i$). Since $\Cc'\subseteq\Cc_i$ is also a canonical cycle of $\Mb^{r_i}_D$, we have $\Cc'=\Cc_i$, i.e. $\Cc_i$ is a lowest canonical cycle of $\Ac_{r_i+1}$.

    For each pair $(l,\Cc'')$ with $l>r_i$ and $\Cc''\subseteq\Cc_i$, we have $\Cc''\subseteq\Mb^{r_i+1}_D\cap\Cc_i=H_{r_i}(\Cc_i)$, so by (F\ref{it.F2}) and Lemma \ref{lem.case2_trim} \eqref{it.trim1}, $\Cc''$ is contained in a canonical cycle of $H_{r_i}(\Cc_i)$, which is also a canonical cycle of $\Mb^{r_i+1}_D$. There are at most $K$ of them, and we choose one of them, denoted by $\Cc_{i+1}$, with $n_{r_i+1}(\Cc_{i+1})\geq(n_k(\Cc_i)-1)/K$. Then $\Cc_{i+1}\subseteq H_{r_i}(\Cc_i)\subsetneq\Cc_i$.

    Each step we lose a factor $K$, so by induction $n(\Cc_i)\geq E/K^i-2$. Since $D=8K^{2K+3}$, we have $E\geq 4K^{2K+2}-1$, hence $n(\Cc_i)\geq1$ for all $i\leq 2K+2$, and we can continue until we obtain $\Cc_{2K+2}$.
\end{proof}

\begin{lemma}\label{lem.case2_convex}
    For $1\leq i\leq 2K+2$, no atom in $\Mb$ has both ancestor and descendant of $\Cc_i$.
\end{lemma}
\begin{proof}
    This property follows for $\Mb^{r_{i-1}+1}_D$ by Proposition \ref{prop_MUMD} \eqref{it.prop_MUMD1}. Note that $\Cc_i$ is a canonical cycle in $\Mb^{r_{i-1}+1}_D$, and $\Mb^{r_{i-1}+1}_D$ is a tree by viewing each canonical cycle as a whole. Hence the property follows for $\Cc_i$.
\end{proof}
We then define the following algorithm.
\begin{definition}[Algorithms \textbf{UP3} and \textbf{DOWN3}]\label{def.case2_updown}
    We define \textbf{UP3} in the same way as \textbf{UP} in Definition \ref{def.alg_up}, with the following modifications for the M-atom $\mf$:
    \begin{enumerate}
        \item we do not apply Definition \ref{def.alg_up} \eqref{it.alg_up_4} to $\mf$, but view $\mf$ as an atom with $l$ particle lines, and cut it as in Definition \ref{def.alg_up} \eqref{it.alg_up_3};
        \item if cutting a $\{3_2\}$ atom $\pf$ makes $\mf$ fully fixed, then we cut $\{\pf,\mf\}$ together (i.e. we drop the restriction ``non-M'' in Definition \ref{def.alg_up} \eqref{it.alg_up_3} for $\mf$).
    \end{enumerate}
    We define \textbf{DOWN3} as the dual of \textbf{UP3}.
\end{definition}
\begin{remark}\label{rem.up3}
    The proof of Proposition \ref{prop.alg_up} \eqref{it.up_proof_1}--\eqref{it.up_proof_2} applies to \textbf{UP3} and \textbf{DOWN3} without change; so the regularity is kept and there is no new bad molecule. Moreover, when an M-atom is cut, it may not be fully fixed, or and may cut with some \{$3_2$\} molecule. We can argue as in the proof of Proposition \ref{prop.alg_up} \eqref{it.up_proof_MO}, and expand part of the M-molecule, to gain at least factor $\eps^{\upsilon}$, as we assume that each two-atom molecule or multi-particle is a good in the proof of Proposition \ref{prop.comb_est_extra}. 
    Hence, Proposition \ref{prop.alg_up} \eqref{it.up_proof_3}--\eqref{it.up_proof_4} also applies to \textbf{UP3} and \textbf{DOWN3}.
\end{remark}

Now we complete the last part of the proof of Proposition \ref{prop.comb_est_extra}. For $1\leq i\leq 2K+1$, denote $\Lc_i:=\Cc_i\setminus\Cc_{i+1}$, and let $\mathrm{Anc}_i$ be the set of atoms in $\Lc_i$ which are ancestors of some atom in $\Cc_{i+1}$. We apply the following cutting sequence:
\begin{enumerate}[(1)]
    \item\label{it.case2_cut1} Cut $\Cc_{2K+2}$ as free, and cut it into elementary molecules by \textbf{UP3}.
    \item\label{it.case2_cut2} For $i=2K+1,2K,\dots,1$, we cut $\Lc_i$ as follows: (a) cut $\mathrm{Anc}_i$ by \textbf{UP3}; (b) cut the fully fixed atoms in $\Lc_i$, until there is no such atom; (c) cut the rest of $\Lc_i$ by \textbf{DOWN3}.
    \item\label{it.case2_cut3} Cut the rest of $\Mb$ by Lemma \ref{lem.case2_rest} with $E=\Cc_1$ and $F=\varnothing$.
\end{enumerate}

We first show that this keeps the regularity and $\#_{\rm bad}=1$. When we start to cut $\Lc_i$, exactly the atoms in $\Cc_{i+1}$ have been cut. An atom in $\mathrm{Anc}_i$ adjacent to an atom $\nf\in\Cc_{i+1}$ is a parent of $\nf$ by Lemma \ref{lem.case2_convex}. So in (a) there are only bottom fixed ends. After (a), an atom which has not been cut and is adjacent to an atom in $\Cc_{i+1}\cup\mathrm{Anc}_i$ is a child of it (a parent would be in $\mathrm{Anc}_i$). An atom cut in (b) has all its top edges fixed, so it only gives top fixed ends to its children. So in (c) there are only top fixed ends. Step \eqref{it.case2_cut1} gives exactly one bad molecule, and each component of $\Lc_i$ is adjacent to $\Cc_{i+1}$ since $\Cc_i$ is connected, so there is no other bad molecule in \eqref{it.case2_cut2}. Finally, by Lemma \ref{lem.case2_convex}, we can apply Lemma \ref{lem.case2_rest} in \eqref{it.case2_cut3}.

Now we count the good molecules. 
We say level $i$ ($1\leq i\leq 2K$) is of \emph{type A}, if some atom in $\Cc_{i+1}$ has a child in $H_{r_i}(\Cc_i)\setminus\Cc_{i+1}$; otherwise we say it is of \emph{type B}.

{\bf Type A.} We claim that we get a good molecule when cutting $\Lc_i$. Let $\nf\in\Cc_{i+1}$ have a child $\mf\in H_{r_i}(\Cc_i)\setminus\Cc_{i+1}$. Then $\mf\in\Lc_i$, and $\mf\notin\mathrm{Anc}_i$ by Lemma \ref{lem.case2_convex}. There is no double bond between $\nf$ and $\mf$, otherwise it forms a cycle in $H_{r_i}(\Cc_i)$ intersecting $\Cc_{i+1}$, so $\mf\in\Cc_{i+1}$ (by viewing each canonical cycle as a whole $H_{r_i}(\Cc_i)$ is a tree).

First, consider another particle line of $\mf$ (not the one of the bond between $\nf$ and $\mf$), and an atom $\pf$ on it above $\mf$. Then $\pf\notin\Cc_{i+1}\cup\mathrm{Anc}_i$ by Lemma \ref{lem.case2_convex}. Going up along this particle line, we reach an atom outside $\Cc_i$ (which is cut in $\Lc_{i'}$ with $i'<i$), or the highest atom of this particle line in $A$. So by induction from the top, no atom of this particle line above $\mf$ is cut before (b) and (c). 

Second, when (c) starts, $\mf$ has only top fixed ends, including the one connecting $\nf$. If $\mf$ has at least two fixed ends, then its component has at least two fixed ends. Otherwise, the bond between $\nf$ and $\mf$ is the only edge of $\mf$ connecting an atom that has been cut. If $\Cc_i$ is not a single box, then by (F\ref{it.F1}) $\mf$ is contained in a cycle in $\Cc_i$. Therefore, for the component containing $\mf$, it either has at least two fixed ends, or has one fixed end and a cycle. This component gives a good molecule by Proposition \ref{prop.alg_up}.  If $\Cc_i$ is a single box $B$, then $H_{r_i}(\Cc_i)$ is in the top part of $B$ (recall the definition of box in Definition \ref{def.box}), and so is $\mf$. We follow a particle line of $\mf$ downward in $B$: it goes through the M-atom of $B$ and reaches the bottom part of $B$.

Therefore, when (c) starts, the component of $\mf$ has no fully fixed atom, and it contains either a cycle and a fixed end, or at least two fixed ends. By Proposition \ref{prop.alg_up} \eqref{it.up_proof_3}--\eqref{it.up_proof_4} (for \textbf{DOWN3}), we get a good molecule in it.

{\bf Type B.} We claim that we get a good molecule when cutting $\Lc_{i+1}$. Let $S=\Mb^{r_i+1}_D$ be the set at the end of Round $r_i$. Since $\Cc_i$ is a lowest node of $\Ac_{r_i+1}$, $S$ contains no child of an atom in $H_{r_i}(\Cc_i)$ outside $\Cc_i$.

First, $\Cc_{i+1}$ is not a single box: such a box is a lowest atom of $S$ (as level $i$ is of type B), so it would be treated in \eqref{it.defMj_a} or \eqref{it.defMj_b}. Let $\nf$ be a lowest atom of $\Cc_{i+1}$. As level $i$ is of type B, $\nf$ has no child in $S$, so $\nf$ is a lowest atom of $S$ contained in a cycle in $S$ (F\ref{it.F1}). Since \eqref{it.defMj_a} and \eqref{it.defMj_b} do not apply, $S\setminus\{\nf\}$ is not connected or loses a particle line; in particular $\nf$ is not a binary atom (a binary lowest atom in a cycle has both parents in this cycle), and $\nf$ is not a box. If $\nf\notin\Cc_{i+2}$, let $\mf:=\nf$ and $S_*:=S$. If $\nf\in\Cc_{i+2}$, then in Round $r_{i+1}$, $\nf$ is a lowest atom contained in a cycle in the final set $S$, in the lowest canonical cycle $\Cc_{i+1}$; if \eqref{it.defMj_c} had not been used for $\Cc_{i+1}$, it would apply to $\nf$. So \eqref{it.defMj_c} is used, and we let $\mf:=\qf$ be the atom removed in it, and $S_*$ be the set $S$ at this moment. In both cases, $\mf\in\Lc_{i+1}$ is a multi-particle C-atom, it is a lowest atom of $S_*$, $\Cc_{i+2}\subseteq S_*$, no atom in $\Cc_{i+1}\setminus S_*$ is a parent of an atom in $S_*$, and $S_*\setminus\{\mf\}$ is not connected or loses a particle line (otherwise $\mf$ is removed in \eqref{it.defMj_a}).

Now we show that when (c) for $\Lc_{i+1}$ starts, $\mf$ has not been cut and is not fully fixed. Since $\mf$ is a lowest atom of $S_*\supseteq\Cc_{i+2}$, we have $\mf\notin\mathrm{Anc}_{i+1}$, so no child of $\mf$ is in $\Cc_{i+2}\cup\mathrm{Anc}_{i+1}$. As an atom cut in (b) has all its parents cut, no child of $\mf$ is cut before (c). It remains to find a top edge of $\mf$ which is not fixed before (c); recall that atoms outside $\Cc_{i+1}$ are cut in a later level or in \eqref{it.case2_cut3}.
\begin{itemize}
    \item If $S_*\setminus\{\mf\}$ loses a particle line, then the top edge of $\mf$ on it is a free end or goes outside $\Cc_{i+1}$.
    \item Otherwise $S_*\setminus\{\mf\}$ is not connected. Since each atom of $\mathrm{Anc}_{i+1}$ reaches $\Cc_{i+2}$ by a path of children, $\Cc_{i+2}\cup\mathrm{Anc}_{i+1}$ lies in one component. Let $P$ be another component; it contains a parent $\pf$ of $\mf$. Each parent of an atom in $P$ is in $P$ or outside $\Cc_{i+1}$, so by induction from the top, no atom of $P$ is cut before (c). In particular the edge between $\mf$ and $\pf$ is not fixed.
\end{itemize}
Hence, when (c) starts, $\mf$ is a multi-particle atom which is not fully fixed. If no good molecule is obtained before $\mf$ is cut, then by the claim above no atom becomes fully fixed, so $\mf$ is cut when it is not fully fixed, and it is a good molecule.
 
Finally, for each level $1\leq i\leq 2K$, we get a good molecule in $\Lc_i$ (type A) or in $\Lc_{i+1}$ (type B). The sets $\Lc_i$ are disjoint, and each $\Lc_m$ is used by at most two levels ($i=m$ and $i=m-1$), so $\#_{\rm good}\geq K$. Together with $\#_{\rm bad}=1$ and $K=(2d+1)/\upsilon$, we get $(\upsilon/2)\#_{\rm good}-d\#_{\rm bad}\geq\upsilon/2$, which is \eqref{eq.extra_case_1}. This completes the proof of Case 2.

\subsection{Algorithms in each induction round}
In this subsection we introduce some algorithms and prove their properties that are used in subsection \ref{subsec_MUMD}.
     

\begin{definition}[Algorithm \textbf{TRANSDN1}]\label{def.alg_transup1}
   Let $A$ be constructed as in the proof of Proposition \ref{prop.comb_est_extra}, and let $\Mb_D$, $\Mb_U$ be constructed Definition \ref{def.construct_MD}. We define the similar algorithm as in  \ref{def.alg_transup}, but with the following modifications:
    \begin{enumerate}
        \item\label{it.alg_transup_1_1}  Throughout all steps, \emph{we treat each canonical cycle in $\Mb_D$ as a whole and regard them as atoms}, but treat each atom in $\Mb_U$ individually. If $\Mb_D$ contains any fully fixed atom $\nf$, then cut it as free, and repeat until there is no such atom left.
        \item\label{it.alg_transup_1_2} Choose a \emph{highest} atom $\nf$ in the set of all atoms with fixed end in $\Mb_D$ (or a highest atom in $\Mb_D$, if $\Mb_D$ does not contain such atom). Let $Z_\nf$ be the set of ancestors of $\nf$ in $\Mb_D$.
        \item \label{it.alg_transup_1_3} Starting from $\nf$, choose a \emph{lowest} atom $\mf$ in $Z_\nf$ that has not been cut. If $\mf$ is a canonical cycle, then cut it as free, and then cut it into elementary molecules by {\bf UP} or {\bf DOWN}. If $\mf$ is a $\{3_2\}$ atom, and has an ancestor $\pf\in \overline{X}^+_0(\Mb_D)$ that becomes fully fixed after cutting $\mf$, then cut $\{\mf,\pf\}$ as free. In all remaining cases, cut only $\mf$ as free. Repeat until all atoms in $Z_\nf$ have been cut. Then go to (\ref{it.alg_transup_1_1}).
        \item \label{it.alg_transup_1_4} Repeat (\ref{it.alg_transup_1})--(\ref{it.alg_transup_3}) until all atoms in $\Mb_D$ have been cut. Then cut the remaining part of $\Mb_U$ as free and cut it into elementary components using \textbf{UP}.
    \end{enumerate}
\end{definition}
\begin{proposition}\label{coro.alg_transup1}
    For the algorithm  \textbf{TRANDN1}, we assume that each canonical cycle in $\Mb_D$ has at most $N$ particles connecting $\Mb_U$, and the number of canonical cycle in $\Mb_D$ is $\leq K$. Then after \textbf{TRANDN1},  $\#_{\{3_23_2B\}}=\#_{\{4_2 (2l)_l\}}=0$ and $\#_{\mathrm{bad}}=1$, and $\#_{\mathrm{good}}\geq \frac{|\overline{X}^+_0(\Mb_D)|}{2}- N\cdot K - (d+1) \#_{\rm mp}(\Mb_D) - 1$, under the assumption that all two-atom molecules are good. Here $\#_{\rm mp}(\Mb_D)$ denotes the number of multi-particle atom in $\Mb_D$ outside any canonical cycle.
\end{proposition}
\begin{proof}
    We have $\#_{\rm bad} = 1$, due to the fact that $\Mb_U$ and $\Mb_D$ are all connected. For each atom in $\overline{X}^+_0(\Mb_D)$, with at most $N\cdot K + (d+1) \#_{\rm mp}$ exceptions, all its children in $\Mb_D$ must be binary atoms outside any canonical cycle. Moreover, for each atom outside canonical cycles, it has at most one fixed end when it is cut, by the same counting argument as in the proof of Proposition \ref{prop.alg_transup}, but viewing each canonical cycle as a whole (the resultant $\Mb_D$ is a tree). Then, by the same argument in the proof of Proposition \ref{prop.alg_transup}, each of such atom in $\overline{X}^+_0(\Mb_D)$ either belong to a good molecule, or there exists a good molecule cut before it. This leads to at least $\frac{|\overline{X}^+_0(\Mb_D)| -  N\cdot K- (d+1) \#_{\rm mp}}{2} - 1$ good molecules, which meets our requirement.
\end{proof}

 \begin{definition}[The function \textbf{SELECT3}]
\label{def.func_select3} Let $\Mb$ be a connected UD molecule with only C-atoms, no bottom fixed end in $\Mb_D$ and no top fixed end in $\Mb_U$, and no fully fixed atoms. We assume $\overline{X}^+_0(\Mb_D) = \emptyset$ (Definition \ref{def.overline_X0}). Let $W$ be the set of atoms in $\Mb_D$ that are connected to the same atom in $\Mb_U  $. Let $Y$ be the set of all canonical cycles in $\Mb_D$, let $Z$ be the set of atoms in $\Mb_D \setminus Y$, with fixed ends, or are  multi-particle atoms, or are involved in weak degeneracies. We define the function $\textbf{SELECT3}=\textbf{SELECT3}(\Mb_D,W,Z,Y)$ as follows.
\begin{enumerate}
\item\label{it.func_select3_1} Consider all the components of $W\cup Z\cup Y$ in $\Mb_D$, which is a finite collection of disjoint subsets of $\Mb_D$.
\item\label{it.func_select3_2} If any two of the subsets in (\ref{it.func_select3_1}), say $U$ and $V$, have the shortest distance (within $\Mb_D$) which is at most 4, then we choose one shortest path between an atom in $U$ and an atom in $V$, and let the atoms on this path be $\nf_j$. Then replace the two sets $U$ and $V$ by a single set which is $U\cup V\cup \{\nf_j\}$.
\item\label{it.func_select3_5} When no scenario in (\ref{it.func_select3_2}) occurs, we output $S:=\textbf{SELECT3}(\Mb_D,W,Z,Y)$ as the union of all the current sets.
\end{enumerate}
\end{definition}
\begin{proposition}\label{prop.select3}
    After applying {\bf SELECT3}, the molecule $\Mb_D$ becomes a proper forest after cutting $S$ as free. Consider all the bonds in $\Mb$ that contains at least one atom in $S$ and one atom in $\Mb_U$. Then the number of these bonds does not exceed $ 2|P_Y|\cdot |Y|+ 10d \cdot (|W|+|Z| +|Y|)$, where $|P_Y|$ denotes the maximal number of particles that connects to $\Mb_U$ of the canonical cycles.
\end{proposition}
\begin{proof}
    We proof the upper bound for the bonds connecting $\Mb_U$ and $S$ in function {\bf SELECT3}.
        Simply note that each single step in Definition \ref{def.func_select3} \eqref{it.func_select3_2}  adds at most $3$ atoms to $S$. Each single step in Definition \ref{def.func_select2} \eqref{it.func_select2_2} decrease the number of subsets by $1$. This means that the total number of steps in Definition \ref{def.func_select3} \eqref{it.func_select3_2} is bounded by $|W|+|Y|+|Z|.$
        Moreover, each component in $Y$ has at most $|P_Y|$ bonds connecting $\Mb_U$ and $S$, and each atom in $W$ and $Z$ has at most $d+1$ such bonds, and an atom in $S \setminus(W\cup  Y \cup Z)$ has at most $2$ such bonds. Hence, the total bonds connecting $\Mb_U$ and $S$ is no more than $2 |P_Y| \cdot |Y| + 10d \cdot (|W|+|Z|+|Y|)$.
\end{proof}
\begin{definition}[Algorithm \textbf{MAINTRUP2}] 
\label{def.alg_maintrup2}
Suppose $\Mb$ is a connected UD molecule, such that $\Mb_U^j$ and $\Mb_D^j$ are constructed in proof of \eqref{prop.comb_est_extra}. Let $Y$ be atoms belong to a canonical cycle in $\Mb_D^j$ and assume each canonical cycle has at most $|P_Y|$ particle lines connecting $\Mb_U$. Let $Y_1$ be those atoms that does not belong to $Y$ and involves in these weak degeneracies, or is a multi-particle atom. We define the following algorithm, which contains three Options that we can choose at the beginning. Let $W$ be the set of atoms \emph{outside} those canonical cycles in $\Mb_D^j$, that connects to the same atom in $\Mb_U\backslash \overline{X}^+_0(\Mb_D^j)$. Denote $\#^+_{\mathrm{conn}}$ by the number of bonds between $\Mb_U$ and $\Mb_D^j$.

In \textbf{Option 1} we do the followings:
\begin{enumerate}
\item Cut $\Mb_D^j$ as free, then cut it into elementary components using \textbf{UP}.
\item Then cut $\Mb_U$ as free and cut it into elementary components using \textbf{UP}.
\end{enumerate}

In \textbf{Option 2} we do the followings:
\begin{enumerate}
    \item  Denote the set of atoms in $\Mb_U$ that connects at least one atom in $W$ by $W'$\footnote{Note in Definition \ref{def.alg_maintrup}, each atom in $W'$ must have descendants in $\overline{X}^+_0(A)$, as $W' \in X(A)$ there. But here, atoms in $W'$ may have no descendants in $\overline{X}^+_0(\Mb_D^j)$.}. If at least half of atoms in $W'$ has descendants in $\overline{X}^+_0(\Mb_D^j)$, we cut all atoms in $\overline{X}^+_0(\Mb_D^j)$ as free, then we cut $\Mb_U$ into elementary components using \textbf{UP}, but only cut one atom each time. Then cut each component in $\Mb_D^j$ by \textbf{DOWN}.
    
    \item\label{prop.chain} If there exists an lowest atom $\nf $ in $\Mb_U$, such that it has at least $2 d^2$ ancestors in $W'$. Then we cut $\Mb_U$ as free, and then apply {\bf UP1}, but first choose the atom $\nf$. After that, we cut $\Mb_D^j$ into elementary molecules by \textbf{DOWN} .

    \item Otherwise, if each atom in $\Mb_U$ has at most $d^2$ ancestors in $W'$, and at least half of atoms in $W'$ has no descendants in $\overline{X}^+_0(\Mb^j_D)$. We first cut $\Mb_D^j$ as free, then cut all fully fixed atoms in $\Mb_U$. After that, we apply \textbf{UP} to $\Mb_U$ and $\Mb_D^j$.
\end{enumerate}

In \textbf{Option 3} we do the followings:
\begin{enumerate}
\item\label{it.alg_maintrup2_0} For each dropped piece $P$ created in Rounds $0,\dots,j-1$ (see the construction of $\Mb^{j+1}_D$ in Section \ref{subsec_MUMD}), we cut $P$ as free, and cut it into elementary molecules by \textbf{UP}.
\item\label{it.alg_maintrup2_1} Cut all atoms in $\overline{X}^+_0(\Mb_D^j)$ as free. If any atom in $\Mb_D^j$ becomes fully fixed, also cut it as free until no such atom exists in $\Mb_D^j$.
\item\label{it.alg_maintrup2_2} If $\Mb_D^j$ remains connected after the above step, let $Z$ be the set of atoms with fixed ends in $\Mb_D^j$, and the atoms in $Y\cup Y_1$ that has not been cut. Define $S=\textbf{SELECT3}(\Mb_D^j,W,Z,Y)$ as in Definition \ref{def.func_select3}, then cut $S$ as free and cut it into elementary components using \textbf{DOWN}. If $\Mb_D^j$ is not connected, apply this step to each connected component of $\Mb_D^j$. After this, there exist only binary atom in $\Mb_D^j$.
\item\label{it.alg_maintrup2_3} Apply the steps (\ref{it.alg_maincr_2})--(\ref{it.alg_maincr_8}) in Definition \ref{def.alg_maincr} to $\Mb_U$, where the set $X_0\cup X_1$ is replaced by $\overline{X}^+_1(\Mb^j_D)$, i.e. the atoms in $\Mb_U$ with a bond to $\Mb^j_D$ (if there is no remaining atom in $\Mb_U$ that belongs to $\overline{X}^+_1(\Mb^j_D)$ or has a fixed end, we choose any lowest atom). Then repeat the step in Definition \ref{def.alg_maincr} (\ref{it.alg_maincr_9}) until all atoms in $\Mb_D^j$ have been cut.
\item\label{it.alg_maintrup2_4} Cut the rest of $\Ac_0\setminus\Mb^j_D$ by \textbf{DOWN}, and then cut the rest of $\Mb$ by Lemma \ref{lem.case2_rest} with $E=\Ac_0$ and $F=\Mb_U$.
\end{enumerate}

We define the dual algorithm \textbf{MAINTRDN2} in the same way (so $\#_{\mathrm{conn}}^+$ is replaced by $\#_{\mathrm{conn}}^-$ etc.).
\end{definition}

\begin{proposition}
\label{prop.alg_maintrup2} Assume all two-atom molecules and multi-particle atoms are good. In \textbf{Option 1} of \textbf{MAINTRUP2} in Definition \ref{def.alg_maintrup2} (same for \textbf{MAINTRDN2}), we have
\begin{equation}\label{eq.alg_maintrup2_1}
\#_{\mathrm{bad}}=1;\qquad \#_{\mathrm{good}}\geq  |Y_1|/10.
\end{equation}
In \textbf{Option 2} of \textbf{MAINTRUP2} in Definition \ref{def.alg_maintrup2} (same for \textbf{MAINTRDN2}), for the first case, we have
\begin{equation}\label{eq.alg_maintrup2_2}
\#_{\mathrm{bad}}=\overline{X}^+_0(\Mb_D^j);\qquad \#_{\mathrm{good}}\geq  \frac{|W|}{2(d+1)} - \overline{X}^+_0(\Mb_D^j).
\end{equation}
For the second case, there exists $d^2$ good multi-particle atoms, and $1$ bad molecule. For the third case, we have
\begin{equation}
   \#_{\mathrm{good}} \geq \frac{|W|}{10d^3}  ,\qquad  \#_{\mathrm{bad}} = 1.
\end{equation}

 In \textbf{Option 3} of \textbf{MAINTRUP2} in Definition \ref{def.alg_maintrup2}, we have all \{$3_23_2$B\} molecules are good, and moreover
\begin{equation}\label{eq.alg_maintrup2_3}
\begin{aligned}
\#_{\{3_23_2A\}}+\#_{\{3_23_2B\}}&\geq\frac{1}{10}\cdot\big(\#_{\mathrm{conn}}^+ - 2|P_Y|\cdot |Y| - 10^3 d^3(|Y_1| + |\overline{X}^+_0(\Mb^j_D)| +|W|+K)  \big),\\
\#_{\mathrm{bad}}&\leq |\overline{X}^+_0(\Mb^j_D)| + |Y_1| + |W| + K + 3(d+1) j \cdot K,
\end{aligned}
\end{equation}
\end{proposition} 
\begin{proof}
    First, \eqref{eq.alg_maintrup2_1} follows by similar arguments in proof of \eqref{eq.alg_maintrup_1}. By viewing each canonical cycle as a whole, and counting the quantity $\sigma$ as before, each atom in $Y_1$ has at most one fixed end when it is cut. In the first case of \textbf{Option 2} of \textbf{MAINTRUP2}, if at least half of atoms in $W'$ has descendants in $\overline{X}^+_0(\Mb_D)$. After cutting atoms in $\overline{X}^+_0(\Mb_D)$ as free, and cut the subsequent fully fixed atom in $\Mb_U$, at least half of atoms in $W'$ has descendants with bottom fixed ends (note that no atoms in $W'$ have been cut at this stage). After applying \textbf{UP}, those atoms in $W'$ must have only bottom fixed end, with at least two free ends. Therefore, by Proposition \ref{prop.est_single} \eqref{it.est_kl}, they are all good molecules. We have 
    \begin{equation*}
        \#_{\mathrm{good}} \geq \frac{|W'|}{2} -| \overline{X}^+_0(\Mb_D)| \geq  \frac{|W|}{2(d+1)} - |\overline{X}^+_0(\Mb_D^j) |,\qquad  \#_{\mathrm{bad}}= \overline{X}^+_0(\Mb_D^j).
    \end{equation*}
    
    For the second case in \textbf{Option 2}, after cutting $\nf$, each atom in $Z_\nf$ will have only bottom fixed ends when it is cut. Hence, for each atom in $W' \cap Z_\nf$, it must be a multi-particle atom with at least two bottom free end. By Proposition \ref{prop.est_single} \eqref{it.est_kl}, they are all good molecules, and there are $2 d^2$ of them. 
    
    For the last case in \textbf{Option 2}, after cutting $\Mb_D$ as free and the subsequent fully fixed atoms in $\Mb_U$, at least half atoms in $W'$ have not been cut. Whenever we cut a such atom $\pf\in W'$, it has no top fixed end (as it has at least two bottom fixed ends). If it is cut with some atom $\pf$, then $\{\pf,\qf\}$ must belong to a good molecule. Moreover, this cutting may turn some atom $\qf$, which is a descendant of some atom in $W'$, or some atom in $W'$, into fully fixed; and cutting $\pf$ as free lead to some atom in $W'$ cut as fully fixed. But the cost of such atoms in $W'$ is at most $2 d^2$ atoms, as each atom in $\Mb_U$ has at most $2d^2$ ancestors in $W'$. If $\pf$ is cut by itself, and it is fully fixed when it is cut, then there must exists a good molecule before cutting $\pf$, again the cost of atoms in $W'$ for each such good molecule is at most $2 d^2$, as each atom of $\Mb_U$ has at most $2d^2$ ancestors in $W'$ by assumption. 

    As a result, we have   
      \begin{equation*}
        \#_{\mathrm{good}} \geq \frac{|W'|}{4d^2}  \geq \frac{|W|}{10d^3}  ,\qquad  \#_{\mathrm{bad}} = 1.
    \end{equation*}
    
    As for the proof of \textbf{Option 3}, here $\Mb_D=\Mb^j_D$ is constructed in Section \ref{subsec_MUMD}, and we assume $|\Sigma_i|<(d+1)K$ for all $1\leq i\leq j$ (otherwise we use Lemma \ref{lem.case2_sigma}). We first bound $\#_{\rm bad}$.

    In Definition \ref{def.alg_maintrup2} \eqref{it.alg_maintrup2_0}, by Proposition \ref{prop_MUMD} \eqref{it.prop_MUMD6}, two dropped pieces have no bond with each other, and a dropped piece $P$ has no bond with $\Mb^j_D$. A parent of an atom in $P$ which is not in $P$ is not in $\Ac_0$, so it is in $A^+$, and it gets a bottom fixed end. A child of an atom in $P$ which is not in $P$ is in $\Ac_0\setminus\Mb^j_D$, $A\setminus\Ac_0$ or $A^-$, and it gets a top fixed end; it does not get a bottom fixed end in this step, as it is not a parent of an atom in another dropped piece. So no fixed end is created in $\Mb^j_D$, and the atoms in $\Ac_0\setminus\Mb^j_D$ only have top fixed ends, also in the later steps (they are below $\Mb^j_D\cup\Mb_U$). Hence the regularity is kept, and in \eqref{it.alg_maintrup2_4} we can apply Lemma \ref{lem.case2_rest}, which creates no bad molecule. As no atom adjacent to $P$ has been cut before, each $P$ gives at most one bad molecule, and there are at most $dK$ dropped pieces in each round (Proposition \ref{prop_MUMD} \eqref{it.prop_MUMD6}). So this step gives at most $dKj$ bad molecules. The step \eqref{it.alg_maintrup2_1} gives at most $|\overline{X}^+_0(\Mb^j_D)|$ bad molecules. In \eqref{it.alg_maintrup2_2}, a component of $S$ can be full only if it comes from a set in $W$, in $Y_1$, or a canonical cycle (there are at most $K$ of them), which gives at most $|W|+|Y_1|+K$ bad molecules. The step in Definition \ref{def.alg_maincr} (\ref{it.alg_maincr_9}) gives no bad molecule, by Part 2 of the proof of Proposition \ref{prop.alg_maincr}.

    It remains to consider the steps (\ref{it.alg_maincr_2})--(\ref{it.alg_maincr_8}) in Definition \ref{def.alg_maincr} in \eqref{it.alg_maintrup2_3}. A bad molecule is created only when we choose an atom $\pf$ in Definition \ref{def.alg_maincr} \eqref{it.alg_maincr_3} without fixed end. Let $\mf$ be a lowest descendant of $\pf$ in $\Mb_U$. No descendant of $\pf$ in $\Mb_U$ has been cut: otherwise we will choose its descendant instead of $\pf$. So no child of $\mf$ has been cut. By Proposition \ref{prop_MUMD} \eqref{it.prop_MUMD2}, all children of $\mf$ are in $\Ac_0$, and they are not all in $\Mb^j_D$ (otherwise $\mf\in\overline{X}^+_0(\Mb^j_D)$ and have been cut in \eqref{it.alg_maintrup2_1}). After $\pf$ is chosen, all atoms in $S_\pf \cap \Mb_U$ are cut, including $\mf$, so each ancestor of $\mf$ chosen later has a fixed end. Hence different bad molecules give different $\mf$. Let $i$ be the smallest number such that $\mf$ has a child in $R_i$ (recall that $R_i = \Mb^{i-1}_D \setminus \Mb^i_D$); then $\mf\in\overline{X}^+_0(\Mb^{i-1}_D)$. By Proposition \ref{prop_MUMD} \eqref{it.prop_MUMD6}, each child of $\mf$ in $R_i$ is either in a component of $\Mb^{i-1}_D\setminus\Ac_i$, where it does not belong to any canonical cycle (otherwise such a canonical cycle would have a bond with $\Mb_U$, contradicts to the choice of $\Ac_i$), or it is a lowest atom removed in \eqref{it.defMj_c} or \eqref{it.type2} in Definition \ref{def.construct_MD}. It is not removed in \eqref{it.defMj_a} or \eqref{it.defMj_b} (see Proposition \ref{prop_MUMD} (\ref{it.prop_MUMD6}b)), and it is not in a dropped piece (otherwise it is cut in \eqref{it.alg_maintrup2_0}, and $\mf$ has a fixed end). So by Proposition \ref{prop_MUMD} \eqref{it.prop_MUMD6} either $\mf\in\Sigma_i$, which gives at most $(d+1)K$ choices for each $i$, or $\mf$ is a parent of some $\qf$ or $\pf$ in Definition \ref{def.construct_MD}, which gives at most $(d+1)K$ choices for each $i$. Summing up, we get the bound on $\#_{\rm bad}$ in \eqref{eq.alg_maintrup2_3}.

    The remaining part is similar as the corresponding part in Proposition \ref{prop.alg_maintrup}. Each \{$3_23_2$B\} belongs to $\Mb_D$. For the lower bound on $\{3_2 3_2\}$, due to the same reason as in proof Proposition \ref{prop.alg_maintrup}, we can cut $S$ by {\bf DOWN} as it has no bottom fixed end.
    
    Denote $Z_0$ be atoms in $\Mb_D \setminus Y $ cut in Definition \ref{def.alg_maintrup2} \eqref{it.alg_maintrup2_1}, and $Z_1$ be atoms in $\Mb_D$ with fixed end after Definition \ref{def.alg_maintrup2} \eqref{it.alg_maintrup2_1}. Then $|Z_1| \leq d(|\overline{X}^+_0(A)| +Z_0)$ due to the same reason as above. Moreover, we have $\rho(Z_0) = 0$ and $\rho(Z_0) \geq |Z_0| -  d|\overline{X}^+_0(\Mb_D)|$. We know $|Z_0| \leq d |\overline{X}^+_0(\Mb_D)|$. Thus, $|Z_1| \leq 2d^2 |\overline{X}^+_0(A)|$. By Proposition \ref{prop.func_select2}, the number of bonds connecting $S$ and $\Mb_U$ is no more than 
    \begin{equation}
        2 |P_Y|\cdot |Y| + 10^2 d^3 ( |Y_1| +|\overline{X}^+_0(\Mb_D)|+ \#_{\mathrm{mp}}(\Mb_D) + |W| +|Y|+ K).
    \end{equation}  
    After cutting $\overline{X}^+_0(A) \cup Z_0 \cup S$ as free, the number of bonds connecting an atom in $\Mb_D$ to an atom in $\Mb_U$ is still at least
    \begin{equation}
        \#_{\mathrm{conn}}^+ -  2 |P_Y|\cdot |Y| - 10^3 d^3 ( |Y_1| + |\overline{X}^+_0(\Mb_D)|+ \#_{\mathrm{mp}}(\Mb_D) + |W| +|Y|+ K),
    \end{equation}
    and we can obtain the lower bound on $\{3_2 3_2\}$ molecule.
\end{proof}

\appendix
\section{Some auxiliary results}\label{app.aux}
In this appendix, we state some useful results. We start with a result concerning the preservation of the norm $\mathrm{Bol}^\beta$ used in Theorem \ref{th.main}, which was proved in Proposition A.1 \cite{DHM24}.
\begin{proposition}
\label{prop.pres_of_decay} Recall the norm $\mathrm{Bol}^\beta$ defined in (\ref{eq.boltzmann_decay_2}). Suppose $n=n(t,x,v)$ is a positive solution to the Boltzmann equation (\ref{eq.boltzmann}), and let $\beta=\beta(t)$ be a fixed, strictly positive, strictly decreasing function with $\beta(0)\geq \beta>0$ and $\beta'(t)\leq -c<0$ for some constant $c$. Assume $n$ satisfies that
\begin{equation}\label{bound1}\sup_{t\in[0,T]}\big\|e^{\beta(t)|v|^2}n(t)\big\|_{L^\infty}\leq A<\infty,
\end{equation} and that
\begin{equation}\label{bound2}\|n(0)\|_{\mathrm{Bol}^{\beta(0)}}+\|\nabla_xn(0)\|_{\mathrm{Bol}^{\beta(0)}}\leq B<\infty,
\end{equation} then for any $t\in[0,T]$ we have
\begin{equation}\label{bound3}\|n(t)\|_{\mathrm{Bol}^{\beta(t)}}+\|\nabla_xn(t)\|_{\mathrm{Bol}^{\beta(t)}}\lesssim B,
\end{equation} where the constants in $\lesssim$ may depend on $(A,\beta,c,T)$.
\end{proposition}
\begin{proof}
    The proof is exactly the same as the hard sphere case, by writing the Boltzmann equation in terms of $\nu$ as in \eqref{eq.boltzmann}.
\end{proof}

\begin{proposition}[\cite{GST14}]\label{prop.coll.kernel}
Define $\Jc_0 = \frac{|(v_1 -v_2) \wedge \nu_{12} | }{|v_1 -v_2|}$ and $\Ec_0 := |v_1 -v_2 |^2$. Define 
$$r_* = r_*(\Ec_0,\Jc_0) := \max \Big( r \in (0,1) \ \big|\  \frac{\Ec_0 \Jc_0^2 }{r^2} + 4 \Phi(r) = \Ec_0\Big).$$
Then the deflection angle (i.e. the angle between $v_1 -v_2$ and $v_3 -v_4$) is equal to $  \pi - 2 \Theta$ where 
\begin{equation}\label{eqn.Theta}
    \Theta = \Theta(\Ec_0,\Jc_0):= \arcsin \Jc_0 + \Jc_0 \int_{r_*}^1 \frac{\d r}{ r^2 \sqrt{1 - \frac{4 \Phi(r) }{\Ec_0} - \frac{\Jc_0^2}{r^2}} }.
\end{equation}
Moreover, under Assumptions \ref{ass.phi1}, \ref{ass.phi2} and \ref{ass.phi3}, for all $\mathcal{E}_0>0$, the function $\mathcal{J}_0\mapsto \Theta(\mathcal{E}_0,\mathcal{J}_0)\in[0,\pi/2]$ satisfies $\Theta(\mathcal{E}_0,0)=0$ and is strictly monotonic: $\partial_{\mathcal{J}_0}\Theta>0$ for all $\mathcal{J}_0\in(0,1)$. Moreover, it satisfies
$$\lim_{\mathcal{J}_0\to 0}\partial_{\mathcal{J}_0}\Theta\in(0,\infty) \quad\text{and}\quad \lim_{\mathcal{J}_0\to 1}\partial_{\mathcal{J}_0}\Theta=0  \ \text{ if }\Phi'(1^-)=0,\qquad \lim_{\mathcal{J}_0\to 1}\partial_{\mathcal{J}_0}\Theta=+\infty\ \text{ if }\Phi'(1^-)<0.$$
Let $\omega_{12}$ be the apsidal line, then $|(v_1-v_2)\cdot\omega_{12}|=|v_1-v_2|\cos\Theta$. Moreover, we have
\begin{equation}\label{eqn.defb}
    b(|v_1 -v_2|, \omega_{12}) = |v_1 -v_2 |\Jc_0^{d-2} \p_\Theta \Jc_0 (\sin \Theta)^{2-d}, \quad [(v_1 -v_2)\cdot \nu_{12}]_- \d \nu_{12} = b(|v_1 -v_2|, \omega_{12}) \d \omega_{12}.
\end{equation}
\end{proposition}

\begin{proof}
Let an energy $\mathcal{E}_0>0$ being fixed, the limiting values $\Theta(\mathcal{E}_0,0)=0$ and $\Theta(\mathcal{E}_0,1)=\pi/2$ are found by direct computation. To prove monotonicity, we use the change of variable
\[
\sin^2\varphi:=\frac{4\Phi(r)}{\mathcal{E}_0}+\frac{\mathcal{J}_0^2}{r^2},
\]
which yields
\[
\Theta(\mathcal{E}_0,\mathcal{J}_0)=\arcsin\mathcal{J}_0
+\int_{\arcsin\mathcal{J}_0}^{\frac{\pi}{2}}
\frac{\sin\varphi}{\frac{\mathcal{J}_0}{r}-\frac{2r^2\Phi'(r)}{\mathcal{E}_0\mathcal{J}_0}}\,d\varphi .
\]
Computing the derivative of this expression with respect to $\mathcal{J}_0$ gives
\begin{equation}\label{eqn.p_theta}
    \begin{aligned}
    \frac{\partial\Theta}{\partial\mathcal{J}_0}(\mathcal{E}_0,\mathcal{J}_0)
    &=\frac{1}{\sqrt{1-\mathcal{J}_0^2}}
    \big(1-\frac{\mathcal{E}_0\mathcal{J}_0^2}{\mathcal{E}_0\mathcal{J}_0^2-2\Phi'(1)}\big)\\
    &\quad +\int_{\arcsin\mathcal{J}_0}^{\frac{\pi}{2}}
    \frac{2\mathcal{E}_0^2\mathcal{J}_0^2r^4\sin\varphi}
    {\bigl(\mathcal{J}_0^2\mathcal{E}_0-2r^3\Phi'(r)\bigr)^3}
    \big(
    r\Phi''(r)+2\Phi'(r)
    +\frac{2r^3}{\mathcal{E}_0\mathcal{J}_0^2}\bigl(\Phi'(r)\bigr)^2
    \big)d\varphi ,
    \end{aligned}
\end{equation}
where $\varphi$ is defined by
\[
\sin^2\varphi=\frac{\mathcal{J}_0^2}{r^2}+\frac{4\Phi(r)}{\mathcal{E}_0}.
\]

In view of the formula giving $\partial_{\mathcal{J}_0}\Theta$, it turns out Assumption \ref{ass.phi2}  implies $\partial_{\mathcal{J}_0}\Theta>0$ for all $\mathcal{J}_0\in(0,1)$, and also the limits
\[
\lim_{\mathcal{J}_0\to 0}\partial_{\mathcal{J}_0}\Theta\in(0,\infty)
\quad\text{and}\quad
\lim_{\mathcal{J}_0\to 1}\partial_{\mathcal{J}_0}\Theta=0
\]
as long as $\Phi'(1)=0$ (if not then $\lim_{\mathcal{J}_0\to 1}\partial_{\mathcal{J}_0}\Theta=\infty$).

Finally, the deflection angle $\pi-2\Theta$ is the angle between $u=v_1-v_2$ and $u'=u-2(u\cdot\omega_{12})\omega_{12}$, we have 
$$1-2(\hat u\cdot\omega_{12})^2=\hat u\cdot\hat u'=\cos(\pi-2\Theta)=1-2\cos^2\Theta$$
where $\hat u=u/|u|$ denotes the unit vector in the direction of $u$. Therefore $|u\cdot\omega_{12}|=|u|\cos\Theta$ and the result follows. 
\end{proof}

\begin{proposition}\label{prop.est.b}
Under Assumptions \ref{ass.phi1}--\ref{ass.phi3}, let $u=v_1-v_2\neq0$, with $|v_1|,|v_2|\leq|\log\eps|^{C^*}$. Recall $C_\Phi=(3C^0_\Phi+2)(6C_\Phi^0+d+3)$ from Definition \ref{def.notation} \eqref{it.defnot3}. Then
\begin{enumerate}
\item\label{it.est.b1}   For every measurable set $E\subseteq\{(v_2,\omega)\in\Rb^d\times\Sb^{d-1}:|v_2|\leq|\log\eps|^{C^*}\}$,
\[
\int_E b(|v_1-v_2|,\omega)\,\d\omega\,\d v_2\leq |\log\eps|^{C^*}|E|^{1/C_\Phi}.
\]
\item\label{it.est.b2}   If $|u\cdot\omega|\geq\lambda$ for some $0<\lambda\leq1$, then
\[
b(|u|,\omega)\leq |\log\eps|^{C^*}\lambda^{-4(3C^0_\Phi+1)}.
\]
\end{enumerate}
\end{proposition}

\begin{proof} Denote by $\Ec_0=|u|^2\leq|\log\eps|^{C^*}$ and $h:=1-\Jc_0$. The proof has three parts.

{\bf Proof part 1.} We first show that 
\begin{equation}\label{eq.est.b.pw}
b(|u|,\omega)\le C\langle\Ec_0\rangle^4\min(1,\Ec_0)^{-(3C^0_\Phi+1)} h^{-(3C^0_\Phi+1)},
\end{equation}
where $\langle \Ec_0 \rangle = \sqrt{1+ \Ec^2_0}$. By Assumptions \ref{ass.phi1} and \ref{ass.phi3}, there exist constants $c_0,M>0$ such that
\[
c_0(1-r)^{C^0_\Phi}\le|\Phi'(r)|\le M,\qquad \Phi(r)\le M(1-r)\qquad\text{for } r\in[1/2,1)
\]
(if $\Phi'(1^-)\neq0$, take $C^0_\Phi=0$). By \eqref{eqn.Theta}, we have $\Theta\ge\arcsin\Jc_0$ and $\sin\Theta\ge\Jc_0$. Together with \eqref{eqn.defb}, we have $b(u,\omega) \le |u|\cdot (\partial_{\Jc_0}\Theta)^{-1}$.
By \eqref{eqn.p_theta} and Assumption \ref{ass.phi2},
$$\partial_{\Jc_0}\Theta\ge\int_{\arcsin\Jc_0}^{\pi/2}\frac{4\Ec_0r^7\Phi'(r)^2\sin\varphi}{D^3}\,d\varphi,
\qquad \sin^2\varphi=S(r):=\frac{\Jc_0^2}{r^2}+\frac{4\Phi(r)}{\Ec_0},$$
where $D:=\Ec_0\Jc_0^2+2 r^3|\Phi'(r)|$. Since $|S'(r)|=2\Jc_0^2/r^3+4|\Phi'(r)|/\Ec_0\ge4|\Phi'|/\Ec_0$, the map $r\mapsto\varphi$ is a bijection from $[r_*,1]$ onto $[\arcsin\Jc_0,\pi/2]$, and $|d\varphi|=|S'|\,dr/(2\sin\varphi\cos\varphi)$. Hence, using $\cos\varphi\le1$,
\[
\partial_{\Jc_0}\Theta\ge\int_{r_*}^1\frac{8r^7|\Phi'|^3}{D^3} \d r\ge c\int_{\max(r_*,1/2)}^1\min\Big(1,\frac{|\Phi'(r)|^3}{\Ec_0^3}\Big)\d r.
\]
The second inequality is because of following: if $r^3|\Phi'|\ge\Ec_0\Jc_0^2$ then $D\le 3r^3|\Phi'|$ and the integrand is $\ge c$, and otherwise $D\le 3\Ec_0$. In both cases the inequality holds.

Let $\delta':=1-\max(r_*,1/2)$. We claim that $\delta'\ge c\min(1,\Ec_0)\,h .$ Recall that $S(r_*)=1$ (Proposition \ref{prop.coll.kernel}), that $r_*\ge\Jc_0$, and that $S(1)=\Jc_0^2$ since $\Phi(1)=0$. We consider three cases.
\begin{enumerate}[(a)]
\item $h\le1/2$. Then $\Jc_0=1-h\ge1/2$, so $r_*\ge\Jc_0\ge1/2$ and $\delta'=1-r_*$. For $r\in[r_*,1]$ we have $\Jc_0\le r$ and $r\ge1/2$, so $|S'(r)|=2\Jc_0^2/r^3+4|\Phi'(r)|/\Ec_0\le 2/r+4M/\Ec_0\le 8+4M/\Ec_0$. Since $\Jc_0\le1$, we get
\[
h=1-\Jc_0\le1-\Jc_0^2=S(r_*)-S(1)\le(8+4M/\Ec_0)(1-r_*),
\]
hence $\delta'=1-r_*\ge h/(8+4M/\Ec_0)\ge c\min(1,\Ec_0)\,h$.
\item $h>1/2$ and $r_*<3/4$. Then $\max(r_*,1/2)<3/4$, so $\delta'>1/4\ge h/4$.
\item $h>1/2$ and $r_*\ge3/4$. Then $\delta'=1-r_*$, and $\Jc_0=1-h<1/2$, so $\Jc_0/r_*<2/3$. Using $S(r_*)=1$ and $\Phi(r_*)\le M(1-r_*)$, we get
\[
\frac59\le1-\frac{\Jc_0^2}{r_*^2}=\frac{4\Phi(r_*)}{\Ec_0}\le\frac{4M(1-r_*)}{\Ec_0},
\]
hence $\delta'=1-r_*\ge5\Ec_0/(36M)\ge c\min(1,\Ec_0)\,h$.
\end{enumerate}
This proves the claim. Using the lower bound of $|\Phi'|$ on $[1-\delta',1-\delta'/2]$, we obtain
\[
\partial_{\Jc_0}\Theta\ge c\,\delta'\min\big(1,\delta'^{3C^0_\Phi}\Ec_0^{-3}\big)\ge c\langle\Ec_0\rangle^{-3}\min(1,\Ec_0)^{3C^0_\Phi+1}h^{3C^0_\Phi+1},
\]
which, with $|u|\le\langle\Ec_0\rangle^{1/2}$, proves \eqref{eq.est.b.pw}.

{\bf Proof part 2.} Now we prove \eqref{it.est.b2}. By Proposition \ref{prop.coll.kernel}, $|u\cdot\omega|=|u|\cos\Theta$ and $\Theta\geq\arcsin\Jc_0$, so
\[
\frac{\lambda}{|u|}\leq\cos\Theta\leq\frac\pi2-\Theta\leq\arccos\Jc_0\leq\frac\pi2\sqrt{1-\Jc_0^2}\leq\frac{\pi}{2}\sqrt{2h},
\]
hence $h\geq\lambda^2/(9|u|^2)$. Moreover $\Ec_0\geq|u\cdot\omega|^2\geq\lambda^2$, so $\min(1,\Ec_0)^{-(3C_\Phi^0+1)}\leq\lambda^{-2(3C^0_\Phi+1)}$. Inserting these into \eqref{eq.est.b.pw}, and using $\Ec_0\leq|\log\eps|^{C^*}$, gives $b\leq|\log\eps|^{C^*}\lambda^{-4(3C^0_\Phi+1)}$.

{\bf Proof part 3.} Now we prove \eqref{it.est.b1}. If $|E|>1$, the bound follows from $\int b\,\d\omega=c_d|v_1-v_2|\le|\log\eps|^{C^*}$ (Lemma \ref{lem.liouville}), so we assume $|E|\le1$.

First fix $v_2$ and a measurable $A\subseteq\Sb^{d-1}$. We show that
\begin{equation}\label{eq.est.b.sphere}
\int_A b(|u|,\omega)\,\d\omega\le C\langle\Ec_0\rangle^4\min(1,\Ec_0)^{-(3C^0_\Phi+1)}|A|^{1/(3C^0_\Phi+2)}.
\end{equation}

Let $\theta_0\in[0,\pi/2]$ be the angle between $\nu$ and $-\hat u$, so that $h=1-\sin\theta_0$ and $[u\cdot\nu]_-=|u|\cos\theta_0$. Let $\eta\in(0,1/2]$. If $h<\eta$, then $|\nu\cdot\hat u|=\cos\theta_0)$ and $\sin\theta_0>1-\eta$ give $\pi/2-\theta_0<\arccos(1-\eta)\le2\sqrt\eta$, hence $|\nu\cdot\hat u|\le\pi/2-\theta_0<2\sqrt\eta$. On the band $\{\nu:|\nu\cdot\hat u|<2\sqrt\eta\}$ we have $[u\cdot\nu]_-\le2\sqrt\eta\,|u|$, and this band has measure at most $C\sqrt\eta$. Since $b\,\d\omega=[u\cdot\nu]_-\,\d\nu$, the set $\{h<\eta\}$ contributes at most $C|u|\,\eta$.   

On $\{h\ge\eta\}$, \eqref{eq.est.b.pw} gives $b\le C\langle\Ec_0\rangle^4\min(1,\Ec_0)^{-(3C^0_\Phi+1)}\eta^{-(3C^0_\Phi+1)}$, so this set contributes at most this bound times $|A|$. If $|A|^{1/(3C^0_\Phi+2)}\le1/2$, we take $\eta=|A|^{1/(3C^0_\Phi+2)}$; then $\eta^{-(3C^0_\Phi+1)}|A|=|A|^{1/(3C^0_\Phi+2)}$, and since $|u|\le\langle\Ec_0\rangle^{1/2}$, both contributions are bounded by the right hand side of \eqref{eq.est.b.sphere}. If $|A|^{1/(3C^0_\Phi+2)}>1/2$, then \eqref{eq.est.b.sphere} follows from $\int b\,\d\omega=c_d|u|\le2c_d|u|\,|A|^{1/(3C^0_\Phi+2)}$.

Now let $0<\delta\le1$, and write $E_{v_2}:=\{\omega:(v_2,\omega)\in E\}$, so that $\int|E_{v_2}|\,\d v_2=|E|$. The part of $E$ with $|v_1-v_2|\le\delta$ contributes at most $c_d\int_{|v_1-v_2|\le\delta}|v_1-v_2|\,\d v_2\le C\delta^{d+1}$. On the part with $|v_1-v_2|\ge\delta$, we have $\Ec_0=|v_1-v_2|^2\in[\delta^2,|\log\eps|^{C^*}]$, so the constant in \eqref{eq.est.b.sphere} is at most $|\log\eps|^{C^*}\delta^{-(6C^0_\Phi+2)}$. Applying \eqref{eq.est.b.sphere} to $A=E_{v_2}$ and H\"older's inequality in $v_2$ over $\{|v_2|\le|\log\eps|^{C^*}\}$, this part contributes at most
\[
|\log\eps|^{C^*}\delta^{-(6C^0_\Phi+2)}\int_{|v_2|\le|\log\eps|^{C^*}}|E_{v_2}|^{1/(3C^0_\Phi+2)}\d v_2\le|\log\eps|^{C^*}\delta^{-(6C^0_\Phi+2)}|E|^{1/(3C^0_\Phi+2)}.
\]
We choose $\delta=|E|^{1/C_\Phi}\le1$, where $C_\Phi=(3C^0_\Phi+2)(6C^0_\Phi+d+3)$. Then $\delta^{d+1}=\delta^{-(6C^0_\Phi+2)}|E|^{1/(3C^0_\Phi+2)}=|E|^{(d+1)/C_\Phi}\le|E|^{1/C_\Phi}$, which proves \eqref{it.est.b1}.
\end{proof}
 
\begin{proposition}[Upper bound on the collision duration]
\label{prop.upperbond.interaction}
Let $[t^-,t^+]$ be a collision interval involving the particles $\pf_1,\ldots,\pf_n$, which evolve according to the short-range dynamics in Definition~\ref{def.short_range}. Assume that these particles do not interact with any particle outside
$\{\pf_1,\ldots,\pf_n\}$ during $[t^-,t^+]$.

Define
$$ x_c(t):=\frac1n\sum_{i=1}^n x_i(t) \qquad \Ic(t):=\sum_{i=1}^n |x_i(t)-x_c(t)|^2. $$
Suppose that, for some $\eps\lesssim\mu\lesssim1$, $ \Ic(t^-),\ \Ic(t^+)\lesssim \mu^2. $
Let
$$ v_c:=\frac1n\sum_{i=1}^n v_i(t^-),\qquad c_{\rm rel}:=\sum_{i=1}^n |v_i(t^-)-v_c|^2. $$
Then
$$c_{\rm rel}(t^+-t^-)^2\leq C(n)\mu^2. $$
\end{proposition}
\begin{proof}
By conservation of momentum, the velocity $ v_c(t):=\frac1n\sum_{i=1}^n v_i(t)$ is independent of time. Set
$$y_i(t):=x_i(t)-x_c(t),\qquad  w_i(t):=v_i(t)-v_c,\qquad q_{ij}(t):=\frac{x_i(t)-x_j(t)}{\eps}.$$
By the conservation of energy, we have
\begin{equation}\label{eq.relative.energy}
    \sum_{i=1}^n |w_i(t)|^2+\sum_{\substack{1\leq i,j\leq n\\ i\neq j}}\Phi(q_{ij}(t))=c_{\rm rel} +\sum_{\substack{1\leq i,j\leq n\\ i\neq j}}\Phi(q_{ij}(t^-)).
\end{equation}
By differentiating $\Ic(t)$ twice, we obtain
$$\begin{aligned}
    \frac{\d \Ic(t)}{\d t}=2\sum_{i=1}^n y_i(t)\cdot w_i(t), \qquad \frac{\d^2 \Ic'(t)}{\d t^2}=2\sum_{i=1}^n |w_i(t)|^2+2\sum_{i=1}^n y_i(t)\cdot \dot v_i(t).
\end{aligned}$$
Using the symmetry property, we have
$$2\sum_{i=1}^n y_i(t)\cdot\dot v_i(t) = - \sum_{\substack{1\leq i,j\leq n\\ i\neq j}}    q_{ij}(t)\cdot\nabla\Phi(q_{ij}(t)).$$
By Assumption \ref{ass.phi2}, and using the fact that $\Phi(1) =0$, $\Phi'(r)\leq 0$, we have $r \Phi'(r) + \frac{1}{2} \Phi(r) \leq 0$. Therefore, by \eqref{eq.relative.energy}, we obtain 
$$\begin{aligned}
  \Ic''(t) &\geq2\sum_{i=1}^n |w_i(t)|^2+\frac{1}{2}\sum_{\substack{1\leq i,j\leq n\\ i\neq j}}\Phi(q_{ij}(t))  \\
    &\geq \frac{1}{2}\Big(\sum_{i=1}^n |w_i(t)|^2+ \sum_{\substack{1\leq i,j\leq n\\ i\neq j}}  \Phi(q_{ij}(t)) \Big) \geq \frac{1}{2}\Big(c_{\rm rel} +\sum_{\substack{1\leq i,j\leq n\\ i\neq j}}\Phi(q_{ij}(t^-))\Big):= c_0.
\end{aligned}$$ 
This means that $\Ic(t)$ is a convex function. Let $ t_m:=\frac{t^-+t^+}{2},$ and define the function $$F(t):= \Ic(t) - \frac{c_0  }{2}(t-t_m)^2.$$ Then by elementary calculation, $F(t)$ is also a convex function. Therefore, $F(t_m)\leq\frac{F(t^-) + F(t^+)}{2}$, which further implies
$$  \Ic(t_m)  \leq  \frac{\Ic(t^-)+\Ic(t^+)}{2} - \frac{c_0}{8} (t^+-t^-)^2.$$
Since $\Ic(t_m)\geq0$, it follows that
$$\frac{c_0}{8}(t^+-t^-)^2 \leq \frac{\Ic(t^-)+\Ic(t^+)}{2}\leq C_0\mu^2.$$
This proves the desired estimate.
\end{proof}

\begin{lemma}\label{lem.tree-count}
Let all constants below be integers. Let $\{B_{m,x,n}\}_{m,x,n\geq0}$ be a family of nonnegative numbers, satisfying
\begin{equation}\label{eqn.tree.B}
B_{m,x,n}\leq\sum_{l=3 }^{n}\sum_{x_0 =1}^{[l/2]} \sum_{\substack{m_1+\cdots+m_l =m-1\\x_1+\cdots+x_l =x-x_0\\n_1+\cdots+n_l =n-l }}C^l \prod_{j=1}^{l }Q_{m_j,x_j,n_j},
\end{equation}
 and
\begin{equation}\label{eqn.tree.Q}
Q_{m,x,n}=\sum_{i=0}^{m}\sum_{r=0}^{x}\sum_{k=0}^{n}(5i+\mathsf{N} \cdot r + 1)\cdot (5(m-i) + \mathsf{N}\cdot (x-r) + 1)\cdot B_{i,r,k}B_{m-i,x-r,n-k}.
\end{equation}
Assume that
\begin{equation}\label{eqn.B000}
\begin{aligned}
         &B_{0,0,0}=1,\qquad B_{m,0,n}=0\quad\textrm{for}\quad\! n\geq1,\qquad B_{0,x,n}=0\quad\textrm{for}\quad \!(x,n)\neq(0,0),\\
         &  B_{m,x,n} = 0 \quad\textrm{if}\quad\! 1\leq  n < 2 x
\end{aligned}
\end{equation}
and $B_{m,0,0}\leq C_0^m.$ Moreover, we assume $m\leq |\log\eps|^{C^*},$ then $$ B_{m,x,n} \leq C^{m+n}_1|\log\eps|^{C^*x}$$ for some constant $C_1>0$ depending on $\mathsf{N}$ but independent of $\eps$.
\end{lemma}
\begin{proof}
Define the generating functions $$ \mathcal B(z,u,v) := \sum_{m,x,n\geq0}B_{m,x,n}z^mu^xv^n, \qquad \mathcal Q(z,u,v) := \sum_{m,x,n\geq0}Q_{m,x,n}z^mu^xv^n.$$
Let $D=z\partial_z$, $E = u \p_u$ and $L = I + 5D + \mathsf{N} E$. For two different generating functions $\Fc$ and $\Gc$, we write $\Fc\preceq \Gc$ when every coefficient of $\Fc$ is bounded by the corresponding coefficient of $\Gc$.

Since $L(z^m u^x v^n) = (1 + 5m + \mathsf{N}x) z^m u^x v^n$, the formula \eqref{eqn.tree.Q} is equivalent to
\begin{equation}\label{eqn.tree.genQ}
\mathcal Q=(L\mathcal B)^2.
\end{equation}
 Let $$ \mathcal B_0(z) := \sum_{m\geq0}B_{m,0,0}z^m \quad \textrm{and} \quad\! \Gc(z,u,v) := \Bc(z,u,v) -\Bc_0.$$
Summing \eqref{eqn.tree.B} over $m,x,n$, we have 
\begin{equation}\label{eqn.B-B0}
\begin{aligned}
     \Bc(z,u,v) - \Bc_0 \preceq  &\sum_{l=3}^n \sum_{x_0=1}^{[l/2]} \sum_{m,x,n\geq 0} z^m u^x v^n  \sum_{\substack{m_1+\cdots+m_l =m-1\\x_1+\cdots+x_l =x-x_0\\n_1+\cdots+n_l =n-l }}\prod_{j=1}^{l }Q_{m_j,x_j,n_j} \\
    =&\sum_{l =3}^n \sum_{x_0=1}^{[l/2]} \sum_{m_j,x_j,n_j\geq 0} z^{1+\sum_j m_j} u^{x_0+\sum_jx_j} v^{l+\sum_jn_j}  
    C^l\prod_{j=1}^{l }Q_{m_j,x_j,n_j} \\
    =& z \sum_{l=3}^n (C v\cdot\Qc)^l\sum_{x_0=1}^{[l/2]}u^{x_0} ,
\end{aligned}
\end{equation}
Moreover, $ \mathcal B_0(z) \preceq \frac{1}{1-C_0z},$ so $\Bc_0$ is analytic near the origin. Since $E \Bc_0 =0$, we have $$L\Bc = (I+5D)\Bc_0+L\Gc.$$ Using \eqref{eqn.tree.genQ} and  \eqref{eqn.B-B0}, we obtain
\begin{equation}\label{eqn.tree.G}
\begin{aligned}
\Gc &\preceq  z \sum_{l\geq 3}(C v\cdot\Qc)^l\sum_{x_0=1}^{[l/2]}u^{x_0}  = z \Theta(z,u,v,L\Gc(z,u,v)),
\end{aligned}
\end{equation}
where  
\begin{equation}\label{eqn.Capital_theta}
      \Theta(z,u,v,w) := \sum_{l\geq 3}\Big[Cv\big((I+5D)\Bc_0(z)+ w)^2\big)\Big]^l \sum_{x_0=1}^{[l/2]} u^{x_0}. 
\end{equation}

The function $\Theta$ is analytic near the origin. Write
$$ \Theta(z,u,v,w) = \sum_{a,x_0,b,q\geq0}\theta_{a,x_0,b,q}z^a u^{x_0}v^b w^q. $$
By the multivariable Cauchy estimate, there exist constant $K>1$ such that 
\begin{equation}\label{eqn.tree.theta}
0\leq \theta_{a,x_0,b,q} \leq K^{1+a+x_0+b+q}.
\end{equation} 

We now iterate \eqref{eqn.tree.G} at the level of coefficients. Every term contributing to $z^m u^x v^n $ is represented by a rooted ordered tree. Indeed, selecting at a vertex $\vf$ the monomial $\theta_{a_\vf,x_{0,\vf},b_\vf,q_\vf}z^{a_\vf} u^{x_{0,\vf}}v^{b_\vf}(L\Gc)^{q_\vf}$ produces a vertex with $q_\vf$ ordered children. Together with the factor $z$ in \eqref{eqn.tree.G}, this vertex $\vf$ contributes $\theta_{a_\vf,x_{0,\vf},b_\vf,q_\vf}z^{1+a_\vf} u^{x_{0,\vf}}v^{b_\vf}$. Suppose that the tree has $s$ vertices then it has $s-1$ edges. We have 
$$ \sum_{\mathsf v}q_{\mathsf v}=s-1.$$
Matching powers of $z,u,v$, we also have
$$ \sum_{\mathsf v}a_{\mathsf v}=m-s, \qquad \sum_\vf x_{0,v} =x, \qquad \sum_{\mathsf v}b_{\mathsf v}=n.$$
Since $a_{\mathrm v}\ge 0$ and $x_{0,\mathrm v}\ge1$, it follows that
$$s\le m, \qquad s\le x.$$
Moreover, by \eqref{eqn.Capital_theta}, we have $b_\vf = l$ and $x_{0,\vf} \leq [l/2]$ at every vertex. Therefore, one has
$$ n=\sum_{\mathrm v}b_{\mathrm v} \ge 2\sum_{\mathrm v}x_{0,\mathrm v} =2x.$$
Therefore, by \eqref{eqn.tree.theta}, the contribution of the $s$ vertices is 
\begin{equation}\label{eqn.theta_est}
    \prod_{\vf}\theta_{a_\vf,x_{0,\vf},b_\vf,q_\vf} \leq K^{\sum_{\vf} (1 + a_\vf + x_{0,\vf} +b_\vf+q_\vf)} \leq K^{s+ m-s+ x + n+s-1} \leq C^{m+n}.
\end{equation}


We next estimate the factors produced by operator $L$. For each edge, the corresponding application of $L$ contributes the factor
$$ 1+5m +\mathsf{N} x,$$
hence, the total contribution is no more than $(1+5m +\mathsf{N} x)^{s-1}\leq C^s \cdot (1 + m + x)^s$. By the elementary inequality $1+t \leq e^t$ and $s \leq x$, we have
\begin{equation}\label{eqn.log_eps_x}
    (1+5m +\mathsf{N} x)^{s-1}\leq C^s \cdot (1 + m + x)^s \leq C^s (1+m)^s(1+\frac{x}{1+m})^s \leq (Ce (1+m))^s \leq |\log \eps|^{C^* x}.
\end{equation}

It remains to count the possible trees and labels. For fixed $s$, the number of rooted ordered trees is at most \(4^{s-1}\). The numbers of possible choices of ${a_{\mathrm v}}, {x_{0,\mathrm v}}$, and \({b_{\mathrm v}}\) are bounded respectively by
\begin{equation}\label{eqn.binom}
     \binom{m-1}{s-1}\leq 2^m, \quad \binom{x-1}{s-1}\leq 2^x,\qquad \binom{n-1}{s-1}\leq 2^{n}. 
\end{equation}  Combining \eqref{eqn.theta_est}--\eqref{eqn.binom}, and summing over \(1\le s\le\min(m,x)\), we obtain
$$ B_{m,x,n} \leq C^{m+n}|\log \eps|^{C^* x}.  $$
This completes the proof. 
\end{proof}

\section*{Acknowledgements}   Y. D. was supported by NSF DMS-2531437 and a Simons Collaboration Grant on Wave Turbulence. N. J. was supported by the National Natural Science Foundation of China under Grant No. 12371180. 

Shortly before completing the manuscript, OpenAI released a series of papers on October 6, among which Paper NO.364 \cite{OpenAI26} is closely related to ours but appears to take a different approach. We emphasize that our project was conceived independently: it began during the third author's visit to Chicago in the fall of 2025, where the central idea originated.

\section*{AI use disclosure}   
The first complete draft of this paper was written by the authors, with limited AI assistance in checking elementary calculations and with the proof of Lemma \ref{lem.tree-count}. After that, the authors used OpenAI's ChatGPT and Anthropic's Claude to assist with language editing and checking proofs. The authors independently verified all AI-assisted revisions and take full responsibility for the content of this paper.

\bigskip 


\end{document}